\documentclass[a4paper, 12pt]{book}

\usepackage{titling}
\title{Combinatorics of the Permutahedra, Associahedra, and Friends}
\author{Viviane Pons}

\usepackage{amsmath, amsfonts, amssymb, amsthm}
\usepackage{etoolbox} 
\usepackage{dsfont} 

\usepackage{hyperref}

\usepackage{xcolor}

\usepackage{graphicx}

\usepackage{listings}

\usepackage{multirow}

\usepackage{hvfloat}
\usepackage{tikz}
\usetikzlibrary{positioning,shapes,shadows,arrows}
\usetikzlibrary{shapes.multipart}
\usetikzlibrary{automata}
\usetikzlibrary{patterns}
\usetikzlibrary{fit}

\usepackage{vmargin}
\usepackage[Glenn]{fncychap}
\usepackage{fancyhdr}

\usepackage{xargs}

\newlength{\beforechapter}
\newlength{\afterchapter}

\usepackage{geometry}
\usepackage{colortbl}

\usepackage{glossaries}

\usepackage{bookmark}

\makeatletter
  \renewcommand{\DOCH}{%
    \settoheight{\myhi}{\CTV\FmTi{Test}}
    \setlength{\py}{\baselineskip}
    \addtolength{\py}{\RW}
    \addtolength{\py}{\myhi}
    \setlength{\pyy}{\py}
    \addtolength{\pyy}{-1\RW}
    \vskip \beforechapter
    \raggedright
    \CNV\FmN{\@chapapp}\space\CNoV\thechapter
    \hskip 3pt\mghrulefill{\RW}\rule[-1\pyy]{2\RW}{\py}\par\nobreak}

\def\@makeschapterhead#1{%
    \vspace*{50\p@}%
    \vskip \beforechapter
    {\parindent \z@ \raggedright
        \normalfont
        \interlinepenalty\@M
        \DOTIS{#1}
        \vskip -40\p@
        \vskip \afterchapter}
    }

\makeatother

\makeatletter
\renewenvironment{thebibliography}[1]
     {\subsection*{\bibname}
      \@mkboth{\MakeUppercase\bibname}{\MakeUppercase\bibname}%
      \list{\@biblabel{\@arabic\c@enumiv}}%
           {\settowidth\labelwidth{\@biblabel{#1}}%
            \leftmargin\labelwidth
            \advance\leftmargin\labelsep
            \@openbib@code
            \usecounter{enumiv}%
            \let\p@enumiv\@empty
            \renewcommand\theenumiv{\@arabic\c@enumiv}}%
      \sloppy
      \clubpenalty4000
      \@clubpenalty \clubpenalty
      \widowpenalty4000%
      \sfcode`\.\@m}
     {\def\@noitemerr
       {\@latex@warning{Empty `thebibliography' environment}}%
      \endlist}
\makeatother

\usepackage{multibib}
\newcites{ext}{Other bibliographic references}
\newcites{me}{International journals}
\newcites{meconf}{International conference proceedings}
\newcites{mepre}{Preprints}
\newcites{mesoft}{Software development}
\newcites{memisc}{Other}

\setmarginsrb{3.2cm}{2cm}{3.2cm}{2.5cm}{1cm}{1cm}{1cm}{1cm}

\makeglossaries

\newtheorem{theorem}{Theorem}
\newtheorem{conjecture}[theorem]{Conjecture}

\newcommand{\NN}{\mathbb{N}}
\newcommand{\K}{\mathbf{k}} 

\newcommand{\Sym}[1]{\mathfrak{S}_{#1}} 

\newcommand{\wole}{\preccurlyeq} 
\newcommand{\woless}{\prec} 
\newcommand{\meet}{\wedge} 
\newcommand{\join}{\vee} 
\DeclareMathOperator{\inv}{inv} 

\newcommand{\biseq}[1]{\tiny{\textbf{#1}}}

\definecolor{darkblue}{rgb}{0,0,0.7} 
\newcommand{\darkblue}{\color{darkblue}} 
\newcommand{\defn}[1]{\emph{\darkblue #1}} 

\newcommand{\Sage}{\gls{sage}}

\newcommand{\hopf}{\mathcal{H}}
\newcommand{\hprod}{\cdot}
\newcommand{\hcop}{\Delta}
\newcommand{\tensor}{\bigotimes}
\DeclareMathOperator{\FQSym}{\mathbf{FQSym}}
\DeclareMathOperator{\PBT}{\mathbf{PBT}}
\DeclareMathOperator{\BF}{\mathbb{F}} 
\DeclareMathOperator{\BG}{\mathbb{G}} 
\DeclareMathOperator{\BP}{\mathbb{P}} 
\DeclareMathOperator{\BQ}{\mathbb{Q}} 

\newcommand{\red}[1]{\textbf{\textcolor{red}{#1}}}
\newcommand{\sred}[1]{\textcolor{red}{#1}}
\definecolor{mygreen}{RGB}{23,103,1}
\definecolor{Prune}{RGB}{99,0,60}

\newcommand{\blue}[1]{\textcolor{blue}{#1}}
\newcommand{\green}[1]{\textcolor{mygreen}{#1}}

\tikzstyle{Red} = [color = red]
\tikzstyle{Blue} = [color = blue]
\tikzstyle{alert} = [color=red, line width = 2]
\tikzstyle{bluealert} = [color=blue, line width =1.5]

\tikzstyle{Path} = [line width = 1.2]
\tikzstyle{StrongPath} =  [line width=2.5]
\tikzstyle{DPoint} = [fill, radius=0.1]
\tikzstyle{Line1} = [dashed]
\tikzstyle{Line2} = [dotted, ultra thick]
\tikzstyle{GrayPath} = [color = gray]

\tikzstyle{Point} = [fill, radius=0.08]

\def \interscale{0.5}

\newcommand\dec{F_{\ge}} 
\newcommand\inc{F_{\le}} 
\newcommand{\trprec}{\vartriangleleft} 

\DeclareMathOperator{\BB}{\mathbb{B}} 
\DeclareMathOperator{\BT}{\mathcal{B}}
\DeclareMathOperator{\Bm}{B^{(m)}}
\DeclareMathOperator{\size}{size} 
\DeclareMathOperator{\graftingTree}{\Delta} 
\DeclareMathOperator{\leftbranch}{\phi} 
\DeclareMathOperator{\compl}{\psi} 
\DeclareMathOperator{\expand}{expand} 
\DeclareMathOperator{\contract}{contract} 
\DeclareMathOperator{\contactsV}{\mathbf{C}} 
\DeclareMathOperator{\risesV}{\mathbf{R}} 
\DeclareMathOperator{\contactsP}{\mathcal{C}} 
\DeclareMathOperator{\risesP}{\mathcal{R}} 

\newcommandx{\rel}[1][1=R]{\mathbin{\mathrm{#1}}} 
\newcommandx{\notrel}[1][1=R]{\mathbin{\!\raisebox{.02cm}{$\not$}\hspace{.02cm}\mathrm{#1}\hspace*{.01cm}}} 
\newcommandx{\Inc}[1]{#1^{\mathsf{Inc}}} 
\newcommandx{\Dec}[1]{#1^{\mathsf{Dec}}}
\newcommand{\IRel}{\mathcal{R}} 
\newcommand{\underprod}[2]{{#1}\backslash{#2}} 
\newcommand{\overprod}[2]{{#1}\slash{#2}} 
\newcommand{\FRel}{\BF\IRel} 
\newcommand{\product}{\cdot} 
\newcommand{\coproduct}{\triangle} 
\newcommand{\IAntisym}{\mathcal{A}} 
\newcommand{\ITrans}{\mathcal{T}} 

\newcommand{\WOEP}{\mathsf{WOEP}} 
\newcommand{\WOFP}{\mathsf{WOFP}} 
\newcommand{\WOIP}{\mathsf{WOIP}} 
\newcommand{\TOEP}{\mathsf{TOEP}} 
\newcommand{\TOFP}{\mathsf{TOFP}} 
\newcommand{\TOIP}{\mathsf{TOIP}} 
\newcommandx{\COEP}[1][1=\signature]{\mathsf{COEP}(#1)} 
\newcommandx{\COFP}[1][1=\signature]{\mathsf{COFP}(#1)} 
\newcommandx{\COIP}[1][1=\signature]{\mathsf{COIP}(#1)} 
\newcommand{\BOEP}{\mathsf{BOEP}} 
\newcommand{\BOIP}{\mathsf{BOIP}} 
\newcommandx{\PIP}[1][1=\orientation]{\mathsf{PIP}_{#1}} 
\newcommandx{\PEP}[1][1=\orientation]{\mathsf{PEP}_{#1}} 
\newcommandx{\PFP}[1][1=\orientation]{\mathsf{PFP}_{#1}} 
\newcommand{\orientation}{\mathbb{O}} 

\newcommand{\chapcitation}[2]{\textcolor{mygreen}{\textit{#1} \\ -- \textit{#2}} \vspace{2em}}

\newcommand{\noarrow}{\relbar\mkern-9mu\relbar}
\newcommand{\decoration}{\delta} 
\newcommand{\includeSymbol}[1]{\ensuremath{%
	\mathchoice
		{\raisebox{-.7mm}{\includegraphics[height=2.2ex]{includes/figures/nodes/#1}}}	
		{\raisebox{-.7mm}{\includegraphics[height=2.2ex]{includes/figures/nodes/#1}}}
		{\raisebox{-.6mm}{\includegraphics[height=1.6ex]{includes/figures/nodes/#1}}}
		{\raisebox{-.5mm}{\includegraphics[height=1ex]{includes/figures/nodes/#1}}}
}}
\robustify{\includeSymbol}
\newcommand{\noneCirc}{\includeSymbol{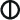}}
\newcommand{\upCirc}{\includeSymbol{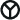}}
\newcommand{\downCirc}{\includeSymbol{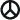}}
\newcommand{\upDownCirc}{\includeSymbol{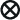}}

\newcommand{\Decorations}{\{\noneCirc{}, \downCirc{}, \upCirc{}, \upDownCirc{}\}} 
\newcommand{\up}[1]{\overline{#1}}

\newcommand{\down}[1]{\underline{#1}}

\newcommand{\updown}[1]{\up{\down{#1}}}

\newcommand{\factCatalan}[1]{\text{\bf C}(#1)} 
\renewcommand{\b}[1]{\mathbf{#1}} 
\newcommandx{\Permutreehedron}[1][1=\decoration]{\mathds{PT}(#1)} 
\newcommandx{\tree}[1][1=T]{\mathrm{#1}} 
\newcommand{\PPT}{\mathbb{P}} 
\newcommand{\automatonU}{\mathbb{U}} 

\DeclareMathOperator{\card}{\#} 
\newcommand{\HH}{\mathcal{H}\mathcal{H}}

\newglossaryentry{sage}{
    name=\texttt{SageMath},
    description={is a free open source mathematical software created in 2005 and mostly written in python. It includes a wide range of mathematics and is built on top of many other open-source packages. Its stated mission is to ``create a viable free open source alternative to Magma, Maple, Mathematica and Matlab''.  It has been developed mostly by researchers, based on free contributions evaluated through code-reviews}
}

\newglossaryentry{jupyter}{
    name=\texttt{Jupyter},
    description={is an open-source software used to create interactive web-based notebooks as computing environments in many languages
    }
}

\newglossaryentry{github}{
	name=\texttt{GitHub},
	description={is a private online platform used to host code and text-based projects using the \texttt{git} version management system}
}

\newglossaryentry{docker}{
	name=\texttt{Docker},
	description={is a system that provides easy to deploy virtual machines called \emph{containers} with specific environments. A \emph{Docker file} is a simple configuration file that describes the needed environment (sofware versions, installation commands, etc.)}
}

\newglossaryentry{zenodo}{
	name=\texttt{Zenodo},
	description={is a general-pupose open repository which allows researcher to deposit research data. It provides a version number and a DOI}
}

\newglossaryentry{PyCon}{
	name=\texttt{PyCon},
	description={is a the name given to general conferences about the Python programming language. The main even is the yearly event PyCon US which often welcomes a few thousand participants. PyCon Fr is also organized every year with a few hundred participants}
}

\newglossaryentry{OpenDreamKit}{
	name=\texttt{OpenDreamKit},
	description={stands for ``Open Digital Research Environments Toolkit for the Advancement of Mathematics''. It is a H2020 European project which funded many different open-source software and associated communities in mathematics. \url{https://opendreamkit.org/}}
}

\newglossaryentry{SageDays}{
    name=\texttt{Sage Days},
    description={are workshop dedicated to the software \texttt{SageMath}. They often include introduction talks and also development time around specific topics.}
}

\newglossaryentry{MathEnJean}{
	name=\texttt{MATh.en.JEANS},
	description={is a French organization which organizes mathematical research projects in schools involving researchers and pupils. The researchers propose mathematical problems on which the pupils work throughout the year. Once a year, pupils present their results in special events.} 
}

\newglossaryentry{PyLadies}{
	name=\texttt{PyLadies},
	description={is an international mentorship program to support women who code in python} 
}

\newglossaryentry{ECCO}{
	name=\texttt{ECCO},
	description={stands for ``Encuentro Colombiano De Combinatoria''. It is a summer school organized in Colombia every other year welcoming both Colombian and international students to follow advanced lectures in combinatorics}
}

\newglossaryentry{CIMPA}{
	name=\texttt{CIMPA},
	description={stands for ``Centre International de Mathématiques Pures et Appliquées''. It is an originally French organization (now funded by France, Spain, Norway, and Switzerland) whose goal is to promote mathematical research in developing countries}
}

\begin{document}

\nociteme{PPTJ23,PP20,CPP19,Pon19,PP18,CP15}

\nocitemeconf{LMP23,PPTJ21,CP19,PP18b,PP17,Pon15,CCP14,CP13}

\nocitemepre{CP23, CP22,Pon22}

\nocitemesoft{PonSage23, PonSage22, PonSage22b, PonSage18, CPSage18, PonSage16, PonSage14, PonSage13}

\nocitememisc{PonMisc22, Myr21, PonMisc19, PonMisc18, PonMisc17, PonMisc17b, PonMisc17c, PonMisc16, PonMisc15}

\frontmatter
\lhead[\oldstylenums \thepage]{\rightmark}
\rhead[\leftmark]{\oldstylenums \thepage}



\begin{titlepage}

	\newgeometry{left=8cm,bottom=2cm, top=3cm, right=1cm}

	\tikz[remember picture, overlay] \node[opacity=1,inner sep=0pt] at (-13mm,-140mm){\includegraphics{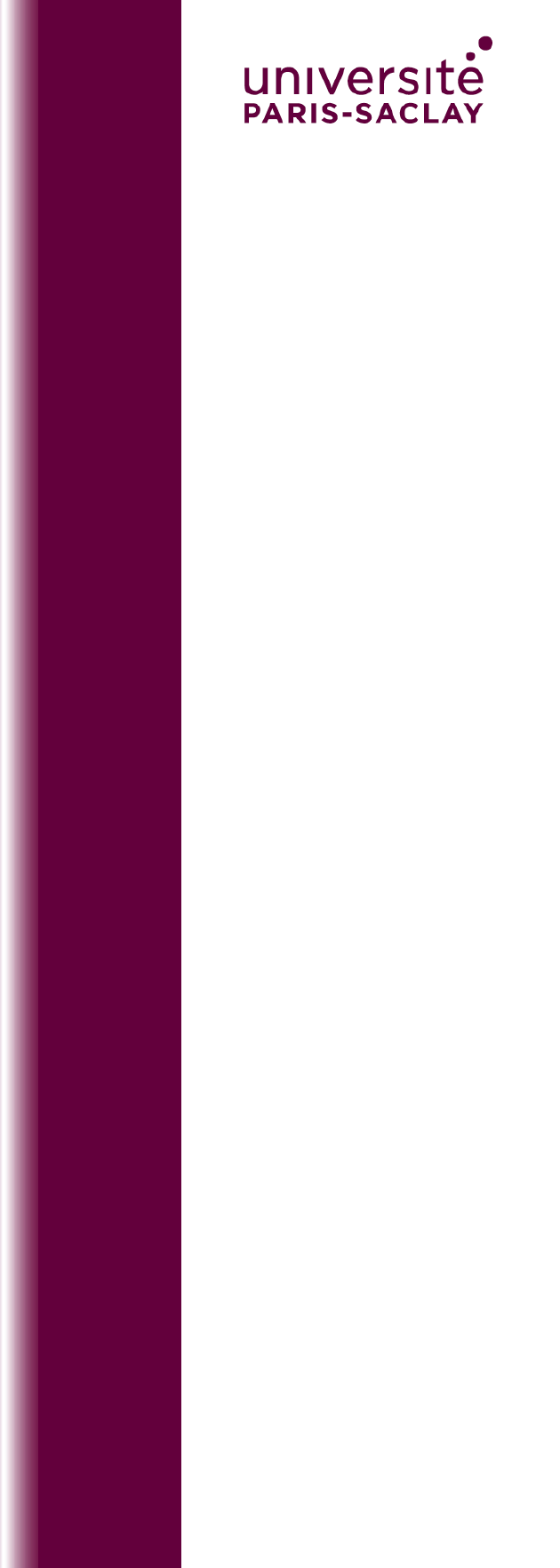}};
	\tikz[remember picture, overlay] \node[opacity=.5, inner sep=0pt] at (1cm, -15cm){\includegraphics[scale=.2]{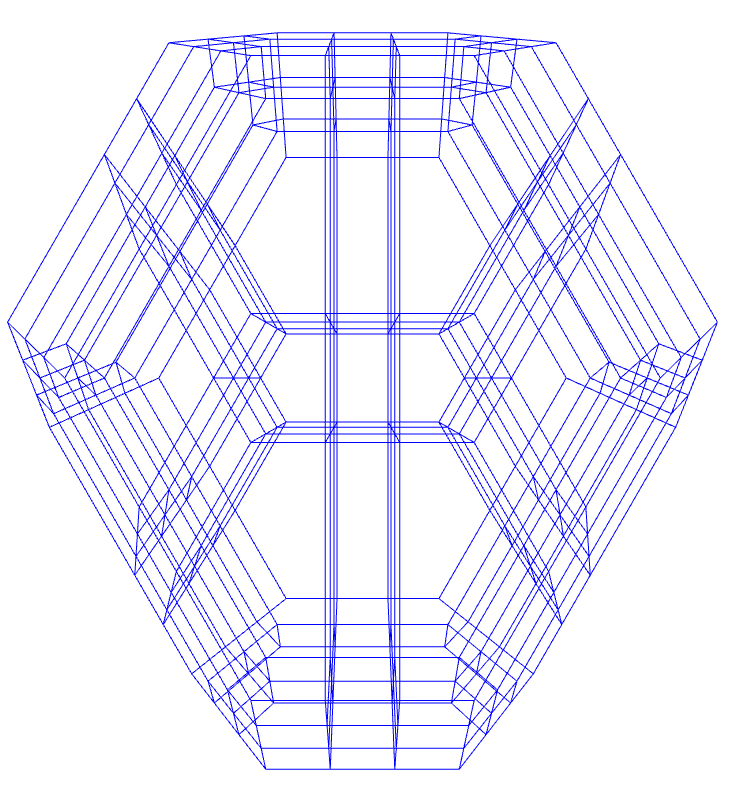}};
	\tikz[remember picture, overlay] \node[opacity=.5, inner sep=0pt] at (8cm, -15cm){\includegraphics[scale=.2]{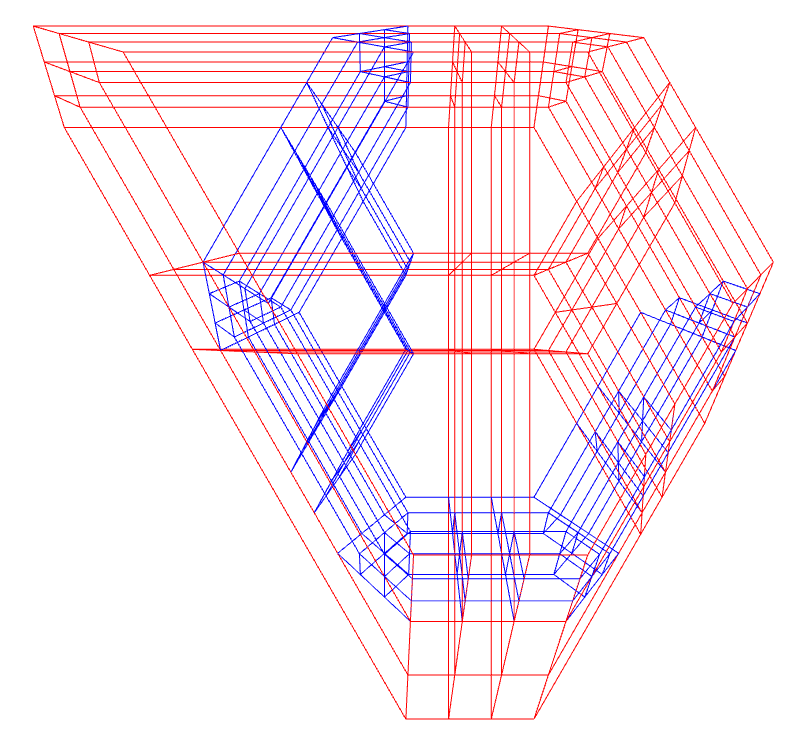}};


	\color{white}

	\begin{picture}(0,0)
		\put(-152,-743){\rotatebox{90}{\Large \textsc{Habilitation à diriger des recherches}}} \\
	\end{picture}


	\flushright{}
	\vspace{10mm} 
	\color{Prune}

	\fontsize{26}{26}\selectfont
	\thetitle\\

	\enskip


	\fontsize{8}{12}\selectfont

	\vspace{1.5cm}

	\normalsize
	\textbf{Habilitation à diriger des recherches}




	\textbf{Présentée et soutenue à l'Université Paris-Saclay le 20/10/2023, par}\\
	\bigskip
	\Large {\color{Prune} \textbf{\theauthor}} 

	\vspace{\fill} 

	\bigskip

	\flushleft{}
	\small {\color{Prune} \textbf{Composition du jury}}\\
	\vspace{2mm}
	\scriptsize
	\begin{tabular}{|p{8.25cm}l}
		\arrayrulecolor{Prune}
		\textbf{Sylvie CORTEEL}                              & Rapportrice \\
		Directrice de recherche, CNRS \& Université Paris-Cité &                           \\

		\textbf{Nathan READING}                                & Rapporteur \\
		Professor, North Carolina State University                                  &                           \\

		\textbf{Victor REINER}                              & Rapporteur              \\
		Professor, University of Minnesota        &                           \\

		\textbf{Mireille BOUSQUET-M\'ELOU}                              & Examinatrice               \\
		Directrice de recherche, CNRS \& Université de Bordeaux                                &                           \\

		\textbf{Florent HIVERT}                      & Examinateur               \\
		Professeur, Université Paris-Saclay       &                           \\

		\textbf{Lionel POURNIN}                               & Examinateur      \\
		Professeur, Université Sorbonne Paris Nord                 &                           \\

		\textbf{Maria RONCO}                               & Examinatrice     \\
		Profesora, Universidad de Talca                    &                           \\
	\end{tabular}

\end{titlepage}
\section*{Abstract}

I present an overview of the research I have conducted for the past ten years in algebraic, bijective, enumerative, and geometric combinatorics. The two main objects I have studied are the \defn{permutahedron} and the \defn{associahedron} as well as the two partial orders they are related to: the weak order on permutations and the Tamari lattice. This document contains a general introduction (Chapters~1 and~2) on those objects which requires very little previous knowledge and should be accessible to non-specialist such as master students. Chapters 3 to 8 present the research I have conducted and its general context. You will find:
\begin{itemize}
\item a presentation of the current knowledge on \defn{Tamari interval} and a precise description of the family of \defn{Tamari interval-posets} which I have introduced along with the \defn{rise-contact involution} to prove the symmetry of the \emph{rises} and the \emph{contacts} in Tamari intervals;
\item my most recent results concerning $q,t$-enumeration of Catalan objects and Tamari intervals in relation with \defn{triangular partitions};
\item the descriptions of the \defn{integer poset lattice} and \defn{integer poset Hopf algebra} and their relations to well known structures in algebraic combinatorics;
\item the construction of the \defn{permutree lattice}, the \defn{permutree Hopf algebra} and \defn{permutreehedron};
\item the construction of the \defn{$s$-weak order} and \defn{$s$-permutahedron} along with the \defn{$s$-Tamari lattice} and \defn{$s$-associahedron}.
\end{itemize}

Chapter~9 is dedicated to the \defn{experimental method} in combinatorics research especially related to the {\tt SageMath} software. Chapter~10 describes the outreach efforts I have participated in and some of my approach towards mathematical knowledge and inclusion.

\lhead[\oldstylenums \thepage]{}
\rhead[Remerciements]{\oldstylenums \thepage}

\textit{Pour mon fils Octave Monnier Pons et mon ``demi-fils'' Ilane Bullet}
\section*{Remerciements}

\chapcitation{For most of history, Anonymous was a woman.}{Virginia Woolf}

Tout d'abord, je souhaite remercier ma rapportrice Sylvie Corteel et mes deux rapporteurs Nathan Reading et Vic Reiner. Merci pour le temps que vous avez consacré à mon travail et pour vos rapports qui me rendent fière de ce que j'ai accompli. J'étends mes remerciements à l'ensemble de mon jury, c'est-à-dire à Mireille-Bousquet Mélou, Florent Hivert, Lionel Pournin et Maria Ronco. C'est un honneur pour moi d'avoir un jury d'une telle qualité, rassemblant des chercheurs et chercheuses de renommées mondiales et dont j'admire le travail. Je remercie aussi François Bergeron et Frédéric Chapoton qui, après avoir été mes rapporteurs de thèse il y a 10 ans, m'ont écrit les lettres de recommandations pour l'HDR. Merci aussi à Nicole Bidoit et Christine Paulin pour m'avoir accompagnée tout au long du processus.

Depuis que j'ai soutenu ma thèse en 2013 et obtenu mon poste actuel en~2014, j'ai bénéficié d'un environnement de travail d'une grande qualité grâce aux collè\-gues, collaborateurs et collaboratrices qui m'ont entourée et soutenue. Ce sont des personnes avec qui j'ai partagé non seulement le plaisir de faire de la science mais aussi des valeurs telles que la science ouverte et inclusive, l'intégrité scientifique, la démocratisation des connaissances, etc. Je remercie ainsi chaleureusement Vincent Pilaud, mon coauteur et collaborateur de longue date avec qui j'ai aussi co-dirigé la thèse de Daniel Tamayo, Nicolas Thiéry avec qui j'ai collaboré sur le projet {\tt OpenDreamKit} et les enseignements de {\tt ProgImpérative}, Florent Hivert avec qui j'ai co-dirigé la thèse de Hugo Mlodecki et qui parraine mon habilitation, et Cesar Ceballos, mon coauteur et qui mène avec moi le projet {\tt PAGCAP}. J'étends ces remerciements à l'ensemble de l'équipe GALaC du LISN ainsi qu'à l'équipe Combi du LIX que j'ai pu croisée en particulier lors de notre séminaire commun et enfin aux membres du projet {\tt PAGCAP}. J'en profite aussi pour remercier mon tout premier coauteur, Grégory Châtel, qui sera content de voir que nos intervalles-posets vivent toujours ! Enfin, je remercie toutes les personnes avec qui j'ai eu des échanges professionnels et scientifiques au cours de ces 10 ans et en particulier les équipes de l'IRIF et de Marne-la-Vallée.

Comme je m’apprête à être habilitée à diriger des thèses, je voudrais remercier les 3 (ex-)doctorants qui ont déjà choisi de me faire confiance : Hugo Mlodecki qui a soutenu en 2022, Daniel Tamayo qui soutient sa thèse la même semaine que mon habilitation et Loïc Le Mogne qui entre dans sa deuxième année de doctorat. Ça a été (et c'est toujours) un grand plaisir de travailler avec vous et je vous souhaite beaucoup de réussite dans vos projets futurs. 

Au cours de ces années, j'ai aussi croisé de très nombreux étudiants et étudi\-an\-tes qui ont fait vivre nos équipes de recherche par leurs idées, leur vivacité, leurs exposés, leurs initiatives et leurs visions. Au risque d'en oublier, je pense à Jean-Baptiste Priez, Aladin Virmaux, Thibault Manneville, Joël Gay, Mathias Lepoutre, Justine Falque, Julian Ritter, Camille Combe, Doriann Albertin, Noémie Cartier, Balthazar Charles, Monica Garcia, Germain Poullot, \'Eva Philippe, Chiara Mantovani, Clément Chenevière. Aram Dermenjian a aussi été un de ces étudiants et c'est ainsi qu'il est devenu un de mes amis les plus chers, merci à lui pour sa présence.

Lors de mon année à Montréal, j'ai été accueillie très chaleureusement par les équipes du LACIM et de L'IRL ce qui m'a permis de passer une année très riche aussi bien sur le plan personnel que scientifique. Je remercie plus particulièrement François Bergeron, Christophe Hohlweg, Olivier Lafitte, Sakina Benhima et Claire Guerrier, ainsi que Srecko Belek pour son chalet sous la neige, toute l'équipe du LACIM pour son accueil et tous les membres de l'IRL et du CRM qui m'ont aidée à naviguer à travers l'administration et l'immigration.

Je souhaite aussi rendre hommage à trois personnes qui ont, chacun à leur façon, contribué à ma construction en tant qu'enseignante-chercheuse et qui nous manquent : Alain Lascoux, Yanis Manoussakis et Marc Zipstein. 

Enfin je souhaite remercier mes proches, ma famille qui est un soutien sans faille : ma mère Sylvie, Rébecca, Alice, et tous les enfants. Merci à Elyah pour son intelligence et sa vivacité, à Maël de nous faire rire avec ses premiers pas, à Ilane pour sa sensibilité et ses pensées si riches, à Octave pour sa tendresse, ses câlins et ses interprétations musicales. Et enfin, merci à mon compagnon Sébastien avec qui je partage toute cette vie. Merci d'être là et de m'épauler, de me comprendre, de comprendre ma vocation et mes engagements. \textit{Ca fait 20 ans et des poussières, qu'on fait face au vent d'hiver, ensemble on a peur de rien.}

\cleardoublepage

\lhead[\oldstylenums \thepage]{Contents}
\rhead[Contents]{\oldstylenums \thepage}
\tableofcontents
\cleardoublepage

\lhead[\oldstylenums \thepage]{List of Figures}
\rhead[Figures]{\oldstylenums \thepage}
\listoffigures
\cleardoublepage

\mainmatter

\lhead[\oldstylenums \thepage]{Introduction}
\rhead[Introduction]{\oldstylenums \thepage}
\chapter*{Introduction}
\phantomsection
\addcontentsline{toc}{chapter}{Introduction}

\chapcitation{It is impossible to be a mathematician without being a poet in soul.}{Sofia Kovalevskaya, Sónya Kovalévsky: Her Recollections of Childhood}

Ten years ago, I defended my PhD thesis on ``orders on permutations''. I studied especially the relation between the \defn{weak order} and the \defn{Tamari lattice}. My approach was mostly combinatorial. The weak order is a partial order on permutations known as the ``bubble sort order''. The minimal element is the identity permutation $1 \dots n$ while the maximal element is $n (n-1) \dots 1$. Each cover relation switches two values in consecutive positions. In particular, any descending path from the maximal element to the identity defines a \defn{sorting network} corresponding to a specific implementation of the bubble sort algorithm (see Section~\ref{sec:perm-order} for more details).

The Tamari lattice~\citeext{Tam62} is a partial order that can be defined especially on binary trees or on any objects counted by the \defn{Catalan numbers}. The cover relation is the well known \defn{binary tree rotation} used in particular in sorting algorithms to balance binary trees~\citeext{AVL62}. As I explain in Chapter~\ref{chap:asso}, it is a \defn{quotient lattice} of the weak order. The construction uses the \defn{binary search tree insertion}: each binary tree is associated with an \defn{interval} of permutations which corresponds to a \defn{lattice congruence class}. 

The combinatorial relation between the weak order and the Tamari lattice extends to geometry and algebra. The Hasse diagram of the weak order is the skeleton of a polytope, the \defn{permutahedron}, which can be constructed as the convex hull of the permutations of size $n$ seen as points in $\NN^n$. Similarly, the Tamari lattice corresponds to the \defn{associahedron}. This polytope was described originally by Tamari himself in his thesis~\citeext{Tam51} and later independently by Stasheff~\citeext{Sta63}. Nowadays, many realizations are known~\citeext{Lod04, HLT11} and it has become a major object in discrete geometry and combinatorics~\citeext{PSZ23}.

In Section~\ref{sec:asso-faces}, I explain how the combinatorial construction of the congruence classes also translate into a geometrical one by ``removing'' faces of permutahedron to obtain the associahedron. Besides, in Sections~\ref{sec:perm-hopf} and~\ref{sec:asso-hopf}, I show the connection between the combinatorial aspects and some \defn{Hopf algebraic} structures on permutations and binary trees~\citeext{MR95, DHT02, LR98, HNT05}.

This combinatorial, algorithmic, geometric, and algebraic correspondence has motivated my work for the past years. I have formed collaborations to deepen my understanding especially of the geometrical aspects and get a broader view of the questions at stake. Indeed, even if we can \emph{explain} the correspondence in many different cases, it still amazes us and leaves many open questions. For example, partial orders structures appear when looking at the skeletons of polytopes. In \emph{some} cases, this partial orders are lattices. Then, they often also correspond to algebraic structures such as Hopf algebras. But we do not yet understand \emph{what} the lattice property brings to the geometric or the algebraic aspects. I do not think these questions have easy answers. They are part of the general pursuit of mathematicians and computer scientists to create connections between distant topics, a pursuit where I believe combinatorics has an essential role to play.

In my work, I have tried to create some of these connections using a combinatorial and algorithmic approach. I have introduced new objects and structures to allow for new points of view: the Tamari interval-posets~\citeme{CP15, Pon19}, the lattice and Hopf algebra of integer posets~\citeme{CPP19, PP20}, the permutrees~\citeme{PP18, PPTJ23}, the $s$-weak order~\citemeconf{CP19} \citemepre{CP22,CP23}. In this document, I present an overview of this work and I attempt to explain the context and mathematical journey that led to each of these constructions and results. I have tried to avoid technical details
to keep a broader view and show the ideas behind the objects, especially through examples. Each chapter ends with some open questions and perspectives.
Besides, some chapters are completed by \Sage{} worksheets with computed examples~\citemesoft{PonSage23}.

\section*{Part~\ref{part:perm-ass}: Permutahedra and Associahedra}

This part is here is to introduce the two main characters of this story and give the notions needed to understand my work. I have tried to make it easy to follow even for non-specialist.

\subsection*{Chapter~\ref{chap:perm}: Permutahedra} I present the permutahedron and the weak order and especially the different points of view (geometrical, combinatorial, and more) through which I have studied them.

\subsection*{Chapter~\ref{chap:asso}: Associahedra} I present the associahedron and the Tamari lattice in their aspects more relevant to my work: mostly through their relations with the permutahedron and the weak order.

\section*{Part~\ref{part:bij}: Bijections on Tamari Intervals and More}

This part presents the work I have conducted in relations to intervals of the Tamari lattice in a \emph{bijective combinatorial} approach.

\subsection*{Chapter~\ref{chap:tamari-intervals}: Intervals of the Tamari Lattice} I give a general overview of the state of knowledge in enumerative and bijective combinatorics on Tamari intervals, especially by listing the different families I know to be in bijection with intervals of the Tamari lattice. I present the functional equations satisfied by those families giving a general explanation on its combinatorial meaning. I then define three families of combinatorial objects, in bijection with Tamari intervals, that were introduced in my work: Tamari interval-posets~\citeme{CP15}, closed flows~\citemeconf{CCP14}, and grafting trees~\citeme{Pon19}. The main result is given in Theorem~\ref{thm:interval-posets}: Tamari interval-posets are in bijection with intervals of the Tamari lattice.

\subsection*{Chapter~\ref{chap:tam-stat}: Statistics and Bijections on Tamari Intervals} I present the results of~\citeme{Pon19} where I define the \defn{rise-contact} involution on Tamari intervals, answering a conjecture of~\citeext{PR12} and~\citeext{BMFPR11}. I explain how to read certain interesting statistics on Tamari intervals and show how natural involutions of Catalan objects extend in this context. The main results are given in Theorems~\ref{thm:rise-contact} and~\ref{thm:m-rise-contact} which express the properties of the rise-contact involution and its generalization to $m$-Tamari intervals. 

\subsection*{Chapter~\ref{chap:qt}: $q,t$-Catalan} I present the results of~\citemepre{Pon22} and~\citemeconf{LMP23} in relations to $q,t$-Catalan numbers. I start by giving the general background on this topic, in relation with representation theory and symmetric functions, especially its link to the enumeration of intervals of the Tamari lattice. I explain the interpretation of the $\zeta$-function I obtained in~\citemepre{Pon22}. The main results of this chapter have been obtained in collaboration with my student Loïc Le Mogne on the $q,t$ enumeration of triangular Dyck paths in~\citemeconf{LMP23}. We have especially Theorem~\ref{thm:qt-lattice-2triangular} and Conjecture~\ref{conj:top-down-sim-sym} which express the relations between triangular $q,t$ enumerations and certain intervals of the $\nu$-Tamari lattices.

\section*{Part~\ref{part:quotient}: Quotients and Sublattices}

This part explores a certain framework to understand the weak order and Tamari lattice as part of bigger structures, in particular using the technologies of \defn{subposets}, \defn{sublattices} and \defn{lattice quotients}.

\subsection*{Chapter~\ref{chap:integer-posets}: Integer Posets} I present the work of~\citeme{CPP19} and~\citeme{PP20} where we introduce a new family of objects called the \defn{integer posets}. I explain how many combinatorial objects can be understood as integer posets: permutations, binary trees, binary sequences, intervals of the weak order, faces of the permutahedron and associahedron and more. Our main result is then given in Theorem~\ref{thm:integer-poset-lattice}: we can define a \defn{weak order} on integer poset and it is a lattice. The classical weak order as well as the Tamari lattice are sublattices of this lattice. We also find many lattice structures like a lattice of weak order intervals or the facial order as subposets. Besides, integer posets can also be endowed with a Hopf algebraic structure and we construct new Hopf algebras using this technology as in Theorem~\ref{thm:integer-poset-hopf}.

\subsection*{Chapter~\ref{chap:permutrees}: Permutrees} I present the results of~\citeme{PP18} and~\citeme{PPTJ23}. We study certain lattice quotients of the weak order which we call the permutree lattices. The permutrees interpolate between permutations, binary trees and binary sequences. We give their combinatorial description and many properties. In particular, we build polytopes called the permutreehedra in Theorem~\ref{thm:permutreehedron} and Hopf algebras in Theorem~\ref{thm:permutree-hopf}. Theorem~\ref{thm:permutree-automaton} is a first step towards the characterization of permutrees in all Coxeter groups using an automaton on reduced words of the symmetric group.

\section*{Part~\ref{part:generalizations}: Generalizations}

\subsection*{Chapter~\ref{chap:s-weak}: The $s$-Weak Order and $s$-Permutahedra} I present results which were initially published in~\citemeconf{CP19} and are being fully written in~\citemepre{CP22} and~\citemepre{CP23}. We define the $s$-weak order, which is a generalization of the weak order on certain decreasing trees. This is motivated by the definition of the $\nu$-Tamari lattices. Our main results are expressed in Theorem~\ref{thm:s-weak-lattice}: the $s$-weak order is a lattice, and Theorem~\ref{thm:s-tam}: the $\nu$-Tamari lattice is a sublattice of the $s$-weak order. Besides, we obtain a certain complex, the $s$-permutahedron, which we believe to be a polyhedral complex. Conjecture~\ref{conj:s-permutahedron} expresses that a realization as a polyhedral subdivision of the permutahedron should exists in all dimensions. 

\section*{Part~\ref{part:exp}: Experimental Approach}

This part explores some attendant aspects of my research which I consider to be essential, especially the experimental approach of my research methodology and the general philosophy in which I conduct my work.

\subsection*{Chapter~\ref{chap:exp}: Epistemology of the Experimental Approach} I give a general overview of my research methodology using computer exploration and experimentation. This is nourished especially by the collaboration I have started with the mathematician historian and philosopher Emmylou Haffner. I explain what code I have written, why, and how I make it available for further research.

\subsection*{Chapter~\ref{chap:outreach}: Outreach} I give an overview of the outreach activities I have conducted these past ten years in particular regarding experimental mathematics and support of free software and I explain what this commitment has meant for me as a scientist.

\cleardoublepage

\lhead[\oldstylenums \thepage]{\S~\thesection \; --- \; \rightmark}
\rhead[Chapter~\thechapter \; --- \; \leftmark]{\oldstylenums \thepage}

\part{Permutahedra and Associahedra}
\label{part:perm-ass}
\chapter{Permutahedra}
\label{chap:perm}

\chapcitation{It is a narrow mind which cannot look at a subject from various points of view.}{George Eliot, Middlemarch }

\begin{figure}[ht]
\center
\begin{tabular}{cc}
\includegraphics[scale=.2]{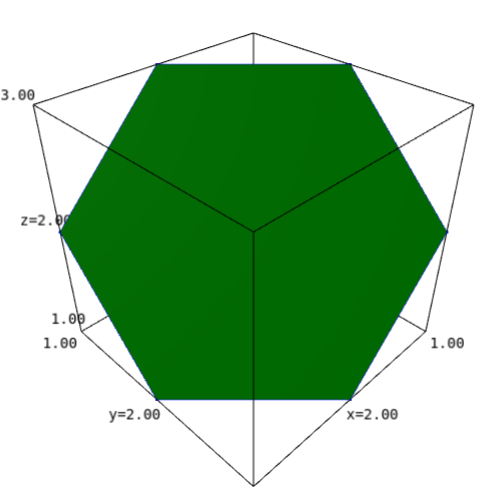}
&
\includegraphics[scale=.2]{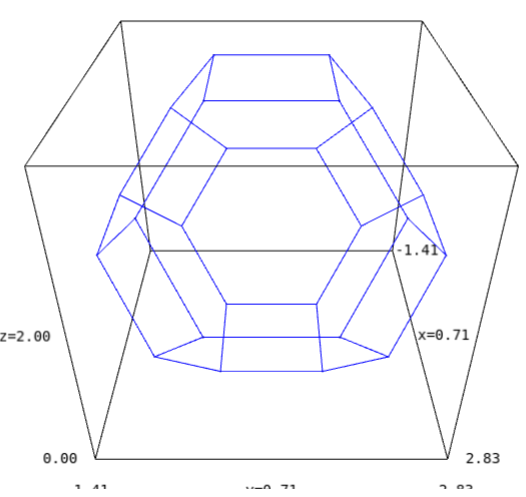}
\end{tabular}
\caption{The permutahedra of dimensions $2$ and $3$ drawn by \Sage}
\label{fig:perm_sage}
\end{figure}

Geometrically, the permutations of $1,\dots,n$ form a set of $n!$ points in $\NN^n$. This set actually lies in a hyperplane of dimension $n-1$ as the sum of coordinates is always $\sum_{1}^n i = \frac{n(n+1)}{2}$. Moreover, the points are in \emph{convex position} and the convex hull is a polytope called the \defn{permutahedron}. As an example, the images of Figure~\ref{fig:perm_sage} can be obtained with {\Sage} using the following code.

\begin{lstlisting}
P3 = Polyhedron(list(Permutations(3)))
P3.plot()
\end{lstlisting}

\begin{lstlisting}
P4 = Polyhedron(list(Permutations(4)))
P4.plot()
\end{lstlisting}

The permutahedron is a central object in my work. It has various properties connecting combinatorics, geometry and algebra. In this chapter, we recall some of these properties, highlighting the different point of views that will arise in the rest of the manuscript.

\section{As a reflection group}
\label{sec:perm-refl}

\begin{figure}[ht]
\center
\includegraphics[scale=.5]{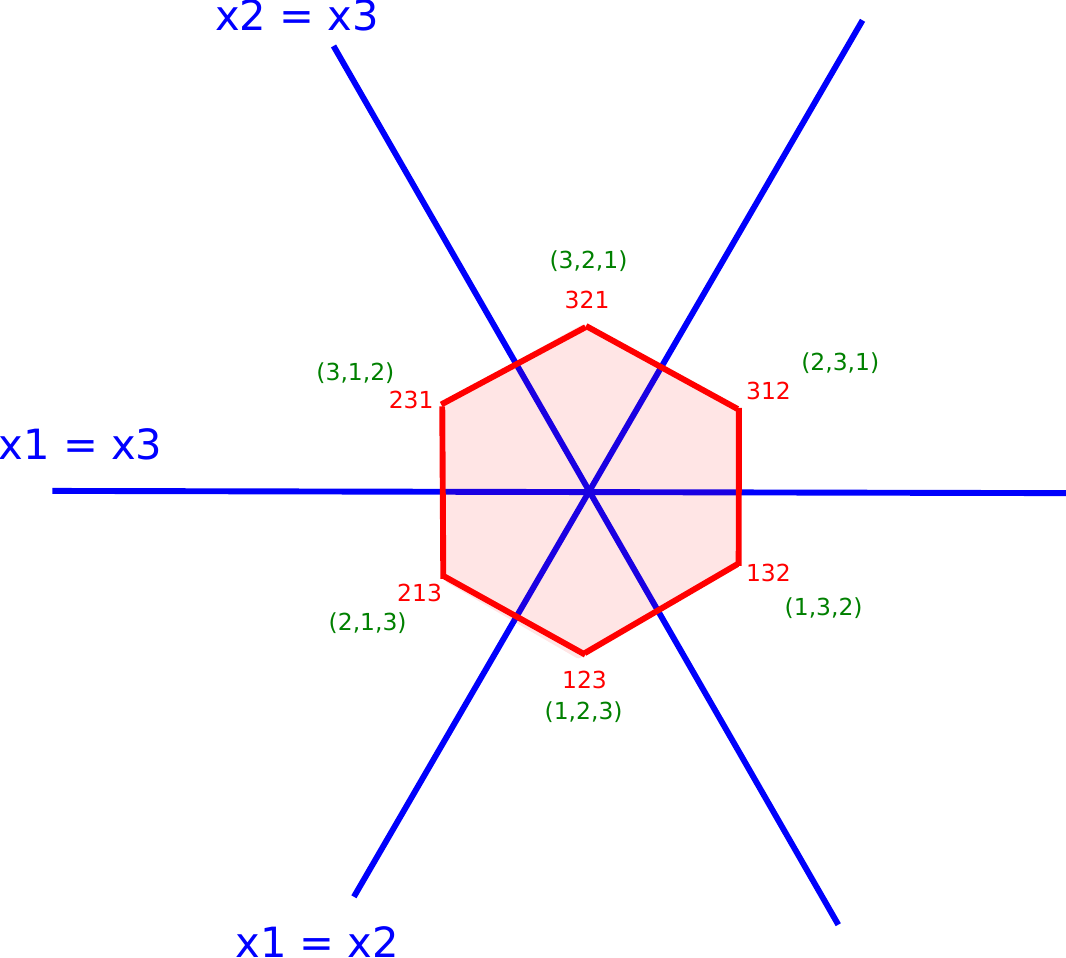}
\caption{The permutahedron as a reflection group.}
\label{fig:refgroup}
\end{figure}

The group of permutations of size $n$ is called the \defn{Symmetric group}. It is a \defn{reflection group}: it is generated by a set of of geometric reflections. The permutahedron is actually the convex hull of the orbit of the identity through the reflections. The reflections are defined by the hyperplanes $x_i = x_j$ for all $i \neq j$. Figure~\ref{fig:refgroup} illustrates the case of dimension $2$ ($n = 3$). We can see that the $3$ reflections corresponding to the $3$ hyperplanes $x_1 = x_2$, $x_1 = x_3$ and $x_2 = x_3$. In green, we write the coordinates of the vertices of the permutahedron. In red, we write the inverse permutation of the green permutation. We will see that it is useful to keep both representations in mind. Each reflection $x_i = x_j$ corresponds to a transposition $(i,j)$ applied either on the positions of the green permutation, or on the values of the red permutation.  

\section{As a partial order and sorting networks}
\label{sec:perm-order}

The \defn{skeleton} of a polytope is the graph structure given by the vertices and edges of the polytope. By choosing a starting point (the identity) and orienting the edges away from the starting point, one defines a partial order, or \defn{poset}, on the vertices. In the case of the permutahedron, the skeleton of the polytope actually corresponds to the \defn{Hasse diagram} of the partial order, \emph{i.e.}, the graph of the \defn{cover relations} of the poset. Recall that a cover relation in a poset is a relation $x \woless z$ such that there is no $y$ with $x \woless y \woless z$. The cover relations are the minimal set of relations which generate the poset through transitivity.

As there are two usual ways of labeling the vertices of the permutahedron (either using the coordinates or the inverse permutation), this defines two isomorphic partial orders on permutations known as the \defn{left and right weak order}. The cover relations correspond to \defn{simple transpositions} $s_i = (i,i+1)$. On the left weak order, these transpositions are applied \emph{to the left} of the permutation \emph{i.e.}, exchanging the values $i$ and $i+1$. On the right weak order, the transpositions are applied \emph{to the right} exchanging the values at positions $i$ and $i+1$. Figure~\ref{fig:weak3} shows the left and right weak order on $\Sym{3}$ with $s_1$ in red and $s_2$ in blue. Figure~\ref{fig:weak4} shows the left and right weak order on $\Sym{4}$ with $s_1$ in blue, $s_2$ in red, and $s_3$ in green.

\begin{figure}[ht]
\center
\begin{tabular}{cc}
\input{includes/figures/perm_gauche3_colors} &
\input{includes/figures/perm_droit3_colors} 
\end{tabular}
\caption{The left and right weak orders on $\Sym{3}$.}
\label{fig:weak3}
\end{figure}

\begin{figure}[ht]
\center
\begin{tabular}{cc}
\input{includes/figures/perm_gauche4_colors} &
\input{includes/figures/perm_droit4_colors} 
\end{tabular}
\caption{The left and right weak orders on $\Sym{4}$.}
\label{fig:weak4}
\end{figure}

\newpage

Each cover relation adds an \defn{inversion} to the permutation, \emph{i.e.}, two values which appear in reverse order. The weak order is \emph{graded} by the number of inversions, with the maximal number of inversions being $\frac{n(n-1)}{2}$ and attained in the maximal permutation $n~n-1 \dots 1$. Moreover, we can give an alternative definition of the weak order using inclusion of inversions. For two permutations $\sigma$ and $\tau$, we define that $\sigma \wole \tau$ if and only if $\inv(\sigma) \subseteq \inv(\tau)$. The definition of $\inv$ depends on the version of the weak order we consider. In the left weak order, we use \defn{position inversions} $\inv_{pos}(\sigma) := \lbrace (i,j) ; i < j \text{ and } \sigma(i) > \sigma(j) \rbrace$. In the right weak order, which is the one we will mostly use in this manuscript, we consider \defn{value inversions} $\inv_{val}(\sigma) := \lbrace (a,b); a < b; \sigma^{-1}(b) < \sigma^{-1}(a) \rbrace$ (in other words, $\sigma(b)$ appears before $\sigma(a)$ in $\sigma$). Note that value inversions are sometimes called \emph{coinversions} in the literature. We find the names position and value inversions less ambiguous especially because we work mostly with the right weak order and the value inversions.

For example, the position inversions of $\sigma = 2314$ are $(1,3)$ and $(2,3)$. It is smaller than $3412$ in the left weak order as the position inversions of $3412$ are $(1,3)$, $(1,4)$, $(2,3)$, and $(2,4)$. On the other hand, the value inversions of $2314$ are $(1,2)$ and $(1,3)$. It is smaller in the right weak order than $2431$ whose value inversions are $(1,2)$, $(1,3)$, $(1,4)$ and $(3,4)$. It is an easy exercise to show that the cover relations of the order defined by inclusion of the inversions are indeed the simple transpositions. 

This also allows for an interpretation of the weak order as a \defn{sorting network}. A sorting network is a series of comparisons and exchange to apply on a list of entries to sort the list. Each descending path from the maximal permutation to the identity ``sorts'' the permutation so each of these paths corresponds to a sorting network. They are actually sorting networks corresponding to different possible implementations of the well known (inefficient) \defn{bubble sort} algorithm. The complexity $O(n^2)$ of the algorithm corresponds to the height of the weak order, \emph{i.e.}, the maximal number of inversions $\frac{n(n-1)}{2}$.

\section{As a lattice}
\label{sec:perm-lattice}

\begin{figure}[ht]
\center
\input{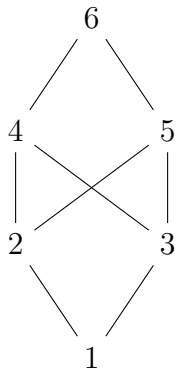}
\caption{An example of a non-lattice}
\label{fig:non-lattice}
\end{figure}

The weak order admits an extra algebraic structure: it is a \defn{lattice}. A lattice can be defined as a partial order which admits two binary operations:
\begin{itemize}
\item the \defn{meet}, written $\meet$, which computes the \emph{greatest lower bound} of two elements $x$ and $y$;
\item the \defn{join}, written $\join$, which computes the \emph{lowest upper bound} of two elements $x$ and $y$.
\end{itemize}
The best way to understand this condition is to look at a counter example. Figure~\ref{fig:non-lattice} shows the Hasse diagram of a poset which is not a lattice. Look for example at elements $4$ and $5$. We take the intersection of the ideals generated by the two elements, \emph{i.e.}, the elements of the poset which are at the same time smaller than or equal to $4$ and $5$. This set is $\lbrace 1, 2, 3 \rbrace$ and does not have a maximal element.  This poset is actually isomorphic to the \defn{Bruhat order} on $\Sym{3}$ (or strong order) which is another well known partial order on permutations.

On the contrary, the meet and join of any two elements of the weak order is well defined. This was proved in~\citeext{GR70} in the context of analyzing voting systems (this is also where the name permutahedron was introduced). The first step in the proof is to characterize inversion sets. This is done in~\citeext[Theorem 2, page 26]{GR70}: a set $I$ of tuples $(i,j)$ with $1 \leq i < j \leq n$ is the inversion set of a permutation if and only if $I$ is \defn{transitive} ($(i,j) \in I$ and $(j,k) \in I$ with $i < j < k$ implies $(i,k) \in I$) and the \emph{complementary} of $I$ is transitive. The complementary of $I$ is the set of tuples $(i,j)$ with $1 \leq i < j \leq n$ which are not in $I$, they correspond to the \emph{non-inversions} of the permutation. Inversion sets can then be used to explicitly compute the join of two permutations $\sigma$ and $\mu$. Let $I_\sigma$ and $I_\mu$ be the respective inversion sets, then $I_\sigma \cup I_\mu$ is not necessary transitive and so it is not an inversion set. Nevertheless, we can prove that the \emph{transitive closure} of $I_\sigma \cup I_\mu$  is an inversion set and it gives the permutation $\sigma \join \mu$.  We have not specified here if we were working with values or positions for the inversion sets: it is not necessary as the characterization and proof are the same.

As an example, in the right weak order the meet of $3214$ and $3142$ is $3124$ and the join is $3421$. Looking at values inversions, we have $\inv_{val}(3214) = \lbrace (1,2), (1,3), (2,3) \rbrace$ and $\inv_{val}(3142) = \lbrace (1,3), (2,3), (2,4) \rbrace$. Their union $\lbrace (1,2), \allowbreak (1,3), (2,3), (2,4) \rbrace$ is not a inversion set because it is not transitive. The transitive closure adds the inversion $(1,4)$ and is the inversion set of $3421$. The left and right weak orders are implemented in {\Sage} and we show examples in~\citemesoft{PonSage23}.

\section{As a Coxeter group}
\label{sec:perm-coxeter}

Reflection groups are also part of a larger family of groups, the \defn{Coxeter groups}. A Coxeter group is a group which admits a certain presentation through a finite set of generators $s_1, \dots, s_n$ with relations $(s_i s_j)^{m_{i,j}} = 1$ where $m_{i,i} = 1$ and $2 \leq m_{i,j} \leq \infty$. Finite Coxeter Groups have been classified. The symmetric group $\Sym{n}$ corresponds to the type $A_{n-1}$. Its generators are the simple transpositions and satisfy the \defn{braid relations}

\begin{align}
s_i^2 &= 1; & \\
(s_i s_j)^2 &= 1 & \text{if } i + 1 < j; \\
(s_i s_{i+1})^3 &= 1. &
\end{align}  

Each element of the group corresponds to an infinite set of words on the generators $s_i$. A word is said to be \defn{reduced} if it is of minimal length. For example, in $\Sym{3}$, $s_1$, $s_1^3$, and $s_1 s_2 s_1 s_1 s_2$ are three words corresponding to the permutation $213$ with respective lengths $1$, $3$, and $5$, but only $s_1$ is reduced. We say that $213$ is of length $1$. The  weak order can be seen as an oriented version of the Cayley graph of the Symmetric group using the Coxeter presentation. A path in the permutahedron skeleton graph corresponds to a word of generators in the group. If the path is following the weak order orientation, \emph{i.e.}, if it is only going up in the poset, it is a reduced word. In particular, for each element, the set of reduced words corresponds to the set of oriented paths from the identity to the element. Depending on the version of the weak order we are working with (left of right), the correspondence between paths and words changes slightly. In the right weak order, generators are added to the right of the word while going along the path, while they are added to the left in the left weak order. 

For example, in Figure~\ref{fig:weak3}, we can see that on the left weak order, the permutation $312$ is obtained by a path starting at $123$ and going through a red edge first ($s_1$) and then a blue edge ($s_2$). As we are on the left weak order, the reduced word corresponding to $312$ is $s_2 s_1$. On the right weak order, the path $s_1$ then $s_2$ goes to $231$. Indeed, the reduced word of $231$ is $s_1 s_2$,  it is the inverse of $312 = s_2 s_1$. Similarly, we can read that, in both the right and left weak orders, the maximal permutation $321$ has $2$ reduced words: $s_1 s_2 s_1$ and $s_2 s_1 s_2$.

\begin{figure}[ht]
\center
\includegraphics[scale=.2]{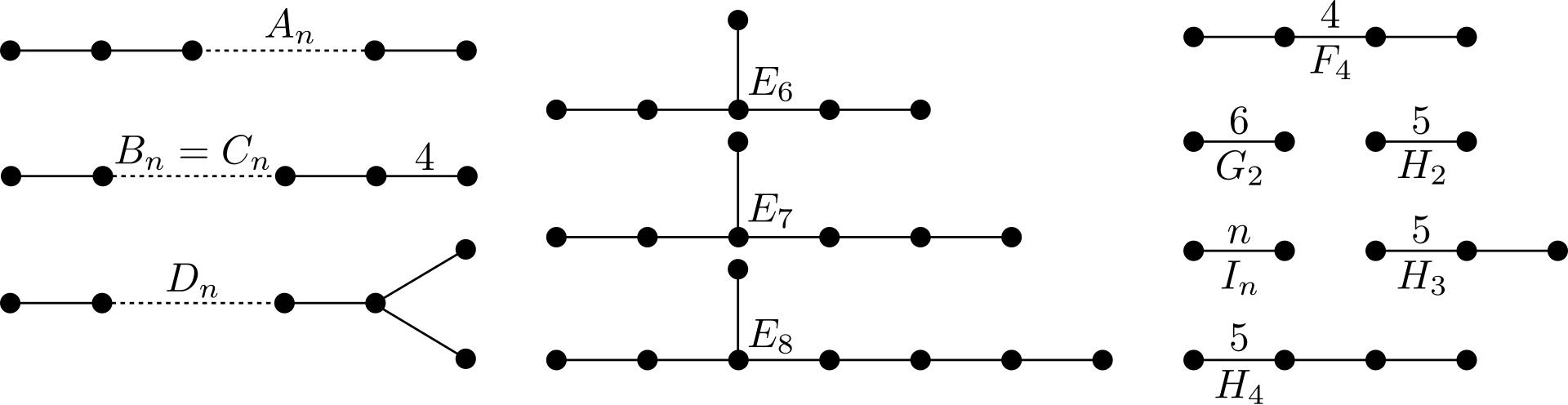}
\caption{The possible connected Dynkin diagrams of finite Coxeter groups}
\small\textsuperscript{ Image under \href{https://creativecommons.org/licenses/by-sa/3.0}{CC BY-SA 3.0} from \href{https://commons.wikimedia.org/wiki/File:Finite_coxeter.svg}{Wikimedia Commons}}
\label{fig:coxeter-dynkin}
\end{figure}

The symmetric group is just one example of a finite Coxeter group. These have been classified~\citeext{Cox35} into 4 infinite families ($A_n$, $B_n$, $D_n$, $I_n$) and $6$ exceptional groups ($E_6$, $E_7$, $E_8$, $F_3$, $H_3$, $H_4$). Other finite Coxeter groups are products of these classified families. We represent Coxeter groups through their \defn{Dynkin diagrams} as in Figure~\ref{fig:coxeter-dynkin}. The generators $s_i$ are the vertices of the diagram. If there is no edge between two vertices $i$ and $j$, then their corresponding generators $s_i$ and $s_j$ commute ($m_{i,j} = 2)$. Otherwise, $m_{i,j}$ is equal to $3$ by default or to the edge label. If the diagram is disconnected, then the group is a product of the connected components. Figure~\ref{fig:coxeter-dynkin} represents the possible connected diagrams of finite Coxeter groups corresponding to the classified families. 

All finite Coxeter groups admit a faithful representation as a reflection group. In particular, as explained in \citeext{Hoh12}, you can construct a group permutahedron by taking the convex hull of the orbit of a generic point. Besides, you can also define a weak order as in type $A$ by interpreting the skeleton of the polytope as the Hasse diagram of a partial order. This order can be defined without any reference to geometry: we say that $\mu \wole \sigma$ in the right (resp. left) weak order if and only if there exist a reduced word of $\mu$ which is prefix (resp. suffix) of a reduced word of $\sigma$. For example, we have $2314 \wole 2431$ (in the right weak order of type $A$) and indeed $2314 = s_1 s_2$ is a prefix of $s_1 s_2 s_3 s_2 = 2431$.  This order has been proved to be a lattice by Bj\"orner~\citeext{Bjo84} for all finite Coxeter groups. 

The work I present here mostly concerns the Symmetric group, otherwise called ``type $A$ case'' in reference to the Coxeter classification. Nevertheless, to study the weak order lattice and the permutahedron, it is important to keep in mind all the aforementioned point of views, and especially this interpretation as a Coxeter group. Besides, the Coxeter approach is a strong motivation of my work. Indeed, the combinatorial work in type $A$ often lays the ground to a more general approach through Coxeter groups. Type $A$ offers many combinatorial tools that are only partially available in other types . Many proofs can be done in a simple elementary way relying on intuition and exploration. They open the door to what \emph{might} be true for other Coxeter groups and sometimes to more generic proofs.

\section{As a Hopf algebra}
\label{sec:perm-hopf}

Hopf algebras are algebraic structures that have raised the interest of the combinatorial community for the past decades. They can be seen as a generalization of the classical notion of \defn{generating functions} studied by Flajolet and Sedgewick~\citeext{FS09} in the context of algorithmic analysis. A generating function can be understood as a formal sum over an infinite combinatorial set where each object of size $n$ is represented by a single monomial $z^n$. The result is a power series where the coefficient of $z^n$ is the number of objects of size $n$. In that sense, generating functions are a way to encode recurrences of integer sequences into an algebraic object. Hopf algebras encode a larger class of structural inductions on combinatorial objects. 

Formally, a combinatorial Hopf algebra is defined as a vector space over a set of combinatorial objects. So instead of working with polynomials or power series, we directly sum the objects themselves and work with formal linear combinations of combinatorial objects. A Hopf algebra $\hopf$ is a special case of a bialgebra. It is endowed with two operations: a product $\hprod : \hopf \tensor \hopf \rightarrow \hopf$ and a coproduct $\hcop : \hopf \rightarrow \hopf \tensor \hopf$. These two operations are respectively \defn{associative} and \defn{co-associative} and admits respectively a \defn{unit} and a \defn{co-unit}. Besides, the coproduct needs to be an algebra morphism. In other words, for any two elements $x$ and $y$ in~$\hopf$, then $\hcop (x \hprod y) = \hcop(x) \hprod \hcop(y)$. We don't go into details into the axioms satisfied by bialgebras and Hopf algebras as they are not completely relevant to our work. We refer the reader to~\citeext{Hiv06} for a more comprehensive presentation.

The initial motivation for the study of Hopf algebras in combinatorics comes from the study of symmetric functions and can be traced back to MacMa-\linebreak
hon~\citeext{Mac15} and later to Rota~\citeext{JR79}. The connection with the permutahedron comes from the Malvenuto-Reutenauer Hopf algebra on permutations~\citeext{MR95} which we call $\FQSym$ after~\citeext{DHT02} who redefine the algebra using the notion of \emph{polynomial realization} with words. We write $\BF_\sigma$ the elements of $\FQSym$ using the basis $\BF$. The product is defined by the \defn{shifted shuffle product} on permutation. We show an example below:

\begin{align}
\BF_{\red{231}} \hprod \BF_{12} &= \BF_{\red{231}45} + \BF_{\red{23}4\red{1}} + \BF_{\red{2}4\red{31}5} + \BF_{\red{23}45\red{1}} + \BF_{4\red{231}5} \\
&+ \BF_{\red{2}4\red{3}5\red{1}} + \BF_{4\red{23}5\red{1}} + \BF_{\red{2}45\red{31}} + \BF_{4\red{2}5\red{31}} + \BF_{45\red{231}} \nonumber \\
&= \sum_{\red{231}45 \wole_R \mu \wole_R 45\red{231}} \BF_\mu.
\end{align}

The letters of the second permutation are \emph{shifted} then the sum runs over all permutations such that the order of the small letters (here $1$, $2$, and $3$ in red) corresponds to the first permutation and the order of the big letters (here $4$ and~$5$) corresponds to the second permutation as in the shuffle of two decks of cards. This is actually a sum over an interval of the right weak order between the element where the letters of the second permutation are in the rightmost positions (here $\red{231}45$) and the element where they are in the leftmost positions (here $45\red{231}$). So, in some sense, the weak order \emph{encodes} the structure of the Hopf algebra.

As $\FQSym$ is a Hopf algebra, the basis $\BF$ also admits a coproduct. It actually corresponds to a product in the \emph{dual basis} $\BG$. The product in $\BG$ is called the \defn{convolution product} on permutations and corresponds to a sum in the left weak order.

\begin{align}
\BG_{\red{312}} \hprod \BG_{12} &= \BG_{\red{312}45} + \BG_{\red{412}35} + \BG_{\red{413}25} + \BG_{\red{512}34} + \BG_{\red{423}15} \\
&+ \BG_{\red{513}24} + \BG_{\red{523}14} + \BG_{\red{514}23} + \BG_{\red{524}13} + \BG_{\red{534}12} \nonumber \\
&= \sum_{\red{312}45 \wole_L \mu \wole_L \red{534}12} \BG_\mu.
\end{align}

In the convolution product, the role of the values and positions are switched: here the orders of the letters placed in the first $3$ positions is given by the first permutation $312$ whereas the order of the letters in the last $2$ positions is given by the second permutation $12$. The basis $\BG$ is isomorphic to $\BF$ using the inversion of permutation just like between the right and left weak orders. This makes $\FQSym$ an auto-dual Hopf algebra whose structure is given by the weak order. Besides, the weak order can be used to define \emph{multiplicative basis} hence proving that $\FQSym$ is free. 

\section{Faces of the permutahedron}
\label{sec:perm-faces}

\begin{figure}[ht]
\center
\begin{tabular}{cc}
\includegraphics[scale=.3]{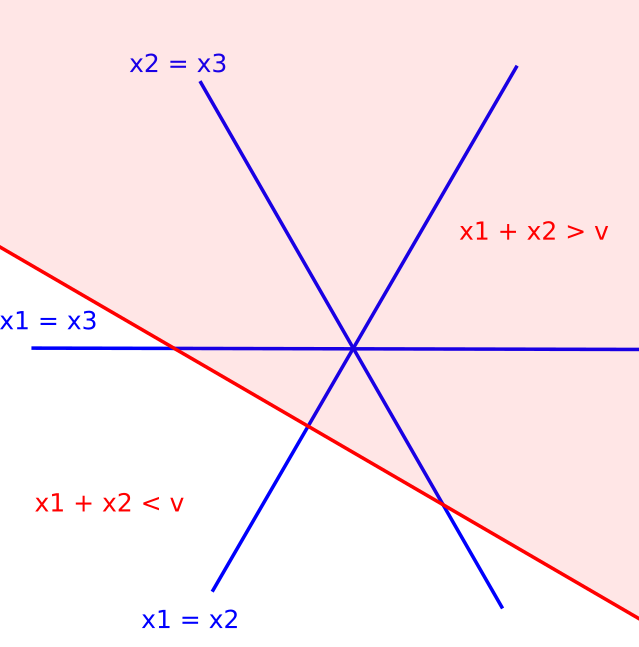} &
\includegraphics[scale=.3]{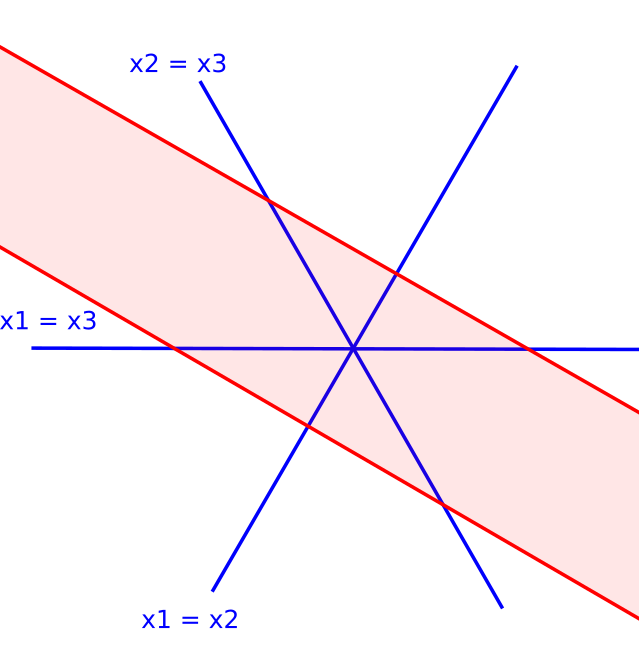}
\end{tabular}

\caption{Inequalities in the plane defined by $x_1 + x_2 + x_3 = 6$.}
\label{fig:permieq}
\end{figure}

The permutahedron can be defined by its vertices, the permutations, as we have seen. As a polytope, it can also be defined through a list of \emph{inequalities}. For example, on the left of Figure~\ref{fig:permieq}, we cut the plane defined by $x_1 + x_2 + x_3 = 6$ (where the permutahedron of dimension $2$ lives) into two half planes using the value of $x_1 + x_2$.  Looking again at Figure~\ref{fig:refgroup}, it is clear that permutahedron of dimension $2$ is defined by $6$ such inequalities.

In dimension $2$, this is an easy exercise to find those inequalities. For example, we always have $x_1 + x_2 = 6 - x_3$ and as $1 \leq x_3 \leq 3$, we find $ 3 \leq x_1 + x_2 \leq 5$ which corresponds to the right image of Figure~\ref{fig:permieq}. In general, for the permutahedron of dimension $n-1$ we find $2^n - 2$ inequalities. They correspond to all non-empty proper subsets of $[n] := \lbrace 1, 2, \dots, n \rbrace$. For each subset $J$, the corresponding inequality is given by $\sum_{j \in J} x_j \geq \binom{|J| + 1}{2}$~\citeext{Rad52}. 

\begin{figure}[ht]
\center
\begin{tabular}{cc}
\input{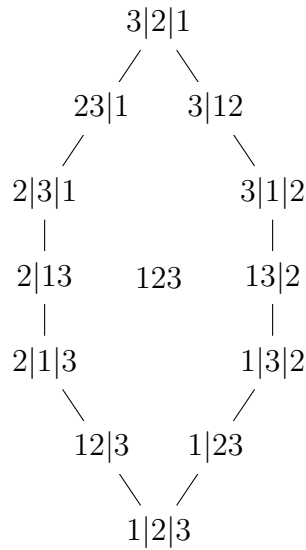} 
\end{tabular}
\caption{Ordered partitions as faces of the permutahedron.}
\label{fig:permfaces}
\end{figure}

The extremal points of the polytope are the points where some inequalities are equalities. They correspond to the \defn{faces} of the polytope and are polytopes themselves. For example, in dimension $2$, the permutahedron has $13$ faces : $1$ face of dimension $2$ (the full polytope), $6$ faces of dimension $1$ (the $6$ edges corresponding to the $6$ inequalities), and $6$ faces of dimension $0$ (the $6$ vertices, a.k.a the permutations). They are counted by the Fubini numbers~\citeext{OEISA000670} and have a nice combinatorial interpretation as \defn{ordered partitions}. The general idea is the following: when you are on the edge between permutation $123$ and permutation $213$, the order between $1$ and $2$ is not defined but they are both before $3$, this is the ordered partition $\lbrace 1, 2 \rbrace, \lbrace 3 \rbrace$ which we write $12|3$. The dimension of the face is given by $n - \ell$ where $\ell$ is the number of parts in the ordered partitions. We illustrate this in Figure~\ref{fig:permfaces}.

This combinatorial interpretation allows for extra structures that extend the ones on permutations. For example, one can define a lattice structure, the \defn{facial weak order} such that the weak order is a sublattice. It was shown in~\citeext{KLNPS01} for type $A$ and in~\citeext{DHP16} for all Coxeter groups. A Hopf algebra structure can also be defined on faces that extend the Malvenuto-Reutenauer one on permutations~\citeext{Cha00}.

\chapter{Associahedra}
\label{chap:asso}

\chapcitation{L'algèbre n’est qu’une géométrie écrite, la géométrie n'est qu'une algèbre figurée.}{Sophie Germain, Oeuvres philosophiques.}

In his thesis~\citeext{Tam51} and later on in~\citeext{Tam62}, Dov Tamari introduces a certain partial order, the ``associativity posets'', on the different ways to parenthesize a word of a given length $n$. He is studying how to replace the usual associativity of a binary operation $(ab)c = a(bc)$ by a \emph{substitution rule} setting some early foundations for what is known nowadays as \emph{term rewriting systems}. For example, there are $5$ ways to parenthesize a word of length $4$:

\begin{tabular}{ccccc}
$((ab)c)d$ & $(a(bc))d$ & $(ab)(cd)$ & $a((bc)d)$ & $a(b(cd))$. \\
\end{tabular}

\begin{figure}[ht]
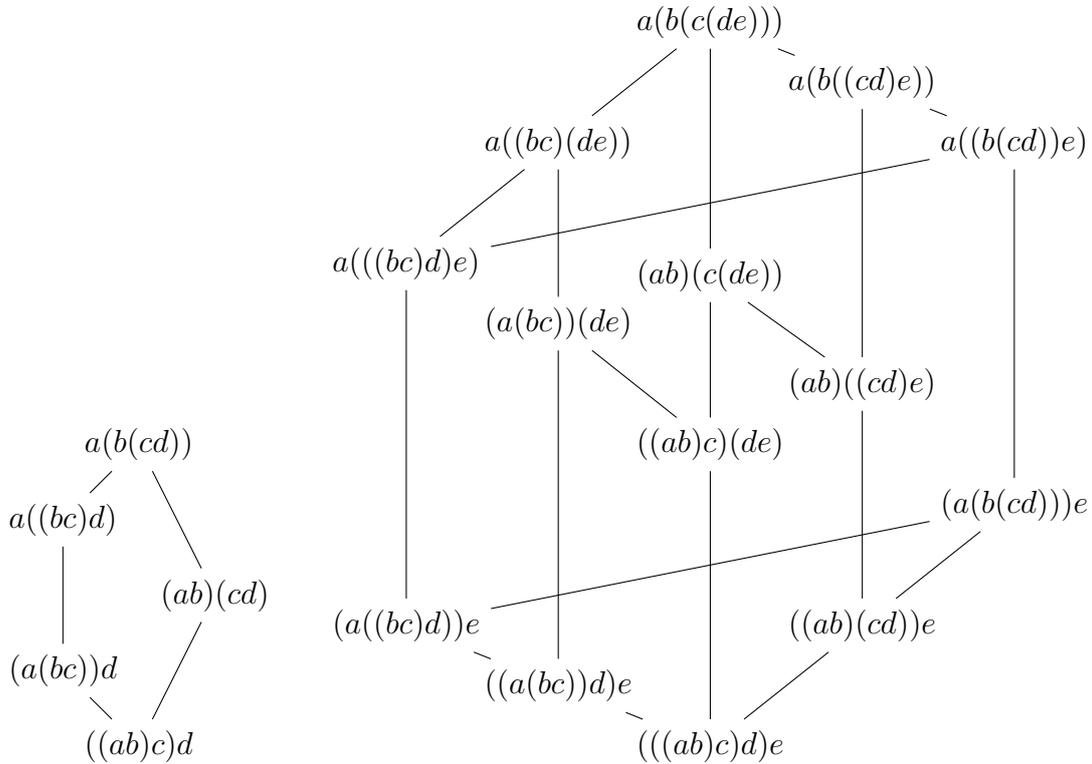

\begin{tabular}{cc}
\input{includes/figures/tamari_words3} &
\input{includes/figures/tamari_words4}
\end{tabular}
\caption{The Tamari lattice on parenthesized words}
\label{fig:tamari-words}
\end{figure}

In general, the number of configurations is counted by the famous Catalan numbers~\citeext{OEIS000108}
\begin{align}
\label{eq:catalan}
C_n = \frac{1}{n+1} \binom{n}{2n}.
\end{align}
By orienting the rewriting such that $(ab)c \woless a(bc)$ as shown in Figure~\ref{fig:tamari-words}, we obtain the partial order known as the \defn{Tamari lattice}. Indeed, it was conjectured by Tamari in his thesis to be a lattice and he later gave several proofs~\citeext{FT67, HT72} with his students Friedman and Huang.

According to Stasheff~\citeext{Sta12}, in his thesis, Tamari also describes the \emph{geometry} of the lattice and gives the first known depiction of the \defn{associahedra}. Just like the weak order, the Tamari lattice Hasse diagram corresponds to the oriented skeleton of a polyotpe. This polytope is the associahedron. It was also introduced independently by Stasheff~\citeext{Sta63} about 10 years after Tamari's thesis and is sometimes known as the \defn{Stasheff polytope}. 

The Tamari lattice and associahedron are now major objects of interest in the combinatorics community. See for example the the collection of papers gathered in~\citeext{MHPS12} to honor the memory of Tamari and the recent survey on Loday's realization of the associahedron~\citeext{PSZ23}. The relations between the Tamari lattice and the weak order has been nourishing many of the questions I considered these last 10 years. In this chapter, I present an overview of the \emph{combinatorial}, \emph{geometrical}, and \emph{algebraical} links between these two structures. This gives the ground motivation of the work I present afterwards.

\section{As a sublattice}
\label{sec:asso-sublattice}

We say that a permutation $\sigma$ of size $n$ contains a pattern $\mu$ of size $k \leq n$ if there is a subword  $w$ of $\sigma$ with $|w| = k$ such that the relative order between the letters of $w$ is given by $\mu$. In other words, the \defn{standardized} permutation of $w$ is $\mu$. For example, the permutation $4213$ contains the pattern $312$ because the subword $413$ is such that the highest value appears first ($4$), then the smallest one $(1)$, then the intermediate one $(3)$, which corresponds to the permutation $312$. In~\citeext{Knu98}, Knuth shows that the so-called \defn{stack sortable} permutations are the permutations \emph{avoiding} the pattern $231$ if we read the values from left to right and $312$ if we read the values from right to left. Both families are counted by the Catalan numbers. Actually, the set of permutations avoiding any pattern of size $3$ is always counted by the Catalan numbers (see~\citeext{Mac15} for the permutations avoiding $123$ and the introduction of~\citeext{SS85} for a more general explanation). 

Now, if we select right-stack sortable permutations in the right weak order, \emph{i.e.}, the permutations avoiding $312$, the induced partial order is actually a \defn{sublattice} of the weak order. 

\begin{figure}[ht]
\center
\input{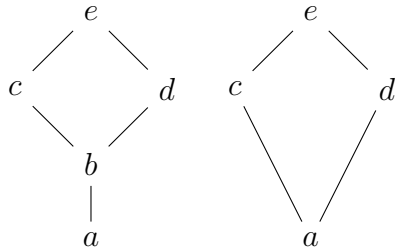}
\caption{Example of a lattice which is not a sublattice.}
\label{fig:non-sub-lattice}
\end{figure}

If $L$ is a lattice, then $I \subset L$ is a sublattice if for any two elements $x$ and $y$ in $I$, then $x \meet_L y$ and $x \join_L y$ are also in $I$. In particular, note that even if $I$ is a lattice (using the order induced by $L$), it might not be a sublattice of $L$. See for example the subposet on Figure~\ref{fig:non-sub-lattice}: it is indeed a lattice but not a sublattice of the lattice on the left because $c \meet d = b$ is not part of the subposet.

\begin{figure}[ht]
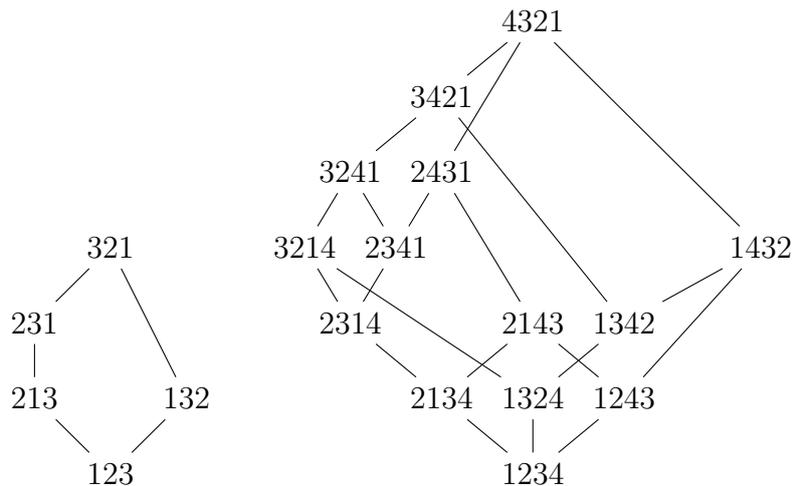

\center
\begin{tabular}{cc}
\input{includes/figures/tamari_perms3}
&
\input{includes/figures/tamari_perms4}
\end{tabular}
\caption{The Tamari lattice as a sublattice of the right weak order.}
\label{fig:tamari-perms}
\end{figure}

In the case of $312$-avoiding permutations, we can show that they form a sublattice of the right weak order using inversion sets as in Section~\ref{sec:perm-lattice}. A permutation avoids $312$ if for all $a < c$ where $(a,c)$ is a (value) inversion, then $(a,b)$ is also a value inversion for all $b$ with $a < b < c$. This property is actually stable when taking the transitive closure of the union of two inversion sets. This shows that the join of two $312$ avoiding permutations is also a $312$ avoiding permutation. A similar proof can be done for the meet. We thus obtain a lattice where the number of elements is given by the Catalan numbers: this is actually the Tamari lattice as we explain in Section~\ref{sec:asso-quotient}. See the example on Figure~\ref{fig:tamari-perms}.

\section{As a quotient lattice}
\label{sec:asso-quotient}

\begin{figure}[ht]
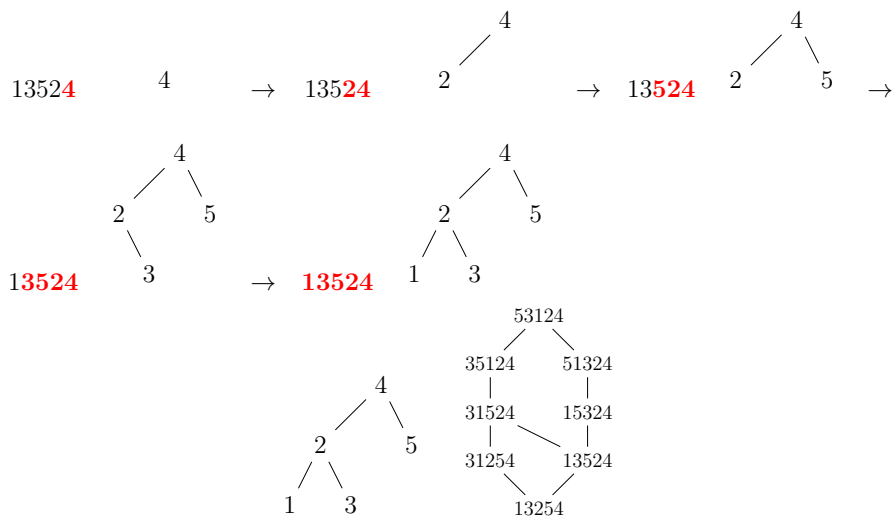

\center
\scalebox{.8}{\input{includes/figures/tree-insertion}}
\scalebox{.8}{\input{includes/figures/bst-linear-extensions}}
\caption{Binary search tree insertion and linear extensions}
\label{fig:bst-insertion}
\end{figure}

The classical algorithm of \defn{binary search tree insertion} gives another connection between the weak order and the Tamari lattice. We consider here \defn{planar rooted binary trees}, which are counted by the Catalan numbers. Recursively, a binary tree is either a leaf (an empty tree) or a non-empty binary tree with a left and a right subtree. We usually do not draw the empty children and the size is given by the number of non-empty nodes. A binary tree whose non-empty nodes are labeled with numbers is said to be a \defn{binary search tree} if for all nodes $\lambda$ labeled by $x$, the nodes in the left subtree of $\lambda$ are labeled with values smaller than or equal to $x$ while the nodes on the right subtree of $\lambda$ are labeled with values greater than $x$. This structure is used for sorting algorithms through binary search tree insertion: a number is inserted through the root and ``falls down'' following branches depending on its value and the label of the node (left if it is smaller, right if it is bigger) until it reaches a leaf, \emph{i.e.}, an empty spot. In Figure~\ref{fig:bst-insertion}, we show the step by step binary search tree insertion of the values of a permutation (from right to left). 

\begin{figure}[ht]
\center
\input{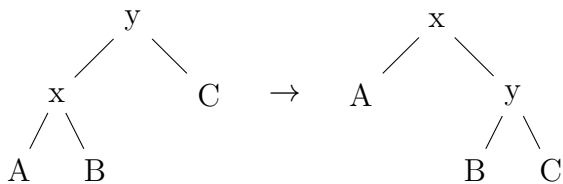}
\caption{Binary tree rotation}
\label{fig:tree-rotation}
\end{figure}

This algorithm is used in conjunction with binary tree rotation (see Figure~\ref{fig:tree-rotation}) in~\citeext{AVL62} to obtain efficient sorting algorithms. Besides, it actually defines a surjection between permutations and binary trees: the image of a permutation is the tree obtained by inserting the values from right to left. Technically, what we obtain is a labeled tree but there is only one way to label a given tree with $1, \dots, n$ such that it is a binary search tree. It consists of labeling the nodes in \defn{in-order}, we call it a \defn{standard binary search tree}. 

It was shown in\citeext{BW91} that the pre-image of a given binary tree is an interval in the right weak order which corresponds to the \defn{linear extensions} of the standard binary search tree: all possible ways to ``read'' a permutation such that a node label is read before its parent. For example, the permutation $12354$ is \emph{not} a linear extension of the tree of Figure~\ref{fig:bst-insertion} because $3$ appears after its parent $2$. In~\citeext{HNT05}, the authors call this pre-image the \defn{sylvester class} of the tree (from the french word ``sylvestre''), it is the transitive closure of the following rewriting rule on permutations:

\begin{equation}
\label{eq:sylvester}
\dots ac \dots b \ldots \equiv \dots ca \dots b \dots
\end{equation}
with $a < b < c$. 

\begin{figure}[ht]
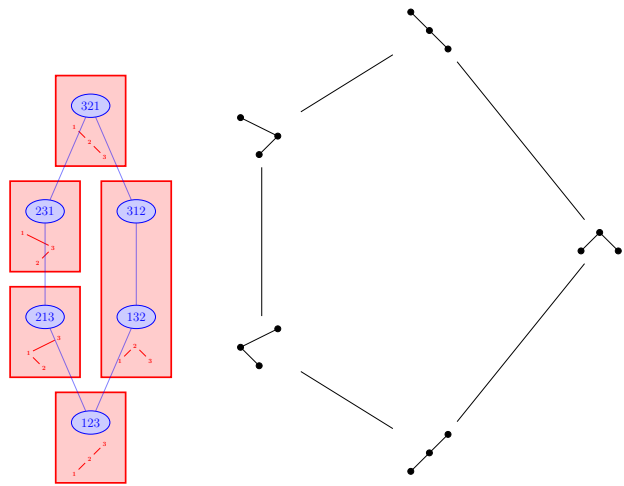

\center
\begin{tabular}{cc}
\scalebox{.4}{\input{includes/figures/tamari_quotient3}} &
\scalebox{.7}{\input{includes/figures/tamari_trees-3}}
\end{tabular}
\caption{The Tamari lattice as a quotient of the right weak order ($n = 3$)}
\label{fig:tamari-quotient3}
\end{figure}

\begin{figure}[h!]
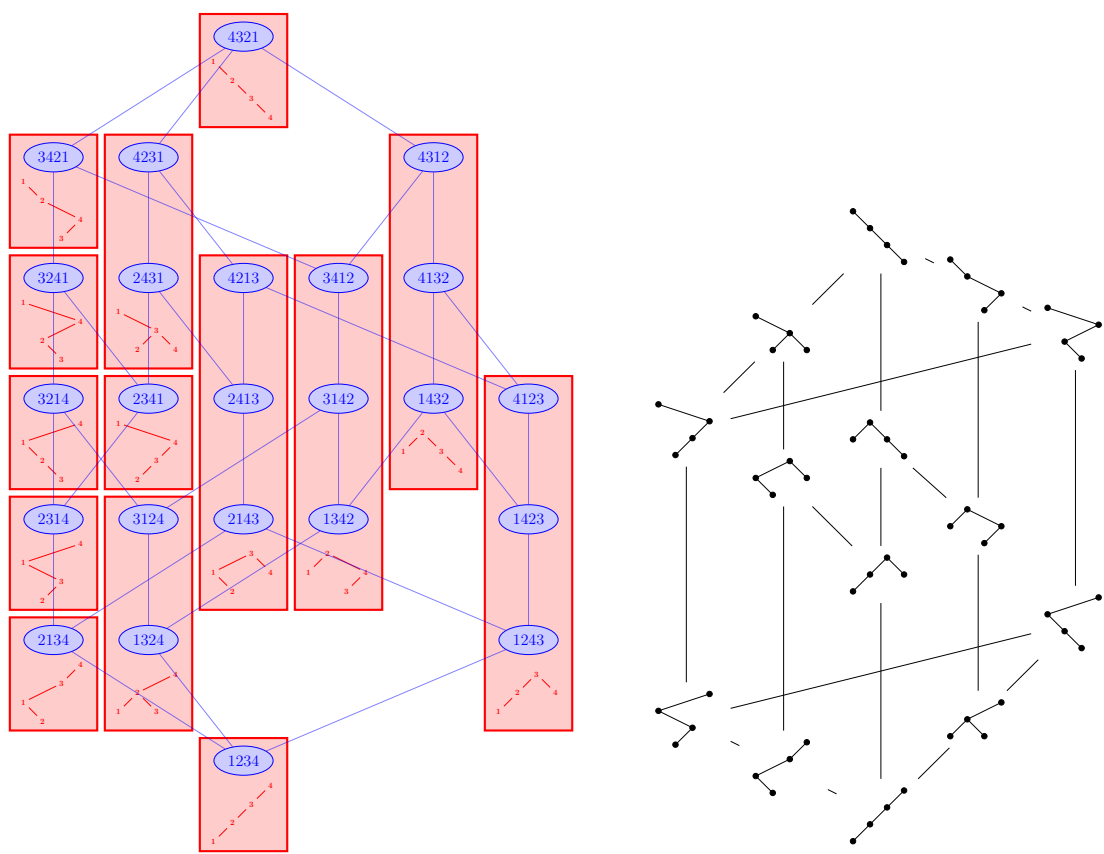

\center
\begin{tabular}{cc}
\scalebox{.5}{\input{includes/figures/tamari_quotient4}} &
\scalebox{.8}{\input{includes/figures/tamari_trees-4}}
\end{tabular}
\caption{The Tamari lattice as a quotient of the right weak order ($n = 4$)}
\label{fig:tamari-quotient4}
\end{figure}

A congruence relation on a lattice is said to be a \defn{lattice congruence} if $x \equiv x'$ and $y \equiv y'$ implies that $x \meet y \equiv x' \meet y'$ and $x \join y \equiv x' \join y'$. This induces a lattice structure on the congruence classes which is called a lattice quotient of $L$. Lattice quotients have been studied in a combinatorial context by Reading~\citeext{Rea04}. In particular, if the congruence classes are intervals and satisfy some extra order-preserving conditions, then they form a lattice quotient. The sylvester classes satisfy the conditions and the quotient lattice we obtain is actually the Tamari lattice as illustrated on Figures~\ref{fig:tamari-quotient3} and~\ref{fig:tamari-quotient4}. Classes are represented by binary trees and the cover relation is the binary tree rotation depicted on Figure~\ref{fig:tree-rotation}. This description of the lattice is actually very close to the original one from Tamari depicted on Figure~\ref{fig:tamari-words}. Indeed the binary tree is just the symbolic expression tree following the associative rule of the parenthesis. In particular, the substitution rule $(ab)c \rightarrow a(bc)$ is exactly the binary tree rotation of Figure~\ref{fig:tree-rotation}.

To connect with the previous section, note that the $312$-avoiding permutations are the minimal elements of the sylvester classes (the maximal elements are the $132$-avoiding permtuations). The construction by quotient means in particular that the subposet of $312$ avoiding permutations is indeed the Tamari lattice. Note that it does not imply that it forms a sublattice. To conclude, a similar construction can be made on the left weak order using the bijection between permutations and \defn{decreasing binary trees} (labeled binary trees where labels are strictly decreasing from root to leaves) as a surjective map between permutations and binary trees defining a lattice congruence.

\section{Other lattice descriptions}
\label{sec:asso-other}

\begin{figure}[ht]
\input{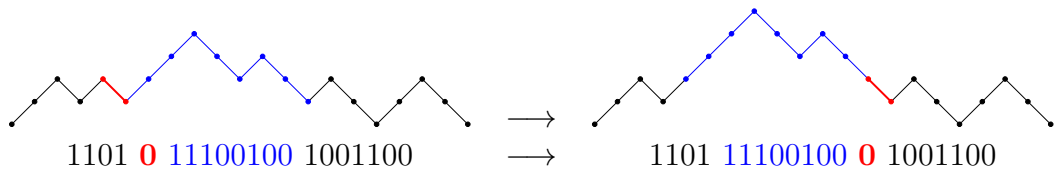}
\caption{Tamari rotation on Dyck paths}
\label{fig:rotation_dyck}
\end{figure}

\begin{figure}[ht]
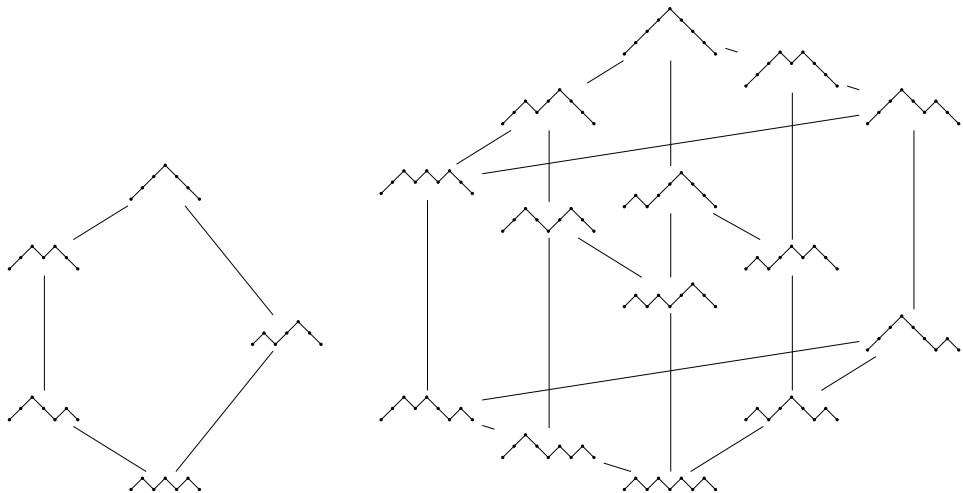

\center
\begin{tabular}{cc}
\scalebox{.5}{\input{includes/figures/tamari_dyck-3}} &
\scalebox{.5}{\input{includes/figures/tamari_dyck-4}}
\end{tabular}
\caption{The Tamari lattice on Dyck paths}
\label{fig:tamari-dyck}
\end{figure}

The Tamari lattice is defined on Catalan objects. As such, there are at least as many descriptions of the lattice as there are objects counted by the Catalan numbers, which is to say \emph{a lot}. For example, it is very common to see the lattice described on \defn{triangulations of a regular polygon} where the cover relations is a \defn{flip} of the triangulation. In my work, I often use the description in terms of \defn{Dyck paths}. A Dyck path is a path formed by two types of steps, ``up'' and ``down'', such that the path starts and finishes on the $x$ axis but never goes below. The cover relation that gives the Tamari lattice is the \defn{Dyck path Tamari rotation} which consist of \emph{lifting} the portion of the path following a down step (the portion is taken until the path reaches again the level of the down step). We show an example on Figure~\ref{fig:rotation_dyck} where we also represent the Dyck path as a binary word where $1$ means ``up'' and~$0$ ``down''. The Tamari lattices for $n=3$ and $4$ on Dyck paths are shown on Figure~\ref{fig:tamari-dyck}.

\section{Loday's associahedron}
\label{sec:asso-loday}

\begin{figure}[ht]
\center
\begin{tabular}{cc}
\includegraphics[scale=.6]{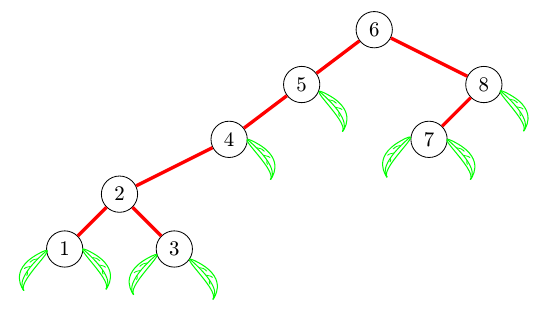}
&
\scalebox{.6}{\input{includes/figures/tamari_coordinates-3}} \\
$ (1,4,1,4,5,18,1,2)$ &
\end{tabular}
\caption{Coordinates of binary trees in Loday's associahedron.}
\label{fig:tree-loday}
\end{figure}

Loday gave the first explicit realization of the associahedron such that every vertex of the Tamari lattice has integer coordinates and corresponds to an extrema of the polytope~\citeext{Lod04,PSZ23}. Each tree of size $n$ corresponds to a point $(v_1, \dots v_n)$ in $\NN^n$ where $v_i$ is equal to the product of the number of leaves in the left subtree times the number of leaves in the right subtree of the node labeled by $i$ in the standard binary search tree labeling. We give an example on Figure~\ref{fig:tree-loday} where we have made the leaves (empty trees) explicit.

\begin{figure}
\center
\begin{tabular}{cc}
\includegraphics[scale=.2]{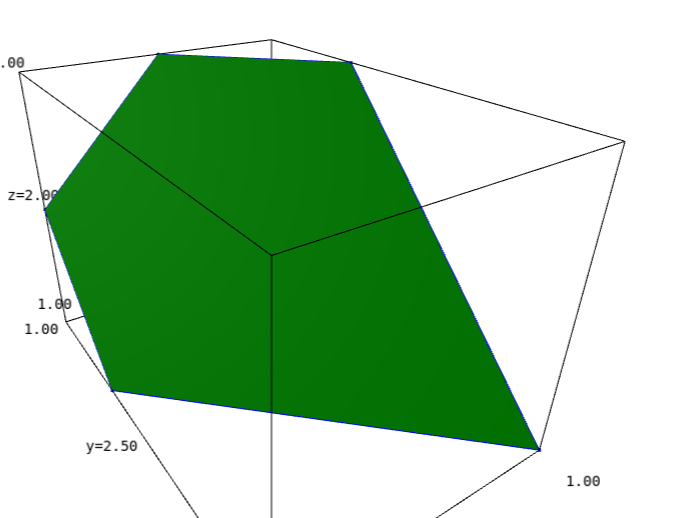} &
\includegraphics[scale=.2]{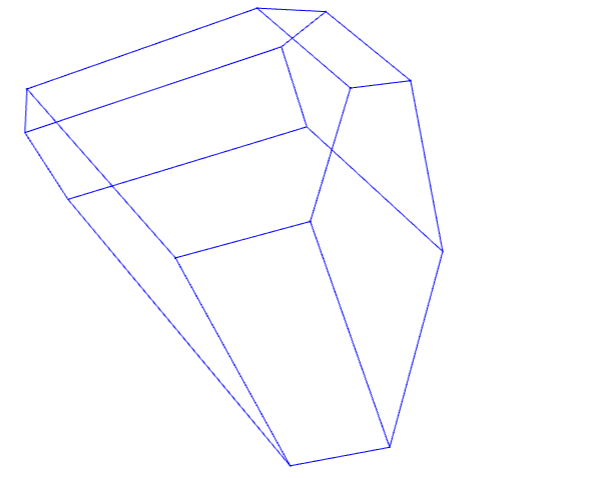}
\end{tabular}
\caption{The associahedra of dimensions 2 and 3 drawn by \Sage}
\label{fig:asso-sage}
\end{figure}

On the right of Figure~\ref{fig:tree-loday}, we show all the coordinates for the Tamari lattice on size $3$. You can check that the sum of the coordinates for each point is equal to $6$ like in the permutahedron. Actually, the four points on the left have permutation coordinates: the inverse permutation of their unique linear extension. They are vertices of both the associahedron and the permutahedron. This is true for all \defn{singletons}: binary trees with a single linear extension. Figure~\ref{fig:asso-sage} presents the associahedra of dimensions 2 and 3 drawn with \Sage{}, the code is available here~\citemesoft{PonSage23}.

\section{As a removahedron}
\label{sec:asso-faces}

\begin{figure}[ht]
\center
\begin{tabular}{cc}
\includegraphics[scale=.2]{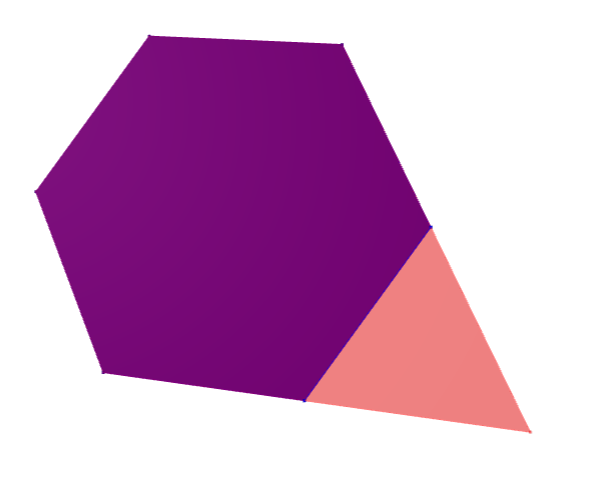} &
\includegraphics[scale=.2]{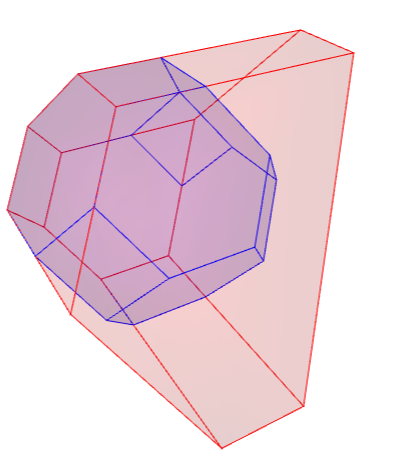}
\end{tabular}
\caption{Permutahedra inside of associahedra}
\label{fig:perm_asso}
\end{figure}

As we have seen, certain vertices of the associahedron correspond to vertices of the permutahedron. As we show on Figure~\ref{fig:perm_asso}, we can actually \emph{embed} the permutahedron inside the associahedron. Geometrically, this means that the list of inequalities defining the associahedron are a subset of the inequalities of the permutahedron we saw in Section~\ref{sec:perm-faces}. The associahedron is then a special case of \defn{removahedron} as defined by Pilaud~\citeext{Pil17}. 

\begin{figure}[ht]
\center
\includegraphics[scale=1]{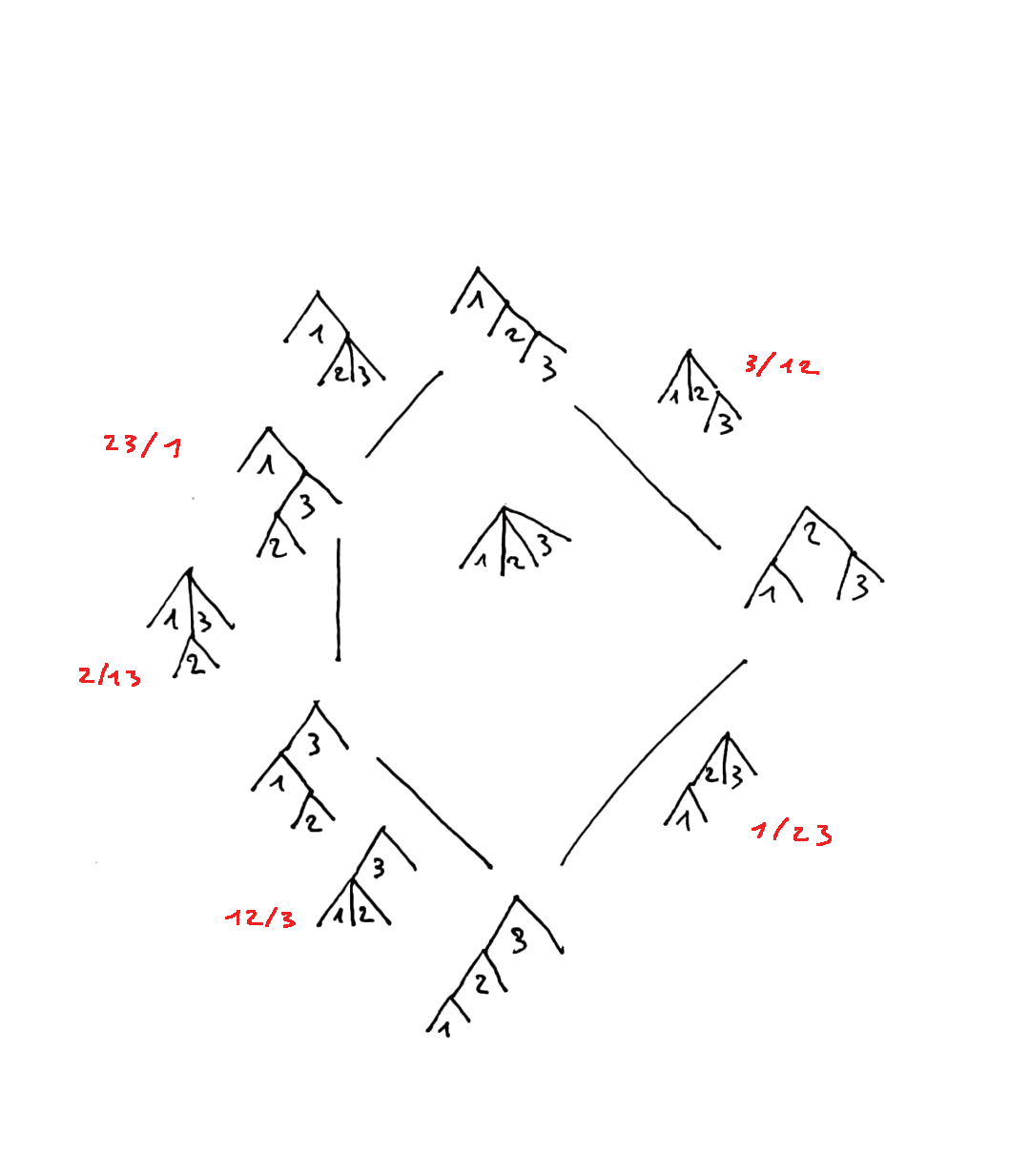}
\caption{Faces of the associahedron as Schr\"oder trees.}
\label{fig:shroder}
\end{figure}

In dimension $2$, the inequality which is ``removed'' is $x_1 + x_3 \geq 3$ corresponding to the ordered partition $13|2$. In general, all inequalities corresponding to faces where no vertex is a singleton binary tree are removed~\citeext{HLT11}. So faces of the associahedron are a subset of faces of the permutahedron. They also admit a combinatorial interpretation in terms of \defn{Schr\"oder trees}. A Schr\"oder tree is a rooted planar tree such that each node has at least two children. We label the ``valleys'' of the tree, \emph{i.e.}, the gap between two subtrees, with increasing labels from left to right. This way, labeled binary Schr\"oder trees actually correspond to binary search tree (the label of the node is the label of the unique valley underneath). They correspond to the vertices of the associahedron, \emph{i.e.}, faces of dimension 0. The other (non binary) trees are the faces of dimensions greater than $0$. More precisely, the dimension is the number of valleys minus the number of internal nodes. Just like ordered partitions are ``in between'' permutations in the weak order, Schr\"oder trees are ``in between'' binary trees. They correspond to partial rotations where a node is merged with its parent and not yet rotated to its right. We illustrate this on Figure~\ref{fig:shroder}.

The inequalities defining the associahedron are given by the Schr\"oder trees with exactly two internal nodes. They correspond to certain two-parts ordered partitions (in red on Figure~\ref{fig:shroder}). Recall that any non-empty proper subset $J$ of $[n]$ defines a permutahedron inequality corresponding to a two-parts ordered partition. For the associahedron, we only keep the subsets that are \emph{intervals}. Indeed, in dimension $2$, the only subset that is not kept is $\lbrace 1, 3 \rbrace$ (ordered partition $13|2$). In Schr\"oder trees, the subset $J$ correspond to the part that is bellow the second internal node. There is no way to get a two-nodes Sch\"roder tree with the bellow part being $\lbrace 1, 3 \rbrace$. Indeed, it would have to be bellow the label $2$ but could not be to its left (because $3 > 2$) nor to its right (because $1 < 2$). 

\section{As a Cambrian lattice}
\label{sec:asso-cambrian}

\begin{figure}[ht]
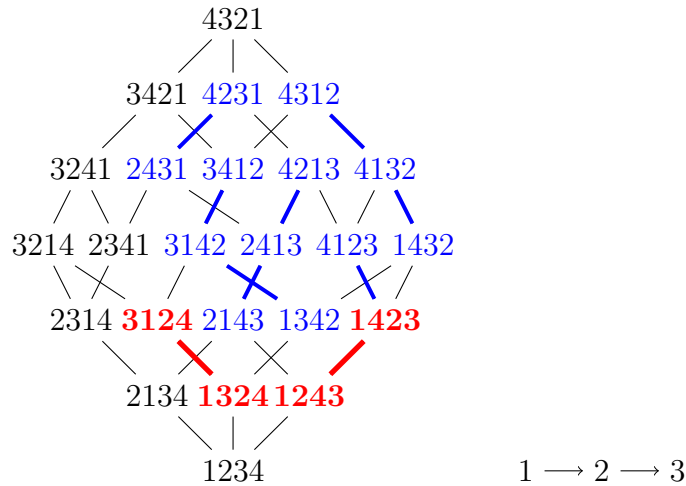

\center
\begin{tabular}{cc}
\input{includes/figures/perm_droit4_congruence_AA} &
\input{includes/figures/oriented_A3_Tamari}
\end{tabular}
\caption{Initial merges (in bold red) generating the sylvester congruence (in bold red and blue) and corresponding orientation of the Dynkin diagram.}
\label{fig:weak-merge}
\end{figure}

We have seen in Section~\ref{sec:perm-coxeter} that the weak order and the permutahedron can be defined for any Coxeter group. This is also true for the Tamari lattice, which then becomes a special case of the Cambrian lattices defined by Reading~\citeext{Rea06}. There are different ways to define and understand Cambrian lattices. We can see them as a special class of quotient lattices of the weak order. Indeed, the sylvester congruence, which defines the Tamari lattice, is characterized by two initial \defn{merge} as we illustrate on Figure~\ref{fig:weak-merge}.

A lattice congruence can be described by a list of merged edges in the Hasse diagram: two elements are in the same class if they are connected by a path of merged edges. For example, in Figure~\ref{fig:weak-merge}, $1342$, $3142$, and $3412$ are in the same class. On the other hand, $3214$ is alone in its class. You can check that these are the same classes as in Figure~\ref{fig:tamari-quotient4}. This representation is possible because all congruence classes form intervals of the lattice. A natural question is then: if I merge a single edge, what other edges need to be merged so that the final congruence is a lattice congruence? We say that the other necessary edges are \defn{forced} by the original one. This has been answered in general by Reading for the weak order using \defn{arc diagrams}~\citeext{Rea04}. The sylvester congruence is then the smallest lattice congruence containing the two initial edges in red on Figure~\ref{fig:weak-merge}: $1324$ -- $3124$, and $1243$ -- $1423$. 

These are the right sides of the two initial hexagons of the lattice. In type $A$, the weak order of size $n$ always starts with $n-2$ hexagons which correspond to the braid relations $s_i s_{i+1} s_i = s_{i+1} s_i s_{i+1}$. The sylvester congruence is then generated by $n-2$ merges: $s_{i+1}$ -- $s_{i+1}s_i$ for $1 \leq i < n- 1$. For $n=4$, this gives indeed the edge between $s_2 = 1324$ and $s_2 s_1 = 3124$, and the edge between $s_3 = 1243$ and $s_3 s_2 = 1423$. This easily generalizes to all finite Coxeter groups and can be encoded by an orientation of the Dynkin diagram. Indeed, each edge between $i$ and $j$ of the diagram is labeled with $k \geq 3$ and corresponds to the relation $(s_i s_j)^k = 1$. In the weak order, it then corresponds to a polygon with $2k$ edges. If $k$ is odd, one side is given by 

\begin{equation*}
s_i \relbar s_i s_j \relbar \dots \relbar (s_i s_j)^{\frac{k-1}{2}} \relbar (s_i s_j)^{\frac{k-1}{2}} s_i,
\end{equation*}
where $\relbar$ represents an edge in the weak order. The other side is 

\begin{equation*}
s_j \relbar s_j s_i \relbar \dots \relbar (s_j s_i)^{\frac{k-1}{2}} \relbar (s_j s_i)^{\frac{k-1}{2}} s_j.
\end{equation*}

\begin{figure}[ht]
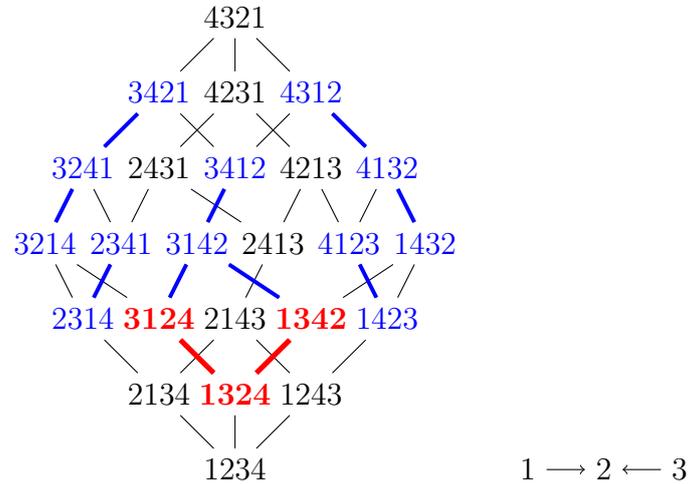

\center
\begin{tabular}{cc}
\input{includes/figures/perm_droit4_congruence_AY} &
\input{includes/figures/oriented_A3_AY}
\end{tabular}
\caption{Another Cambrian congruence on $A_3$ with corresponding orientation}
\label{fig:weak-merge-cambrian}
\end{figure}

Then an orientation of the Dynkin diagram defines a lattice congruence by deciding which side of the polygon to merge. If the orientation is $i \rightarrow j$, then all edges between $s_j s_i$ and $(s_i s_j)^{\frac{k-1}{2}}$ (for $k$ odd) are merged. If the orientation is $i \leftarrow j$, then all edges between $s_i s_j$ and $(s_i s_j)^{\frac{k-1}{2}}$ are merged. For type $A$, the polygon is always an hexagon and there is only one edge to merge. The orientation corresponding to the sylvester congruence is written on the right of Figure~\ref{fig:weak-merge}. We present another Cambrian congruence for type $A_3$ on Figure~\ref{fig:weak-merge-cambrian}.

The lattices corresponding to Cambrian congruences are called the Cambrian lattices. Not only do they define some ``Tamari'' lattices for all Coxeter groups, they also define generalizations of the classical Tamari lattice in type $A$ with different orientations of the Dynkin diagram. They share some of its properties. For example, their elements are always counted by the Catalan numbers: you can check that there are $14$ classes on Figure~\ref{fig:weak-merge-cambrian}. The unoriented graphs of the lattice Hasse diagram are isomorphic and they give different realizations of the associahedron~\citeext{HLT11}. 

The orientation of the Dynkin diagram somehow indicates in which order $s_i$ and $s_j$ much appear in the reduced word to avoid certain patterns. This is formalized through the system of \defn{$c$-sorted} words described by Reading~\citeext{Rea06}. Each orientation defines a \defn{Coxeter element} by reading a linear extension of the oriented Dynkin diagram. For example, for the Tamari orientation of Figure~\ref{fig:weak-merge}, this would be $c = s_1 s_2 s_3$. For the orientation of Figure~\ref{fig:weak-merge-cambrian}, this would be $c' = s_1 s_3 s_2 = s_3 s_1 s_2$. By taking an infinite concatenation of the Coxeter element, you get the $c^{\infty}$ word, for respective previous examples,

\begin{align*}
c^{\infty} &= s_1 s_2 s_3 | s_1 s_2 s_3 | s_1 s_2 s_3 | \dots \\
c'^{\infty} &= s_1 s_3 s_2 | s_1 s_3 s_2 | s_1 s_3 s_2 | \dots.
\end{align*}

Now each reduced word can be read as a subword of $c^{\infty}$. Let us call $I_1$ the subset of letters used from the first copy of $c$, $I_2$ from the second copy and so on. A reduced word is said to be \defn{$c$-sortable} if it exists as a subword of $c^{\infty}$ with $I_1 \subset I_ 2 \subset I_3 \dots$. For example, $s_1 s_2 s_1 = 3214$ is both $c$-sortable and $c'$-sortable as we have (for both $c$ and $c'$) $I_1 = \lbrace s_1, s_2 \rbrace$ and $I_2 = \lbrace s_1 \rbrace \subset I_1$. On the other hand, $s_2 s_3 = 1342$  is $c$-sortable but not $c'$-sortable. Indeed, for $c'$, we have $I_1 = \lbrace s_2 \rbrace$ and $I_2 = \lbrace s_3 \rbrace$. Elements of the group which have at least one $c$-sortable reduced word generalize the $312$-avoiding permutations for Cambrian lattices: they correspond to the minimal elements of the Cambrian congruence classes.

\section{As a sub and quotient Hopf algebra}
\label{sec:asso-hopf}

\begin{figure}[ht]
\center
\input{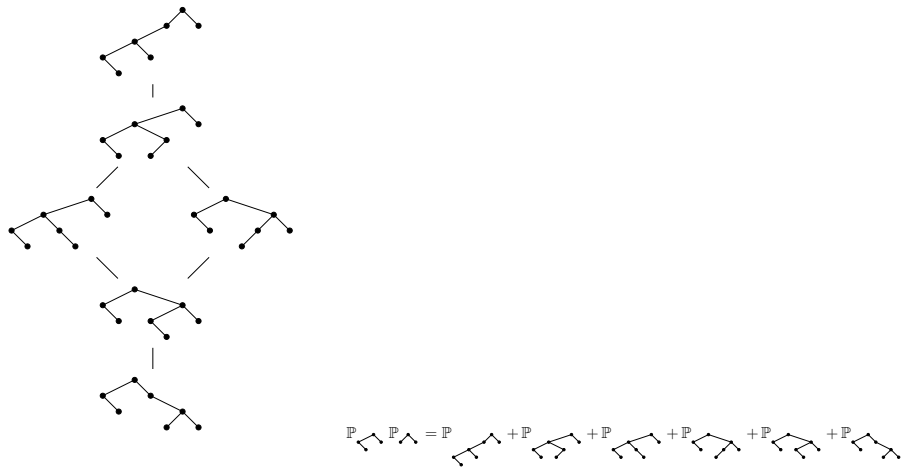}
\caption{Product of $\PBT$ as an interval of the Tamari lattice}
\label{fig:pbt-product}
\end{figure}

In~\citeext{LR98}, Loday and Ronco describe a Hopf algebra on binary trees as a quotient Hopf algebra of the Malvenuto-Reutenauer Hopf algebra on permutations that we described in~\ref{sec:perm-hopf}. They also note the correspondence between the Hopf algebraic relations and the geometrical relations between the permutahedron and associahedron. The Loday-Ronco Hopf algebra is studied in~\citeext{HNT05} where it is called $\PBT$: they show how it can be constructed using sylvester classes and binary search tree insertions. Indeed, we can define a basis $\BP$ indexed by binary trees as a sum over elements from the basis $\BF$ of $\FQSym$: the sum is over all linear extension of the standard binary search tree. See the example below.

\begin{equation}
\BP_{
\scalebox{.5}{\input{includes/figures/trees/T4-7}}
}  = \BF_{2143} + \BF_{2413} + \BF_{4213}.
\end{equation}
You can check that the binary search tree insertion of Section~\ref{sec:asso-quotient} gives the binary tree for all three permutations. The basis $\BP$ forms a sub Hopf algebra of $\FQSym$. As the product between two elements of $\BF$ is an interval in the weak order, the product of two elements of $\BP$ gives a sum over sylvester classes which is an interval of the Tamari lattice, see the example on Figure~\ref{fig:pbt-product}.

Besides, the basis $\BP$ admits a dual basis $\BQ$ which forms a quotient Hopf algebra of $\FQSym$. This is the original Loday-Ronco algebra which can be described through decreasing binary trees mimicking the relation between the left weak order and the Tamari lattice.

\section{The cube}
\label{sec:asso-cube}

\begin{figure}[ht]
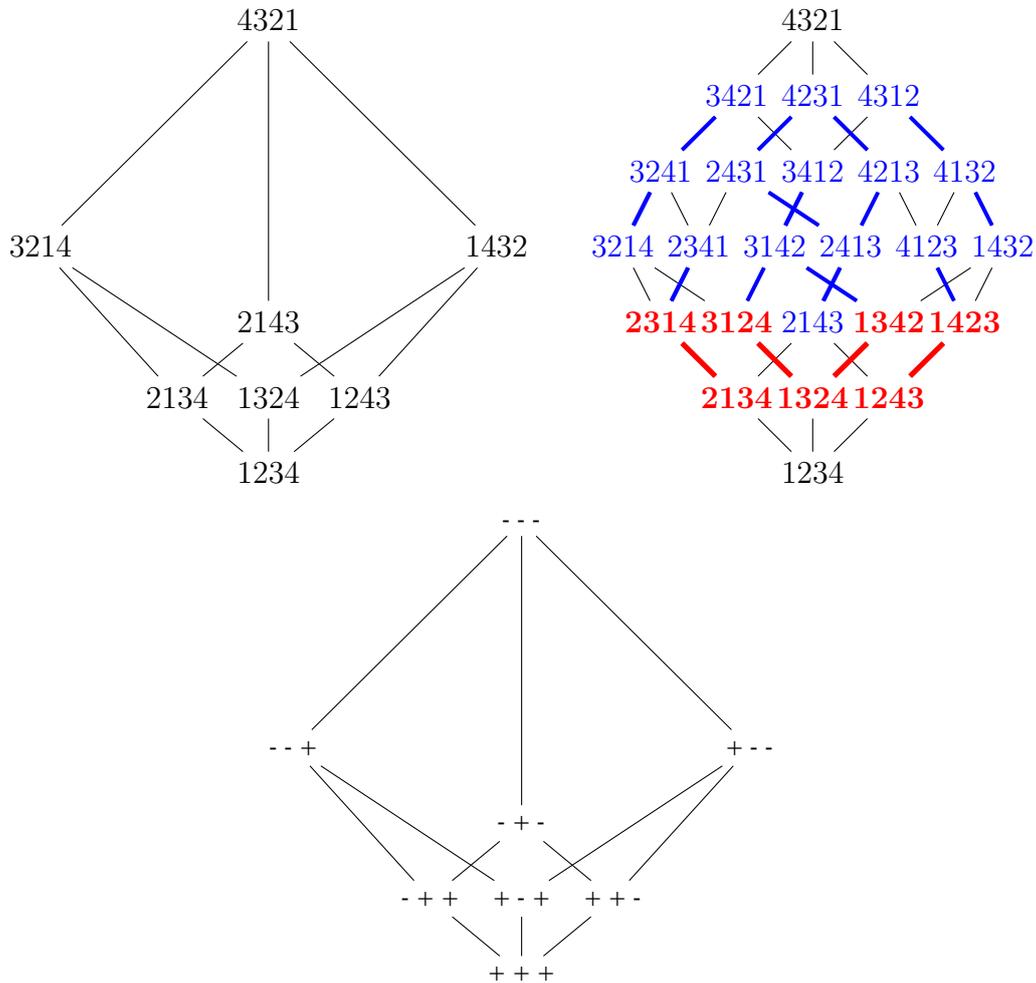

\center
\begin{tabular}{cc}
\input{includes/figures/cube_perms4}
&
\input{includes/figures/perm_droit4_congruence_XX}
\end{tabular}
\input{includes/figures/cube_binaryseq}
\caption{The boolean lattice as a sub and quotient lattice of the weak order.}
\label{fig:boolean}
\end{figure}

The associahedron can be seen as an \emph{intermediate step} between the permutahedron and the cube. Indeed, the cube corresponds to the boolean lattice which is also a sublattice and a quotient lattice of both the weak order and the Tamari lattice. The cube appears as a sublattice if you take permutations avoiding both patterns $312$ and $213$, see Figure~\ref{fig:boolean} for an example. In terms of quotient, the cube is given by the minimal lattice congruence obtained by merging \emph{both sides} of all initial hexagons of the weak order. This congruence is actually a surjective map between permutations and \defn{binary sequences}  based on the \defn{recoils} of the permutation. To a permutation of size $n$, we associate a binary sequence $v_1 \dots v_{n-1}$ with $v_i = +$ if $i+1$ appears after $i$ in the permutation and $v_i = -$ otherwise ($i$ and $i+1$ form a recoil). A compatible map also exists on binary trees using the \defn{canopy} of the binary tree. The congruence classes of the recoils map are then unions of sylvester classes. Besides, the cube is also a removahedron which can be constructed from the associahedron by removing more inequalities. And finally, it also corresponds to the Hopf algebra of recoils~\citeext{GKL+95} which is a sub hopf algebra of $\FQSym$ and $\PBT$. 

\section{$m$-Tamari and $\nu$-Tamari}
\label{sec:nu-tam}

\begin{figure}[ht]
\begin{tabular}{cc}
\scalebox{.5}{\input{includes/figures/mTamari-3-3-paths}}
&
\includegraphics[scale=1.5]{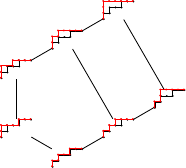}
\end{tabular}
\caption{The $m$-Tamari lattice and a $\nu$-Tamari lattice}
\label{fig:m-tam-nu-tam}
\end{figure}

We finish this chapter with two generalization of the Tamari lattice that we study in our work. The $m$-Tamari lattice appears in~\citeext{BPR12}. It is best described on $m$-ballot paths: they are paths between $(0,0)$ and $(mn,n)$ made of vertical and horizontal steps and which stay above the line $y = \frac{x}{m}$. In particular, an $m$-ballot path of size $n$ has $n$ vertical steps and $n \times m$ horizontal steps. For $m = 1$, they are exactly Dyck paths (where ``up'' is vertical and ``down'' is horizontal). A rotation operation can be defined for all $m$ which is very similar to the Dyck path rotation and also gives a lattice. See the example for $m=3$ and $n=3$ on Figure~\ref{fig:m-tam-nu-tam}. Actually, by replacing every vertical step by a series of $m$ up steps, we obtain an upper ideal of the classical Tamari lattice of size $n \times m$.  

Another way to understand $m$-ballot paths is to say that they are all the paths above the path represented by $(10^m)^n$ (where $1$ means vertical and $0$ horizontal). In~\citeext{PRV17}, Préville-Ratelle and Viennot show that you can define a lattice on paths above \emph{any} initial down path. This is what we call the $\nu$-Tamari lattice (where $\nu$ represents the initial down path). The rotation is a generalization of the Dyck path rotation. At each point, we compute the number of horizontal steps that can be taken to the right until crossing the initial down path. This is the distance between the point and the original path. A \defn{rotation} switches an horizontal step with the portion of path directly following the step, where the path is taken up to next point at a similar distance. We show an example on Figure~\ref{fig:m-tam-nu-tam}. Just like the $m$-Tamari lattices, the $\nu$-Tamari lattices can also be seen as a ``piece'' of a larger Tamari lattice. They also have interesting properties. For example, Ceballos, Padrol, and Sarminento provide a geometric realization of the $\nu$-Associahedron~\citeext{CPS16}.  

\part{Bijections on Tamari Intervals and More}
\label{part:bij}

\chapter{Intervals of the Tamari Lattice}
\label{chap:tamari-intervals}

\chapcitation{There are times when I feel like I'm in a big forest and don't know where I'm going. But then somehow I come to the top of a hill and can see everything more clearly.}{Maryam Mirzakhani, Brillant 10, 2005.}

In~\citeext{Cha05}, Chapoton proves that the number of \defn{intervals} in the Tamari lattice satisfies a surprisingly nice formula:

\begin{equation}
\label{eq:tam-intervals}
I_n = \frac{2}{n(n+1)} \binom{4n+1}{n-1} = \frac{2 (4n+1)!}{(n+1)!(3n+2)!}
\end{equation}
where $I_n$ is the number of couples of elements $x$, $y$, in the Tamari lattice on size~$n$, with $x \wole y$. I write ``nice formula'' because not only is it a closed product formula, but it is almost a rather simple binomial coefficient. It suggests that a \defn{bijective proof} is within reach. I write ``surprisingly'' because it is not at all common for combinatorial lattices to admit such a formula. For example, there is no such formula (and no hope for a product formula) for the intervals of the weak order. Similarly, other weak order lattice quotients, especially other type $A$ Cambrian lattices do not have such a formula. 

We write below the first numbers of the sequence~\citeext{OEIS000260},

\begin{equation}
1, 1, 3, 13, 68, 399, 2530, 16965, 118668.
\end{equation}
As referenced in the \emph{OEIS}, these numbers also count the number of simple, rooted, planar, triangular maps. The formula can indeed be found in~\citeext[formula~1.5]{Tut62}. Again, this is a surprising link between two structures that come from very different backgrounds. It turns out that these numbers are found in many different contexts. Along with my co-authors, I have introduced some objects in bijection with Tamari intervals: Tamari interval-posets~\citeme{CP15}, closed flows~\citemeconf{CCP14} and grafting trees~\citeme{Pon19}. In this chapter, I give a general overview of the current knowledge on objects counted by the Tamari interval numbers and thus place my own contributions in the global context to extract some interesting questions to explore.

\section{Generating functions, Catalan numbers and Tamari intervals}
\label{sec:gen-series}

As we have seen in Chapter~\ref{chap:asso}, the number of elements in the Tamari lattice is counted by the Catalan numbers~\eqref{eq:catalan}. This is probably one of the most famous sequences in combinatorics. In his book~\citeext{Sta15}, Stanley gives 214 objects counted by the Catalan numbers. One reason to find the Catalan numbers ``everywhere'' is that they encode the simplest possible binary structure. Indeed, let $F(z)$ be the generating functions of Catalan numbers,

\begin{equation}
F(z) := \sum_{n \geq 0} C_n z^n,
\end{equation}
Where $C_n$ is equal to~\eqref{eq:catalan}. Then, $F$ satisfies the following functional equation

\begin{equation}
\label{eq:catalan-eq}
F(z) = 1 + z F(z)^2.
\end{equation}
From there, it is immediate to find a recursive formula for $C_n$. To find the closed product formula, one needs to solve the quadratic equation and expand the solution using Newton binomial theorem. I have actually used this several times in talks for high schoolers and general audience as it is a very nice and understandable proof mixing combinatorics and analysis. In particular, this means that all Catalan objects can be \emph{decomposed} into two smaller objects just like a binary tree can be decomposed into a left and a right subtree. 

The formula~\eqref{eq:tam-intervals} looks very much like the Catalan formula. It is actually the solution of another functional equation which requires an extra parameter $x$. Let

\begin{equation}
\phi(z) := \sum_{n \geq 0} I_n z^n,
\end{equation}
where $I_n$ is~\eqref{eq:tam-intervals}. Then, $\phi(z) = \Phi(z,1)$ where 

\begin{equation}
\label{eq:tam-interval-eq}
\Phi(z,x) = 1 + zx \Phi(z,x) \frac{x\Phi(z,x) - \Phi(z,1)}{x - 1}.
\end{equation}
We give this exact functional equation in~\citeme[Theorem 3.2]{CP15} and this was already proven by Chapoton in~\citeext{Cha05} where he also solves it to obtain \eqref{eq:tam-intervals}. The extra parameter is called the \defn{catalytic parameter}: it is necessary to keep this extra parameter to express the functional equation. Just like the Catalan functional equation~\eqref{eq:catalan-eq}, this is quadratic, but the right part of the product is now transformed through a divided difference using the catalytic parameter. 

To understand what it means in terms of combinatorics, we can compare a simple example on both Catalan objects and Tamari intervals. In the generating function of Catalan, an object of size $k$ corresponds to a monomial $z^k$. The product in the functional equations means that it can be combined to an object of size $n-1-k$ to form a new object corresponding to a monomial $z \times z^k \times z^{n-1-k} = z^n$. In the case of Tamari intervals, two objects of respective size $k$ and $n-1-k$ are represented by two monomials $x^{v_1} z^k$ and $x^{v_2}z^{n-1-k}$. But the \emph{right} monomial is transformed through the divided difference into a sum $z^{n-1-k}(1 + x + x^2 + \dots + x^{v_2})$. It means that the combinations of two smaller objects now gives a set of $v_2 + 1$ objects of size $n$ which all carry an extra parameter. One way to prove that a combinatorial family is counted by $I_n$ is then to exhibit the underlying Catalan binary decomposition and the extra catalytic parameter.

\section{A map of maps using maps and trees}
\label{sec:map-map}

\begin{figure}[ht]
\scalebox{.7}{\input{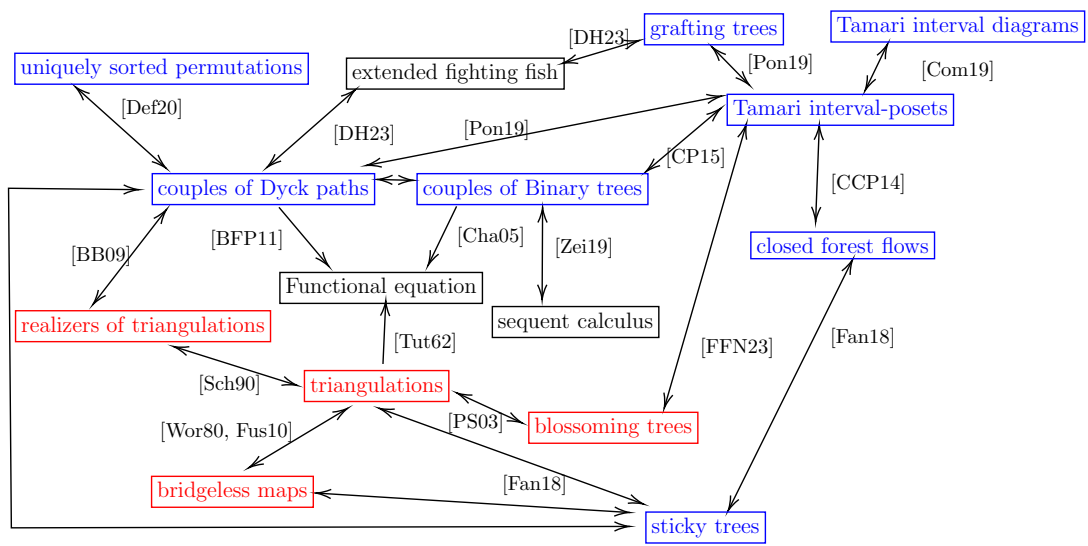}}
\caption{Landscape of bijections between objects counted by Tamari intervals}
\label{fig:tam-bijections}
\end{figure}

On Figure~\ref{fig:tam-bijections}, I present an non-exhaustive landscape of objects counted by $I_n$ with their bijections. I have written in blue the objects whose definitions are directly related to the Tamari lattice itself and in red the ones coming from planar maps. The study of planar maps, \emph{i.e.} graphs embed in the plane, is a field by itself which is connected to the study of the Tamari lattice through this surprising connection between Tamari intervals and triangulations (simple, rooted, planar, triangular maps). Coming from the Tamari side, my knowledge of the map aspect is limited and the landscape concentrates more on the Tamari related aspects. Besides, this landscape is restricted for now to families counted by $I_n$ thus excluding bijections on subsets of intervals or generalizations which are also numerous. 

A first comment on this landscape is the great variety of objects that appear. On top of the red (maps) and blue (Tamari) ones, we see sequent calculus from the work of Noam Zeilberger~\citeext{Zei19} and extended fighting fish~\citeext{DH23}, which represent certain branching surfaces or, equivalently, walks in the quarter plane. In Chapter~\ref{chap:tam-stat}, we explore some specific interesting statistics on intervals and see how they relate to some of these objects. Besides, this landscape concentrates on combinatorics: the $I_n$ numbers also appear in some algebraic settings which we explore in Chapter~\ref{chap:qt}.

As we have said, Chapoton is the first to count Tamari intervals~\citeext{Cha05}. He proves that they satisfy the functional equation~\eqref{eq:tam-interval-eq} and notes the connection with maps but does not provide a direct bijection.  The proof of Chapoton is based on binary trees and the catalytic parameter is the number of nodes on the left branch of the maximal tree in the interval. A bijection between triangulations and intervals is later given by Bernardi and Bonichon~\citeext{BB09} using \defn{Schnyder woods} or \defn{realizers}: this is a way to encode a triangulation with three intertwined planar trees. The bijection uses Dyck paths. A triangulation can be represented by multiple set of realizers: Tamari intervals then correspond to \defn{minimal realizers} while general realizers correspond to couples of Dyck paths $(P,Q)$ where $P$ is \emph{below} $Q$. 

The functional equation can also be expressed in terms of Dyck paths. The catalytic parameter then corresponds to some statistic of the lower Dyck path called the number of \defn{contacts}. This is what is done in~\citeext{BMFPR11} where the authors solve a more a general functional equation corresponding to intervals in the $m$-Tamari lattices. We discuss it in Section~\ref{sec:intervals-m-tam}. In this paper, they also notice that another statistic called the \defn{initial rise} could serve as a catalytic parameter. We have explored this question in~\citeme{Pon19} and discuss it in Chapter~\ref{chap:tam-stat}.

My main contribution to the landscape of Figure~\ref{fig:tam-bijections} is the introduction of \defn{Tamari interval posets}. At the end of my thesis, I started a collaboration on Tamari intervals with Grégory Ch\^atel, who was also a PhD student in my institution. We introduced Tamari interval-posets in a conference paper~\citemeconf{CP13}. In the full version~\citeme{CP15}, we extend this definition to $m$-Tamari intervals. In~\citemeconf{CCP14} and~\citeme{Pon19}, I present other useful objects in bijection with Tamari interval-posets. Nevertheless, even though  they give interesting insights about intervals, these objects did not permit to obtain a direct bijective proof of Chapoton's formula~\eqref{eq:tam-interval-eq}. This could be done indirectly through the bijection with maps and their interpretation with \defn{blossoming trees} but was not completely satisfactory. This is now solved by a very recent work of Fang, Fusy and Nadeau~\citeext{FFN23} which is actually still in writing: this is a direct bijection between blossoming trees (which can be counted directly by~\eqref{eq:tam-intervals}) and Tamari interval posets.

\section{Tamari interval-posets}
\label{sec:tam-interval-posets}

\begin{figure}[ht]
\center
\input{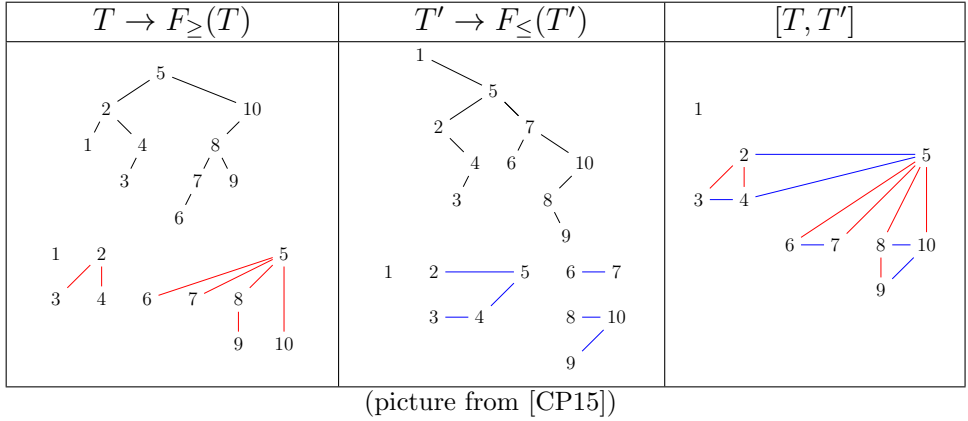}

\footnotesize{(picture from~\citeme{CP15})}
\caption{Construction of Tamari interval-posets from binary trees }
\label{fig:interval-poset}
\end{figure}

Tamari interval-posets of size $n$ are partial orders on $[n]$ satisfying certain conditions. We write $a \trprec b$ to say that $a$ is smaller than $b$ in the poset (to differentiate from the symbol $\woless$ used for the Tamari lattice and the weak order) and $a < c$ if $a$ is smaller than $c$ as integers. The definition stated in~\citeme[Definition 2.7]{CP15} says that a partial order on $[n]$ is a Tamari interval-poset if for all $a < b < c$, 

\begin{align}
\label{eq:interval-poset-bc}
a \trprec c & \Rightarrow b \trprec c; \\
\label{eq:interval-poset-ba}
c \trprec a & \Rightarrow b \trprec a.
\end{align}
We give an example on the right of Figure~\ref{fig:interval-poset}. We see that we have $1 \trprec 5$ which implies $b \trprec 5$ for $1 < b < 5$. Similarly, we have $10 \trprec 5$ which implies  $b \trprec 5$ for $5 < b < 10$. Besides, Tamari interval-posets have been integrated into \Sage. The examples of this Section can be found in~\citemesoft{PonSage23}. The main result of~\citeme{CP15} is

\begin{theorem}[Theorem 2.8 from~\citeme{CP15}]
\label{thm:interval-posets}
Tamari interval-posets are in bijection with intervals of the Tamari lattice.
\end{theorem}

The bijection is built from binary trees using the connection between the Tamari lattice and the weak order. Recall that each binary tree can be seen itself as a poset whose linear extensions are the permutations of the sylvester class (see Figures~\ref{fig:bst-insertion}, \ref{fig:tamari-quotient3}, and~\ref{fig:tamari-quotient4}). Let $[T,T']$ be a Tamari interval expressed as a couple of binary trees. We construct a new poset called the \defn{decreasing forest} of~$T$, $\dec(T)$, by keeping all \defn{decreasing relations} ($c \trprec a$ with $a < c$) of the binary tree seen as a poset. For example, on Figure~\ref{fig:interval-poset}, you see that $6$, $7$, $8$, $9$, and $10$ are below $5$ in the binary tree: they are decreasing relations and are kept in~$\dec(T)$. On the other hand, the increasing relation $2 \trprec_T 5$ is not part of $\dec(T)$. Similarly, we construct the \defn{increasing forest}, $\inc(T')$, of the maximal tree in the interval $T'$ by keeping its \defn{increasing relations}. The interval-poset is constructed by taking relations from both the decreasing forest of $T$ and the increasing forest of $T'$. 

We usually draw it as in Figure~\ref{fig:interval-poset} as the union of the Hasse diagrams of increasing relations (in blue, oriented from left to right) and decreasing relations (in red, oriented from bottom to top). For example, $9 \trprec 10$ is drawn on Figure~\ref{fig:interval-poset} even though it is not a cover relation as we have $9 \trprec 8 \trprec 10$, because it is a cover relation of the increasing forest. This representation is very useful for many constructions on interval-posets.

The linear extensions of a Tamari interval-poset correspond to an interval of the weak order. More precisely, the linear extensions of the decreasing forest of a binary tree $T$ correspond to the interval between the minimal permutation of its sylvester class and the maximal permutation of the weak order. For example, for the tree $T$ of Figure~\ref{fig:interval-poset}, the linear extensions of $\dec(T)$ are the interval $[1~3~4~2~6~7~9~8~10~5, 10~9~8~7~6~5~4~3~2~1]$. Similarly, the linear extensions of the increasing forest of a binary tree $T'$ correspond to the interval between the identity and the maximal permutation of the sylvester class of $T'$. For $T'$ as in Figure~\ref{fig:interval-poset}, we get the interval $[1~2~3~4~5~6~7~8~9~10, 9~8~10~6~7~3~4~2~5~1]$. The linear extensions of the Tamari interval-poset are then the intersections between the two intervals, here $[1~3~4~2~6~7~9~8~10~5, 9~8~10~6~7~3~4~2~5~1]$. 

\begin{figure}
\input{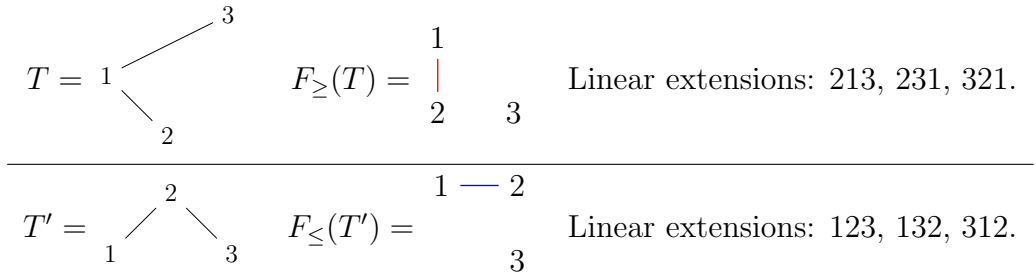}
\caption{Decreasing and increasing forests of two incompatible trees.}
\label{fig:incompatible-forests}
\end{figure}

In particular, if you take two trees $T$ and $T'$ which are not comparable in the Tamari lattice, then $\dec(T)$ and $\inc(T')$ have no linear extensions in common: the intersection of the intervals is empty. In this case, you cannot construct the interval-poset because the relations of $\dec(T)$ and $\inc(T')$ are incompatible. For example, on Figure~\ref{fig:incompatible-forests}, we see that $2 \trprec 1$ in $\dec(T)$ while $1 \trprec 2$ in $\inc(T')$: we cannot construct a poset with both relations.

\begin{figure}[ht]
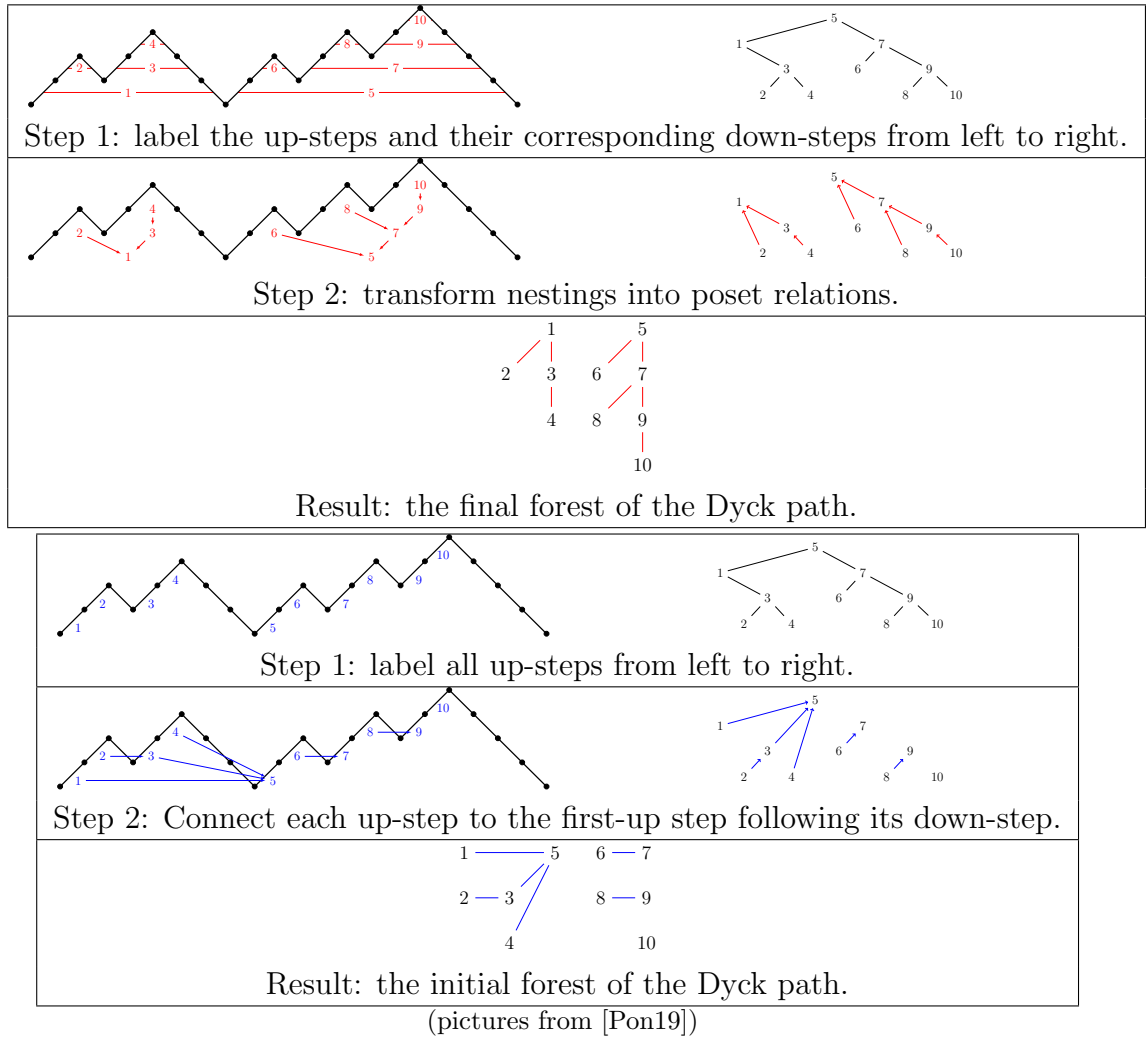

\center
\input{includes/figures/dyck-dec-forest}

\input{includes/figures/dyck-inc-forest}

\footnotesize{(pictures from~\citeme{Pon19})}
\caption{Bijection between Dyck paths and decreasing and increasing forests}
\label{fig:dyck-dec-inc}
\end{figure}

The interval-poset conditions~\eqref{eq:interval-poset-bc} and~\eqref{eq:interval-poset-ba} directly come from the increasing and decreasing forests construction. In particular, if an increasing forest is compatible with a decreasing forest (\emph{i.e.}, they share a common linear extension), then the result of their union is always an interval-poset. Besides, the construction of the forests is actually a very classical bijection between binary trees and planar forests. The decreasing forest is the ``left sibling -- right child'' while the increasing forest is the symmetric ``right sibling -- left child'': for each node, its left (resp. right) subtree becomes its left (resp. right) sibling while its right (resp. left) subtree becomes its child. In particular, the nodes on the left (resp. right) branch become the roots of the trees in the forest. This bijection easily translates to Dyck paths as illustrated on Figure~\ref{fig:dyck-dec-inc}.

To recover the functional equation~\eqref{eq:tam-interval-eq}, we describe a composition of Tamari interval-posets that translates into the product and divided difference. This is illustrated on Figure~\ref{fig:composition-interval-poset}. The catalytic parameter is the number of connected components in the decreasing forest: \emph{i.e.} the poset where you keep only the decreasing relations (in red). Following the bijection with binary trees, this corresponds to the number of nodes on the left branch of the minimal tree. On the example of Figure~\ref{fig:composition-interval-poset}, it is equal to $2$ for the interval on the left and to $3$ for the interval on the right. The composition adds an extra node (the node $4$) and creates a set of interval-posets where nodes $1$ to $3$ correspond to the left interval while nodes $5$ to $8$ correspond to the right one. Increasing relations have been added between all nodes from the left interval to the node $4$. The set is over all possible ways to add decreasing relations between the right interval and the new node $4$. In the generating function, the left interval would be the monomial $x^2 z^3$ while the right is $x^3 z^4$. The composition gives the polynomial

\begin{align}
x^2 z^3 xz \frac{x^4 z^4 - z^4}{x - 1} &= x^3 z^8 (1 + x + x^2 + x^3) \\
&= x^3 z^8 + x^4 z^8 + x^5 z^8 + x^6 z^8. 
\end{align} 

\begin{figure}
\center
\input{includes/figures/composition}

\footnotesize{(picture from~\citeme{CP15})}
\caption{Composition of interval-posets}
\label{fig:composition-interval-poset}
\end{figure}

This composition (and associated decomposition) allows us to recover the functional equation~\eqref{eq:tam-interval-eq}. It also gives a natural way to write the generating function $\Phi(z,x)$ of Tamari intervals as a sum over binary trees

\begin{equation}
\Phi(z,x) = \sum_{T} z^{\size(T)} \BT_T(x)
\end{equation}
where the sum is over binary trees and $\BT_T(x)$ is a polynomial defined recursively such that $\BT_T(1)$ is the number of trees smaller than or equal to $T$ in the Tamari lattice. The recursive definition (Definition 1.1 of~\citeme{CP15}) is given by $\BT_T(x) = 1$ if $T$ is empty and otherwise

\begin{equation}
\label{eq:tam-bt}
\BT_T(x) = x \BT_L(x)  \frac{x \BT_R(x) - \BT_R(1)}{x - 1}.
\end{equation}

As we see on Figure~\ref{fig:tam-bijections}, Tamari interval-posets have allowed us and others to create new bijective connections between Tamari intervals and other combinatorial families. They also give a nice characterization of some special families of intervals: \defn{exceptional intervals} and \defn{modern intervals} which have appeared in the work of Rognerud~\citeext{Rog20, Rog21}.

\section{Closed flows}
\label{sec:tam-flows}

In~\citeext{Cha14}, Chapoton studies certain polynomials appearing in the context of the Pre-Lie operad. They can be interpreted using \defn{flows on planar forests}. It appeared that some of these polynomials were related to the Tamari lattice and actually to the Tamari polynomials that we defined in~\citeme{CP15}. This lead us to define a simple bijection between closed flows and Tamari interval-posets~\citemeconf{CCP14}. 

\begin{figure}
\center
\scalebox{.7}{\input{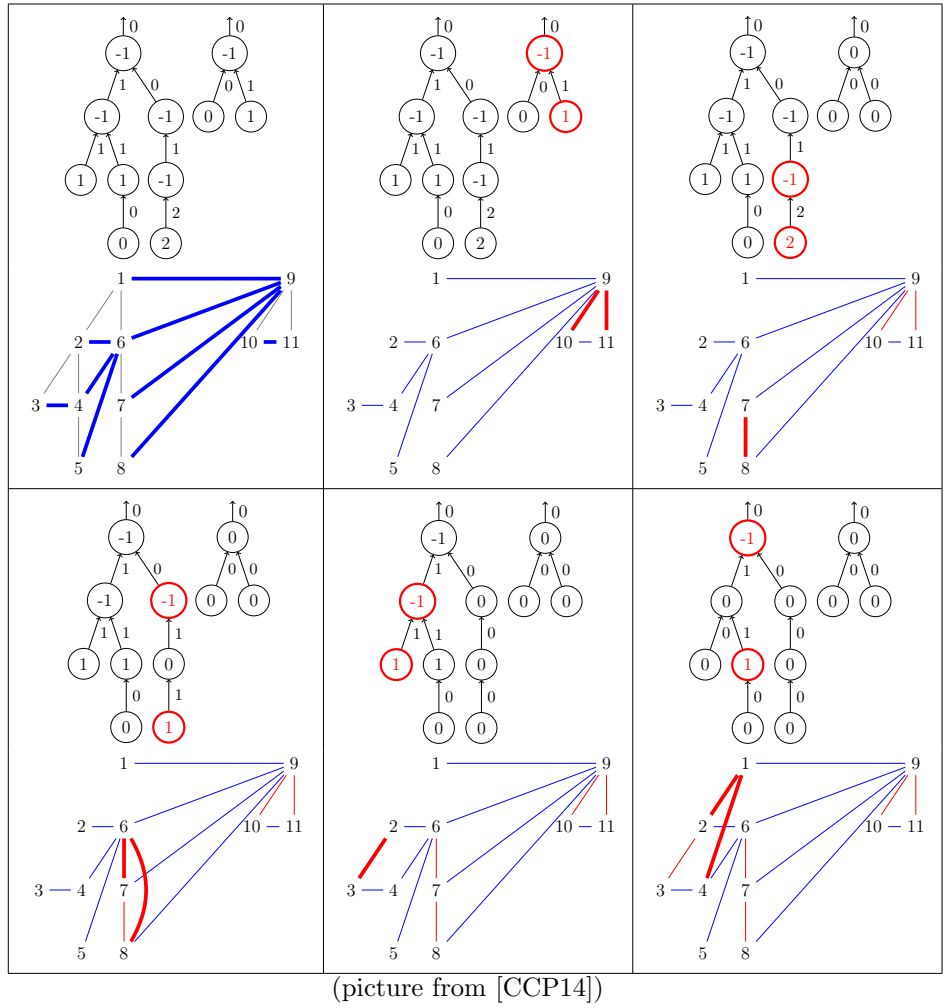}}

\footnotesize{(picture from~\citemeconf{CCP14})}
\caption{Bijection between closed flows and Tamari interval-posets}
\label{fig:flows-tamari}
\end{figure}

A \defn{flow} can be seen as some fluid going up a tree. Each node of the tree can be a \defn{source} or a \defn{leak} by respectively adding or subtracting to its incoming flow. On each edge, the flow must be weakly positive. We show an example on Figure~\ref{fig:flows-tamari}. The value added /removed to the incoming flow is written on each vertex. By definition, we state that the only possible value for a leak is $-1$ but sources do not have an upper limit. The \defn{exit rate} is the sum of the out flow for all trees. We say that a flow on a forest is \defn{closed} if the exit rate is $0$. In~\citemeconf{CCP14}, we present a bijection between closed flows and Tamari interval-posets and obtain the following as a consequence.

\begin{theorem}[Theorem~4.1 from~\citemeconf{CCP14}]
The number of closed flows on a given forest is the number of trees smaller than or equal to a certain tree in the Tamari lattice.
\end{theorem}

We illustrate the bijection on Figure~\ref{fig:flows-tamari}. The general idea is the following: the flow forest gives the increasing forest of the interval-poset through a classical bijection. Then each vertex of the flow corresponds to a vertex of the Tamari interval-posets. We add the decreasing relations by reading the leaks in decreasing order. Each leak $i$ is connected to the closest possible source $k$ and all relations $j \trprec i$ for $i < j \leq k$ are added to the interval-poset.

\section{Grafting trees}
\label{sec:grafting}

\begin{figure}[ht]
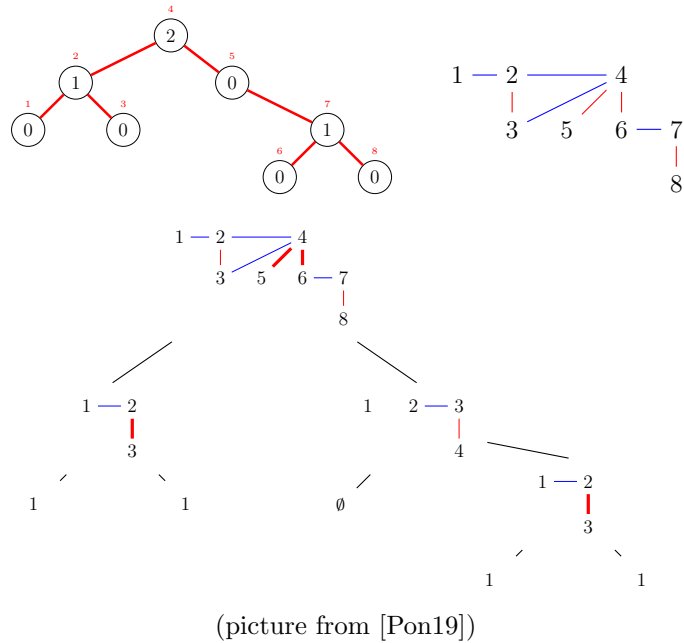

\center
\begin{tabular}{cc}
\scalebox{0.6}{
\input{includes/figures/grafting_tree_example}
}
&
\scalebox{0.8}{\input{includes/figures/interval-posets/I8-ex3}}
\end{tabular}

\input{includes/figures/grafting-tree-construct}

\footnotesize{(picture from~\citeme{Pon19})}

\caption{Example of grafting tree with corresponding interval-poset decomposition.}
\label{fig:grafting-tree}
\end{figure}

The composition of Tamari interval-posets of Figure~\ref{fig:composition-interval-poset} can be seen as a bijection between interval-posets of size $n$ and triples $(I,J,r)$ such that $I$ and $J$ are Tamari interval-posets with $\size(I) + \size(J) = n$ and $r$ is an extra parameter between $0$ and the so-called \defn{contacts} of $J$. The \emph{contacts} are the catalytic parameter, they correspond to the number of nodes on the left branch of the lower binary tree, or equivalently to the number of non-initial contacts between the lower Dyck path and the origin, or the number of connected components in the decreasing forest of $J$. This gives a recursive decomposition of Tamari interval-posets which can be encoded in a labeled binary tree where the label is the parameter $r$. We show this on Figure~\ref{fig:grafting-tree}.

In~\citeme{Pon19}, I call these trees \defn{grafting trees}. They are a direct combinatorial interpretation of the functional equation~\eqref{eq:tam-interval-eq}. The tree itself is the upper tree of the interval. The label of a node $i$ is the number of children of $i$ in the decreasing forest of the interval-poset. They can be interpreted as a list of rotations to apply to the minimal tree of the Tamari lattice to reach the lower binary tree of the interval. In particular, if all the labels are $0$, then the lower binary tree of the interval is the minimal element of the Tamari lattice.

Grafting trees are all labeled binary trees such that for each node $v$, the label $\ell(v)$ satisfies $0 \leq \ell(v) \leq \size(T_R(v)) - \sum_{w \in T_R(v)} \ell(w)$ where $T_R(v)$ is the right subtree of $v$. For example, on Figure~\ref{fig:grafting-tree}, the label of the root must be comprised between $0$ and $3$ because the right subtree has $4$ nodes and labels summing to $1$. If labels all have maximal values, then it corresponds to the interval $[T,T]$ made of a single tree.

As they directly translate the functional equation, grafting trees are very natural objects to consider when working with Tamari intervals. I use them in~\citeme{Pon19} to describe certain involutions as I explain in Chapter~\ref{chap:tam-stat}. They also make a connection to maps as there is an easy bijection between grafting trees and $(1,1)$-description trees found in~\citeext{CS03}.

\section{Interval-posets in $m$-Tamari}
\label{sec:intervals-m-tam}

In~\citeext{BPR12}, Bergeron and Préville-Ratelle conjecture that the number of intervals of the $m$-Tamari lattices of Section~\ref{sec:nu-tam} is given by 

\begin{equation}
\label{eq:mtam-intervals}
I_{m,n} := \frac{m+1}{n(mn+1)}\binom{(m+1)^2n +m}{n-1}.
\end{equation}
This is proved in~\citeext{BMFPR11} using an $m$-generalization of the functional equation~\eqref{eq:tam-interval-eq}

\begin{equation}
\label{eq:mtam-intervals-eq}
\Phi^{(m)}(z,x) = 1 + \Bm(\Phi, \Phi, \dots, \Phi)
\end{equation}
where $\Bm$ is an $(m+1)$ linear form defined by 
\begin{equation}
\Bm(f,g_1,\dots, g_m) = zx f \Delta( g_1 \Delta( \dots \Delta(g_m)) \dots ),
\end{equation}
and 
\begin{equation}
\Delta(g) = \frac{xg(x) - g(1)}{x - 1}.
\end{equation}
In other words, the product and divided difference are applied $m$ times. The authors of~\citeext{BMFPR11} prove that the intervals of the $m$-Tamari lattices satisfy the functional equation using $m$-Dyck paths and the \defn{contact} statistic as a catalytic parameter. Then they solve the equation with a guess-and-check approach to obtain~\eqref{eq:mtam-intervals}. Note that unlike the classical case, this is not a simple equation and requires advanced techniques.

\begin{figure}
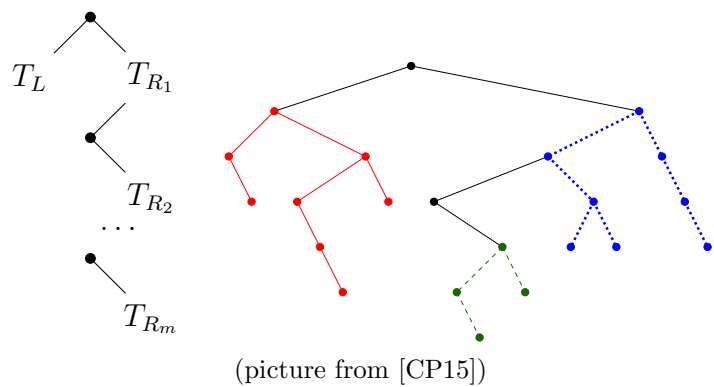

\center
\begin{tabular}{cc}
\input{includes/figures/m-binary-trees}&
\scalebox{.6}{\input{includes/figures/mbinary-example}}
\end{tabular}

\footnotesize{(picture from~\citeme{CP15})}
\caption{Decomposition of $m$-binary trees and example for $m=2$.}
\label{fig:m-binary}
\end{figure}

As the $m$-Tamari lattices can be seen as a special interval of the Tamari lattice, we can express intervals of the $m$-Tamari lattice as some subfamily of Tamari interval-posets. This is what we do in~\citeme{CP15}. This requires to characterize the binary trees which appear in the $m$-Tamari lattice and see how they are in bijection with $(m+1)$-ary trees. We introduce the notion of $m$-binary trees and show how each tree can be decomposed into a set of $m$ ``roots'' and $m+1$ subtrees as we illustrate on Figure~\ref{fig:m-binary}.

This gives the following bijection.

\begin{theorem}[Theorem~4.6 of \citeme{CP15}]
Intervals of the $m$-Tamari lattice in size~$n$ are in bijection with Tamari interval-posets of size $n \times m$ such that, for all $1 \leq i \leq n$, we have

\begin{equation}
im \trprec im - 1 \trprec \dots \trprec im - (m-1).
\end{equation}
\end{theorem}

\begin{figure}[ht]
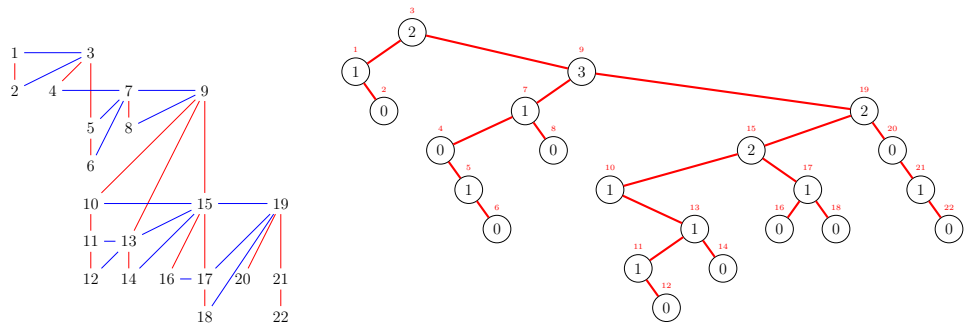

\center
\begin{tabular}{cc}
\scalebox{.5}{\input{includes/figures/interval-posets/I22-ex1}}
&
\scalebox{.5}{\input{includes/figures/m-grafting-tree}}
\end{tabular}
\caption{An $m$-Tamari interval-poset and its corresponding grafting tree.}
\label{fig:m-interval-example}
\end{figure}

See an example on Figure~\ref{fig:m-interval-example} for $m=2$ and $n=11$: you can check that we have $2i \trprec 2i - 1$ for all $i$.
Using this, we recover the functional equation~\eqref{eq:mtam-intervals-eq} and, for each element of the lattice, find a recursive formula generalizing~\eqref{eq:tam-bt}. In~\citeme{Pon19}, we also explain how this translates to grafting trees. An $m$-Tamari interval poset corresponds to a grafting tree such that $\ell(v_i) \geq 1$ if $i \not\equiv 0$ mod $m$. See the example on Figure~\ref{fig:m-interval-example} for $m = 2$: the labels of all odd value nodes are greater than or equal to $1$. In particular, this implies that the grafting tree is an $m$-binary tree. We call these trees the $m$-grafting trees.

\section{Open questions}
\label{tam:questions}

\begin{figure}[ht]
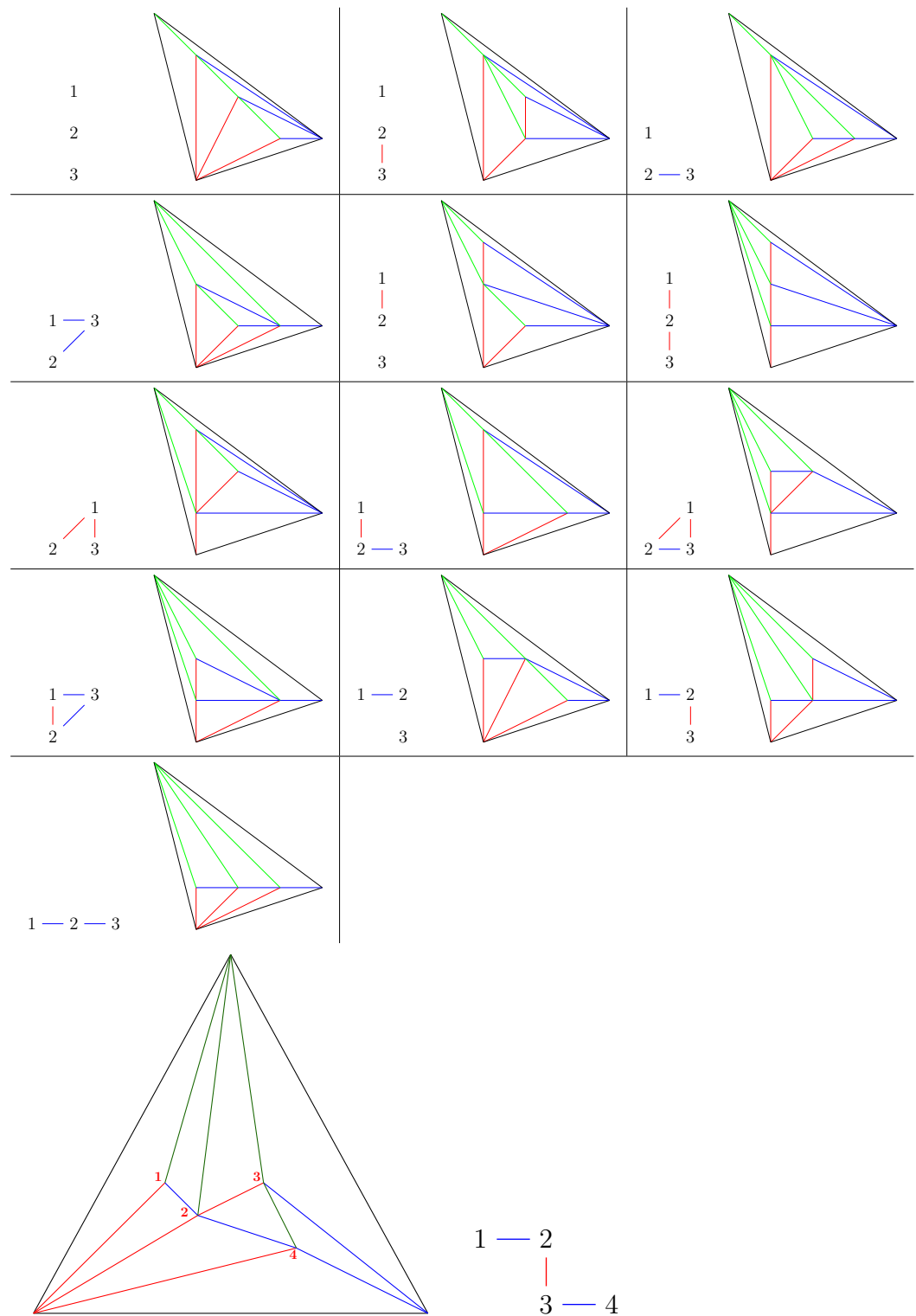

\scalebox{.8}{\input{includes/figures/triangulation_bijection-3}}

\input{includes/figures/example_triangulation_bijection}
\caption{Bijection between interval-posets and triangulations}
\label{fig:tamari-intervals-triangulations}
\end{figure}

On Figure~\ref{fig:tamari-intervals-triangulations}, we present the bijection of Bernadi and Bonichon~\citeext{BB09} directly between Tamari interval-posets and triangulations. We show the $13$ intervals for size $3$ and one example in size $4$. There seems to be a direct link between the decreasing and increasing forests of Tamari interval-posets on one side and two of the trees of the realizers of the triangulation on the other side. On the last example, we see that the red tree of the triangulation (rooted on the bottom left) is exactly the decreasing forest of the interval-poset. The blue tree (rooted on the bottom right) is a bit different from the increasing forest but each node $v$ has the same number of children in both the realizer and the interval-poset. 

I never had to time to fully formalize and prove this description and to explore it. This leaves some questions open. Triangulations have a clear symmetry of order $3$, what is this symmetry on Tamari interval-posets? Basically, what is the ``third tree'' (in green, rooted at the top) on the Tamari interval-posets? What is the connection between this bijection and the new bijection of~\citeext{FFN23}? Beside, neither~\citeext{BB09} nor~\citeext{FFN23} (to our knowledge) generalize to $m$-Tamari. The main open question then remains: is there a ``map'' interpretation of the intervals of the $m$-Tamari lattices and a bijective proof of~\eqref{eq:mtam-intervals}?

\chapter{Statistics and Bijections on Tamari Intervals}
\label{chap:tam-stat}

\chapcitation{I counted everything. I counted the steps to the road, the steps up to church, the number of dishes and silverware I washed … anything that could be counted, I did.}{Katherine Johnson}

Having multiple combinatorial families to represent Tamari intervals is useful when we try to understand their inner structure. Indeed, each family has natural statistics and symmetries. By combining them, we sometimes reveal more properties. When solving the functional equation~\eqref{eq:mtam-intervals-eq}, the authors of~\citeext{BMFPR11} noticed that there was a symmetry between the \defn{contacts} of the lower path and the \defn{initial rises} of the upper path. This was the result of an analytic proof (on intervals of $m$-Tamari) and they did not have a combinatorial explanation, even for $m=1$. As the contacts are the catalytic parameter of the functional equation, this suggests that another decomposition of Tamari intervals is possible based on the initial rises. In~\citemeconf{CCP14}, we provide this decomposition and note that it actually gives an involution on Tamari intervals. But it is only in~\citeme{Pon19} that I provide a ``natural way'' to see this involution using grafting trees. This also led me to define a very interesting transformation on grafting trees related to $m$-Tamari intervals. It is the key element to prove combinatorially the symmetry of the two statistics in the $m$-Tamari case as it was stated in~\citeext{BMFPR11} and to solve a more general conjecture of Préville-Ratelle~\citeext[Conjecture 17]{PR12}.

In this chapter, I give an overview of my results in~\citeme{Pon19}. I show the two natural involutions on Tamari intervals that led to the \defn{rise-contact involution}. I explain their effects on different interesting statistics and then show how to generalize the involution to $m$-Tamari intervals. Note that the rise-contact involution has been implemented into \Sage. The computations originally provided with the paper~\citeme{Pon19} are here~\citemesoft{PonSage18}, the specific examples of this Chapter are available here~\citemesoft{PonSage23}.

\section{Contacts, descents, and complement involution}

We have mentioned many times that the functional equation of the Tamari intervals requires a catalytic parameter. Different choices are possible for the statistic it refers to. In~\citeext{BMFPR11} and~\citeext{CP15}, the number of \defn{non initial contacts} of the lower path is taken, \emph{i.e.}, the number of time the path comes back to the $x$ axis. For short, we call it the \defn{contact value} of the interval. On binary trees, this is the number of nodes on the left branch of the lower tree. On Tamari interval-posets, this is the number of components in the decreasing forest (in red).

\begin{figure}[ht]
\input{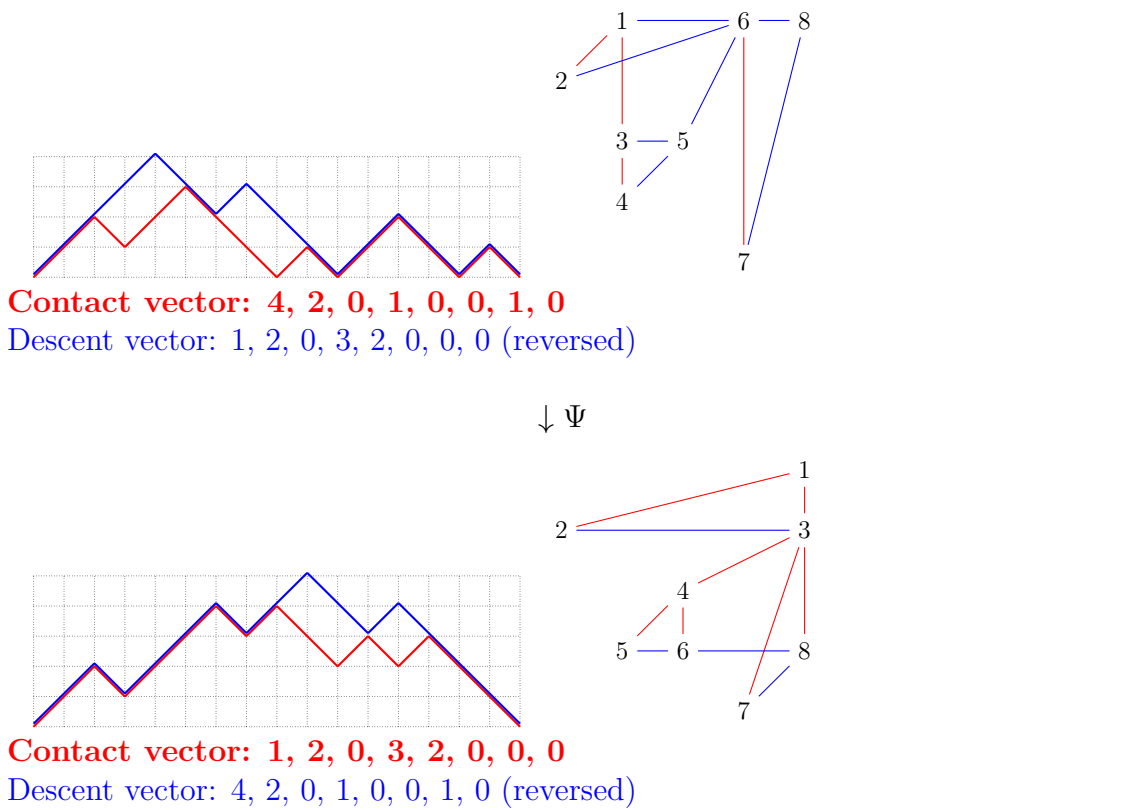}
\caption{The complement involution on Tamari intervals}
\label{fig:tam-complement}
\end{figure}

We can actually do more and construct a \defn{contact vector} $c_0, c_1, \dots c_{n-1}$, that completely characterizes the lower path of the interval. The first value of the vector is the \emph{contact value}. Then note that each up step of the path defines a smaller Dyck path by taking the portion of the path following the up step until it reaches its corresponding down step. For example, on the first interval of Figure~\ref{fig:tam-complement}, the subpath following the first up step of the lower Dyck path is given by $101100$. It has two non-initial contacts which gives $c_1 = 2$. The second up step precedes a peak, so it defines an empty subpath and we have $c_2 = 0$ and so on. 

This is even easier to read on the Tamari interval-poset. As we said, $c_0$ is the number of connected components of the decreasing forest (the graph in red). In the first interval of Figure~\ref{fig:tam-complement}, it is indeed $4$, with the components being $\lbrace 1,2, 3, 4 \rbrace$, $\lbrace 5 \rbrace$, $\lbrace 6, 7 \rbrace$ and $\lbrace 8 \rbrace$. Then $c_i$ is the number of children in the decreasing forest of the vertex $i$. We can see in the example, that $1$ has two children ($2$ and $3$), giving $c_1 = 2$, while $2$ has no children, giving $c_2 = 0$, and so on. 

The Tamari lattice is trivially isomorphic to its revered order. Indeed , if two binary trees are such that $T \wole Q$, it is immediate that their mirror images $T'$ and $Q'$ where the left and right subtrees have been recursively exchanged are such that $Q' \wole T'$. This is a trivial involution on binary trees but not so natural on Dyck paths. On the other hand, it is easy to define on Tamari interval-posets. Sending the interval $[T, Q]$ to $[Q', T']$ corresponds to relabeling the poset such that~$i$ becomes $n+1-i$. This ``switches'' the increasing and decreasing relations of the poset. We show the involution on Figure~\ref{fig:tam-complement}: see that we had $7 \trprec 8$ and $6 \trprec 8$ on the the first interval, which corresponds to the relations $2 \trprec 1$ and $3 \trprec 1$ on its image. We call this involution the \defn{complement} of the interval-poset and write it~$\Psi$.

As the involution transforms the decreasing forest in red into the increasing forest in blue on the second interval, it is clear what the image of the contact vector is. Where we looked at components and children in decreasing forest, we now look at the increasing forest, reading the vertices in decreasing order. This gives another obvious choice for the catalytic parameter: the number of components of the increasing forest, which is the number of nodes on the right branch of the upper tree in the interval. A little observation indicates that this is the length of the final descent in the upper path. More precisely, the complement involution exchanges the contact vector of the lower path with the reversed descent vector of the upper path (the number of down steps following each up step from right to left). See for example that in the second interval of Figure~\ref{fig:tam-complement}, the final up-step of the upper path is followed by $4$ down steps, the previous one by $2$ down steps and so on.

\section{Rises, descents, and left-branch involution}

The surprising result of~\citeext{BMFPR11} is that the \defn{initial rise} of the upper path could serve as a catalytic parameter (and not only the final descent). The initial rise is the number of initial up steps. Of course, reversing the path gives a natural involution on Dyck paths which sends the final descent to the initial rise but it is \emph{not} compatible with the Tamari order. In other words, if $D \wole D'$ in the Tamari lattice, it is not the case in general for their reversed version. Nevertheless, using the representation of Tamari intervals as grafting trees, we find that the equivalent of ``reversing the path'' is actually easy to define and exchanges the descent vector with the rise vector. The general idea is to reverse the upper Dyck path and let the lower Dyck path naturally ``follow''.

\begin{figure}
\input{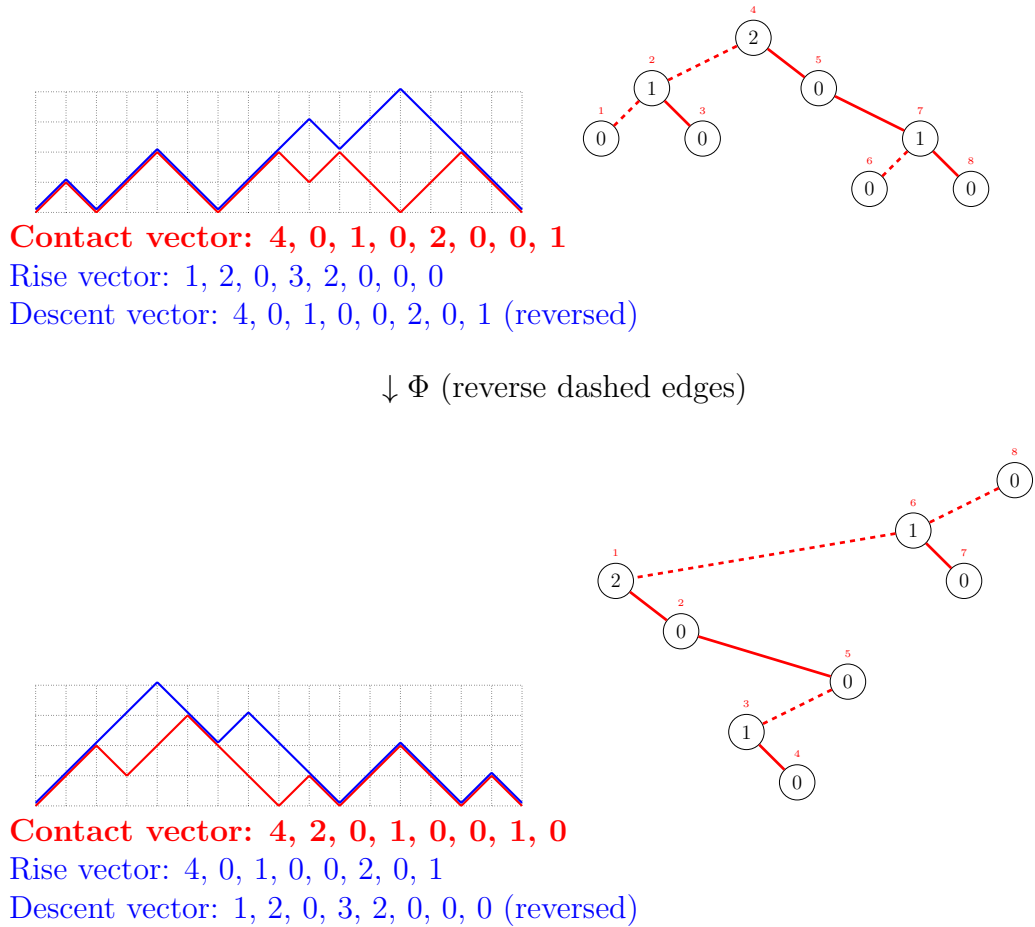}
\caption{The left-branch involution on Tamari intervals}
\label{fig:tam-left-branch}
\end{figure}

Remember that the shape of the grafting tree corresponds to the upper binary tree while the labels give the lower element of the interval. Dyck pack reversal easily translates to binary trees through an involution known as the ``left-branch involution'' which consists of reversing all nodes of every ``left'' branch of the tree. We illustrate this on Figure~\ref{fig:tam-left-branch}: we have written the left branches with dashed lines to make it easier to follow. The labels move along with their nodes and the result is again a grafting tree corresponding to a Tamari interval where the descent and rise vectors have been exchanged. We call this the \defn{left-branch involution} on Tamari intervals and write it $\Phi$.

The natural question is: what happens to the contact vector? To answer, we need to understand how the contact vector can be read directly on the grafting tree. This is actually quite easy. The \emph{contact value} is given by the size of the tree minus the sum of the labels. On the example of Figure~\ref{fig:tam-left-branch}, this is $8 - 4 = 4$. Then, the value $c_i$ is the label on the node $i$. Reading the nodes from left to right on the first tree of Figure~\ref{fig:tam-left-branch}, we obtain indeed $0, 1, 0, 2, 0, 0, 1$. You can check that it corresponds to the contact vector of the lower Dyck path. Following the involution $\Phi$, it is clear that the \emph{contact value} $c_0$ stays the same. The rest of the contact vector is permuted as the nodes have changed their positions but the values themselves stay the same.

\section{Distance and Tamari inversions}

Another interesting statistic appears in the work of Préville-Ratelle~\citeext{PR12}: the \defn{distance}. In Chapter~\ref{chap:qt}, we discuss in more details the algebraic context in which the statistic appears, in relations with $qt$-Catalan numbers. For now, we only give its combinatorial definition and how it can be computed on different objects. The \defn{distance} of a Tamari interval is the length of the \defn{maximal chain} in the lattice between its lower element and upper element. For example, if the interval is reduced to a single element, then the \emph{distance} is $0$, if it corresponds to an edge, the \emph{distance} is $1$. If the interval is a pentagon (with one size of length $3$ and another of length $2$), then the \emph{distance} is $3$. On a Tamari interval, seen as a couple of Dyck paths or binary trees, this statistic is not so easy to compute as one need to find this maximal chain and prove that it is indeed maximal. It turns out that it has a very simple interpretation on Tamari interval-posets as well as on grafting trees. Besides, it also appear in the study of \defn{extended fighting fish} where the \emph{distance} is the \defn{area} of the final surface~\citeext{DH23}.

We call a \defn{Tamari inversion} on an interval-poset, a couple $(a,c)$ with $1 \leq a < c \leq n$ such that: there is no $b$ with $a \leq b < c$ and $c \trprec b$, and no $b$ with $a < b \leq c$ with $a \trprec b$. We show an example on Figure~\ref{fig:tam-inversions}. In particular, $(i,i+1)$ is always a Tamari inversion if you do not have $i \trprec i+1$ nor $i+1 \trprec i$. We prove in~\citeme[Proposition 28]{Pon19} that the number of Tamari inversions give the \emph{distance} of the interval. The intuition is that Tamari inversions correspond to rotations that can be applied one by one from the lower element to the upper element.

\begin{figure}[ht]
\begin{center}
\input{includes/figures/distance-example}
\end{center}

Tamari inversions: $(4, 7), (5, 6), (5, 7)$.

\caption{Tamari inversions on an interval-poset and grafting tree}
\label{fig:tam-inversions}
\end{figure}

The \emph{distance} can also easily be read on the grafting tree. It somehow computes a certain \emph{deficit} between the labels and the maximum value they could take. Remember that a label on a node needs to be smaller than or equal to the size of its right subtree minus the sum of the labels. For example, on the root node $v_4$ of the grafting tree of Figure~\ref{fig:tam-inversions}, we need $\ell(v_4) \leq 4 -1$. The actual label is $2$, so there is a ``deficit'' of $1$. The deficit $d_i$ is then the size of the right subtree of $v_i$ minus the sum of its labels minus the label $\ell(v_i)$. The only other non-zero value besides $d_4$ on Figure~\ref{fig:tam-inversions} is $d_5 = 3 - 1 - 0 = 2$. Then the \emph{distance} is the sum of all deficits. Actually each value $d_i$ gives the number of Tamari inversions $(i,*)$.

Using this characterization, it is easy to see that the two previous involutions $\Psi$ and $\Phi$ preserve the \emph{distance} statistic. For the complement $\Psi$: the definition of Tamari inversions is symmetric towards the increasing and decreasing relations. Each Tamari inversion $(a,c)$ in $I$ corresponds to a Tamari inversion $(n+1-c, n+1-a)$ on $\Psi(I)$. For the \emph{left-branch} involution $\Phi$, we see that even though the shape of tree is changed, the size of the right subtrees and their labels are not. 

\section{The Rise-contact involution}

Using $\Phi$ and $\Psi$, we obtain the following result.

\begin{theorem}[Theorem 54 of \citeme{Pon19}]
\label{thm:rise-contact}
Let $\beta$ be the \defn{rise-contact involution} defined by $\beta := \Phi \circ \Psi \circ \Phi$. Then $\beta$ is an involution on Tamari intervals which exchanges the \emph{contact value} and the \emph{initial rise} (these are the first values of the contact and rise vectors) and the values of the contact vector and rise vector (the order is not kept) while keeping the \emph{distance}. 
\end{theorem}

\begin{figure}[ht]
\center
\input{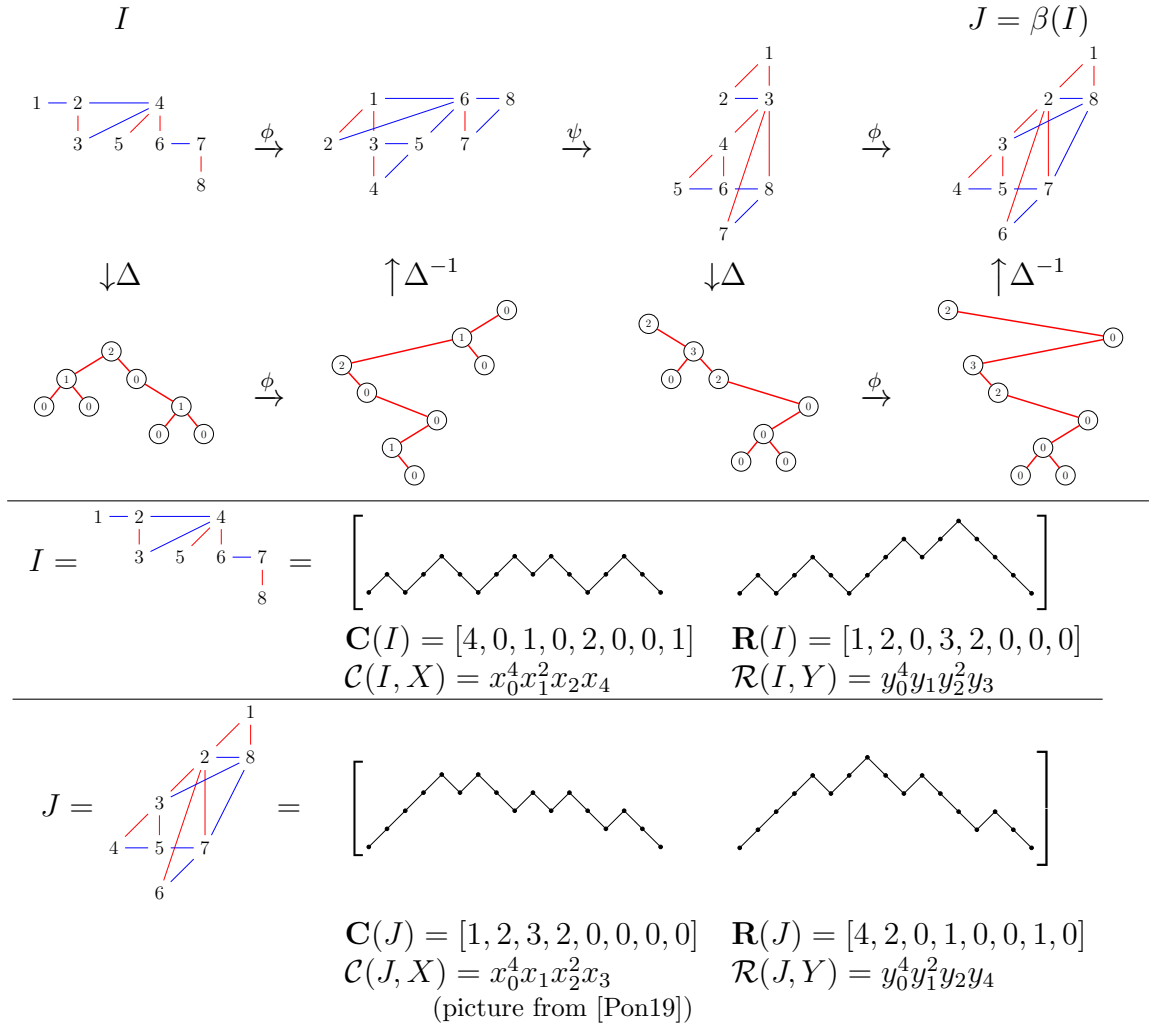}

\footnotesize{(picture from~\citeme{Pon19})}
\caption{The rise-contact involution on Tamari intervals}
\label{fig:rise-contact}
\end{figure}

The \emph{rise-contact} $\beta$ is the conjugate of the $\Psi$ involution by $\Phi$ and so is trivially an involution. This proves \citeext[Conjecture 17]{PR12} for the case $m=1$. We show an example of the complete involution on Figure~\ref{fig:rise-contact}.

\section{The expand-contract operation on $m$-Tamari intervals}

It is possible to define generalizations of the contact and rise vectors on intervals of the $m$-Tamari lattice. In~\citeext[Conjecture 17]{PR12}, Préville-Ratelle conjectures that they are equi-distributed and symmetric (when forgetting the order of the values in the vector). This generalizes the symmetry found in~\citeext{BMFPR11}. This suggests that a rise-contact involution also exists on $m$-Tamari intervals. But the generalization is not direct. Indeed, even though $m$-Tamari intervals can be trivially interpreted as intervals of a Tamari lattice on size $n \times m$, this set is \emph{not} stable by the rise-contact involution we defined previously. The idea behind the generalization is to find \emph{another} interpretation of the $m$-Tamari intervals as a subset of intervals in $n\times m$ stable through the rise-contact involution.

\begin{figure}[ht]
\center
\input{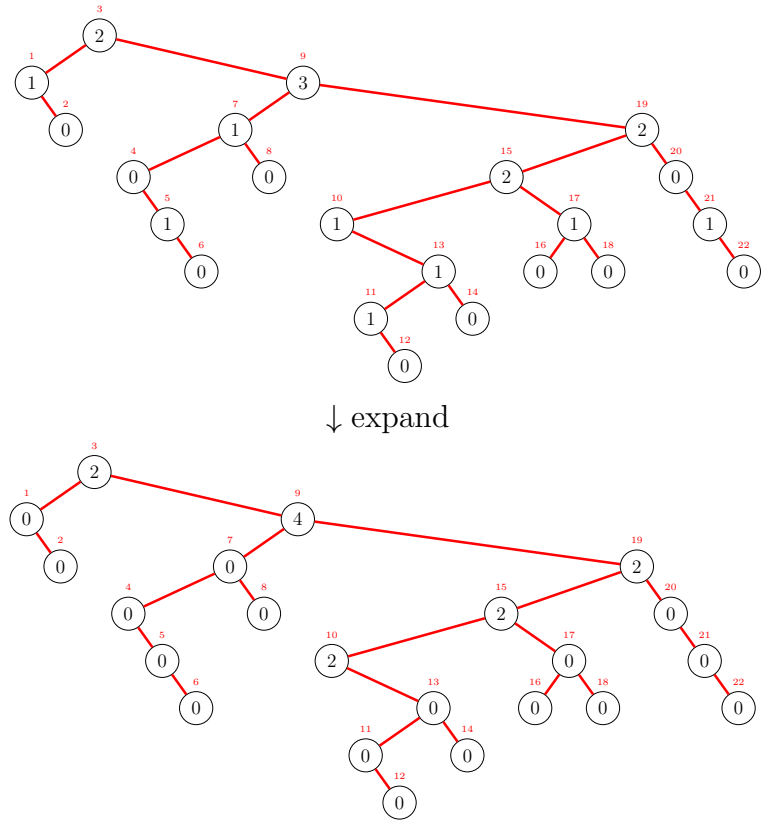}

\caption{The expand operation on $m$-grafting trees}
\label{fig:m-tam-expand}
\end{figure}

To do that, we use the characterization of $m$-grafting trees of Section~\ref{sec:intervals-m-tam}. We have the following issue: the rise vector of an $m$-grafting tree is $m$-divisible (all values are multiples of $m$) but not the contact vector. This dissymmetry comes from the interpretation of $m$-ballot paths as Dyck path of size $n \times m$ by replacing each up step by $m$ up-steps. This actually changes the rise vector but not the contact vector. To solve this, we apply a simple transformation on the labels of $m$-grafting tree: $\ell'(v_i) = m \ell(v_i)$ if $i \equiv 0$ mod $m$, and $\ell'(v_i) = m(\ell(vi) - 1)$ otherwise. This is the \defn{expand} operation. The inverse operation is called the \defn{contract}. We illustrate this on Figure~\ref{fig:m-tam-expand}. In~\citeme{Pon19}, we prove that the \emph{expand / contract} operations define a bijection between $m$-grafting trees and grafting trees of size $n \times m$ whose contact and rise vectors are $m$-divisible. In particular, this new set of intervals is now stable through the classical rise-contact involution $\beta$ which gives the following result.

\begin{theorem}[Theorem 74 of~\citeme{Pon19}]
\label{thm:m-rise-contact}
The \defn{$m$-rise-contact} $\beta_m := \expand \circ \beta \circ \contract$ is an involution on $m$-Tamari intervals exchanging the $m$-contacts and the $m$-rises and proving Conjecture~17 of \citeext{PR12}.
\end{theorem}

\section{Open questions}

In~\citeme{Pon19}, the ideas are more important than the proofs. The definition of the \emph{expand} operation especially is very simple. The proof that it is a bijection (especially that the result is a grafting tree) is only a slightly technical proof by induction. But somehow, we miss some combinatorial insight. What ``are'' the expanded $m$-grafting trees? Can they be used for other purposes? (For example, to obtain a bijective proof of the formula counting intervals). I would also like to study and understand the meaning of these different involutions in the $\nu$-Tamari lattices. It seems possible to define an \emph{expand} operation there too and some sort of rise-contact involution.

\chapter{$q,t$-Catalan} 
\label{chap:qt}

\chapcitation{You must give me leave to judge for myself, and pay me the compliment of believing what I say.}{Jane Austen, Pride and Prejudice}

The enumeration of Tamari intervals is related to \defn{diagonal harmonic polynomials}. This algebraic background may give us a hint on \emph{why} we find such nice formulas as they also appear as the dimensions of certain spaces and can be obtained through symmetric functions. In particular, the study of those polynomials motivated the definition of the $m$-Tamari lattices in~\citeext{BPR12}. In his thesis~\citeext{PR12}, Préville-Ratelle gives a good introduction to the algebraic questions and their relations to combinatorics, stating many problems and conjectures.

The combinatorial aspects of these questions have been an important motivation for my work. In the year 2021--2022, I obtained a \emph{délégation CNRS} (sabbatical) to work in Montreal with François Bergeron in particular. He introduced me to his recent work on \defn{triangular partitions}~\citeext{BM22} which gives an interesting generalization of Dyck paths with many open questions. I spent most of the year working on these questions. In this chapter I present the general background of $q,t$-Catalan combinatorics and their relation to Tamari intervals and triangular partitions. I give an overview of two papers~\citemepre{Pon22} and~\citemeconf{LMP23}. The first one is a short note on a description of the Zeta function (from~\citeext{Hag08}) which I had known for some years and put in writing recently. The second is the result of my collaboration on triangular partitions with Loïc Le Mogne who joined me in Montreal for his Master internship and is now doing a PhD under my supervision.

\section{Background on diagonal harmonic polynomials}

\subsection{The $q,t$-Catalan and the Tamari intervals}

Let us work on a ring of multivariate polynomials in multiple set of variables. This ring is naturally graded by the degrees of the polynomials in the different sets  of variables. The action of the symmetric group $\Sym{n}$ which permutes simultaneously all sets of variable is called the \defn{diagonal action} and is compatible with the degree graduation. Invariant polynomials under this action are called \defn{diagonal polynomials}. The coinvariant spaces are the graded quotient of the ring by the ideal of invariant polynomials. They can be realized as \defn{diagonal harmonic polynomials}. We are interested in the computation of the dimensions of those spaces.

When working with a single set of variables, this is the classical theory of harmonic polynomials. This was studied in the 1950 by Shephard-Todd,
Weyl, and Chevalley in connection with invariant theory and reflection groups~\citeext{Hum90, Bro10}. The extension of the theory to multi sets of variables has been studied in the past years in different fields: combinatorics, but also quantum mechanic, representation theory, algebraic geometry.

In two sets of variables, the dimensions of the coinvariant spaces are given by $(n+1)^{(n-1)}$ which is the number of \defn{parking functions}. They can be understood as some Dyck paths with increasing labels on up steps. The connection with combinatorics is explained in~\citeext{Hag08, Hag20}. It comes from the work of Bergeron, Garsia and Haiman using advanced algebraic tools such as representation theory. The symmetric functions are often used to express the characters of the symmetric groups through the Schur functions. In the study of coinvariant spaces also, symmetric functions play an essential role through a basis of \defn{modified Macdonald polynomials}. This allows to compute the dimensions of the spaces in the form of polynomials in $q$ and $t$ such that the values for $q=t=1$ are the Catalan numbers. These are symmetric, Schur-positive polynomials. 

This last sentence has many repercussions for combinatorics. The $q,t$-Catalan enumeration means that the polynomials could be obtained by sending each combinatorial object of a Catalan family to a monomial $q^it^j$ where $i$ and $j$ are specific statistics on the object. These statistics have been proved to be the \defn{area} and \defn{dinv}. The \emph{area} of a Dyck path is a very natural statistic counting the number of squares below the path. The \emph{dinv} (and its equivalent the \emph{bounce}) is a more mysterious statistic: we give different interpretations in Section~\ref{sec:dinv}. The representation theory behind the definition tells us that the two statistics are symmetric but there is, to this day, no direct combinatorial proof. 

A simple proof would not only be satisfactory for the mind, it could also lead to a better understanding of those polynomials and give better computation techniques. Indeed, as we explain in Section~\ref{sec:triang}, there are generalizations of the $q,t$ symmetry and the polynomials are very difficult to compute. This is where the ``Schur positivity'' comes into play and gives interesting combinatorial questions. It tells us that the $q^it^j$ monomials can be ``put together'' to form Schur polynomials. Each Catalan object belongs to a certain subset corresponding to a specific Schur monomial but we do not know \emph{how} to partition the objects into such sets and which Schur monomials appear in the $q,t$ enumeration of a Catalan number. Obtaining a combinatorial interpretation (or an algorithm) to compute those monomials would be a huge progress.

In~\citeext{BPR12}, Bergeron and Préville-Ratelle work with $3$ sets of variables instead of $2$ and connect the theory of coinvariant spaces with the Tamari lattice. Indeed, now the dimensions of the spaces are given by polynomials in $3$ symmetric variables $q$, $t$, and $r$. The values of the polynomials when $q=t=r=1$ are the number of intervals of the Tamari lattice (and $m$-Tamari lattice for a certain $m$ generalization of the coinvariant spaces). The $3$ variables then correspond to $3$ statistics on Tamari intervals. One of them is the \defn{distance} of the interval which we studied in Chapter~\ref{chap:tam-stat}, it somehow generalizes the notion of \emph{area} of a Dyck path: the \emph{distance} of the interval between the minimal element of the lattice and a Dyck path $P$ is the \emph{area} of $P$. The second statistic is the \emph{dinv} of the maximal element of the interval. The third statistic has no interpretation so far. 

Besides, the polynomials are also Schur positive. Moreover, they correspond to ``the same'' Schur functions as the $2$ variables case. Indeed, Schur functions are formal objects that form a basis of symmetric functions in an infinite number of variables. They are indexed by integer partitions and can be expanded into any number of variables. If the integer partition has $k$ parts, expanding the Schur function in less than $k$ variables gives $0$. The Schur polynomials appearing in the $q,t$-Catalan enumeration are then the expansion in $2$ variables of certain Schur functions $s_\lambda$  where $\lambda$ has, at most, two parts. The Schur polynomials in the $q,t,r$ enumeration of Tamari intervals are these same Schur functions expanded in $3$ variables plus some extra terms $s_\lambda$ where $\lambda$ has $3$ parts.

\subsection{Triangular partitions and Dyck paths}
\label{sec:triang}

\begin{figure}[ht]
\center
\input{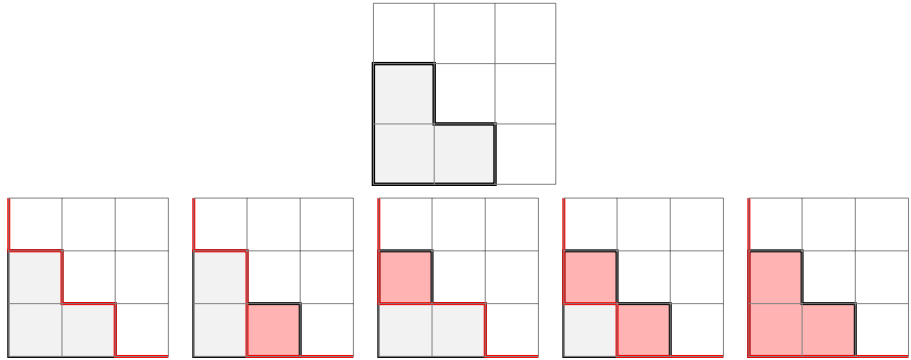}

\caption{Dyck paths of size $3$ as sub-partitions of $(2,1)$}
\label{fig:dyck-paths-partitions}
\end{figure}

An \defn{integer partition} $\lambda$ of size $n$ is a list of positive integers $\lambda_1 \geq \lambda_2 \geq \dots \geq \lambda_k$ with $\sum_i \lambda_i = n$ where $k$ is called the \defn{length} of the partition. For convenience, we often consider partitions ending with an infinite number of $0$ parts ($\lambda_j := 0$ for $j > k$). A \defn{sub-partition} $\mu$ of $\lambda$ is a partition such that $\mu_i \leq \lambda_i$ for all $i \geq 0$. Dyck paths of size $n$ are trivially in bijection with sub-partitions of the \defn{staircase} integer partition $(n,n-1, \dots, 1)$. Indeed, we often represent partitions using their \defn{Ferrers diagram} as on Figure~\ref{fig:dyck-paths-partitions} : each line of the diagram contains $\lambda_i$ boxes. The top partition of Figure~\ref{fig:dyck-paths-partitions} is then $(2,1)$. If $\mu$ is a sub-partition of $\lambda$, its Ferrers diagram is \emph{contained} in the Ferrers diagram of $\lambda$. The frontier between the two partitions draws a \emph{path} below the initial partition (in red on Figure~\ref{fig:dyck-paths-partitions}). When $\lambda$ is the staircase partition, we obtain a Dyck path by a simple rotation. See for example the $5$ sub-partitions of $(2,1)$ corresponding to the $5$ Dyck paths of size $3$. The Dyck paths of size $4$ are the sub-partitions of $(3,2,1)$.

This suggests a natural generalization of Dyck paths: changing the partition~$\lambda$. This includes in particular the case of \defn{rational Dyck paths}. Let $a$ and $b$ be two co-prime integers. Bizley~\citeext{Biz54} gives a nice enumeration for paths made of north and east steps, above the line $y = \frac{a}{b} x$ and ending in $(a,b)$,

\begin{equation}
\frac{1}{a+b} \binom{a+b}{a,b}.
\end{equation}
These are called the \defn{rational Catalan numbers}. Indeed, for $a = n$ and $b = n+1$, we recover the Catalan numbers. In our settings, this is counting the partitions whose Ferrers diagrams fit below the line between $(0,b)$ and $(a,0)$. They are all sub-partitions of the maximal fitting partition which we call a \defn{rational partition}. Combinatorial properties of rational Dyck paths are known for being related to representation theory and coinvariant spaces (See for example~\citeext{ALW16, ALW15, TW18}). In particular, there exists a $q,t$ enumeration that generalizes the $q,t$-Catalan case with similar open questions. In~\citeext{Ber17}, Bergeron summarizes some of these questions (both algebraic and combinatorial) in the \defn{rectangular} case (where the co-prime condition for $a$ and $b$ is dropped). 

\begin{figure}[ht]
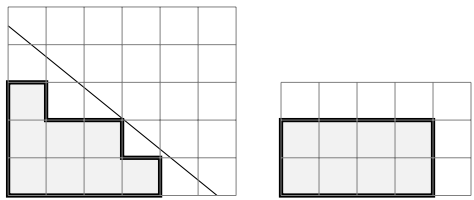

\center
\begin{tabular}{cc}
\input{includes/figures/triangular_partition}
&
\input{includes/figures/non_triangular}
\end{tabular}

\caption{The triangular partition $(4,3,1)$ and the non-triangular partition~$(4,4)$.}
\label{fig:triang}
\end{figure}

The work of~\citeext{BHMPS23} suggests that the connection between combinatorial enumeration, representation theory and symmetric functions exists beyond the rectangular and rational cases: by taking paths below any line. This motivated the work of Bergeron and Mazin on \defn{triangular partitions}~\citeext{BM22}. A integer partition is said to be \defn{triangular} if it is the maximal partition fitting below the line joining $(r,0)$ and $(0,s)$ for $r$ and $s$ non-negative real numbers. In particular, $r$ and $s$ are not required to be integers. We show an example on Figure~\ref{fig:triang} next to a partition which is not triangular.

In~\citeext{GHSR20}, the authors use a formula from Gorsky and Negut~\citeext{GN15} to obtain a $q,t$-enumeration of the sub-partitions of any partition $\lambda$, thus generalizing the $q,t$-Catalan numbers. But the signification of each $q,t$ monomial is not clear and coefficients are not always positive.  Conjecture 7.1.1 from \citeext{BHMPS23} states that the coefficients might be positive if the partition lies under a certain convex curve. In the case where the partition is triangular, the authors obtain an explicit formula using symmetric functions and prove a combinatorial interpretation of the statistics corresponding to $q$ and $t$ (generalizing the \emph{area} and \emph{dinv}). In~\citeext[Conjecture 1]{BM22}, Bergeron and Mazin further conjecture that this $q,t$-enumeration is Schur-positive.

\section{Understanding the ``dinv''}
\label{sec:dinv}

\subsection{Zeta function and area sequences}

In the $q,t$-Catalan enumeration, the combinatorial interpretation for the $q$ exponent comes naturally: it is the \emph{area} of the Dyck path. If we see Dyck paths as a sub-partition $\mu$ of a partition $\lambda$, this is the number of boxes which are part of~$\lambda$ but not of $\mu$: the red boxes on Figure~\ref{fig:dyck-paths-partitions} and~\ref{fig:triangular-area}. Even in the classical case, the interpretation of the second statistic is more difficult. In~\citeext{Hag03}, Haglund gives a first conjectural interpretation of the exponent of $t$ as the \defn{bounce} statistic. At the same period, Haiman finds another interpretation that we now call the \defn{dinv}. This leads Haglund to come up with a bijection from Dyck paths of size~$n$ to themselves sending the \emph{area} and \emph{bounce} to the \emph{dinv} and \emph{area} respectively. This is the \defn{zeta map} $\zeta$~\citeext[page 50]{Hag08}. In particular, it proves that the \emph{area}, the \emph{dinv}, and the \emph{bounce} are equi-distributed (but not the $q,t$-symmetry). As $q,t$-enumerations arose in more general settings, some generalization of the $\zeta$ map also appeared. In~\citeext{ALW15}, the author describe a general $\zeta$ map for all rational Dyck paths. They conjecture that it is bijective but it appears a difficult question. It was eventually solved in~\citeext{TW18}.

\begin{figure}[ht]
\center
\scalebox{.8}{\input{includes/figures/area_sequence_dyckpath}}

Area sequence: $01211123301101221$

Area: $22$

\caption{A Dyck path (as a sub-partition) with its area and area sequence.}
\label{fig:triangular-area}
\end{figure}

I discovered the $q,t$-enumeration of Catalan number at some time during my thesis. I found the problem of finding a combinatorial proof of the symmetry fascinating and tried to solve it. I did not succeed but I came up with a bijective map that I soon discovered was the $\zeta$ map (more precisely $\zeta^{-1}$). As it was not a major result, I did not care to write it at the time. I finally did last year in~\citemepre{Pon22} after discussing it with Hugh Thomas as he had actually never seen this description anywhere. As I explain below, it is based on the \defn{area sequence} of Dyck paths. We read the sequence value by value and make an \emph{insertion} in a new area sequence that we build iteratively. The process is very simple and can be easily implemented. I actually provide the full implementation with the paper in~\citemesoft{PonSage22}.

The \defn{area sequence} of a Dyck path is a very common way to represent Dyck paths and to compute the area. It consists of ``counting the boxes'': when represented as in~\ref{fig:triangular-area}, we read the number of red boxes, line by line, from top to bottom. This gives a bijection between Dyck paths of size $n$ and sequences of integers $a_1, \dots, a_n$ where $a_1 = 0$ and $0 \leq a_{i+1} \leq a_i +1$. The \emph{area} is of course the sum of the $a_i$. The \emph{dinv} is also easy to compute from the area sequence: it is the number of couples $(a_i, a_j)$ with $i < j$ and $a_j = a_i$ or $a_j = a_i - 1$. On Figure~\ref{fig:triangular-area}, the \emph{dinv} is 65.

I define an operation called the \defn{insertion} on the area sequence. If $a_1 \dots a_n$ is an area sequence, then inserting in position $i$ gives the new sequence \linebreak ${a_1 \dots a_i (a_i+1) a_{i+1} \dots a_n}$. By convention, inserting in position $0$ adds a $0$ at the beginning of the sequence. I then explain how we can control the increase of the \emph{dinv} through the insertion and I define specific \defn{insertion points} such that if the last inserted value increased the \emph{dinv} by $v_i$, then the new sequence has $v_i + 1$ insertion points corresponding to potential increases of the \emph{dinv} of $0, 1, \dots v_i +1$. From an area sequence $a_1 \dots a_n$, we can then construct a new sequence by inserting at the insertion point in the image which increases the \emph{dinv} by exactly $a_i$. 

\begin{theorem}[Theorem~6 of~\citemepre{Pon22}]
The insertion process described previously defines a bijection on area sequences such that the \emph{dinv} of the image is the \emph{area} of the sequence.
\end{theorem}

We show an an example on Figure~\ref{fig:zeta-bij} where insertion points are made visible. We prove that this map is actually $\zeta^{-1}$. In particular, we give a direct proof that the \emph{area} of the new sequence is the \emph{bounce} of the original one.

\begin{figure}[ht]
\center
\begin{tabular}{ll}
$0$        & $\cdot_1 0 \cdot_0$ \\
\hline
$01$        & $\cdot_2 0 \cdot_1 0 \cdot_0$ \\
 \hline
$012$        & $ \cdot_3 0 \cdot_2 0 \cdot_1 0 \cdot_0$ \\
\hline
$0121$        & $0 ~ 0 \cdot_2 1 \cdot_0 0 \cdot_1$ \\
\hline
$01211$        & $0 ~ 0 ~ 1 \cdot_1 0 \cdot_2 1 \cdot_0$ \\
\hline
      $012111$   & $0 ~ 0 ~ 1 \cdot_2 2 \cdot_0 0 ~ 1 \cdot_1$ \\
\hline
      $0121112$ & $0 ~ 0 ~ 1 \cdot_3 2 \cdot_1 2 \cdot_0 0 ~ 1 \cdot_2$ \\
\hline
      $01211123$ & $0 ~ 0 ~ 1 \cdot_4 2 \cdot_2 2 \cdot_1 2 \cdot_0 0 ~ 1 \cdot_3$ \\
\hline
      $012111233$ & $0 ~ 0 ~ 1 ~ 2 \cdot_3 2 \cdot_2 2 \cdot_1 0 ~ 1 \cdot_4 2 \cdot_0$ \\
\hline
      $0121112330$ & $0 ~ 0 ~ 1 ~ 2 ~ 2 ~ 2 ~ 0 ~ 1 ~ 2 \cdot_1 3 \cdot_0$ \\
\hline
      $01211123301$ & $0 ~ 0 ~ 1 ~ 2 ~ 2 ~ 2 ~ 0 ~ 1 ~ 2 \cdot_2 3 \cdot_1 3 \cdot_0$ \\
\hline
      $012111233011$ & $0 ~ 0 ~ 1 ~ 2 ~ 2 ~ 2 ~ 0 ~ 1 ~ 2 ~ 3 \cdot_2 4 \cdot_0 3 \cdot_1$ \\
\hline
      $0121112330110$ & $0 ~ 0 ~ 1 ~ 2 ~ 2 ~ 2 ~ 0 ~ 1 ~ 2 ~ 3 ~ 4 \cdot_1 5 \cdot_0 3$ \\
 \hline
      $01211123301101$ & $0 ~ 0 ~ 1 ~ 2 ~ 2 ~ 2 ~ 0 ~ 1 ~ 2 ~ 3 ~ 4 \cdot_2 5 \cdot_1 5 \cdot_0 3$ \\
\hline
      $012111233011012$ & $0 ~ 0 ~ 1 ~ 2 ~ 2 ~ 2 ~ 0 ~ 1 ~ 2 ~ 3 ~ 4 \cdot_3 5 \cdot_2 5 \cdot_1 5 \cdot_0 3$ \\
\hline
      $0121112330110122$ & $0 ~ 0 ~ 1 ~ 2 ~ 2 ~ 2 ~ 0 ~ 1 ~ 2 ~ 3 ~ 4 ~ 5 \cdot_3 6 \cdot_0 5 \cdot_2 5 \cdot_1 3$ \\
\hline
      $01211123301101221$ & $0 ~ 0 ~ 1 ~ 2 ~ 2 ~ 2 ~ 0 ~ 1 ~ 2 ~ 3 ~ 4 ~ 5 ~ 6 \cdot_1 5 ~ 5 \cdot_2 6 \cdot_0 3$ \\
\hline
\end{tabular}

\caption{An example of the inverse zeta map bijection}
\label{fig:zeta-bij}
\end{figure}

\subsection{The deficit in triangular Dyck paths}
\label{sec:deficit}

\begin{figure}[ht]
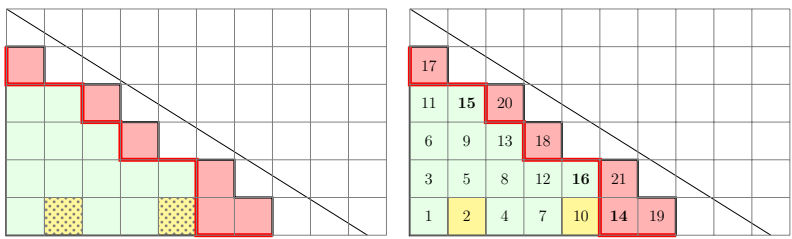

\center

\scalebox{.5}{
\begin{tabular}{cc}
\input{includes/figures/sim}
&
\input{includes/figures/deficit}
\end{tabular}
}

\caption{Similar cells, Triangular young tableau and deficit}
\label{fig:sim-deficit}
\end{figure}

In~\citeext{BHMPS23}, the authors generalize the \emph{dinv} statistic to all triangular Dyck paths by counting the \defn{$p$-balanced hooks} in the Ferrers diagram of the sub-partition. This is reformulated in~\citeext[Section 4.1]{BM22} and given a new name: the \defn{sim} statistic. Remember that a triangular partition is the maximal partition lying under a given line. Actually, for any given triangular partition, there is an infinite number of lines to choose from. The authors of~\citeext{BM22} provide a very simple criteria for deciding if a a given partition is triangular using the \defn{hooks} of the different cells. They also define a certain notion of \defn{mean slope} as an explicit line choice. A cell of a sub-partition is said to be \defn{similar} if its hook is ``compatible'' in a certain sens with the mean slope of the triangular partition. The number of similar cells is called the \defn{sim}. It corresponds to the known \emph{dinv} statistic on regular and rational Dyck paths. We show an example on the left of Figure~\ref{fig:sim-deficit}: the similar cells are colored in green while the non similar cells are in yellow with dotted pattern.

In~\citemeconf{LMP23}, we present a new combinatorial interpretation of the similar cells using a new statistic called the \defn{deficit}. The cells of the Ferrers diagram of a partition $\lambda$ can be filled with numbers to form a \defn{tableau}. It is called a \defn{standard Young tableau} if it contains exactly the numbers $1, \dots, n$ such that the numbers are increasing on each line and each column. Let $\theta$ be a standard Young tableau of shape $\lambda$ where $\lambda$ is a triangular partition. We say that there is a $\theta$-inversion in a sub-partition $\mu$ of $\lambda$ if there is a cell $c$ inside $\mu$ with a higher value in $\theta$ than another cell $c'$ outside of $\mu$. For example, on the right of Figure~\ref{fig:sim-deficit}, we see that~$15$ is inside the sub-partition while $14$ is outside. Then the cell which lies at the ``hook'' of $c$ and $c'$ is said to be a \defn{deficit cell}. The number of deficit cells gives the \emph{deficit} of the triangular Dyck path with tableau $\theta$. 

In the example of Figure~\ref{fig:sim-deficit}, we see two deficit cells: the cell number $2$ which is at the hook of $15$ and $14$, and the cell number $10$ which is at the hook of $16$ and~$14$. These are the two cells which were identified as \emph{non similar} cells on the left. Indeed, we prove that on a certain tableau that we call the \defn{triangular tableau}, the deficit cells are exactly the non similar cells. This tableau is constructed using the mean slope defined by Bergeron and Mazin (or more precisely, a slight irrational deformation): we move the slope towards the origin and number the cells in decreasing order as we touch them. Using an irrational slope ensures that slope touches the cells one by one. We show an example on Figure~\ref{fig:triangular-tableau}.

\begin{figure}[ht]
\center

\begin{tabular}{cccc}
\input{includes/figures/triangular_tableau}
&
\input{includes/figures/triangular_tableau2}
&
\input{includes/figures/triangular_tableau3}
&
\input{includes/figures/triangular_tableau4}
\end{tabular}

\footnotesize{(pictures from~\citemeconf{LMP23})}
\caption{Construction of the triangular tableau}
\label{fig:triangular-tableau}
\end{figure}

This interpretation has been very useful and it also opens new horizons. Indeed, we might wonder what are the $q,t$-enumerations of other standard Young tableau of shape $\lambda$, using the tableau deficit to compute the power of $t$. Most of the time, it is not symmetric but it actually defines a certain family of tableaux, which we call the \defn{similar-symmmetric} tableaux, or \defn{sim-sym}, on which we recover the original $q,t$-enumeration. Characterizing those tableaux seems to be a difficult question: we find some of them as ``slope tableaux'' using a similar construction as for the mean slope, but not all of them. We were nevertheless able to characterize them in the case of triangular partitions of length $2$. 

\section{Lattice interval enumeration and Schur-po\-si\-tivity}
\label{sec:triangular-intervals}

We have briefly explained the relations between the $q,t$-enumeration of Dyck paths, the symmetric functions and representation theory. In the case of triangular partitions, these relations are still conjectural (especially concerning representation theory). When I met him in Montreal, François Bergeron was working on this question. In particular, he was able to compute, or partially compute, a certain Schur-positive symmetric function for each triangular partition. The expansion in two variables of those functions are the $q,t$ enumeration of triangular Dyck paths. The question he asked me was: what happens in $3$ variables? For classical Dyck paths and $m$-Dyck paths, we know we obtain an enumeration on Tamari and $m$-Tamari intervals but what happen on the other cases? 

The most natural choice is to look at the intervals of the $\nu$-Tamari lattice that we mentioned in Section~\ref{sec:nu-tam}. In certain cases, it gives the wanted enumeration but not in all cases. For some partitions, we have the right number of intervals but loose the $q,t$-symmetry or the Schur-positivity. In other cases, the number of intervals does not correspond to the expansion of the symmetric function. I spent lots of time trying to find a lattice structure that would give the proper $q,t$-interval enumeration. Even though I was able to compute some examples (I actually explain that in Chapter~\ref{chap:exp}), I could not find a proper general definition. 

Around this time, I was joined by Loïc Le Mogne who started his internship and was interested in working on these questions. I asked him to look at the case of triangular partitions of length $2$. Indeed, this case is much easier than the general case but still largely non-trivial. As we came up with the definition of \emph{sim-sym} tableaux, it appeared to us that a certain maximal chain in the $\nu$-Tamari lattice defines a Young tableau, which we call the \defn{top-down tableau}, which \emph{sometimes} is \emph{sim-sym}. In particular, in the length $2$ case, the top-down tableau is always \defn{sim-sym} but does not necessarily corresponds to the triangular tableau. We then state the following theorem.

\begin{theorem}[Theorem~5 of \citemeconf{LMP23}]
\label{thm:qt-lattice-2triangular}
The $q,t$-enumeration of intervals of the $\nu$-Tamari lattice on triangular $2$-partitions is symmetric and Schur-positive.
\end{theorem}

The statistic corresponding to $q$ is the distance while the statistic corresponding to $t$ is the number of similar cells using the top-down tableau. The proof of the theorem has been obtained by Loïc Le Mogne based on case-by-case analysis using the structure of the lattice. It will appear in the long version of the paper.

Using this result and some extra computation, we express a new conjecture on $q,t$-enumerations of intervals.

\begin{conjecture}[Conjecture~4 of~\citemeconf{LMP23}]
\label{conj:top-down-sim-sym}
Let $\lambda$ be a triangular partition such that the top-down tableau of $\lambda$ is \emph{sim-sym}, then the $q,t$-enumeration of $\nu$-Tamari intervals is symmetric and Schur-positive.
\end{conjecture}

\section{Open questions}

We have seen already that the general field of $q,t$-Catalan combinatorics has \emph{many} open questions as we can read in~\citeext{PR12} and~\citeext{Ber17}. Some of these questions are very difficult but they often leave a trace of more approachable problems. The study of triangular partitions and our own results are still very new and we have many leads to investigate:

\begin{itemize}
\item We would like to better understand sim-sym tableaux. In particular, there should be a bijection between elements of similar \emph{area} and \emph{sim} using different tableaux. 
\item The question of finding the ``good'' lattice structure is still open. Even if the general case is difficult, we have had new ideas on how to approach it. One would be to further study the length $2$ case. Indeed, we might be able to define a different lattice for each sim-sym tableau and prove that they give a similar enumeration of intervals.
\item I did not have the time yet to investigate the $\zeta$ function on triangular partitions. I do not know if the general definition of sweep maps from~\citeext{ALW15} which was later proven to be bijective in~\citeext{TW18} works for triangular partitions (it is said to work on the rectangular and rational cases). In any case, it needs specific investigation to understand how it is or could be defined and what it tells us on the combinatorics of the triangular case.
\item There are some notions of \emph{contacts} and \emph{rises} on triangular Dyck paths, in relation in particular with \defn{triangular parking functions}. We need to further investigate the expressions that arise in the symmetric functions related to triangular partitions and look for proper combinatorial interpretation. In particular, we need to look at classical bijections and correspondences of the classical case and see how they generalize.
\item Similarly, we need to investigate how to get ``triangular versions'' of some classical Catalan objects, in particular objects that are in relations with Tamari intervals as in Figure~\ref{fig:tam-bijections}.
\end{itemize}

\part{Quotients and Sublattices}
\label{part:quotient}

\chapter{Integer Posets}
\label{chap:integer-posets}

\chapcitation{Definitions belong to the definers, not the defined.}{Toni Morrison, Beloved.}

As I presented in the previous chapters, Tamari interval-posets were the result of a collaboration with Grégory Châtel. This started around 2012 while I was finishing my thesis and Gréory Châtel was a fellow student. Around 2013, we started working with Vincent Pilaud to understand the connection between our work on Tamari intervals and Cambrian lattices~\citeext{Rea06}. We wanted to answer the following questions: can the intervals of the Cambrian lattices be seen as interval-posets? Can we use this characterization to count the intervals? Is it possible to define other algebraic structures, such as Hopf algebras, on intervals? The short answers are: yes, no, yes. In particular, there is a specificity of the Tamari lattice for counting intervals. Nevertheless, as we started exploring those questions, it appeared that many combinatorial objects that were considered very different in nature could be seen as \defn{integer poset}. It became interesting to study this specific family and how it related to previously known algebraic structures, especially lattices and Hopf algebras. This collaboration resulted in $3$ papers~\citeme{PP18,CPP19}, and~\citeme{PP20}. 

This chapter is dedicated to~\citeme{CPP19} where we explore the lattice structure of integer posets, and~\citeme{PP20} which is about Hopf algebraic aspects. The results from both papers were originally presented in a FPSAC conference poster~\citemeconf{PP18b}. The paper~\citeme{PP18} describes some family of objects called the \defn{permutrees}, we present it in details in Chapter~\ref{chap:permutrees}.  Even though it was written before the two other papers, the original motivation came from the study of integer posets. We start here by explaining how many objects can be described in terms of integer posets. We then give the main ideas to construct the lattice and Hopf algebra of integer posets through the notion of integer relations.

\section{Generalizing the Tamari interval-posets}
\label{sec:rel-objects}

We have defined Tamari interval-posets as posets on integers satisfying certain relations. Besides, the linear extensions of a Tamari interval-poset give a certain interval of the weak order. Actually, any interval of the weak order can be encoded as the linear extensions of a certain poset on integer. Tamari interval-posets are just a special case. 

\begin{figure}[ht]
\center
\input{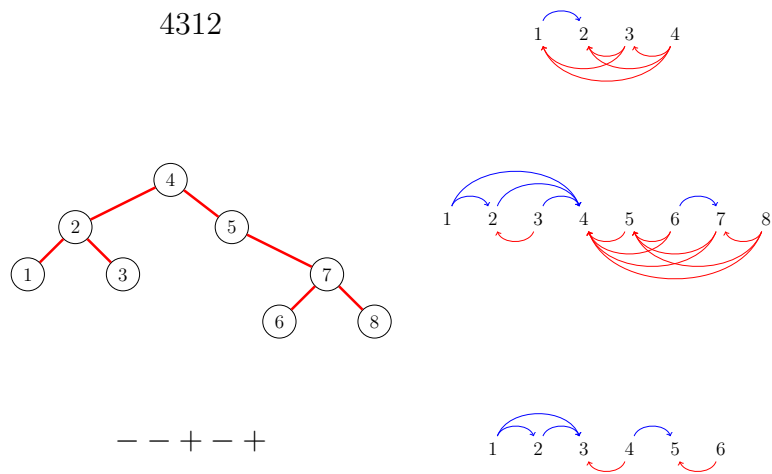}

\caption{Permutations, binary trees, and binary sequences as integer posets}
\label{fig:rel_objects}
\end{figure}

Indeed, let us consider permutations as a linear order on integers: they correspond to a special case of partial orders on integers where the order is total. Similarly, standard binary search trees define a poset by orienting all edges towards the root. The linear extensions are the sylvester class of the tree. This is also the case of binary sequences: we interpret any $-$ sign by an increasing relation and any $+$ sign by a decreasing relation. In this case, linear extensions also give a certain interval of the weak order corresponding to the \defn{recoil class} we mentioned in Section~\ref{sec:asso-cube}. We show examples of these constructions on Figure~\ref{fig:rel_objects}. In this chapter, we always represent integer posets by writing the integers in their natural order on a line and writing \emph{all} relations (not just cover relations) by dividing them intro \defn{increasing relations} (in blue, from left to right, on top) and \defn{decreasing relations} (in red, from right to left, on bottom).

\begin{figure}[ht]
\center
\input{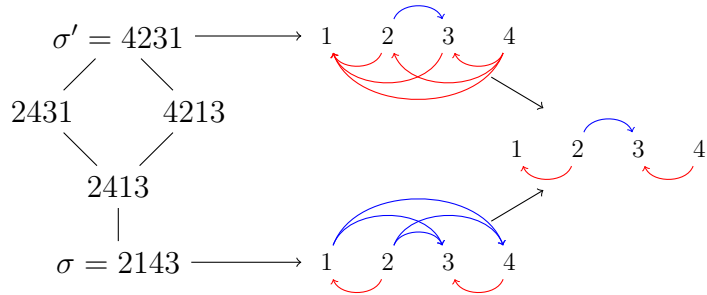}

\caption{Construction of a weak order interval-poset}
\label{fig:woip}
\end{figure}

Now we can construct interval-posets of the corresponding lattices (the weak order, the Tamari lattice, and the boolean lattice). For each interval $e \wole f$, we take the decreasing relations of the smaller element $e$ and the increasing relations of the larger element $f$ as illustrated on Figure~\ref{fig:woip} for permutations. This always give a certain poset on integers and we characterize the obtained families depending on the considered lattice. The conditions that need to be satisfied are always local conditions on triplets $a < b < c$. Basically, the increasing (resp. decreasing) relation $a \trprec c$ (resp. $c \trprec a$) implies \emph{something} on the increasing (resp. decreasing) relations $a \trprec b$ and $b \trprec c$ (resp. $c \trprec b$ and $b \trprec a$). We illustrate the different conditions in Figure~\ref{fig:interval-conditions}.

\begin{figure}[ht]
\center
\input{includes/figures/interval-conditions}

\caption{Conditions on interval-posets}
\center 

\label{fig:interval-conditions}
\end{figure}

We then obtain three families of integer posets representing three families of intervals in the weak order, the Tamari lattice and the boolean lattice respectively. In~\citeme{CPP19}, we call them $\WOIP$ for \defn{weak order interval-posets}, $\TOIP$ for \defn{Tamari order interval-posets} and $\BOIP$ for \defn{boolean order interval-poset}. The elements of those latices (permutations, binary trees, and boolean sequences) are the \defn{maximal posets} (posets with the maximal number of relations) that we can construct satisfying the interval conditions (we call them $\WOEP$, $\TOEP$, and $\BOEP$ where the E stands for \defn{element}). In particular, you can always recover the element from the increasing (resp. decreasing) relations. For a permutation, this means recovering the permutation from the inversions. For a binary tree, this corresponds to the bijection between the tree and its decreasing / increasing forest. Besides, we can also characterize \defn{faces} of the corresponding polytopes as special kind of intervals: the $\WOFP$~\citeme[Proposition 2.10]{CPP19} are the interval-posets corresponding to the faces of the permutahedron, the $\TOFP$~\citeme[Proposition 2.10]{CPP19} correspond to faces of the associahedron, while all $\BOIP$ are actually faces of the cube. 

The nice enumeration of Tamari intervals come from the fact that the increasing and decreasing relations form a forest of rooted trees. It is not the case for the weak order especially and the characterization does not lead to interesting decompositions in terms of enumeration. Nevertheless, all this objects are known to appear in different algebraic structures: lattices and Hopf algebras. This is what motivates the following sections which proposes a general framework to define these structures.

\section{Integer Relations}

\subsection{As a lattice}

Let $S$ be a set of size $n$ and $\rel$ a \defn{binary relation} on $S$, \emph{i.e.}, a set of couples $(x,y)$ with $x \in S$ and $y \in S$. We write $x \rel y$ if $x$ is in relation with $y$ and $x \notrel y$ if not. We consider only \defn{reflexive} relations, \emph{i.e.}, $x \rel x$ for all $x \in S$. But $\rel$ is in general \emph{not} symmetric so $x \rel y$ does not imply $y \rel x$. The \defn{boolean lattice} is a well known structure on binary relations where the order is given by inclusion of the relations. The minimal element is the empty relation, while the maximal is the relation $\rel$ where $x \rel y$ for all $x$ and $y$ in $S$. The meet and join are given by taking respectively the intersection and union of relations.

We now consider the case where $S = 1,\dots, n$ and we say that $\rel$ is an \defn{integer relation}. In this case, the set $S$ is actually endowed with $2$ relations: the relation $\rel$ and the natural order relation between integers. This allows for a new lattice definition which we call the \defn{weak order on integer relations} as it is indeed a generalization of the weak order on permutations as we explain in Section~\ref{sec:rel-induced}. 

We see the integer relation $\rel$ as the union of two subrelations : the increasing relation $\Inc{\rel}$, \emph{i.e.}, the couples $a \rel b$ with $a < b$ and decreasing relation $\Dec{\rel}$, \emph{i.e.}, the couples $b \rel a$ with $a < b$. We say that a relation $\rel$ is smaller than or equal to a relation $\rel[S]$, $\rel \wole \rel[S]$, if and only if $\Dec{\rel} \subseteq \Dec{\rel[S]}$ and $\Inc{\rel[S]} \subseteq \Inc{\rel}$ (Definition 1.1 of~\citeme{CPP19}).

\begin{figure}[ht]
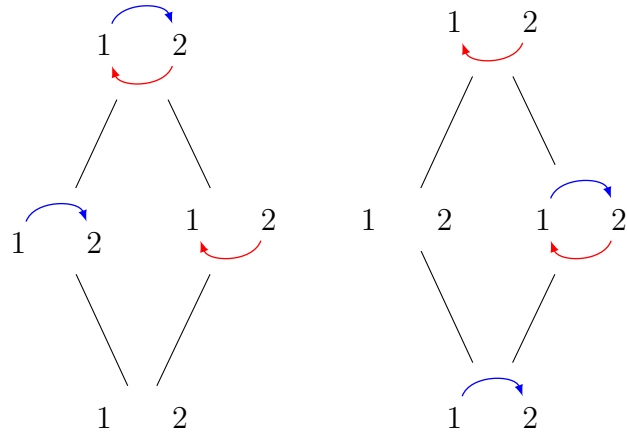

\center
\begin{tabular}{cc}
\input{includes/figures/rel2_boolean}
&
\input{includes/figures/rel2_weak_order}
\end{tabular}
\caption{The boolean lattice (on the left) and the weak order (on the right) on integer relations.}
\label{fig:rel2_wo}
\end{figure}

We present the example for relations of size 2 on Figure~\ref{fig:rel2_wo}. We show both the classical boolean lattice on the left and the weak order on the right. On the boolean lattice, the minimal element is the empty relation and at each step, one relation is added. On the weak order, the minimal element is the relation containing all possible increasing relations and at each step, either one increasing relation is removed or one decreasing is added. In particular, the weak order on integer relations is isomorphic to the classical boolean lattice. Indeed, it is a product of the reversed boolean lattice on increasing relations by the boolean lattice on decreasing relations. It then forms a cube of dimension $n(n-1)$.

Note that we always represent relations as in Figure~\ref{fig:rel2_wo}. Integers $1$ to $n$ are written on a single line. Each increasing relation $a \rel b$ is shown as a blue arc (from left to right) over the vertices while each decreasing relation $b \rel a$ is a red arc (from right ro left) under the vertices.

As it is isomorphic to the boolean lattice, the weak order on integer relation is a lattice. We can easily compute the join in the meet as shown on Figure~\ref{fig:rel2_wo}: the meet is obtained by taking the union of increasing relations and intersection of decreasing relations. The join is the symmetric operation.

\begin{figure}[ht]
\center
\input{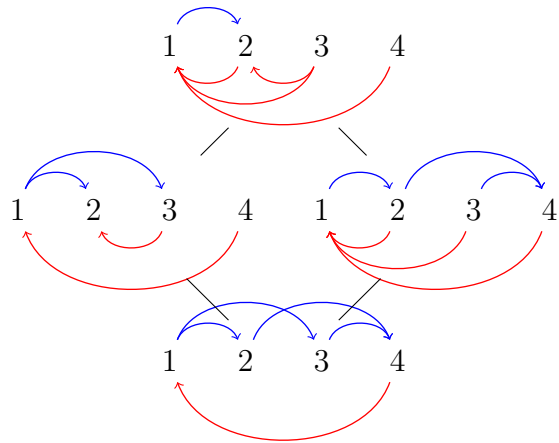}

\caption{Join and meet on integer relations}
\label{rel_inf_sup}
\end{figure}

\subsection{As a Hopf algebra}

Let $\IRel$ be the set of all integer relations, $\IRel_n$ the set of integer relations of size $n$, and $\K$ a field. Then we can define a Hopf algebra on the vector space $\K\IRel := \bigoplus_{n \ge 0} \K\IRel_n$.

The product (on the fundamental basis $\FRel$) is given by what we call the \defn{shifted shuffle} of two relations $\rel$ of size $n$ and $\rel[S]$ of size $m$~\citeme[Definition 3]{PP20}. This is a set containing all relations $T$ such that the subrelation of $T$ on $1, \dots, n$ is $R$ and the subrelation on $n+1, \dots, n+m$ is the shifted version of $\rel[S]$, written $\overline{\rel[S]}$. In other words, an relation $T$ of the shifted shuffle is obtained by concatenating $\rel$ and $\overline{\rel[S]}$ and adding some increasing relations from $\rel$ to $\overline{\rel[S]}$ and some decreasing relations from $\overline{\rel[S]}$ to $\rel$. In particular, the shifted shuffle has cardinality $2^{2mn}$. 

\begin{figure}
\center
\includegraphics[scale=.6]{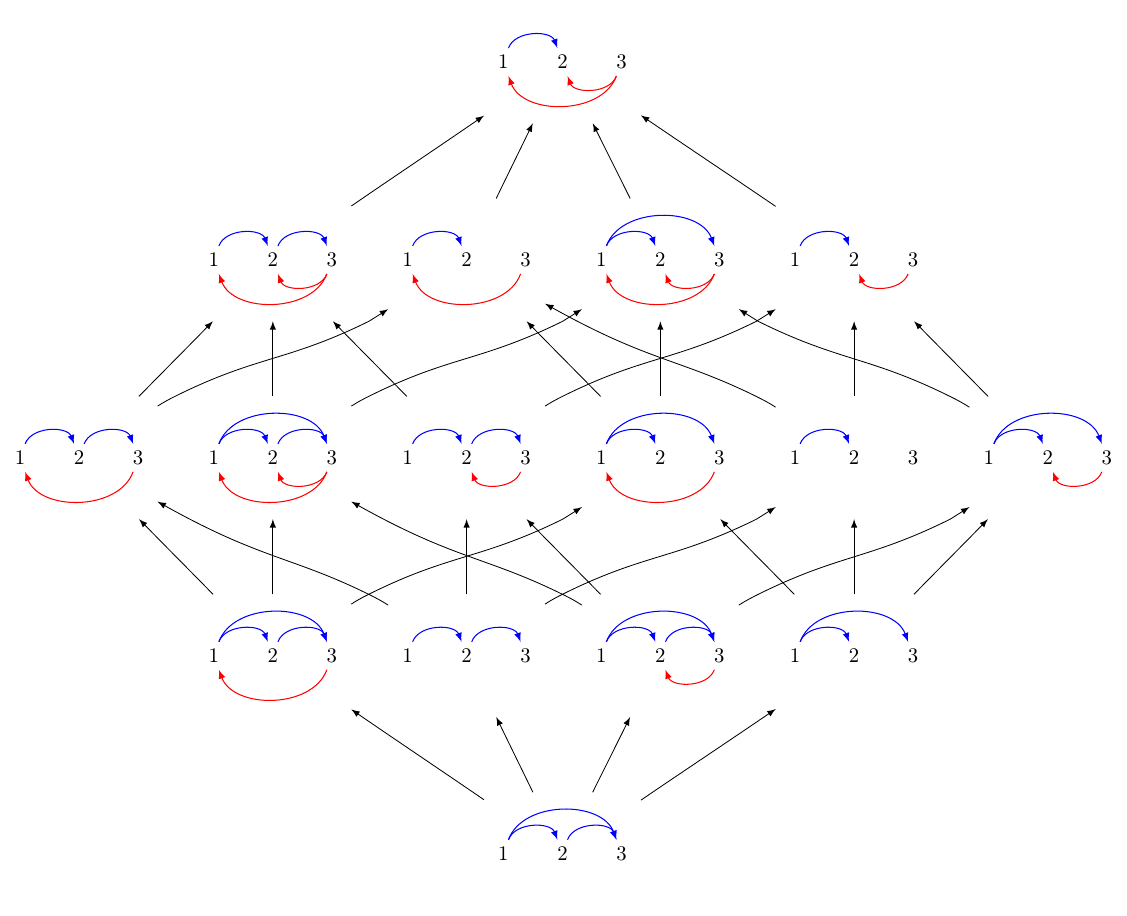}
\input{includes/figures/productRelations}

\footnotesize{(picture from~\citeme{PP20})}
\caption{Product in the Hopf algebra of integer relations as an interval}
\label{fig:rel-product}
\end{figure}

The product between $\rel$ and $\rel[S]$ is the sum over all the elements of the shifted shuffle. This actually gives a sum over an interval of the weak order on integer relations. The minimal element is given by $\underprod{\rel}{\rel[S]}$, the relation where all increasing relations have been added between $\rel$ and $\overline{\rel[S]}$ and no decreasing relations. While the maximal element is $\overprod{\rel}{\rel[S]}$, the relation where all decreasing relations have been added between $\overline{\rel[S]}$ and $\rel$ and no increasing relations. See an example on Figure~\ref{fig:rel-product}.

The coproduct on $\K\IRel$ is defined using the notion of \defn{total cuts} of a relation (Definition~10 of~\citeme{PP20}). A total cut on relation $\rel$ of size $n$ is a partition of $1,\dots,n$ into two sets $X$ and $Y$ such that for all $x \in X$ and $y \in Y$, we have $x \rel y$ and $y \notrel x$. In particular, there are always two trivial total cuts: the case where either $X$ or $Y$ is the empty set. The coproduct is the sum over all the total cuts, where the relation $\rel$ is divided into the two (standardized) subrelations given by the total cut. For example, we have

\begin{align*}
\coproduct \Big(
\FRel_{\scalebox{.5}{\input{includes/figures/relations/r3_12_13_32}}} \! \!
\Big)
=
\FRel_{\scalebox{.5}{\input{includes/figures/relations/r3_12_13_32}}} \! \!  \otimes \FRel_\varnothing
+
\FRel_{\scalebox{.5}{\begin{tikzpicture}[baseline, scale=\interscale]
\node(T1) at (0,0) {1};
\end{tikzpicture}}} \! \!
\otimes
\FRel_{\scalebox{.5}{\begin{tikzpicture}[baseline, scale=\interscale]
\node(T1) at (0,0) {1};
\node(T2) at (1,0) {2};
\draw[->, line width = 0.5, color=red] (T2) edge [bend left=70] (T1);
\draw[->,line width = 0.5, color=white, opacity=0] (T1) edge [bend left=70] (T2);
\draw[->,line width = 0.5, color=white, opacity=0] (T2) edge [bend left=70] (T1);
\end{tikzpicture}}} \! \!
+
\FRel_{\scalebox{.5}{\begin{tikzpicture}[baseline, scale=\interscale]
\node(T1) at (0,0) {1};
\node(T2) at (1,0) {2};
\draw[->, line width = 0.5, color=blue] (T1) edge [bend left=70] (T2);
\draw[->,line width = 0.5, color=white, opacity=0] (T1) edge [bend left=70] (T2);
\draw[->,line width = 0.5, color=white, opacity=0] (T2) edge [bend left=70] (T1);
\end{tikzpicture}}} \! \!
\otimes
\FRel_{\scalebox{.5}{\begin{tikzpicture}[baseline, scale=\interscale]
\node(T1) at (0,0) {1};
\end{tikzpicture}}} \! \!
+
\FRel_\varnothing \otimes  \FRel_{\scalebox{.5}{\input{includes/figures/relations/r3_12_13_32}}}
\end{align*}
where the total cuts are $(\{1,2,3\}, \varnothing)$, $(\{1\}, \{2,3\})$, $(\{1,3\},\{2\})$, and $(\varnothing, \{1,2,3\})$.

We prove in~\citeme[Proposition~15]{PP20} that the product and coproduct indeed define a graded Hopf algebra on $\K\IRel$. The correspondence with the weak order allows us to define \defn{multiplicative bases} and obtain that the Hopf algebra is freely generated by certain elements indexed by \defn{under-indecomposable} relations (relations $\rel[T]$ which cannot be written as $\rel[T] = \underprod{\rel}{\rel[S]}$). 

\section{Integer Posets}

\subsection{As a lattice}

A \defn{poset} is a binary relation which is \defn{reflexive}, \defn{antisymmetric} and \defn{transitive}. As we only consider reflexive integer relations, we are interested in studying the induced subposets of the lattice of integer relations on antisymmetric and transitive integer relations. 

\begin{figure}[ht]
\center
\scalebox{.6}{\input{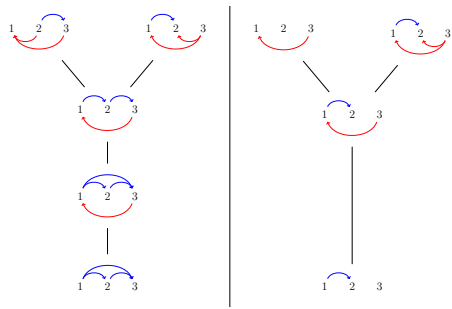}}

\caption{Meet of two transitive relations.}
\label{fig:meet_transitive}
\end{figure}

\begin{figure}[ht]
\center
\includegraphics[scale=.6]{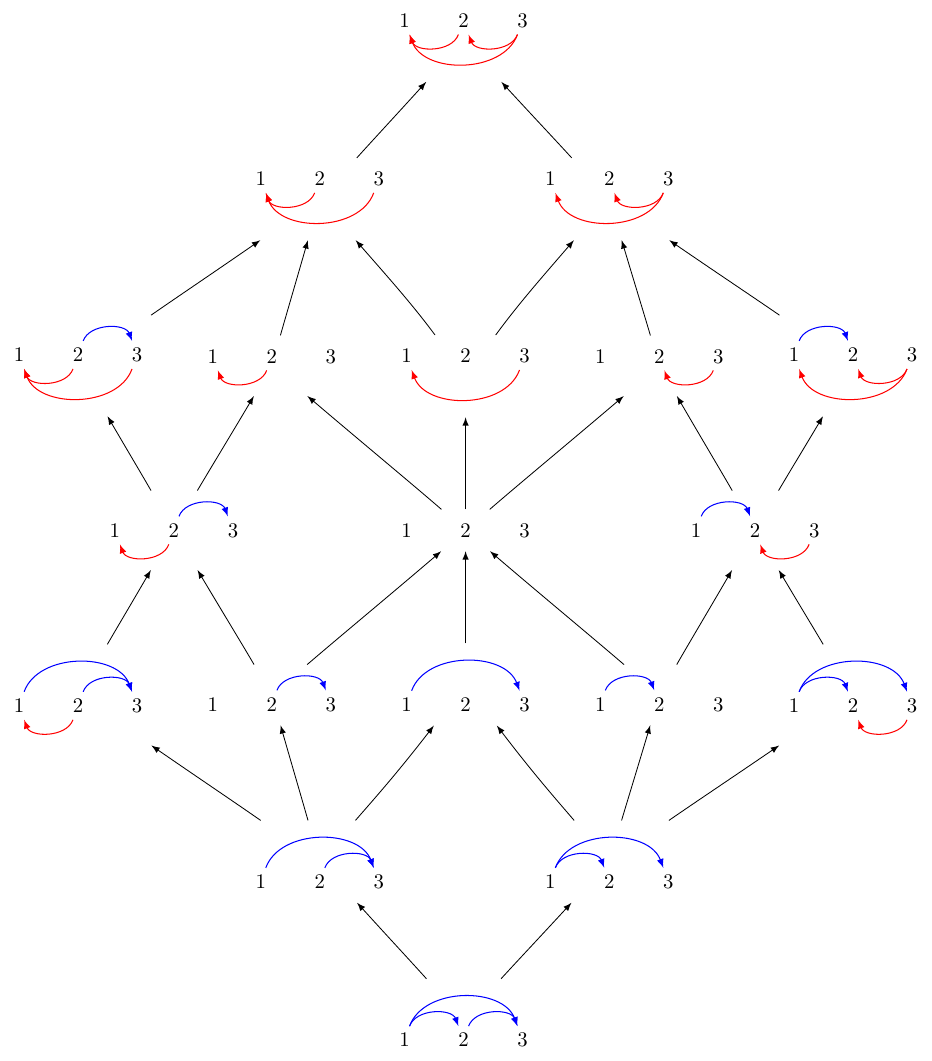}

\footnotesize{(picture from~\citeme{CPP19})}
\caption{The weak order lattice on integer posets}
\label{fig:weak-posets}
\end{figure}

We write $\IAntisym_n$ the set of antisymmetric integer relations of size $n$. It is easy to prove that it forms a sublattice of the weak order on integer relations~\citeme[Proposition~1.4]{CPP19}. Indeed, suppose that $\rel[T] = \rel \meet \rel[S]$ is \emph{not} an antisymmetric relation. This means that we have $x$ and $y$ with $x \rel[T] y$ and $y \rel[T] x$. We suppose $x < y$. The decreasing relations of~$\rel[T]$ are the intersection of the decreasing relations of~$\rel$ and~$\rel[S]$, so we have $y \rel x$ and $y \rel[S] x$. As the increasing relations of $\rel[T]$ are the union of the increasing relations of~$\rel$ and~$\rel[S]$, we have the relation $(x,y)$ in at least one of~$\rel$ or~$\rel[S]$. So one of~$\rel$ or~$\rel[S]$ is also not antisymmetric.

On the other hand, transitive relations $\ITrans_n$ do \emph{not} form a sublattice of the lattice of integer relations. Indeed, it is well known that the union of two transitive relations is not itself transitive. Nevertheless, the weak order on transitive integer relations is still a lattice. We explain a process to construct the meet of two transitive relations as illustrated in~\ref{fig:meet_transitive}. Let $\rel$ and $\rel[S]$ be two transitive relations. Their meet as integer relations is not transitive. We first apply a transitive closure on the \emph{increasing} relations. The result is still not transitive as there might be a chain involving a decreasing relation. In this case, we remove the problematic decreasing relations.

This operation is called the \defn{transitive decreasing deletion} and we prove \linebreak in~\citeme[Proposition~1.14]{CPP19} that it indeed gives the meet of the two relations in the poset of transitive relations. Beside, this is compatible with the antisymmetry condition: the antisymmetic (transitive) relations form a sublattice of the transitive relations. This gives the main theorem of~\citeme{CPP19}.

\begin{theorem}[Theorem~1 from~\citeme{CPP19}]
\label{thm:integer-poset-lattice}
The weak order on integer posets is a lattice.
\end{theorem}

See Figure~\ref{fig:weak-posets} for an example in size $3$. Note that this was later generalized to other finite Coxeter groups in~\citeext{GP20}.

\subsection{As a Hopf algebra}

Looking at integer posets in the Hopf algebra of integer relations, we obtain the two following properties~\citeme[Proposition~30]{PP20}:
\begin{itemize}
\item if at least one integer poset appears in the product of two relations $\rel$ and~$\rel[S]$, then both $\rel$ and $\rel[S]$ are integer posets;
\item the coproduct of an integer poset only contains integer posets.
\end{itemize}

\begin{figure}[ht]
\center
\includegraphics[scale=.5]{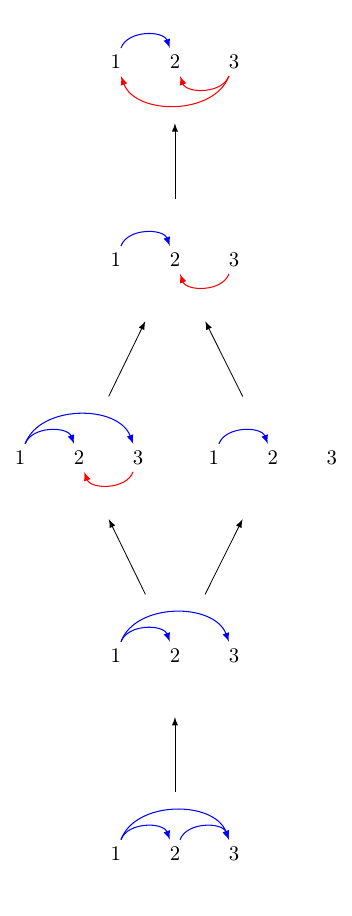}

\footnotesize{(picture from~\citeme{PP20})}

\caption{Interval corresponding to the product of integer posets.}
\label{fig:product-integer-posets}
\end{figure}

This allows us to define a Hopf algebra on integer posets as a quotient of the Hopf algebra of binary relations by the ideal generated by relations which are \emph{not} posets. In other words, we just ``remove'' the non-posets from the result of the product. For example, the result of the product of Figure~\ref{fig:rel-product} in the Hopf algebra of integer-posets is now given by a sum over the interval of Figure~\ref{fig:product-integer-posets}.

\subsection{General audience paper}

\begin{figure}[ht]
\center
\includegraphics[scale=.2]{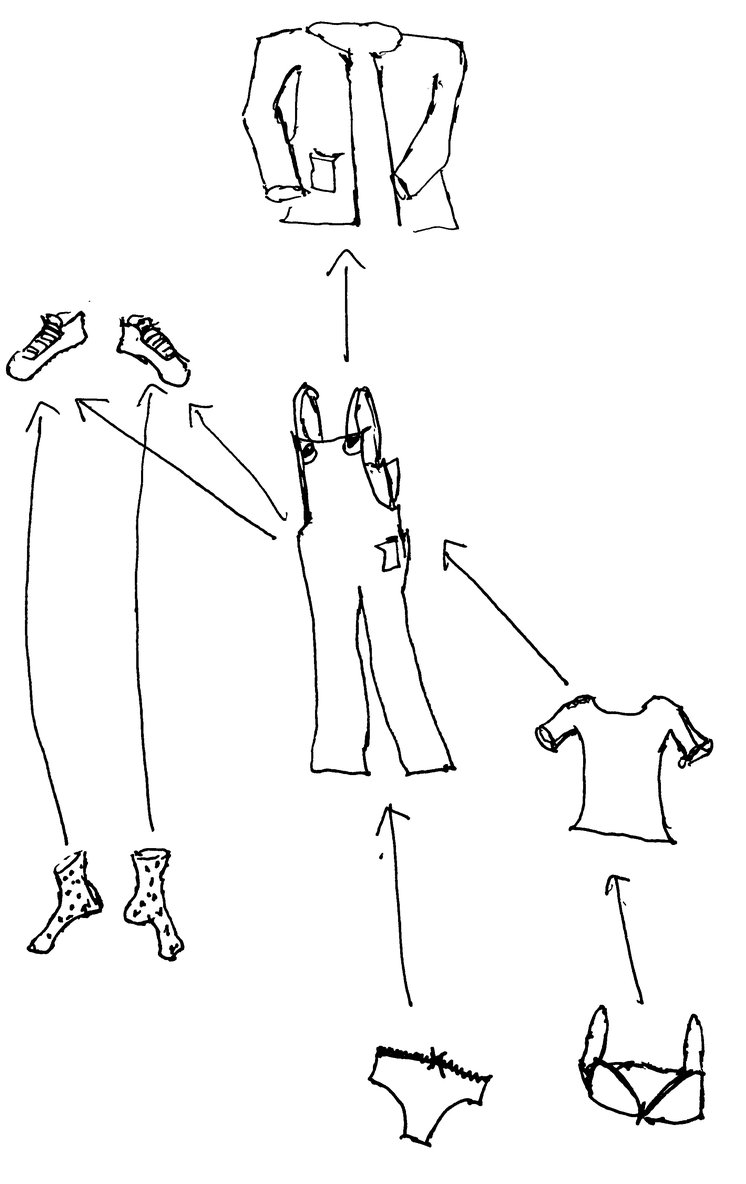}

\caption{The Closet-Poset}
\label{fig:closet-poset}
\end{figure}

The lattice aspects of this work concerning binary integer relations and integer posets have been turned into a general audience paper published in the magazine \emph{Interstices}~\citememisc{PonMisc17}. This was an occasion to present the notion of \defn{poset} using the illustration of the \defn{closet-poset} of Figure~\ref{fig:closet-poset}.

\section{Induced structures on permutations, binary trees and more}
\label{sec:rel-induced}

Permutations seen as $\WOEP$ (see Section~\ref{sec:rel-objects}) form a sublattice of the lattice of integer posets~\citeme[Proposition 2.3]{CPP19}. This is the weak order. Indeed the order is given by the inclusion of the inversions. Moreover, the $\WOIP$ and $\TOIP$ induce some subposets. They are both lattices but only $\TOIP$ forms a  sublattice. This led us to study more precisely in which cases we obtain a sublattice and to describe processes to compute the meet (and the join) when the induced poset is a lattice but not a sublattice. 

This motivated the definition of \defn{permutrees} which we describe precisely in Chapter~\ref{chap:permutrees}. Recall by looking at Figure~\ref{fig:interval-conditions} that, in a weak order interval, for each increasing relation $a \trprec c$ we need to have $a \trprec b$ or $b \trprec c$ for $a < b < c$ (with symmetrical conditions on decreasing relations). In the case of Tamari interval-posets, we ``force'' to always have $b \trprec c$. We could also force to always have $a \trprec b$. In the boolean lattices, we force to have \emph{both}. It is quite natural to ``mix'' the forcing conditions of the different objects. For example, in size $4$, we can decide that each increasing relation $a \trprec c$ with $a < 2 < c$ forces the relation $2 \trprec c$, whereas each increasing relation $a \trprec c$, with $a < 3 < c$ forces the relation $a \trprec 3$. 

We encode the rules into a \defn{decoration}. Each decoration gives a certain family of intervals and elements: the permutrees. In particular, this includes and generalizes the \defn{Cambrian trees} defined by Reading in~\citeext{Rea06}. In~\citeme[Theorems 2.77, 2.82, 2.85, and 2.88]{CPP19}, we study the relations between the lattices on intervals and elements and the lattice of integer posets. We find in particular a sufficient condition on the decoration to obtain sublattices of the lattice of integer posets: it needs to be \defn{covering}: for each $b$, with $1 < b < n$, either $a \trprec b$ or $b \trprec c$ should be explicitly forced (or both). We obtain like this that the Tamari lattice, the Cambrian lattice, and the boolean lattice as well as their lattices of intervals are sublattices of the lattice of integer posets. 

Besides, using the characterization of the faces of the permutahedron as $\WOFP$, of the associahedron as $\TOFP$, we obtain induced subposets of the integer posets which are lattices. For the permutahedron, this is the \defn{facial order} defined\linebreak in~\citeext{KLNPS01} and more generally in~\citeext{DHP16}. The facial order on the faces of the associahedron is found in~\citeext{DHP16, NT06, PR06}. This construction also generalizes to the faces of the permutreehedron (we define this polytope in the next Chapter).

In~\citeme{PP20}, we study the induced Hopf algebraic structures and obtain the following.

\begin{theorem}[from~\citeme{PP20}]
\label{thm:integer-poset-hopf}
The intervals of weak order, Tamari lattice, boolean lattice and more, can be endowed with Hopf algebraic structures.
\end{theorem}

For the intervals of the weak order, we are able to quotient the Hopf algebra of integer posets to obtain a Hopf algebra of weak order intervals ($\WOIP$). The process is very similar to the quotient from integer relations to integer posets: we ``remove'' the non-weak order intervals from the product of two weak order intervals in the Hopf algebra of integer posets. Besides, the $\WOIP$ Hopf algebra can also be obtained as a subalgebra of the integer-posets algebra by constructing certain fibers of a deletion process defined in~\citeme{CPP19}. Using similar processes, we reconstruct the Malvenuto-Reutenauer Hopf algebra on permutations~\citeext{MR95} and the Chapoton Hopf algebra on faces of the permutahedron~\citeext{Cha00}.

We then prove in~\citeme[Proposition 84]{PP20} that the $\TOIP$ (intervals of the Tamari lattice) form a \emph{sub Hopf algebra} of the Hopf algebra of weak order intervals. Similarly, we recover the Loday-Ronco Hopf algebra~\citeext{LR98} and a Hopf algebra on the faces of the associahedron~\citeext{Cha00}.  All these constructions generalize to permutrees.

\section{Open questions}

In~\citeext{Foi13}, Foissy describes a Hopf algebra of \defn{double posets} and a specific case where one of the poset is actually linear. By looking at the product, it is clear that it is not the same as the one we define. It still raises questions on these structures. Indeed, there are different ways to define products and co-products on double relations (and more generally on $k$-tuples of binary relations). One choice gives our Hopf algebra, another gives Foissy Hopf algebra. It would be interesting to classify these different possibilities: which ones are self dual? In case they are not, what is the dual? In which case can we obtain Hopf algebras on posets, on permutations or permutations intervals? What are the structures that we recover or define this way? 

\chapter{Permutrees}
\label{chap:permutrees}

\chapcitation{Mathematical science shows what is. It is the language of unseen relations between things. But to use and apply that language, we must be able to fully to appreciate, to feel, to seize the unseen, the unconscious.}{Ada Lovelace}

In~\citeext{CP17}, Ch\^atel and Pilaud present a Hopf algebra on \defn{Cambrian trees}. Cambrian trees are certain non rooted binary trees which can be used to represent the elements of the \defn{Cambrian lattices} of Reading~\citeext{Rea06}. The Cambrian lattices are lattice quotients of the weak order and each Cambrian tree represent a class. In particular, Cambrian trees can be represented as integer posets (see Chapter~\ref{chap:integer-posets}). While working on integer posets, it appeared that we could define a more general family: the \defn{permutrees}. The name ``permutree'' comes from ``permutation'' and ``tree''. Indeed, permutrees are a family that regroups permutations, rooted binary trees, Cambrian trees, binary sequences, and interpolations between these objects. A permutree comes with a certain \defn{decoration}. Each decoration defines a specific family of permutrees, a lattice quotient of the weak order and a polytope. Besides, we are able to define a Hopf algebra of permutrees (using all decorations) which is a subalgebra of a certain decorated version of the Malvenuto-Reutenauer Hopf algebra on permutations~\citeext{MR95}.

In this chapter, we present the results of two papers: first~\citeme{PP18} (also presented at EUROCOMB in~\citemeconf{PP17}) in collaboration with Vincent Pilaud, \linebreak then~\citeme{PPTJ23} (also presented at FPSAC in~\citemeconf{PPTJ21}) in collaboration with Vincent Pilaud and Daniel Tamayo Jiménez. I consider this second paper as a first step towards defining permutrees for all finite Coxeter groups.


\section{Lattice quotients and decorations}
\label{sec:permutree-congruence}

In Section~\ref{sec:asso-cambrian}, we explain how Cambrian lattices can be constructed as quotients of the weak order defined by a certain orientation of the Dynkin diagram of the group. In type $A$, this means choosing an orientation $i \rightarrow i+1$ or $i \leftarrow i+1$ for each $1 \leq i < n - 1$. If we have $i \rightarrow i+1$, we contract the edge between the permutations $s_{i+1} = 1 \dots i (i+2) (i+1) \dots n$ and $s_{i+1} s_i = 1 \dots (i+2) i (i+1) \dots n$. This is the ``right'' side of the initial hexagon formed by $s_i$ and $s_{i+1}$. Symmetrically, if we have $i \leftarrow i+1$, we contract the edge between the permutations $s_{i} = 1 \dots (i +1) i (i+2) \dots n$ and $s_{i} s_{i+1} = 1 \dots (i+1) (i+2) i \dots n$, the ``left'' side of the hexagon. We then take the lattice congruence generated by these contractions.

\begin{figure}[ht]
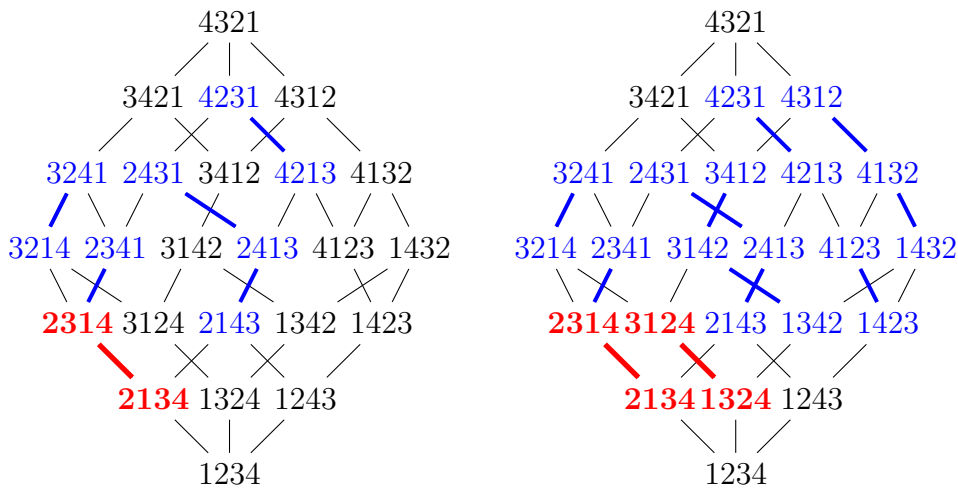

\center
\begin{tabular}{cc}
\input{includes/figures/perm_droit4_congruence_YI.tex}
&
\input{includes/figures/perm_droit4_congruence_XI.tex}
\end{tabular}

\caption{Two examples of permutree congruences.}
\label{fig:permutree-cong}
\end{figure}

For Cambrian congruences, each edge needs to to be oriented in exactly one direction. The \defn{permutree congruences} are the generalizations where we now allow an edge to \emph{not} be oriented or to be oriented in both directions. See on Figure~\ref{fig:permutree-cong} the examples on size $4$ for the orientations $1 \leftarrow 2 \noarrow 3$ and $1 \leftrightarrow 2 \noarrow 3$. The initial contracted edges are in red. Other contracted edges are in blue. Two permutations are in the same class if there is a path of contracted edges between them.

\hvFloat[capPos=b, capAngle=0, capWidth=w, objectAngle=90, capVPos=c, objectPos=c]{figure}
{\begin{minipage}{21cm}\vspace*{-2.1cm}\centerline{\includegraphics[scale=.3]{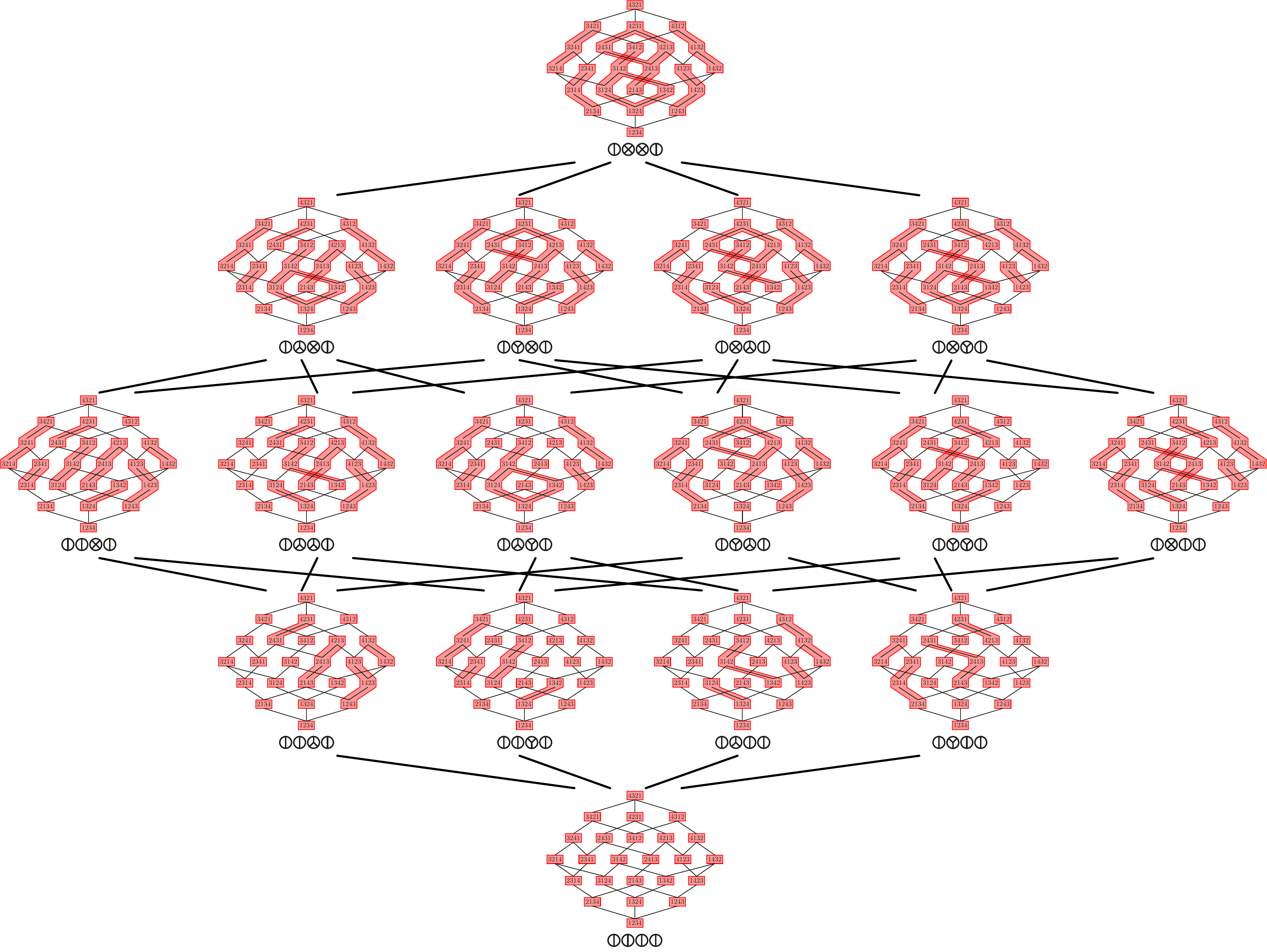}}
\centerline{\footnotesize{(picture from \citeme{PP18})}}
\end{minipage}}
{The permutree congruences, for all decorations~$\decoration \in \noneCirc{} \cdot \Decorations^2 \cdot \noneCirc{}$.}
{fig:fibersPermutreeCongruences}

These congruences form a subfamily with interesting properties of the set of all lattice congruences of the weak order studied by Reading in~\citeext{Rea04}. We show all possible permutree congruences for size $4$ on Figure~\ref{fig:fibersPermutreeCongruences}. We see that they are ordered by refinement: orienting an edge in the Dynkin diagram corresponds to contracting a side of a polygon and to merge some congruence classes. The most refined congruence corresponds to the case where no edge is contracted: this is the trivial congruence where each permutation is its own class and we have $n!$ classes. If all edges are contracted, we obtain the cube congruence we described in Section~\ref{sec:asso-cube} with $2^{n-1}$ classes. If exactly one edge is contracted for each hexagon, this is a Cambrian congruence. If we contract always the same size, this is the sylvester congruence. In both cases, the number of classes is given by the Catalan numbers. In Section~\ref{sec:permutrees-permutrees}, we explain how in general we can associate each class with a tree satisfying certain properties called a \defn{permutree}.

To each orientation of the Dynkin diagram, we associate a certain word $\decoration \in \Decorations^n$ called the \defn{decoration} using the following rules

\begin{tabular}{cc}
$i \noarrow i+1$ & $\decoration_{i+1} = \noneCirc{}$; \\
$i \rightarrow i+1$ & $\decoration_{i+1} = \downCirc{}$; \\
$i \leftarrow i+1$ & $\decoration_{i+1} = \upCirc{}$; \\
$i \leftrightarrow i+1$ & $\decoration_{i+1} = \upDownCirc{}$. 
\end{tabular}

These symbols correspond to certain type of nodes in the permutrees as we explain in Section~\ref{sec:permutrees-permutrees}. The first and last symbols $\decoration_1$ and $\decoration_n$ are set to $\noneCirc$ by convention. We could actually use any symbols: the first and last symbol of the word only play a role for the combinatorial definition but do not change the corresponding Dynkin diagram orientation nor the congruence of the weak order. Using this convention, the orientations $1 \leftarrow 2 \noarrow 3$ and $1 \leftrightarrow 2 \noarrow 3$ of Figure~\ref{fig:permutree-cong} correspond to the words $\noneCirc{} \upCirc{} \noneCirc{} \noneCirc{}$ and  $\noneCirc{} \upDownCirc{} \noneCirc{} \noneCirc{}$ respectively.

We can use the decorations to interpret the minimal elements of the congruence classes in terms of pattern avoidance. Recall that the minimal elements of the sylvester classes are the $312$-avoiding permutation. In case of permutrees, we need to specify the \emph{value} of the middle element of the pattern.  Minimal permutations in the permutree classes are the ones avoiding $cab$ subwords with $a < b < c$ and $\decoration_b \in \lbrace \downCirc{}, \upDownCirc{} \rbrace$ and $bca$ with $a < b < c$ and $\decoration_b \in \lbrace \upCirc{}, \upDownCirc{} \rbrace$. For example, $2413$ is not minimal in its class for $\decoration = \noneCirc{} \upCirc{} \noneCirc{} \noneCirc{}$ because it contains $241$ and $\delta_2 = \upCirc$. This characterization is given in~\citeme[Corollary 2.14]{PP18}. We use it in~\citeme{PPTJ23} to characterize the reduced words of the minimal class elements as explained in Section~\ref{sec:permutree-words}.

\section{Permutree insertion and permutree lattice}
\label{sec:permutrees-permutrees}

\begin{figure}[ht]
\center
\includegraphics{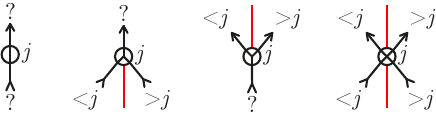}

\caption{The labeling rules on permutrees.}
\label{fig:permutree-rule}
\end{figure}

Permutrees are directed (from bottom to top) planar non rooted trees, labeled by exactly the numbers $1, \dots, n$ and such that each node is of certain type defined by a symbol in $\Decorations$. The image of the symbol represents the incoming / outgoing edges and implies certain rules illustrated on Figure~\ref{fig:permutree-rule}. A node $\noneCirc{}$ has one incoming edge and one outgoing edge. A node $\downCirc{}$ has two incoming edges and one outgoing edge.  A node $\upCirc{}$ has one incoming edge and two outgoing edges. And a node $\upDownCirc{}$ has two incoming edges and two outgoing edges. The nodes coming from or going to the \emph{left} of a given node $i$ must be labeled with values smaller than $i$. Symmetrically, the nodes coming from or going to the \emph{right} of a given node $i$ must be labeled with values bigger than $i$. The \defn{decoration} is the words of node symbols read from $1$ to $n$ (from left to right). We say that the symbols $\downCirc{}$ and $\upDownCirc{}$ are \defn{down symbols} as they have two incoming edges. Similarly, we write that $\upCirc{}$ and $\upDownCirc{}$ are \defn{up symbols}. In particular $\upDownCirc{}$ is both an up and a down symbol while $\noneCirc{}$ is neither.

\begin{figure}[ht]
  \centerline{\includegraphics[scale=.8]{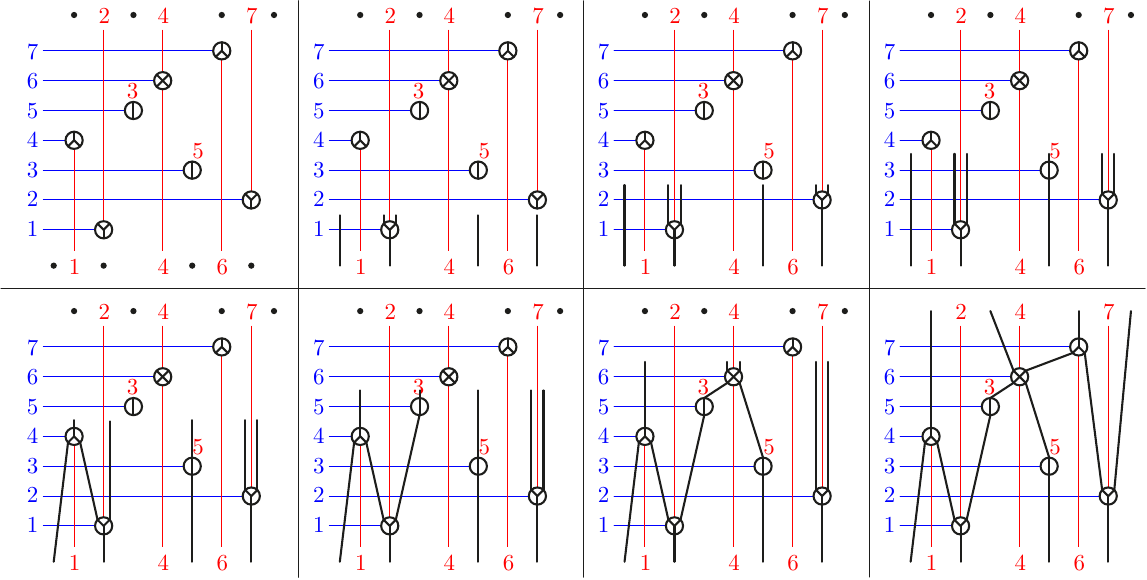}}
  
  \centerline{\footnotesize{(picture from \citeme{PP18})}}
  \caption{The insertion algorithm on the decorated permutation~$\up{2}\up{7}5\down{1}3\updown{4}\down{6}$.}
  \label{fig:insertionAlgorithm}
\end{figure}

Permutrees are actually better understood as the result of a certain \defn{insertion algorithm} from permutations that we illustrate on Figure~\ref{fig:insertionAlgorithm}. The algorithm takes a permutation $\sigma$ of size $n$ and a decoration $\decoration \in \Decorations^n$ and constructs a permutree with decoration $\decoration$. We often write the decoration directly on the values of the permutation for convenience. A value $v$ is underlined if $\decoration_v$ is a down symbol and overlined if $\decoration_v$ is an up symbol. For example the decorated permutation $\up{2}\up{7}5\down{1}3\updown{4}\down{6}$ corresponds to the decoration $\downCirc{}\upCirc{}\noneCirc{}\upDownCirc{}\noneCirc{}\downCirc{}\upCirc{}$.

The algorithm is as follows. For all $1 \leq i \leq n$, we place the symbol $\decoration_i$ at position $(\sigma_i, i)$. The values $1, \dots, n$ on the $y$-axis are called the \defn{levels} of the nodes. We trace a red line below any down symbol and above any up symbol. The \defn{label} of the node is the value $\sigma_i$ in the permutation and is written below (or above) the red line. These lines correspond to certain ``walls'' that cannot be crossed by the edges of the permutree. They enforce the labeling rules attached to the nodes. We start an edge at every interstice between the red lines at the bottom of the picture including the borders. We then ``grow'' the tree by attaching edges to the first available entry point they see without crossing a red line.

\begin{figure}[ht]
  \centerline{\includegraphics[scale=.8]{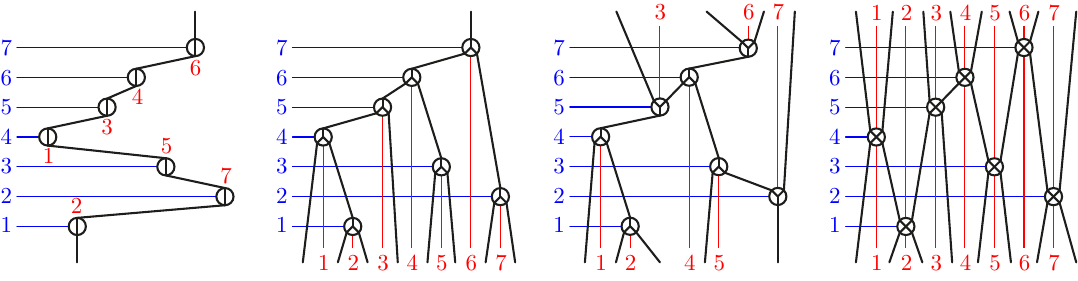}}
  
   \centerline{\footnotesize{(picture from \citeme{PP18})}}
  \caption{Leveled permutrees corresponding to a permutation (left), a leveled binary tree (middle left), a leveled Cambrian tree (middle right), and a leveled binary sequence (right).}
  \label{fig:permutationsBinaryTreesCambrianTreesBinarySequences}
\end{figure}

The result is a \defn{leveled permutree}. Decorated permutations are in bijection with leveled permutrees. The permutree itself is the graph in black oriented from bottom to top. By ``forgetting'' the levels, we obtain a surjection between decorated permutations and permutrees. The fiber of a permutree is the set of the linear extensions of the graph. These fibers correspond to the permutree congruence defined in Section~\ref{sec:permutree-congruence}.  Indeed, if $\decoration = \noneCirc^n$, then there is no red line and the permutree is a linear graph labeled by the permutation (see the first example of Figure~\ref{fig:permutationsBinaryTreesCambrianTreesBinarySequences}). If $\decoration = \downCirc{}^n$, the algorithm is the binary search tree insertion described in Section~\ref{sec:asso-quotient} (see the second example of Figure~\ref{fig:permutationsBinaryTreesCambrianTreesBinarySequences}). If $\decoration \in \lbrace \downCirc{}, \upCirc{} \rbrace^n$, this is the Cambrian tree insertion that can be found in~\citeext{CP17}. If $\decoration = \upDownCirc{}^n$, each node is directly attached to the previous one and is either above or below: this gives a binary sequence as described in Section~\ref{sec:asso-cube}.  

\begin{figure}[ht]
\center
\begin{tabular}{ccccc}
\includegraphics[scale=.8]{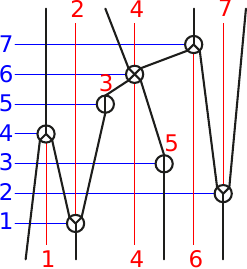}
& &
\includegraphics[scale=.8]{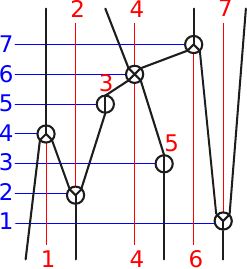}
& &
\includegraphics[scale=.8]{includes/figures/leveledPermutree1} \\
$\up{7}\up{2}5\down{1}3\updown{4}\down{6}$ &
$\equiv$ &
$\up{2}\up{7}5\down{1}3\updown{4}\down{6}$ &
$\equiv$ &
$\up{2}\up{7}\down{1}53\updown{4}\down{6}$
\end{tabular}

\caption{The permutree congruence on decorated permutations.}
\label{fig:permutree-cong-permutations}
\end{figure}

The description of the permutree congruence classes as the linear extensions of permutrees gives the characterization of the minimal (and also maximal) elements which we explained in Section~\ref{sec:permutree-congruence}. Actually, it describes the possible local moves on any decorated permutation within a congruence class. Indeed, the permutree congruence is the transitive closure of two local rules

\begin{align}
\dots ac \dots \down{b} \dots & \equiv \dots ca \dots \down{b} \dots, \\
\dots \up{b} \dots ac \dots & \equiv \dots \up{b} \dots ca \dots,
\end{align}
where $a < b < c$ and the written decoration on $b$ are the \emph{necessary} requirements (there an be extra decorations on $a$, $b$ and $c$). See the example given on Figure~\ref{fig:permutree-cong-permutations}. We have that $\up{7}\up{2}5\down{1}3\updown{4}\down{6} \equiv \up{2}\up{7}5\down{1}3\updown{4}\down{6}$ as we are able to exchange $\up{7}$ and $\up{2}$ because they have a \defn{witness} on the right (either $\updown{4}$ or $\down{6}$) with a down decoration. You can check indeed that the two levels permutrees are different but that the underlying permutrees are similar. Similarly, we have $\up{2}\up{7}5\down{1}3\updown{4}\down{6} \equiv \up{2}\up{7}\down{1}53\updown{4}\down{6}$ because $5$ and $\down{1}$ can be exchanged thanks to the witness $\up{2}$ to the left with an up decoration.

\begin{figure}[ht]
\center
\includegraphics[scale=.5]{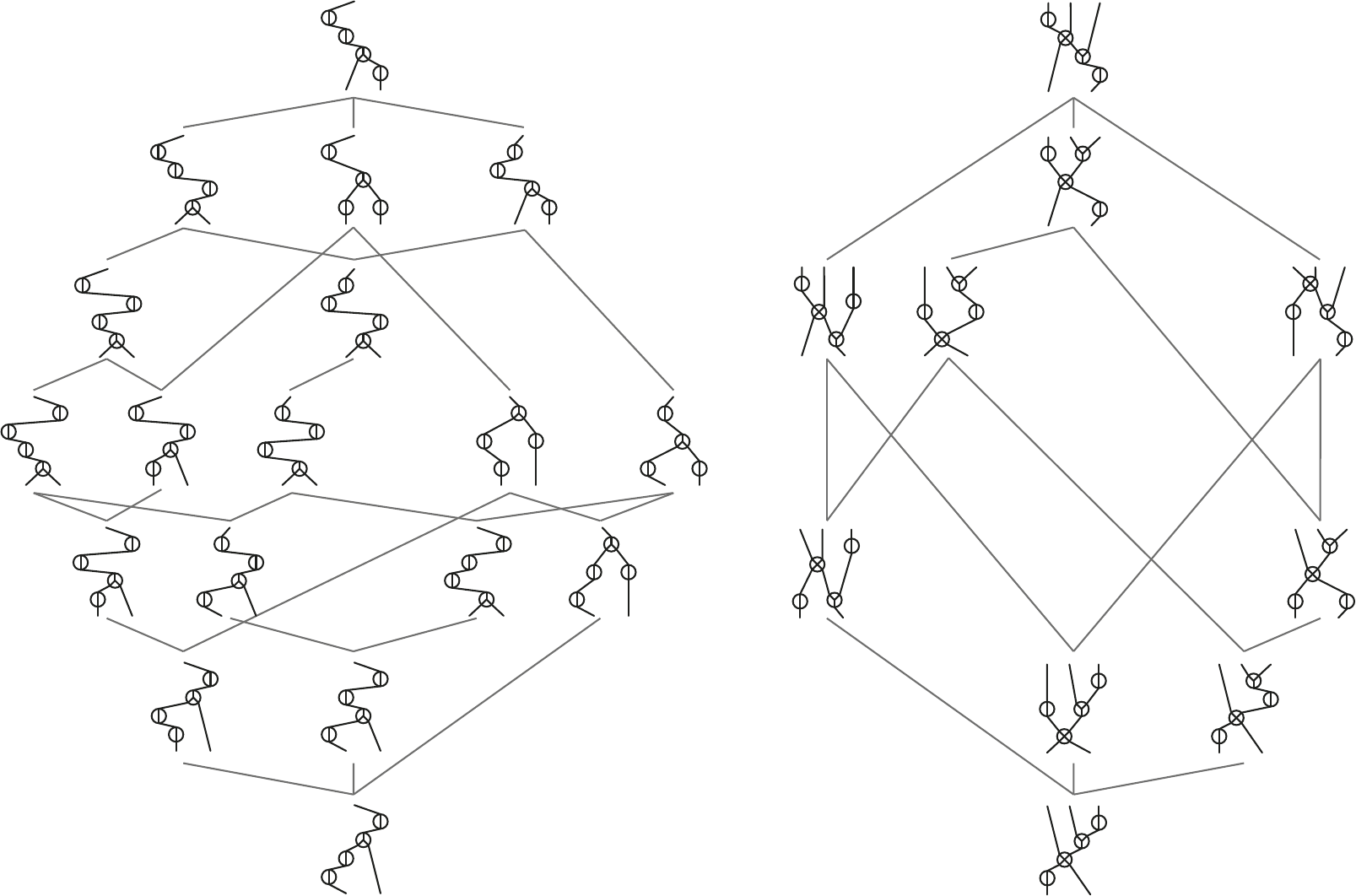}

\footnotesize{(picture from \citeme{PP18})}

\caption{The $\decoration$-permutree lattices, for the decorations $\decoration = \noneCirc{}\noneCirc{}\downCirc{}\noneCirc{}$ (left) and $\decoration = \noneCirc{}\upDownCirc{}\upCirc{}\noneCirc{}$~(right)}
\label{fig:permutreeLattices}
\end{figure}

Each decoration $\decoration$ defines a lattice congruence~\citeme[Proposition 2.13]{PP18}. We thus obtain $\decoration$-permutree lattices (see Figure~\ref{fig:permutreeLattices}). The rotation on permutrees is natural generalization of the rotation on binary trees.

\section{Numerology}

There is a very natural question: given a certain decoration $\decoration$, how many $\decoration$-permutrees are there? This also corresponds to counting the equivalence classes. It is known for some decorations giving specific subfamilies: the factorial numbers for permutations, the Catalan numbers for binary and Cambrian trees, and the powers of $2$ for the binary sequences. In general, we call this the \defn{factorial-Catalan number} $\factCatalan{\decoration}$ and provide a recursive formula which interpolates between the factorial and Catalan numbers. 

Indeed, the Catalan numbers satisfy a well known recursion

\begin{equation}
\label{eq:catalan}
C_n = \sum_{k=0}^{n-1} C_k C_{n-1-k}
\end{equation}
with $C_0 = 1$. On the other hand, the factorial numbers satisfy

\begin{equation}
\label{eq:factorial}
n! = n (n-1)! = \sum_{i=1}^n (n-1)!.
\end{equation}
The idea is then to \emph{mix} those two summations depending on the actual symbols of $\decoration$. The first essential property is given in~\citeme[Corollary 2.21]{PP18}: only the positions of the $\noneCirc{}$ and $\upDownCirc{}$ matter. In particular a $\downCirc$ symbol can be changed into a $\upCirc$ without changing the value of $\factCatalan{\decoration}$. 

Besides, the case of the $\upDownCirc{}$ symbol is actually trivial. If $\decoration_i = \upDownCirc$, we have that $\factCatalan{\decoration} = \factCatalan{\decoration_1 \dots \decoration_{i-1} \noneCirc{}}\factCatalan{\noneCirc{}\decoration_{i+1} \dots \decoration_n}$. For example, for $\decoration = \upDownCirc{} \upDownCirc{} \upDownCirc{} \upDownCirc{}$, we indeed get

\begin{equation}
\factCatalan{\decoration} = \factCatalan{\noneCirc{}\noneCirc{}} \factCatalan{\noneCirc{}\noneCirc{}} \factCatalan{\noneCirc{}\noneCirc{}}  = 2^3.
\end{equation}

Then for $\decoration \in \lbrace \noneCirc{}, \upCirc{}, \downCirc{} \rbrace^n$, we get a recursive formula summing over all possibilities to remove one symbol,

\begin{equation}
\factCatalan{\decoration} = \sum_{\decoration_i = \noneCirc{}} \factCatalan{\decoration_1 \dots \decoration_{i-1}\decoration_{i+1} \dots \decoration_n} + \sum_{\decoration_i \in \lbrace \upCirc{}, \downCirc{} \rbrace} \factCatalan{\decoration_1 \dots \decoration_{i-1}} \factCatalan{\decoration_{i+1} \dots \decoration_n}.
\end{equation}
When the removed symbol is $\noneCirc{}$, we simply sum on the decoration of size $n-1$ as in~\eqref{eq:factorial}. When the removed symbol is $\upCirc{}$ or $\downCirc{}$, the word is split into two parts as in~\eqref{eq:catalan}. Using these rules, we can compute the number of permutrees for the examples of Figure~\ref{fig:permutreeLattices}.

\begin{align}
\factCatalan{\noneCirc{} \noneCirc{} \downCirc{} \noneCirc{}} &= 2 \factCatalan{\noneCirc{} \downCirc{} \noneCirc{}} + \factCatalan{\noneCirc{} \noneCirc{}} \factCatalan{\noneCirc{}} + \factCatalan{\noneCirc{} \noneCirc{} \downCirc{}} \\
&= 2 \times 5 + 2 + 6 = 18
\end{align}

\begin{align}
\factCatalan{\noneCirc{} \upDownCirc{} \upCirc{} \noneCirc{}} &= \factCatalan{\noneCirc{} \noneCirc{}} \factCatalan{\noneCirc{} \upCirc{} \noneCirc{}} \\
&= 2 \times 5 = 10
\end{align}

\section{Permutreehedra}

The weak order, the Tamari lattice, the Cambrian lattices, the boolean lattice: the Hasse diagrams of these lattices all correspond to the skeletons of certain polytopes which are all \defn{generalized permutahedra}~\citeext{Pos09}. A generalized permutahedron is a polytope where every edge follows the direction of a permutahedron edge $e_i - e_j$. When we drew our first permutree lattices, they ``looked nice'' in the sense that it felt like it was the graph of a polytope. For example, in size $n$, any element is adjacent to $n-1$ elements: this is a characteristic of vertices in a simple polytope.

Indeed, we realized that we could generalize Loday's construction of the associahedron (See Section~\ref{sec:asso-loday}) and obtain a geometrical realization of the permutree lattices. More precisely, for every decoration $\decoration$, we define a polytope  $\Permutreehedron$  called the \defn{permutreehedron} as the convex hull of the vertices $\b{a}(\tree)$ for $\tree$ a $\decoration$-permutree. Each coordinate of $\b{a}(\tree)$ is computed as follows:

\[
\b{a}(\tree)_i =
\begin{cases}
1 + d & \text{if } \decoration_i = \noneCirc{}, \\
1 + d + \down{\ell}\down{r} & \text{if } \decoration_i = \downCirc{}, \\
1 + d - \up{\ell}\up{r} & \text{if } \decoration_i = \upCirc{}, \\
1 + d + \down{\ell}\down{r} - \up{\ell}\up{r} & \text{if } \decoration_i = \upDownCirc{},
\end{cases}
\]
where $d$ is the number of descendants of $i$ , $\down{\ell}$ and $\down{r}$ are the number of left and right descendants of $i$, and $\up{\ell}$ and $\up{r}$ are the number of left and right ascendants of $i$. For example, the coordinates of the permutree of Figure~\ref{fig:insertionAlgorithm} are $[7, -4, 3, 8, 1, 12, 1]$.

\begin{theorem}[Theorem~3.4 of~\citeme{PP18}]
\label{thm:permutreehedron}
For each decoration $\decoration$, The permutreehedron $\Permutreehedron$ is the polytope corresponding to the $\decoration$-permutree congruence.
\end{theorem}

The word ``corresponding'' here refers to the \defn{polyhedral fan} of the polytope. In~\citeext{Rea05}, Reading gives the explicit fan corresponding to any lattice quotient of the weak order of a Coxeter group. In particular, this describes the $\decoration$-permutree fan of which the permutahedron fan is a refinement. Our theorem states that the permutreehedron indeed realizes the permutree fan. This includes in particular that the Hasse diagram of the lattice is the skeleton of the polytope.

Beside, we can also construct the permutreehedron as a \defn{removahedron} as in Section~\ref{sec:asso-faces}. Recall that the faces of the permutahedron can be represented by ordered partitions (see Section~\ref{sec:perm-faces}). The facets correspond to partitions of $1,\dots,n$ into two sets $I$ and $J$. For the permutreehedron, we only keep the facets that correspond to an \defn{edge cut} of a permutree. For each oriented edge between two vertices of a permutree, the edge cut partitions $1,\dots, n$ into two sets $I$ and $J$. $I$ contains all the vertices connected to the source of the edge and $J$ all other vertices. For example, on the tree of Figure~\ref{fig:insertionAlgorithm}, the different values for $I$ are $\lbrace 2, 3, 4, 5, 6, 7 \rbrace$, $\lbrace 1, 2 \rbrace$, $\lbrace 1, 2, 3 \rbrace$, $\lbrace 5 \rbrace$, $\lbrace 1,2,3,4,5 \rbrace$, and $\lbrace 7 \rbrace$. These are all the facets containing the tree.

\hvFloat[capPos=b, capAngle=0, capWidth=w, objectAngle=90, capVPos=c, objectPos=c]{figure}
{\begin{minipage}{21cm}\vspace*{-2.1cm}\centerline{\includegraphics[scale=.3]{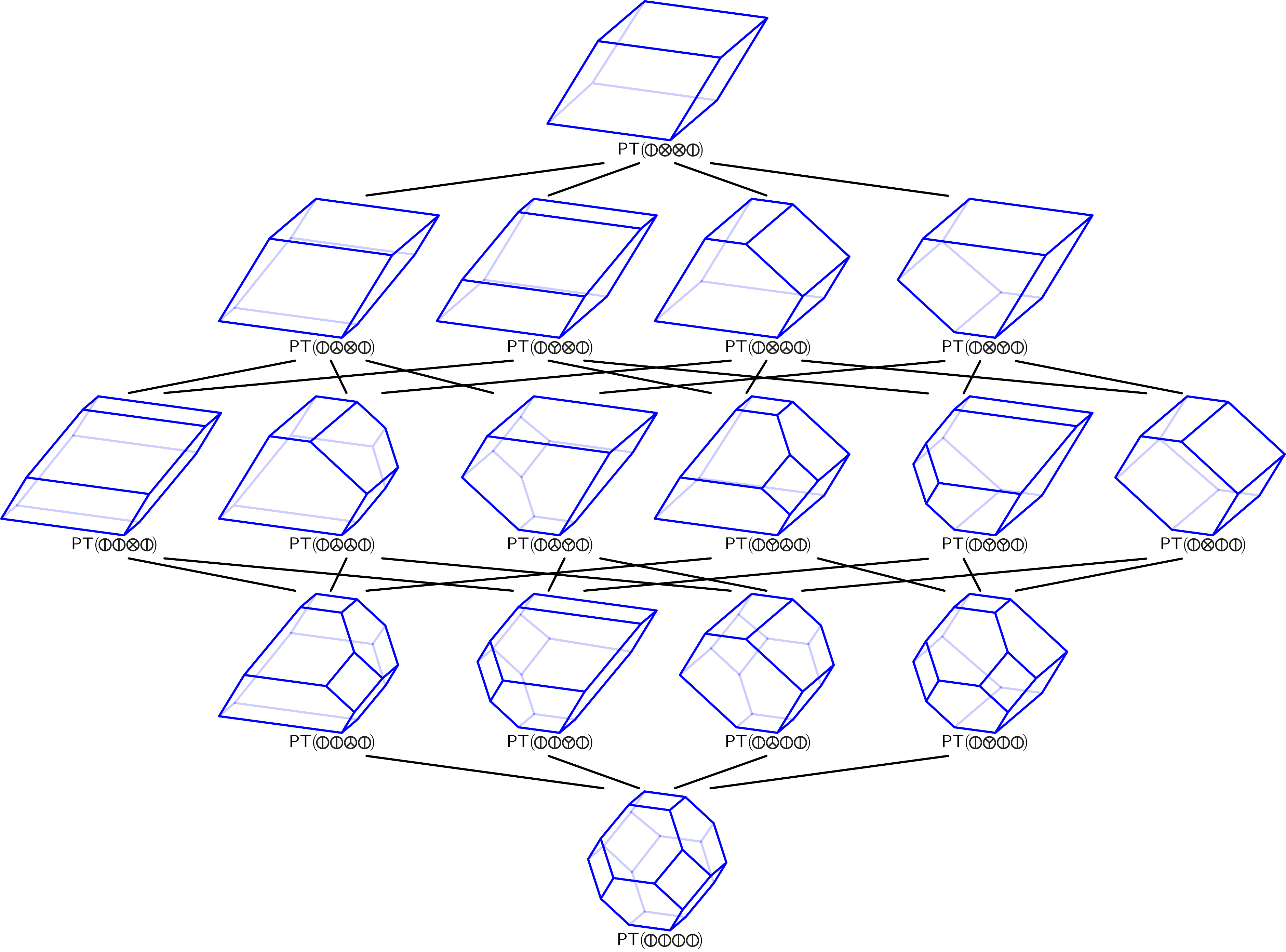}}
\centerline{\footnotesize{(picture from \citeme{PP18})}}
\end{minipage}}
{The permutreehedra, for all decorations~$\decoration \in \noneCirc{} \cdot \Decorations^2 \cdot \noneCirc{}$.}
{fig:permutreehedra}

Figure~\ref{fig:permutreehedra} shows all permutreehedra for decorations of size $4$. You can see the correspondence with Figure~\ref{fig:permutree-cong-permutations}. When we add an orientation to an edge of the Dynkin diagram, this corresponds combinatorially to coarse the lattice congruence, \emph{i.e.}, to merge some congruence classes together. Geometrically, this corresponds to removing some facets so that the polytopes can be seen inside one another. 

In~\citeext{PS19}, Pilaud and Santos proved that all lattice congruence of the weak order can be realized a generalized permutahedra. But their realization is not the same as the one of permutreehedra. In particular, the quotientopes are not removahedra. Actually, it was proven in~\citeext{APR21} that the permutree congruences are exactly the congruences that can be realized as removahedra.

\section{Schr\"oder permutrees}

\begin{figure}[ht]
\center
\includegraphics[scale=.8]{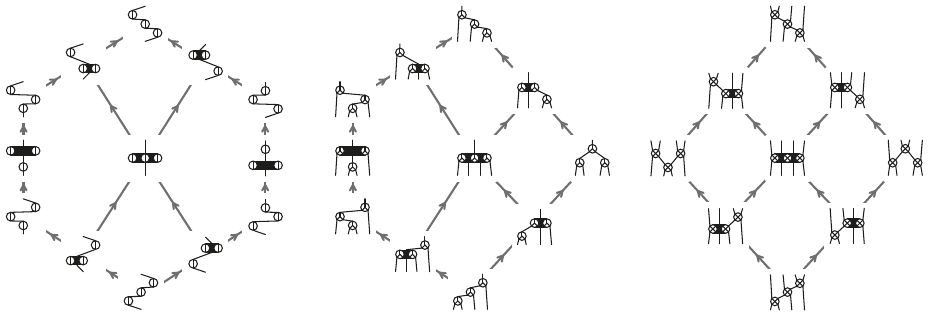}

\footnotesize{(picture from \citeme{PP18})}
\caption{The Schr\"oder permutree lattices for $\noneCirc{}^3$, $\downCirc{}^3$, $\upDownCirc{}^3$}
\label{fig:shroder-permutree-lattice}
\end{figure}

We define a concept of \defn{Schr\"oder permutree} which are the equivalent of Schr\"oder trees for the permutrees (see Section~\ref{sec:asso-faces}). They  correspond to permutrees where some nodes have been merged together. It gives a combinatorial description to the faces of the permutreehedra. Most properties and operations on permutrees easily generalize to Schr\"oder permutrees. We have an insertion algorithm (from ordered partitions), a congruence relation, and a counting formula. The congruence is a lattice congruence from the facial weak order and thus we obtain a facial order on Schr\"oder permutrees as illustrated on Figure~\ref{fig:shroder-permutree-lattice}.

\section{Hopf algebras}

We have seen in Section~\ref{sec:perm-hopf} the definition of the Malvenuto-Reutenauer Hopf algebra ($\FQSym$) and how it relates to the weak order. In particular, the sylvester congruence not only define a lattice congruence but also the Loday-Ronco Hopf algebra on binary trees as a subalgebra of $\FQSym$ (or, dually, a quotient Hopf algebra). The same construction extends to permutrees. For this, we use a decorated version of $\FQSym$. Basically, the product is now defined on decorated permutations and the decoration moves along with the value it is attached to as in the following example. 

\begin{align}
\BF_{\sred{\up{2}3\down{1}}} \hprod \BF_{\updown{1}2} &= \BF_{\sred{\up{2}3\down{1}}\updown{4}5} + \BF_{\sred{\up{2}3}\updown{4}\sred{\down{1}}} + \BF_{\sred{\up{2}}\updown{4}\sred{3\down{1}}5} + \BF_{\sred{\up{2}3}\updown{4}5\sred{\down{1}}} + \BF_{\updown{4}\sred{\up{2}3\down{1}}5} \\
&+ \BF_{\sred{\up{2}}\updown{4}\sred{3}5\sred{\down{1}}} + \BF_{\updown{4}\sred{\up{2}3}5\sred{\down{1}}} + \BF_{\sred{\up{2}}\updown{4}5\sred{3\down{1}}} + \BF_{\updown{4}\sred{\up{2}}5\sred{3\down{1}}} + \BF_{\updown{4}5\sred{\up{2}3\down{1}}}. \nonumber
\end{align}

Conversely, on the co-product (or the product in the dual Hopf algebra), the decoration is attached to its position as in

\begin{align}
\BG_{\sred{\up{3}1\down{2}}} \hprod \BG_{\updown{1}2} &= \BG_{\sred{\up{3}1\down{2}}\updown{4}5} + \BG_{\sred{\up{4}1\down{2}}\updown{3}5} + \BG_{\sred{\up{4}1\down{3}}\updown{2}5} + \BG_{\sred{\up{5}1\down{2}}\updown{3}4} + \BG_{\sred{\up{4}2\down{3}}\updown{1}5} \\
&+ \BG_{\sred{\up{5}1\down{3}}\updown{2}4} + \BG_{\sred{\up{5}2\down{3}}\updown{1}4} + \BG_{\sred{\up{5}1\down{4}}\updown{2}3} + \BG_{\sred{\up{5}2\down{4}}\updown{1}3} + \BG_{\sred{\up{5}3\down{4}}\updown{1}2}. \nonumber
\end{align}

We define a basis indexed by permutrees. The element indexed by a permutree~$T$ is a sum over all decorated permutations $\sigma$ whose insertion algorithm results in~$T$: these are the linear extensions of the tree or, equivalently, all the elements of the permutree class. See the example below for the tree of Figure~\ref{fig:insertionAlgorithm}. 

\[
\PPT_{\!\!\includegraphics{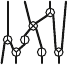}} = \BF_{\up{2}\down{1}35\updown{4}\up{7}\down{6}} + \BF_{\up{2}\down{1}35\up{7}\updown{4}\down{6}} + \BF_{\up{2}\down{1}3\up{7}5\updown{4}\down{6}} + \dots + \BF_{\up{7}5\up{2}3\down{1}\updown{4}\down{6}} + \BF_{\up{7}5\up{2}3\updown{4}\down{1}\down{6}} + \BF_{\up{7}5\up{2}3\updown{4}\down{6}\down{1}}.
\]

\begin{theorem}[Theorem~4.5 of \citeme{PP18}]
\label{thm:permutree-hopf}
The elements $(\PPT_T)$ where $T$ is a permutree form a Hopf subalgebra of the decorated version of $\FQSym$. 
\end{theorem}

\begin{figure}[ht]
\center

\includegraphics[scale=.7]{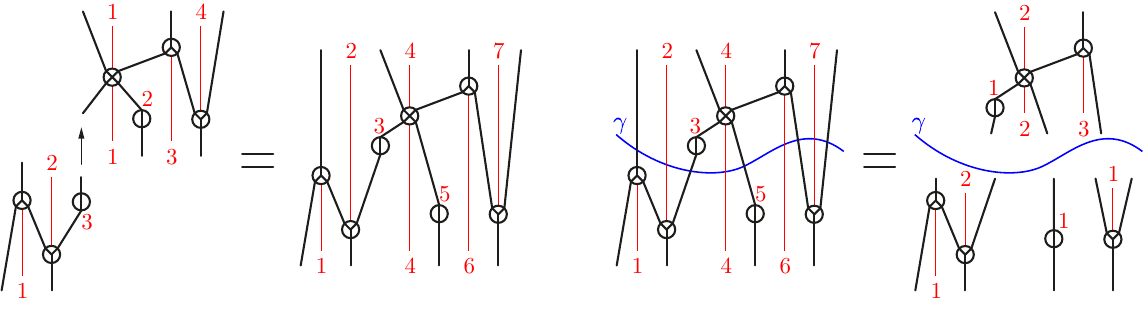}

\footnotesize{(picture from~\citeme{PP18})}
\caption{Grafting two permutrees (left) and cutting a permutree (right).}
\label{fig:permutree-product-coproduct}
\end{figure}

On Figure~\ref{fig:permutree-product-coproduct}, we show two operations on permutrees (grafting and cutting) which give a direct combinatorial descriptions of the product and coproduct. Besides, by duality, the dual Hopf algebra is a quotient of the decorated dual of $\FQSym$ and we also obtain direct descriptions of the operations. 

A similar construction can be made on Schr\"oder permutrees using a decorated version of the Chapoton Hopf algebra on the faces of the permutahedra~\citeext{Cha00}.

\section{Reduced words and other groups}
\label{sec:permutree-words}

Cambrian lattices are defined for all finite Coxeter groups. In particular, the minimal elements of each class are proved to be the \defn{$c$-sortable} elements as we explain in Section~\ref{sec:asso-cambrian}. The concept of permutree congruence can also easily be extended to other finite Coxeter groups but the identification of minimal class elements turns out to be more difficult than for the the Cambrian case.

Indeed, the $c$-sorting algorithm can be seen as a generalization of classical \defn{stack sorting} (which would correspond to the choice of $s_1 \dots s_{n-1}$ as a Coxeter element in type $A$). From an algorithmic point of view, the $c^{\infty}$ word is a sorting network. At each step, we try to apply a certain operation if possible. The list of operations to ``try'' is decided in advance and does not depend of the success / failure of the previous operations. The list of operations actually applied give us a subword of $c^{\infty}$ which identifies if an element is minimal in its Cambrian class or not. The most natural question would then be: is there such a sorting network for all permutree congruences? We give a negative answer in~\citeme[Remark 5.9]{PPTJ23} which indicates that the characterization of the permutree classes in general is a difficult problem.

Our approach in~\citeme{PPTJ23} is to study the reduced words of minimal permutations in the permutree classes in type $A$. The basic case corresponds to orienting exactly one edge of the Dynkin diagram, \emph{i.e.}, to contract the side of exactly one hexagon. We can then construct the $\automatonU(j)$ automaton of Figure~\ref{fig:automaton-uj} which corresponds to orienting the edge $(j-1) \leftarrow j$. 

\begin{figure}[ht]
\center

\scalebox{.8}{\input{includes/figures/uj_automaton}}
\caption{The automaton~$\automatonU(j)$ to.}
\label{fig:automaton-uj}
\end{figure}

Even though we have not explicitly drawn the looping transitions, this is a complete automaton: at each state, any transposition $s_i$ can be read. Only the transitions that result in a change of state are written explicitly. The states of the first and second line are accepting states whereas the states in the third line are rejecting states.  

\begin{theorem}[Theorem~1.1 of~\citeme{PPTJ23}]
\label{thm:permutree-automaton}
A permutation is minimal in its permutree class if and only if it admits a reduced word accepted by the automaton.
\end{theorem}

For example, in size $4$, for $U_2$, corresponding to the orientation $1 \leftarrow 2 \noarrow 3$ of Figure~\ref{fig:permutree-cong} (left), the reduced word $s_2 s_1 s_3$ is accepted which means that the permutation $3142$ is minimal in its class. On the other hand, $s_1 s_2$ (permutation $2314$) and $s_2s_1s_2s_3$ (permutation $3241$) are rejected by the automaton. You can check that it is the case of all their reduced words and they are not minimal in their classes.

We actually provide an algorithm so that, given a permutation, only one reduced word needs to be tested in the automaton. This offers a new characterization of minimal permutree classes elements besides the pattern avoidance, just like stack sortable permutations correspond to $231$-avoiding permutations. By taking products of automata, we generalize our algorithm to all permutree congruences in type~$A$. 

\section{Open questions}

We believe that similar automata can be found for other types of finite Coxeter groups. Daniel Tamayo worked on this subject during his thesis under the supervision of Vincent Pilaud and myself. We were able to find a general rule for type~$B$. For type $D$, we could find some automata that worked in small sizes but we do not have a general rule at this date. From our observations, the general case automata (if they exist) in types different than $A$ and $B$ are much more elaborate. In type $A$ especially, many simplifications arise, giving very simple automata. 

More generally, there are still many areas to explore for permutrees in other Coxeter groups. For example, can we obtain the permutreehedron as Coxeter removahedra and do the results of~\citeext{APR21} can be extended? We would also like to consider permutrees as subsets of root systems like in the work of Gay and Pilaud~\citeext{GP20}.

\part{Generalizations}
\label{part:generalizations}

\chapter{The $s$-Weak Order and $s$-Permutahedra}
\label{chap:s-weak}

\chapcitation{I just knew there were stories I wanted to tell.}{Octavia E. Butler}

As I was working on the $m$-Tamari lattices by the end of my thesis, I wondered what was the equivalent of the weak order in the $m$-Tamari case. This is indeed a very natural question considering the many links between the weak order and the Tamari lattice in combinatorics, geometry and algebra. I studied the \defn{$m$-permutations}, \emph{i.e.}, the permutations with $m$ copies of each letter. I was interested in the lattice properties and at the same time, Novelli and Thibon were working on $m$ versions of the Hopf algebraic structures. In this context, the notion of \defn{metasylvester congruence} was introduced in~\citeext{NT20}. This lead me to define the \defn{metasylvester lattice} in~\citemeconf{Pon15}. 

Even though I am very interested in all the geometric aspects related to my work, my own background and approach is very much on the combinatorial side. In 2014, I could ``see'' that this new lattice of mine was ``nice'' and looked like a permutahedron that had been cut into smaller pieces, but I had no idea how to formalize that property. Around this time, I attended a talk by Cesar Ceballos at \emph{Séminaire Lotharingien de Combinatoire} where he presented the results he obtained with Padrol and Sarmiento on the realization the $\nu$-Tamari lattices~\citeext{CPS16}. I went to talk to him right after to show him my own work and see if we could work on a geometric realization together. This resulted in a fruitful collaboration which is still going on with many fascinating results and open questions. We presented a first version of our work at FPSAC 2019~\citemeconf{CP19}. In 2022, we released the first paper~\citemepre{CP22} of a series of at least two. It covers the combinatorial aspects of~\citeext{CP19} and is a generalization of~\citemeconf{Pon15}. We are working at the moment on the second paper~\citemepre{CP23} where we cover the combinatorial definition of the $s$-permutahedra complex. In this Chapter, I present an overview of all these results as well as the many different leads for future work. Computations and examples are available here~\citemesoft{PonSage22b}. 

\section{The metasylvester lattice}

\begin{figure}[ht]
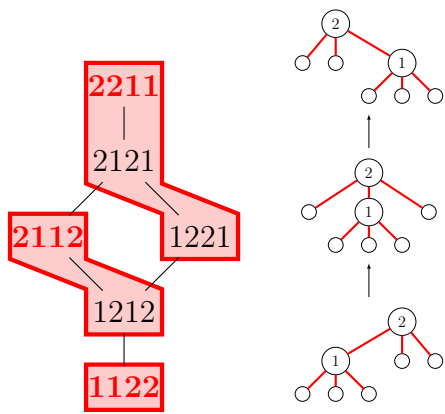

\center

\begin{tabular}{cc}
\input{includes/figures/mperms-2}
&
\scalebox{.5}{\input{includes/figures/sweak_22}}
\end{tabular}

\caption{The lattice of $2$-permutations (left) with metasylvester classes and the corresponding $s$-weak order.}
\label{fig:m-perms}
\end{figure}

As the $m$-Tamari lattice is an upper ideal of the Tamari lattice of size $n \times m$, it is a sublattice and a quotient lattice of a certain upper ideal in the weak order. It is actually more natural to look at the reversed $m$-Tamari which is a sublattice of the weak order on permutations where every letter is repeated $m$ times as in Figure~\ref{fig:m-perms}.

The \defn{metasylvester congruence}~\citeext{NT20} is a monoid congruence defined as the transitive closure of the following rules

\begin{align}
\dots ac \dots a &\equiv \dots ca \dots a \\
\dots b \dots ac \dots b \dots &\equiv \dots b \dots ca \dots b \dots
\end{align}
where $a < b < c$. As an example, the classes on Figure~\ref{fig:m-perms} are the red shapes. As both local moves require to have a letter repeated, when $m = 1$, on the classical weak order, each permutation is its own class. Besides, a \defn{sylvester} class as defined by~\eqref{eq:sylvester} in~Section~\ref{sec:asso-quotient} is a union of metasylvester classes. 

The metasylvester congruence is not a lattice congruence. Indeed, it is not compatible with the meet operation. For example, take $\sigma = 121332$ and $\mu = 131223$. We have that $\sigma \equiv 332112$ and $\mu \equiv 311223$. But $\sigma \meet \mu = 112323$ is not equivalent to $332112 \meet 311223 = 311223$. Nevertheless, the join operation is compatible and actually the maximal elements of the class have a nice characterization. They are \defn{Stirling permutations}, \emph{i.e.}, permutations avoiding the pattern $121$. They form a sublattice of the weak order and are in bijections with $m$-decreasing trees. In~\citemeconf{Pon15} I call this lattice the metasylvester lattice. This is actually a special case of the $s$-weak order which I defined with Ceballos in~\citemeconf{CP19} and~\citemepre{CP22}.

\section{The $s$-weak order}

\begin{figure}[ht]
\center

\input{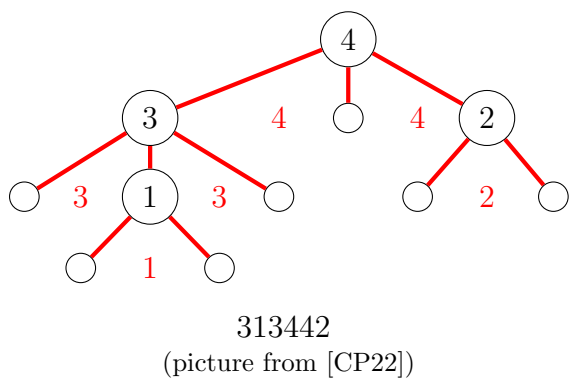}

\footnotesize{(picture from~\citemepre{CP22})}

\caption{An $s$-decreasing tree and the corresponding $s$-permutation.}
\label{fig:sdtree-sper}
\end{figure}

Let $s$ be a sequence of $s$ non negative integers. We define the \defn{$s$-decreasing trees}. They are ordered rooted trees on $n$ nodes labeled $1$ to $n$ such  that the node $i$ has $s(i) + 1$ children and so that the labels are decreasing from the root to the leaves. In particular, there are $(1 + s(n))(1 +s(n) + s(n -1)) \cdots (1 + s(n) + \dots + s(2))$ $s$-decreasing trees for a given sequence $s$. When $s$ does not contain any zeros, then $s$-decreasing trees are in bijection with $s$-permutations: permutations where each letter $i$ appears $s(i)$ times and avoiding the pattern $121$. The bijection is illustrated on Figure~\ref{fig:sdtree-sper}. 

When $s = (m,m\dots, m)$, the $s$-permutations are the maximal elements of the metasylvester classes and so we have a lattice structure. This generalizes easily for all $s$ without zeros: the $s$-permutations form a sublattice of the weak order. This is not the case when $s$ contains zero but we can still extend the definition. We generalize the notion of \defn{tree inversion} which was already introduced in~\citemeconf{Pon15}. Basically, the tree-inversions are a \defn{multiset} where each inversion $(b,a)$ appears with a certain \defn{cardinality} $\card(b,a)$. The cardinality of $(b,a)$ expresses where is $a$ regarding $b$ in the tree. If it is to the left or in the first subtree of $b$, then $\card(b,a) = 0$. If it is to the right of $b$ or in the last subtree, we have $\card(b,a) = s(b)$. The intermediate values correspond to the case where $a$ belongs to a \defn{middle child} of $b$. We show an example on Figure~\ref{fig:dtree}.

\begin{figure}[ht]
\center

$s = (0,0,2,1,3)$

\begin{tabular}{clc}
\multirow{4}{*}{\scalebox{.8}{\input{includes/figures/example_dtree}}}&
$T_0^5 = \begin{aligned}\scalebox{.4}{\input{includes/figures/dtrees/T05}}\end{aligned} $ & (strict) left child \\ \cline{2-3}
 & $T_1^5 = \emptyset $ & \multirow{2}{*}{middle children} \\
 & $T_2^5 = \begin{aligned}\scalebox{.4}{\input{includes/figures/dtrees/T25}}\end{aligned}$ & \\ \cline{2-3}
 & $T_3^5 = \emptyset $ & (strict) right child \\
\end{tabular}

\vspace{.5cm}
\begin{tabular}{llll}
$\card(5,4) = 0$ & $\card(5,3) = 2$ & $\card(5,2) = 2$ & $\card(5,1) = 2$ \\
                 & $\card(4,3) = 1$ & $\card(4,2) = 1$ & $\card(4,1) = 1$ \\
                 &                  & $\card(3,2) = 2$ & $\card(3,1) = 2$ \\
                 &                  &                  & $\card(2,1) = 0$
\end{tabular}

\footnotesize{(picture from~\citemepre{CP22})}

\caption{An $s$-decreasing tree and its tree-inversions.}
\label{fig:dtree}
\end{figure}

In particular, when $s=(1,1,\dots,1)$, then $s$-decreasing trees are decreasing binary trees, in bijection with permutations. In this case, the tree-inversions are the classical value inversions of permutations: each cardinality $\card(b,a)$ can either be $0$ (no inversion, $a$ is before $b$) or $1$ ($b$ is before $a$). Another interesting case is when $s$ contains some zeros: if $s(b) = 0$, then the node $b$ has one child and $\card(b,a) = 0$ for all $a < b$, independently of the position of $a$. 

\begin{figure}[htbp]
\begin{center}
\begin{tabular}{cc}
\includegraphics[height=8cm]{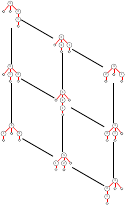} &
\includegraphics[height=8cm]{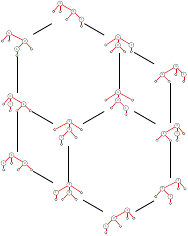} \\
$s = (0,0,2)$ & $s = (0,1,2)$ \\
\includegraphics[height=8cm]{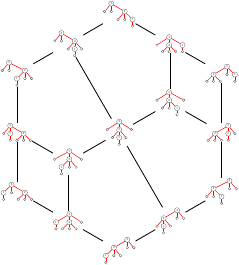} &
\includegraphics[height=8cm]{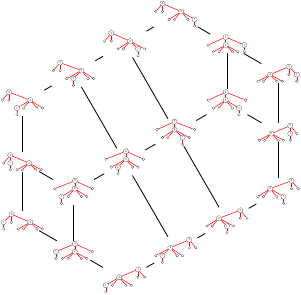} \\
$s = (0,2,2)$ & $s = (0,3,2)$
\end{tabular}

\footnotesize{(picture from~\citemepre{CP22})}

\caption{Examples of $s$-weak lattices.}
\label{fig:sweak_lattice}
\end{center}

\end{figure}

The tree-inversions form a multiset and we can generalize the notion of inclusion of classical sets. For two multisets $T$ and $S$, we have $T \subseteq S$ if, for all $a < b$, $\card_T(b,a) \leq \card_S(b,a)$. This gives a partial order on $s$-decreasing trees which we call the $s$-weak order as it is a direct generalization of the classical weak order. See some examples on Figure~\ref{fig:sweak_lattice}.

\begin{theorem}[Theorem~1.21 of~\citemepre{CP22}]
\label{thm:s-weak-lattice}
The $s$-weak order is a lattice.
\end{theorem}

\begin{figure}[htbp]
\begin{center}
\includegraphics[width = 0.5\textwidth]{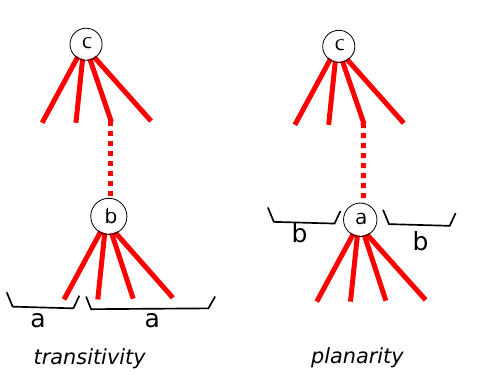}

\footnotesize{(picture from~\citemepre{CP22})}
\caption{Illustration of the transitivity and planarity conditions on $s$-tree inversion sets.}
\label{fig:transitivity_planarity}
\end{center}
\end{figure}

Our proof is based on a generalization of the notion of a \defn{transitive} relation to a multiset of inversions: it is transitive if for all $c > b > a$, $\card(b,a) > 0$ implies that $\card(c,a) \geq \card(c,b)$. In the classical case where cardinalities are either~$0$ or~$1$, this indeed says that if $(b,a)$ and $(c,b)$ are inversions, we have $\card(c,a) = 1$: $(c,a)$ is also an inversion. A multiset of inversions is an \defn{$s$-tree inversion set} if it is \defn{transitive} and \defn{co-transitive} (the non-inversions are also transitive). We actually call this second property the \defn{planarity} of the set, it translates to: for all $c > b > a$ such that $\card(c,a) =i$, either $\card(b,a) = s(b)$ or $\card(c,b) \geq i$. It is easy to see that the tree-inversions of an $s$-decreasing tree satisfy these properties as illustrated on Figure~\ref{fig:transitivity_planarity}. We prove that they are sufficient properties and give an explicit construction to recover the tree from the tree-inversions.

To prove that the $s$-weak order is a lattice, we then generalize the notions of \defn{union} and \defn{transitive closure} to multisets of inversions. The join of two $s$-decreasing trees is obtained by taking the transitive closure of the union of their tree-inversions. We also characterize the \defn{tree-ascents}, \emph{i.e.}, the cover relations of the poset. Basically, a cover relation corresponds to moving a node $i$ to the right. There are some technicalities in \emph{when} and \emph{how} we can do this move: nodes move by taking along their \defn{middle children} but leave behind their \defn{left child}. The nodes are not allowed to move if they have a \defn{strict right child}. For example, on the tree of Figure~\ref{fig:dtree}, $1$ and $2$ are allowed to move to the right of $5$ but not $3$ because $2$ is a strict right child of $3$ (while $1$ is not a strict right child of $2$ because $s(2) = 0$). 

Furthermore, we prove that the $s$-weak order is \defn{polygonal}, \defn{semidistributive} and \defn{congruence uniform}. These are all interesting lattice properties that are defined in particular in~\citeext{Rea16}. It implies some properties on the the lattice congruences of the $s$-weak order. Our proof is based on the notion of~$\HH$ lattices found in~\citeext{CCPBM04}. It also uses some useful properties of the intervals of the $s$-weak order found in~\citeext{Lac22} who studies some topological properties of the lattice based on our initial extended abstract~\citemeconf{CP19}.

\section{Link with the $\nu$-Tamari lattice}

In~\citemeconf{Pon15}, I claimed that the $m$-Tamari lattice was both a sublattice and a quotient lattice of the metasylvester lattice. In~\citemepre{CP22}, we generalize and prove this result on the $s$-weak order. We define the \defn{$s$-Tamari trees}: there are the trees such that $\card(c,a) \leq \card(c,b)$ for all $a < b < c$. This means that a node $a$ cannot move faster to the right than a node $b$ of greater value. 

\begin{theorem}[Theorem~2.2 of~\citemepre{CP22}]
\label{thm:s-tam}
The set of $s$-Tamari trees form a sublattice of the $s$-weak order.
\end{theorem}

\begin{figure}[ht]
\begin{center}
\begin{tabular}{cc}
\scalebox{1.5}{\input{includes/figures/sweak_s023_stam.tex}}
&
\includegraphics[width=0.5\textwidth]{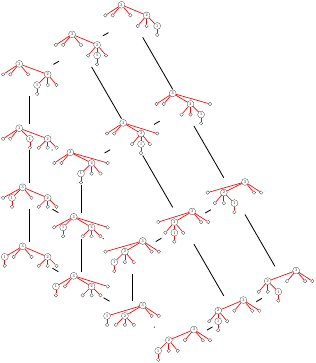}
\end{tabular}

\footnotesize{(picture from~\citemepre{CP22})}
\end{center}
\caption{$s$-decreasing trees and $s$-Tamari trees for $s=(0,2,3)$. On the left picture, the non $s$-Tamari trees are in black.}
\label{fig:s-Tamari_trees023}
\end{figure}

This is illustrated on Figure~\ref{fig:s-Tamari_trees023}. We give explicit descriptions of the cover relations and prove that this is isomorphic to a certain $\nu$-Tamari lattice as defined in Section~\ref{sec:nu-tam} where $\nu$ is simply the reversed sequence $s$.

Besides, when $s$ does not contain any zeros, then the $s$-Tamari lattice is a quotient lattice of the $s$-weak order and we give direct characterization of the congruence classes. In particular, the minimal elements of the classes are the $s$-decreasing trees.

\section{Pure intervals and $s$-permutahedra}

\begin{figure}[ht]
\center
\includegraphics[scale=.5]{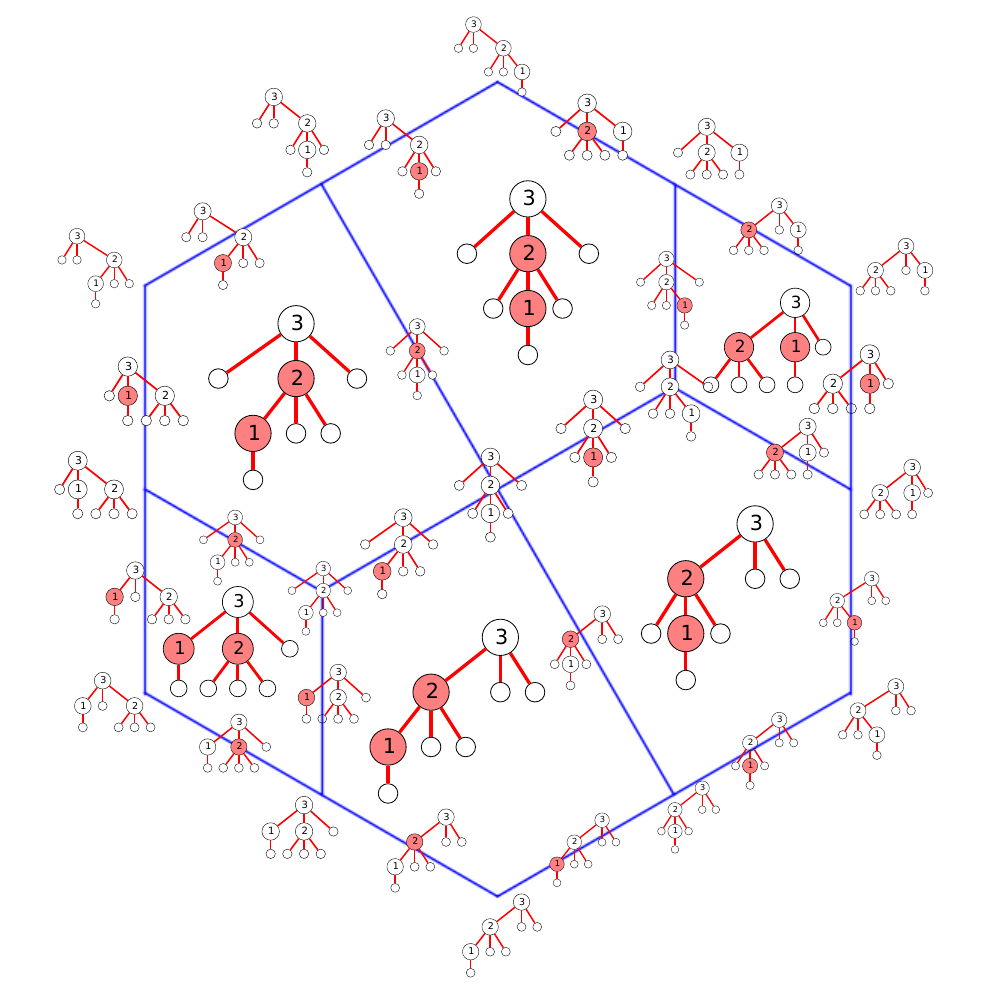}

\caption{The $s$-permutahedron complex for $s = (0,2,2)$.}
\label{fig:s-permutahedron}
\end{figure}

Looking at the examples of Figure~\ref{fig:sweak_lattice}, a geometrical structure seems to appear. Indeed, the lattices can be drawn in a way that ``cuts'' the plane into smaller polygonal areas. These polygons are glued together in a ``nice'' way and form a bigger polygon. All examples of Figure~\ref{fig:sweak_lattice} are for $n = 3$. For bigger $n$, we notice the same phenomenons of cutting a polytope of dimension $n-1$ into smaller polytopal cells. This is what we call a \defn{polytopal complex}. The road from \emph{seeing} to \emph{proving} turned out to be a much longer path than we initially thought and at this date, we still have no final proof that the $s$-weak order defines a polytopal complex. Nevertheless, in~\citemeconf{CP19} and~\citemepre{CP23} we define what we call a \defn{combinatorial complex}: \defn{the $s$-permutahedron} as shown on Figure~\ref{fig:s-permutahedron}.

\begin{figure}[h!]
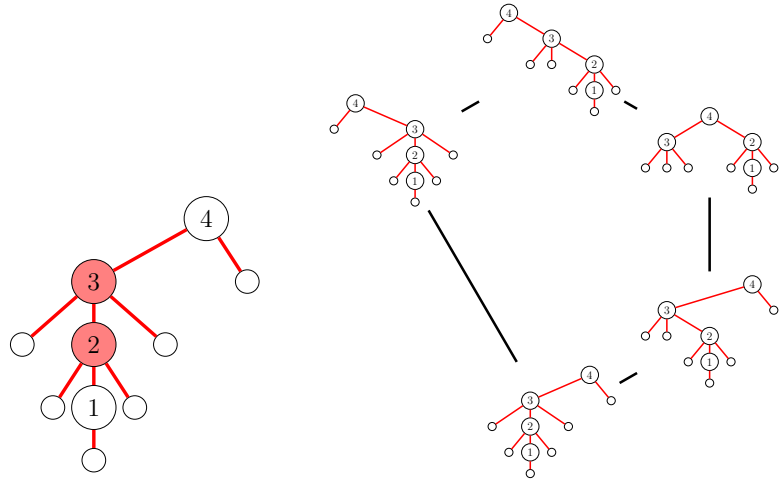

\center
\scalebox{.8}{
\begin{tabular}{cc}
\input{includes/figures/pure_face}
&
\input{includes/figures/example_pure}
\end{tabular}
}
\caption{Example of a pure interval. On the left, the minimal tree with colored tree-ascents and on the right the corresponding interval.}
\label{fig:example-pure}
\end{figure}

In~\citemepre{CP22}, we proved that the $s$-weak order is a \defn{polygonal lattice}. This means that if a tree $T$ is covered by two trees $S$ and $S'$, then the interval between $T$ and $S \join S'$ is a \defn{polygon}: it is formed by two chains that meet only in their extremities. The polygons that appear have been characterized by Lacina~\citeext{Lac22}: there are either squares, hexagons or pentagons. We believe that that this generalizes to other dimensions: if a tree $T$ is covered by $k$ trees $S_1, \dots, S_k$, then the interval between $T$ and $\join_{i=1}^k S_i$ is the skeleton of a $k$ dimensional polytope which is a generalized permutahedra.

This is what we call a \defn{pure interval} of the $s$-weak order. More precisely,  a pure interval is defined by an $s$-decreasing tree $T$  with a set of selected tree-ascents. Each tree-ascent gives a covering element $S$ and the pure interval is the interval between $T$ and the join of the covering elements $S$. We show an example on Figure~\ref{fig:example-pure}. Our paper~\citemepre{CP23} studies the combinatorial properties of the pure intervals. The $s$-permutahedra is the set of pure intervals on which we can define notions of \defn{containment}, \defn{intersection} and \defn{dimension}.

The ``dimension'' of a pure interval is the number of selected tree-ascents. At this stage, this is a formal notion which should later correspond to the actual dimension of the polytope. 

Not all intervals are pure intervals (just like not all intervals of the weak order are faces of the permutahedron). The first question we ask is then: if $T \wole S$, how can we tell if this corresponds to a pure interval? We prove that there is a simple characterization using the \defn{variations} of the interval, \emph{i.e.}, the inversions~$(b,a)$ such that $\card_T(b,a) < \card_S(b,a)$. We define a notion of \defn{essential variation} and show in particular that the selected tree-ascents of the pure interval are the \defn{minimal essential variations} of the interval. This characterization leads to  our main result in this paper.

\begin{theorem}[from~\citemepre{CP23}]
\label{thm:s-weak-pure}
The intersection of two pure intervals is a pure interval.
\end{theorem}

Besides, the dimension of the intersection is always smaller than or equal to the dimensions of the initial pure intervals. We reach equality only if one pure interval is already included in the other. You can check this property on Figure~\ref{fig:s-permutahedron}: the $2$-dimensional faces can intersect only in an edge (pure interval of dimension $1$) or in a tree (pure interval of dimension $0$). To prove that the $s$-permutahedra is indeed a polytopal complex, we need to prove further that each pure interval can be realized as a polytope. We have serious leads for this result which I explain in Section~\ref{sec:sweak-open}. We hope to develop them in a sequel paper.

\begin{figure}[ht]
\center
\includegraphics[scale=.35]{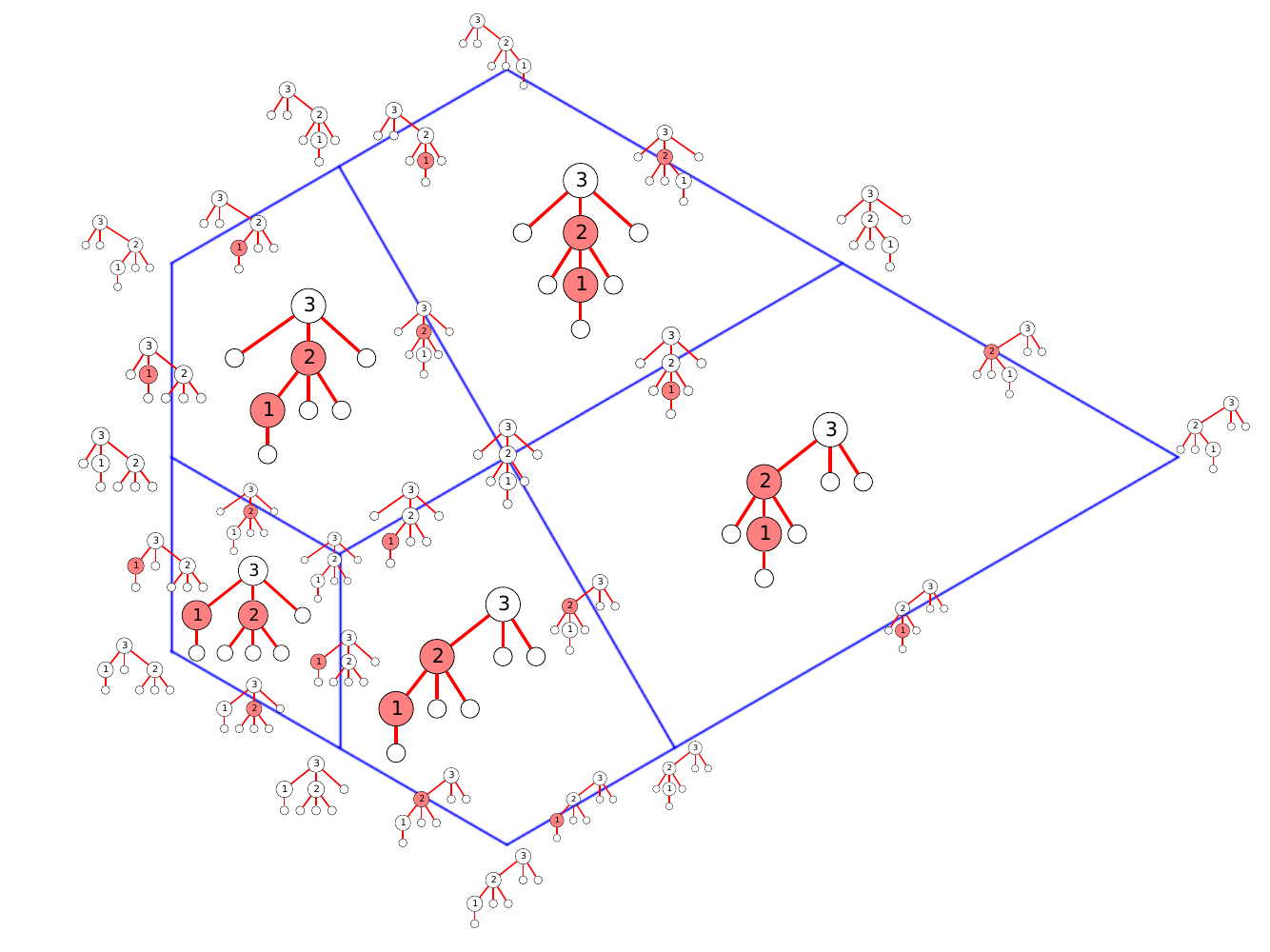}

\caption{The $s$-associahedron complex for $s = (0,2,2)$.}
\label{fig:s-asso}
\end{figure}

We can extend the definitions of pure intervals to the $s$-Tamari lattice and thus define the \defn{$s$-associahedron} (see Figure~\ref{fig:s-asso}). In this case, we can actually prove that it is a polytopal complex through a direct bijection with the $\nu$-associahedron of Ceballos, Padrol and Sarmiento~\citeext{CPS16}.

\begin{figure}[h!]
\begin{center}
\scalebox{.8}{
\begin{tabular}{cccc}
\includegraphics[height=2cm]{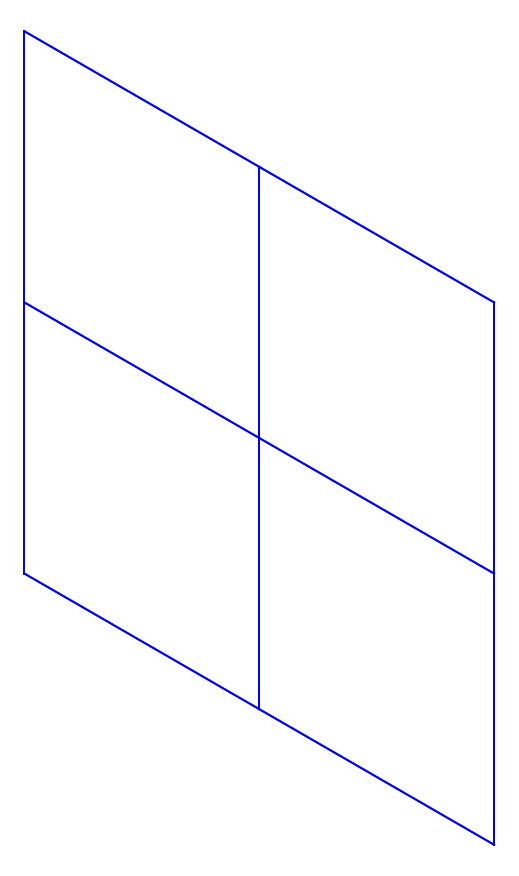} &
\includegraphics[height=2cm]{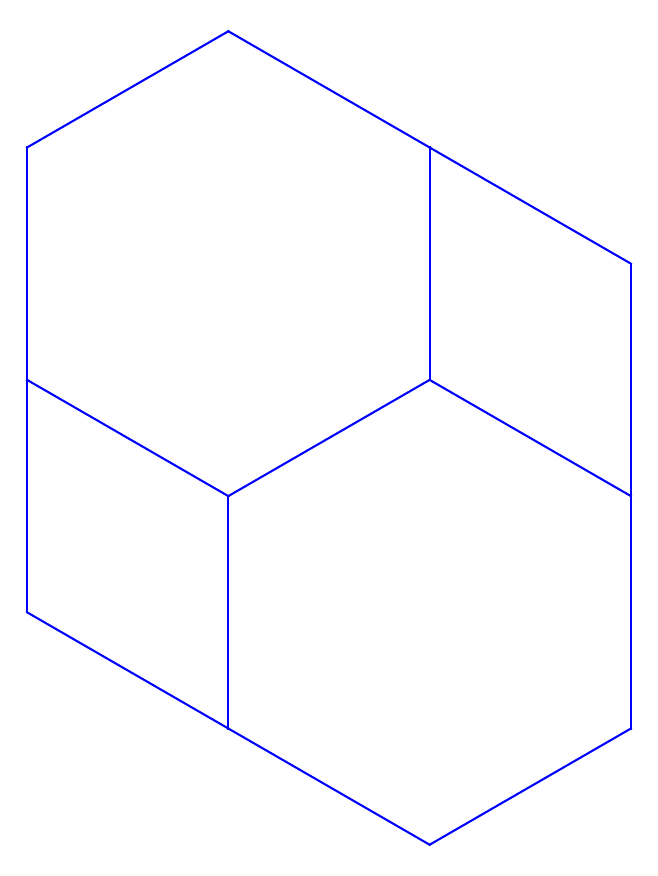} &
\includegraphics[height=2cm]{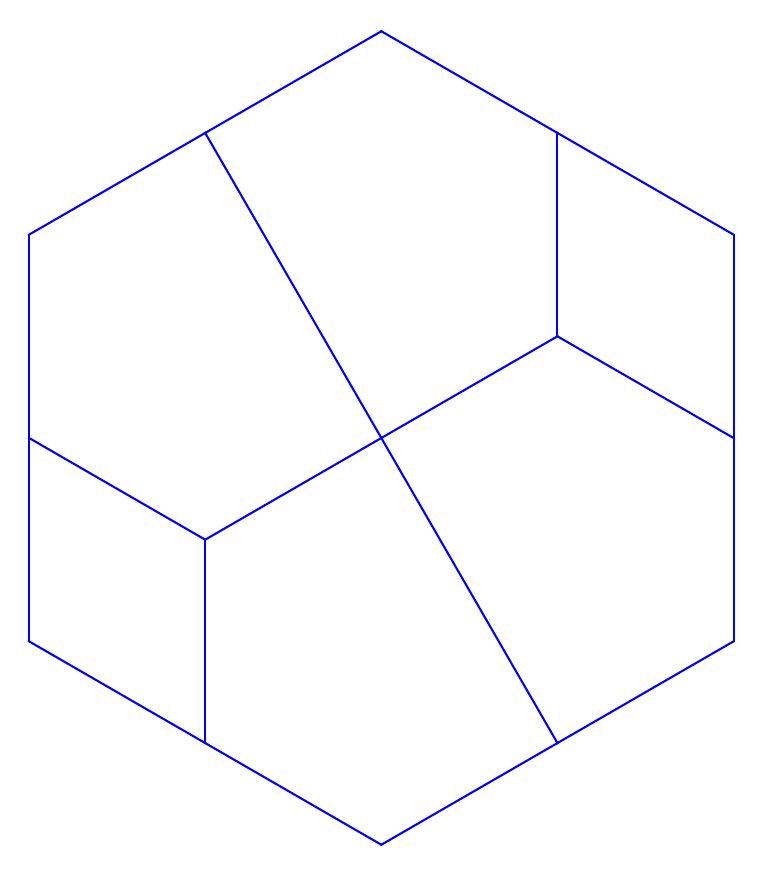} &
\includegraphics[height=2cm]{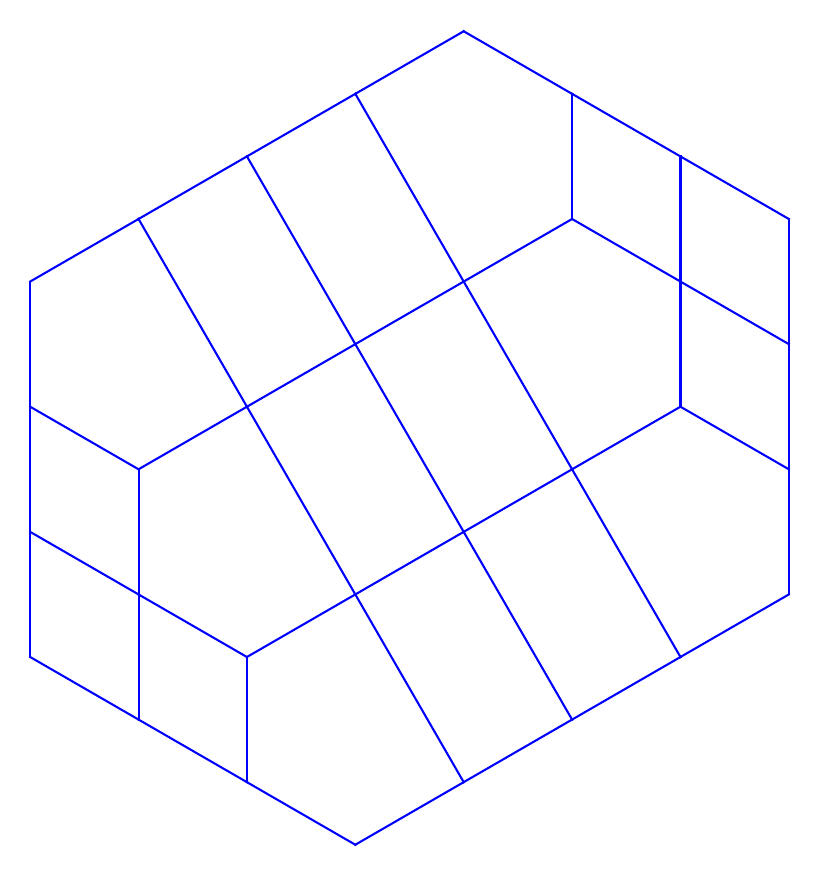} \\
$s = (0,0,2)$ & $s = (0,1,2)$ & $s = (0,2,2)$ & $s = (0,4,3)$   
\end{tabular}
}

\scalebox{.8}{
\begin{tabular}{cc}
\includegraphics[height=1.5cm]{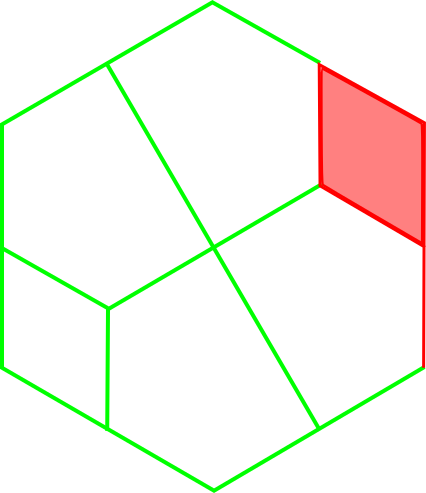} &
\includegraphics[height=1.5cm]{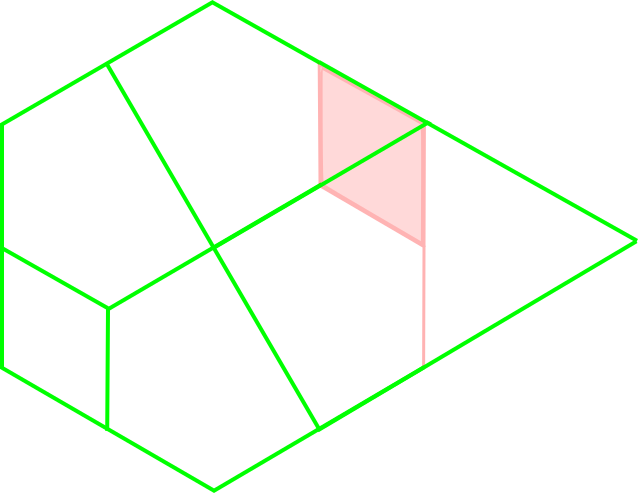} \\
$s = (0,2,2)$ & 
\end{tabular}
}

\scalebox{.8}{
\begin{tabular}{ccc}
\includegraphics[height= 3cm]{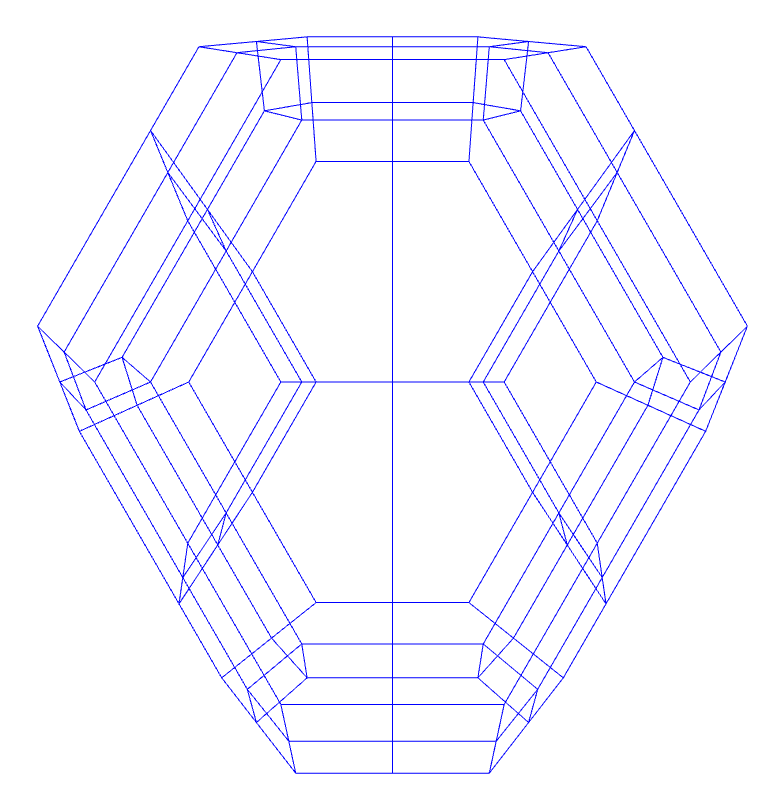}&
\includegraphics[height= 3cm]{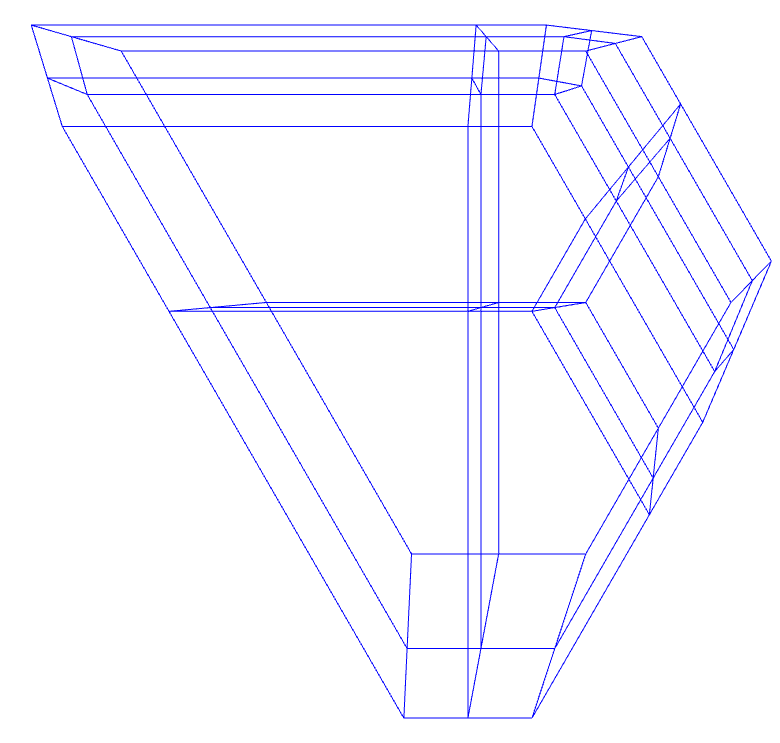}&
\includegraphics[height= 3cm]{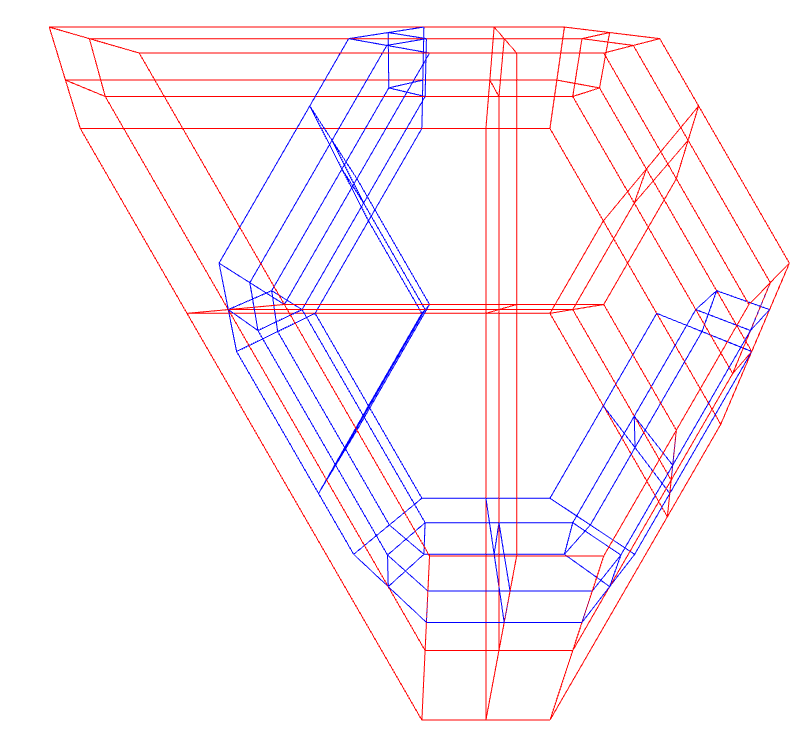} \\
$s = (0,2,2,2)$ & & \\ \hline
\includegraphics[height= 3cm]{includes/figures/s0333}&
\includegraphics[height= 3cm]{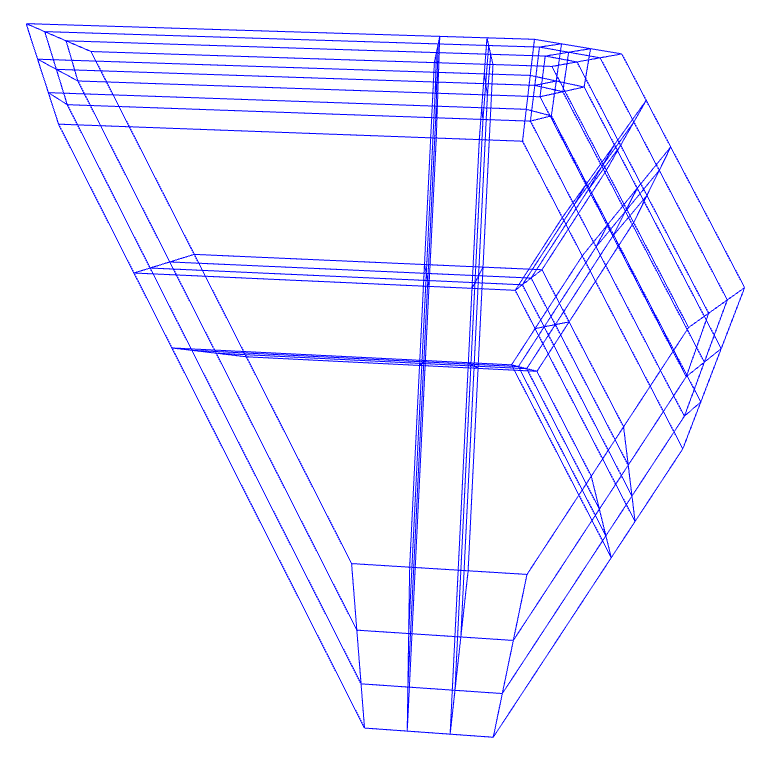}&
\includegraphics[height= 3cm]{includes/figures/s0333_both} \\
$s = (0,3,3,3)$ & & 

\end{tabular}
}

\footnotesize{(pictures from~\citemeconf{CP19})}
\caption{Some realizations for $s$-permutahedra and $s$-associahedra.}
\label{fig:s-permutahedron-geometry}
\end{center}
\end{figure} 

We express a conjecture at the end of~\citemeconf{CP19} which should also appear in~\citemepre{CP23}.

\begin{conjecture}[from~\citemeconf{CP19} and~\citemepre{CP23}]
\label{conj:s-permutahedron}
A geometric realization of the $s$-permutahedron as a polytopal subdivision of permutahedron exists in all dimensions.
\end{conjecture}

By removing some facets, one should obtain a realization of the $s$-associahedron. We have such a realization for dimensions 2 and 3 as we show on Figure~\ref{fig:s-permutahedron-geometry}. More can be found on \href{https://www.lri.fr/~pons/static/spermutahedron/}{this webpage}\footnote{https://www.lri.fr/\~{}pons/static/spermutahedron/}.

\section{Perspectives}
\label{sec:sweak-open}

Our work on the $s$-weak order and $s$-permutahedra is still in progress. We are finishing our second paper at the moment and hope to have at least $2$ more. Besides, our initial extended abstract~\citemeconf{CP19} has raised the interest of the scientific community. For example, Lacina~\citeext{Lac22} explored some topological properties of the lattice. At the moment, the group of~\citeext{GMPT23} is working on a positive solution of our conjecture. I present here the current state of our (and other's) research and the perspectives we would like to explore.

\subsection{Ascentopes}

\begin{figure}[ht]
\center
\begin{tabular}{cc}
\includegraphics[scale=1]{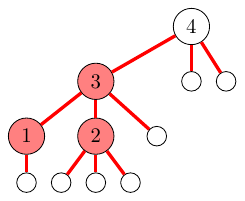}
&
\includegraphics[scale=.5]{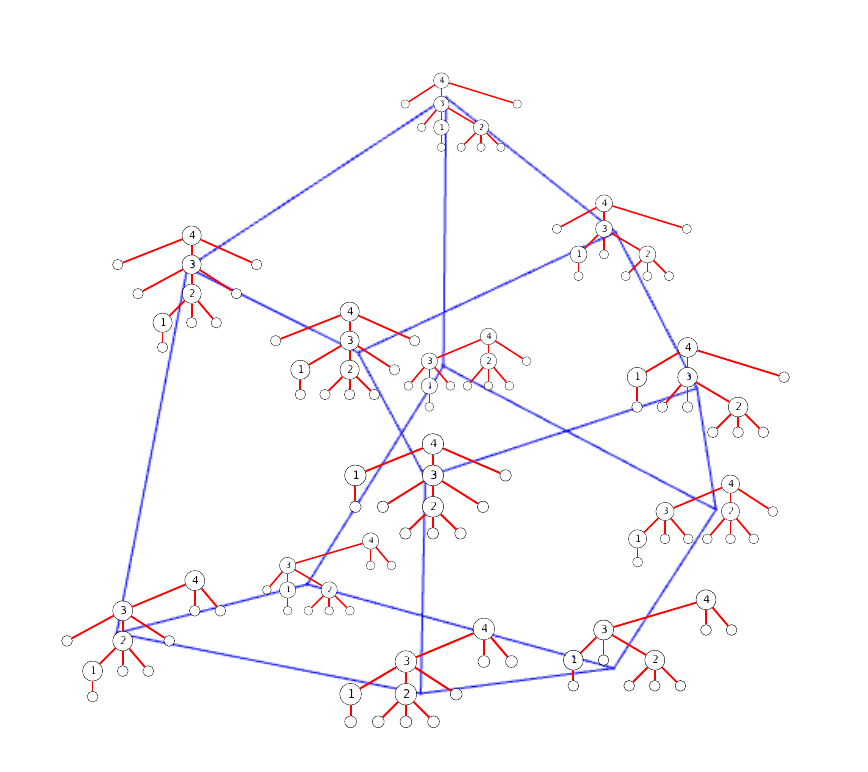}

\end{tabular}

\caption{A pure interval and its corresponding ascentope}
\label{fig:ascentope}
\end{figure}

\begin{figure}[ht]
\center
\includegraphics[scale=.8]{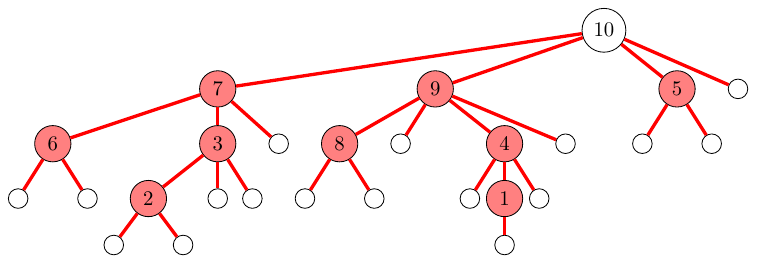}

$f$ vector = $(1, 2178, 9801, 19008, 20790, 14082, 6099, 1680, 282, 26, 1)$

\caption{A pure interval of dimension $9$ with the $f$ vector of its ascentope.}
\label{fig:bigface}
\end{figure}

We do not prove yet in~\citemepre{CP23} that the $s$-permutahedra is a polytopal complex. Nevertheless, we actually have a polytopal construction for each pure interval which we call the \defn{ascentope}. An ascentope is a certain removahedron defined using an $s$-decreasing tree $T$ and a set of selected tree-ascents. Experimentally, our construction gives a geometrical realization of the pure intervals. See Figure~\ref{fig:ascentope} for an example in $3$ dimensions. We can compute the skeleton of the polytope and check that it corresponds to the Hasse diagram of the interval. We also computed $f$-vectors for big examples such as the one of Figure~\ref{fig:bigface}. The $f$ vector gives the number of faces in each dimension. As we understand the faces combinatorially, we can check that the geometric results are consistent with our combinatorial construction: it is the case on all the examples we have computed.

The next step is obviously to prove that this construction is indeed a realization of the pure intervals and to study the properties of this family of polytopes. Note that the ascentope realization does not solve the more general conjecture which requires to construct the polytopal complex as a polytopal subdivision. Indeed, there is no indication that the ascentopes can be ``glued'' together. Nevertheless, it would prove that the $s$-permutahedra is indeed a polytopal complex for all $s$. 

\subsection{Solving the conjecture}

Since 2019 and~\citemeconf{CP19}, I have been mentioning the conjecture of the polytopal subdivision realization a few times. In particular, I had started discussing this with \'Eva Philippe who is a PhD student in Paris (under the supervision of Arnau Padrol). Around the same time, Rafael Gonz\'alez D'Le\'on launched a research project to strengthen the links between France and South American research in combinatorics. As I was traveling to Colombia for the ECCO conference in 2022, I planed a few extra days of research over there. I also presented the conjecture during an ``open problem'' session at the conference. 

Following this, Rafael Gonz\'alez suggested that we look at the problem through the lens of \defn{flow polytopes}. These are certain polytopes defined by \defn{flows} on a graph which can be used to obtain \defn{triangulations} and \defn{polytopal subdivisions}~\citeext{MM19}. This looked very promising and I invited \'Eva Philippe and my PhD student Daniel Tamayo to join the research discussions. Together with an extended group of researchers, they came up with a certain flow polytope which gave a dual realization of the $s$-permutahedra. Even though the proof is not finalized yet, they were able to show me explicit bijections and coordinates to construct the polytopal subdivision. Their paper is in writing~\citeext{GMPT23} and would partially solve the conjecture.

The construction of~\citeext{GMPT23} only works for an~$s$ which do not contain zeros. The next natural step would then be to see if the result can be extended to the zero case. Besides, we also conjectured that the realization should give a realization of the $s$-associahedron. Is it the case? Can it be proved using flows? What is the link with the current realization of Ceballos, Padrol and Sarmiento~\citeext{CPS16}? In both cases, Tropical geometry seem to be necessary to compute the vertices coordinates. Is there a realization which does not rely on this?

Besides, what are the polytopes we obtain for each pure interval? Are they removahedron or even generalized permutahedron? Are they related to the ascentopes that we have been defining with Cesar Ceballos?

Finally, the dual of the $s$-permutahedron seems to have interesting properties. It looks like a certain \defn{Shi-arrangement} where some walls have been removed. We would like to study this construction to better understand its properties.

\subsection{Lattice quotients and geometric realizations}

We have proved that the $s$-Tamari lattice is a lattice quotient of the $s$-weak order when $s$ does not contain any zeros. What are the other lattice quotients? In particular, Tamayo, Philippe and Pilaud have been working on some $s$-Permutrees with interesting combinatorial and geometrical properties. When $s$ contains some zeros, the Tamari case defines another lattice which is not isomorphic to $\nu$-Tamari: what is it? What are the elements counted by? Does it have an interpretation in terms of paths?

\subsection{Other groups and other lattices}

An natural question is also to look at other Coxeter groups. For this, the perspective of the generalized Shi-arrangement would be a good starting point. Do we also find polytopal subdivision? This is the case in type $B$ for the Ceballos, Padrol and Sarmiento realization of the $\nu$-Tamari lattices~\citeext{CPS16}. More generally, our approach through pure intervals works in many other combinatorial lattices with nice geometric realization. Is there a way we can extend our proof and characterizations to these other cases?

\part{Experimental Approach}
\label{part:exp}

\chapter{Epistemology of the Experimental Approach}
\label{chap:exp}

\chapcitation{To me programming is more than an important practical art. It is also a gigantic undertaking in the foundations of knowledge.}{Grace Hopper}

In 2010, I started a research internship at the end of my master degree that would eventually lead to my PhD. My first assignment was to attend a workshop on the open source mathematical software \Sage{}. At the time, I had no experience with \Sage{} nor even with the python language (on which \Sage{} is based). Nevertheless, with a strong background in computer science, I had no difficulties to learn the system. By the end of my internship, I had actually developed the first draft of what is now an external package on bases for multivariate polynomials~\citemesoft{PonSage16}.

Already, at this early stage of my research career, it felt natural to me to approach mathematical knowledge and research through experimentation and code. Indeed, it has since become an essential part of my research methodology. In other words: I have written some code related to all my research and all my papers. Computer exploration and experimentation is as essential to me as having a pen and paper. Moreover, I believe that, to a certain degree, experimentation is part of any mathematical research, whether a computer is used or not. On the other hand, as it is used on early stages of the research, it can often be absent of the final result. This creates a difficult problem for historians and philosophers who try to understand \emph{how} mathematical knowledge is created. For us mathematicians, and especially mathematicians who use computers, this creates another challenge: how do we value this work? In particular, how do we value the necessary technical work to create the tools that are needed?

A first step, in my opinion, is to make this work visible and to explain the experimental process in mathematics. Indeed, it has appeared to me that many people do not associate experimentation with mathematics and do not really understand what we mean by that. And probably, different mathematicians have different practices and mean different things. This has led me to have an introspective approach of my own methods and to discuss it with historians and philosophers. In this context I was invited by Emmylou Haffner to give a talk at a seminar of mathematical history and philosophy in SPHERE~\citememisc{PonMisc22}. Emmylou Haffner studies  ancient mathematican's drafts to understand the creative process~\citeext{Haf22} like in this post where she looks at Elie Cartan's drafts~\citeext{Haf17}. The seminar was coordinated by Arilès Remaki whose thesis focuses on Leibnitz~\citeext{Rem21} and especially on the combinatorial explorations that led to his work. Looking at Leibnitz's drafts with combinatorial pictures and computational tables makes an interesting connection with our own exploratory methods now enhanced by computer.

In this chapter, I use the research example that I presented in~\citememisc{PonMisc22} to explain my research approach and unfold the early stages of research creation through experimentation using both the computer and ``pen and paper''. Then I give an overview of the research code that I have written for the past ten years and of the good practices that I have put in place to share, reuse and valorize this work.

\section{The experimental process}
\label{sec:exp-process}

When I gave this talk at the SPHERE seminar, I was in the middle of a research project that led to what I presented in Chapter~\ref{chap:qt}. In particular, I was looking at a certain lattice structure on triangular Dyck paths which would give a combinatorial interpretation to some Schur expansions computed by François Bergeron. I decided to present this problem as a typical example of research investigation.

\subsection{To understand objects}

\begin{figure}[p]
\center
\begin{tabular}{ccc}
\includegraphics[height=10cm]{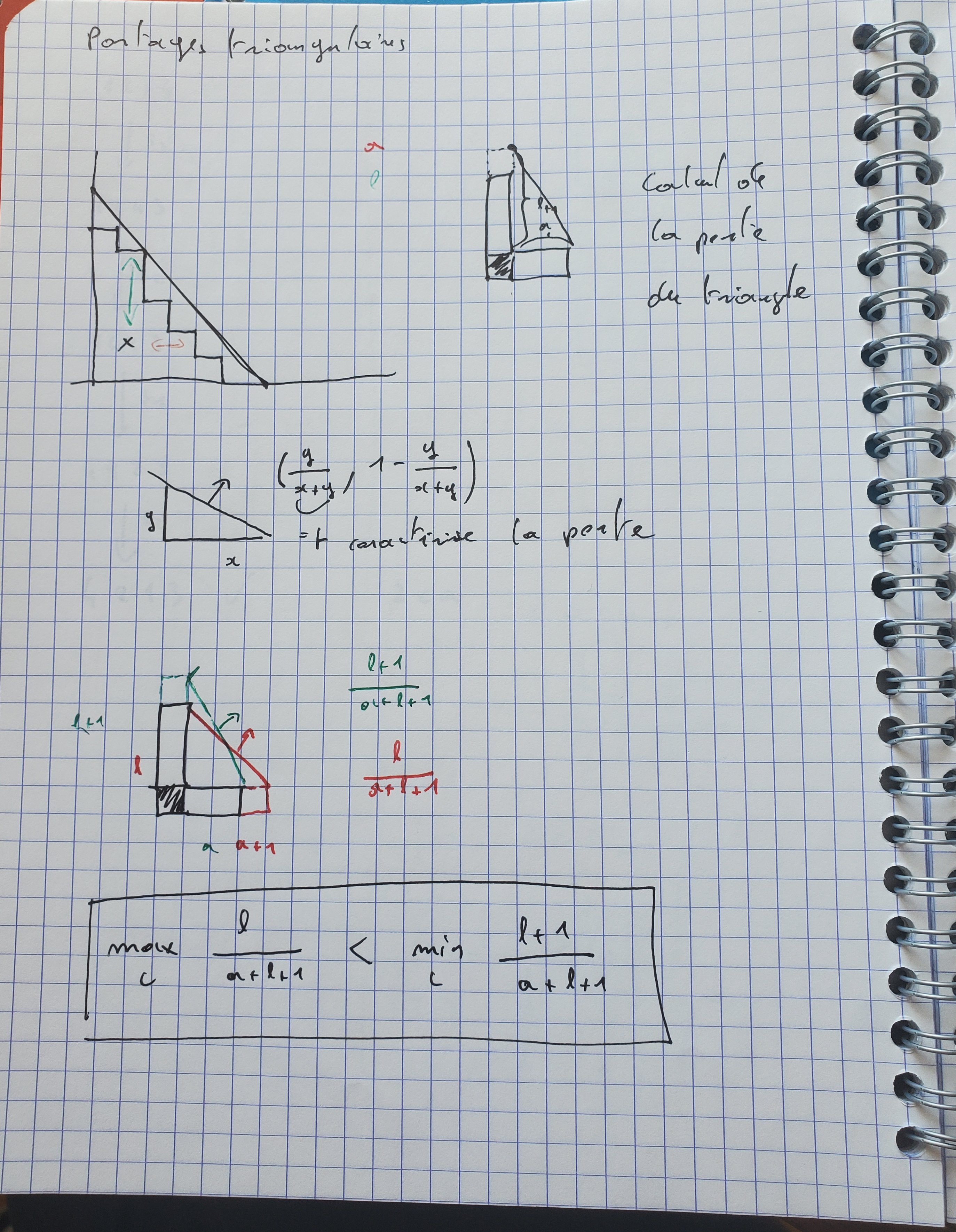}&
\includegraphics[height=10cm]{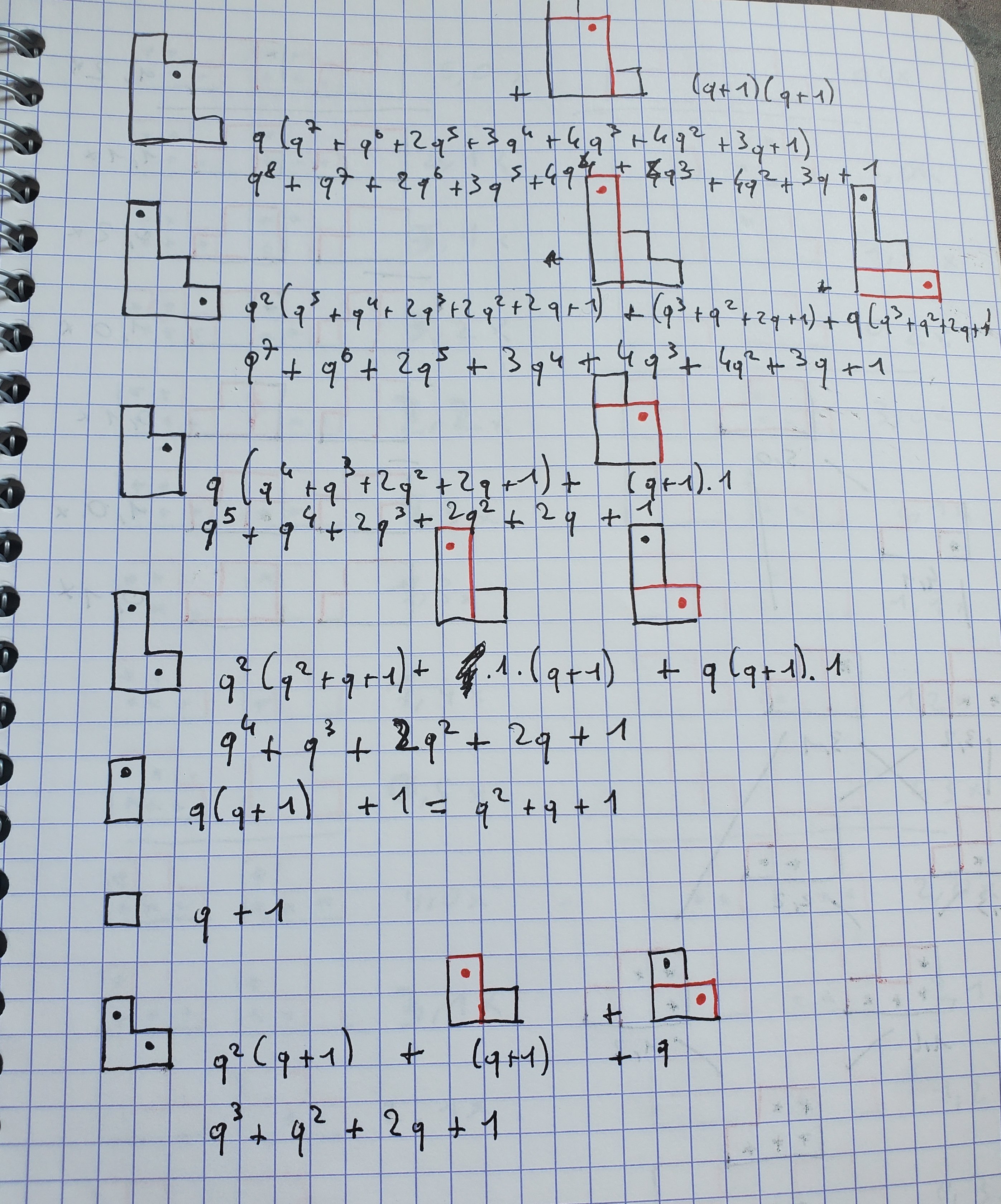}
\end{tabular}

\includegraphics[height=10cm]{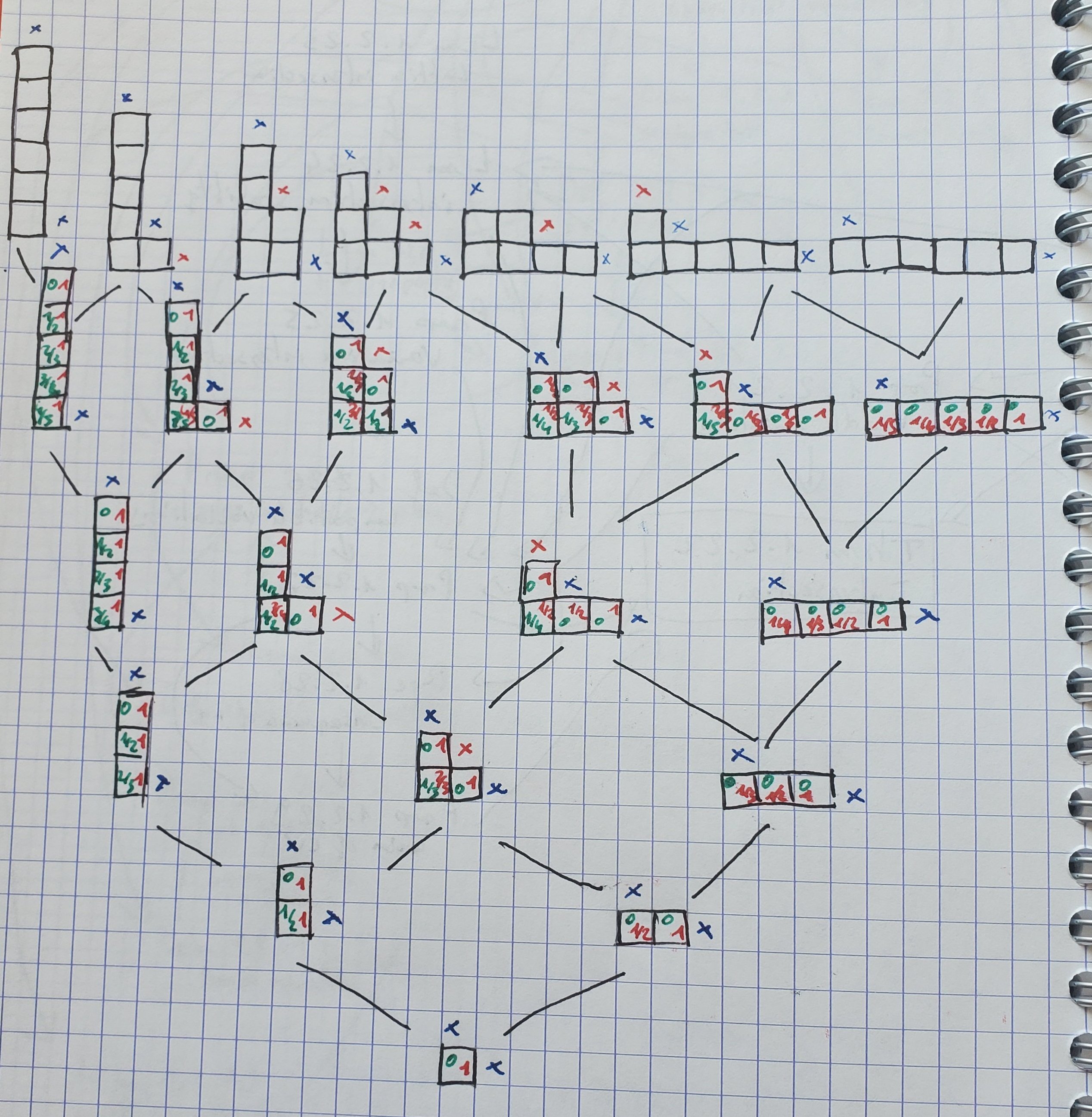}

\caption{Pictures of early drafts of triangular partitions and Dyck paths}
\label{fig:drafts-triangular}
\end{figure}

Looking at my drafts, it became apparent that my first task had been to ``understand'' the paper~\citeext{BM22} that Bergeron presented to me when I arrived in Canada. On Figure~\ref{fig:drafts-triangular}, I show some initial drafts kept in my written notes. We see that I am recomputing some of the examples of~\citeext{BM22} to familiarize myself with the material. I also started to code at this stage. I implemented several classes to represent triangular partitions and triangular Dyck paths along with methods to compute the properties described in the paper. For example, below is a piece of code which is part of the triangular partition class to test if another partition is \defn{similar} in the sense of~\citeext{BM22}

\begin{lstlisting}
    def is_sim_cell(self,i,j,tau):
        return self.cell_min_slope(i,j) < tau.mean_slope() \
        and tau.mean_slope() <= self.cell_max_slope(i,j)

    def sim_cells(self,tau):
        for i,j in self.cells():
            if self.is_sim_cell(i,j,tau):
                yield (i,j)

    def is_similar(self,tau):
        return all(self.is_sim_cell(i,j,tau) \
         for i,j in self.cells())
\end{lstlisting}

The constant back and forth between my manual computations and examples and the automated ones is what ensures me that I am progressing on solid grounds. The manual computation and often manual enumeration of objects gives me the time to understand the processes to the point where I can code and automate the task. Doing so assures me that I have not overlooked a crucial detail or misunderstood a definition. Besides, sometimes the error is not mine. In the early paper draft printed to me by Bergeron, a typo had replaced a ``$\leq$'' sign by a ``$<$'' sign: I was able to spot it when I tested the definition.

\subsection{To eliminate false statements and generate counter example}

\begin{figure}[ht]
\center

\includegraphics[height=10cm]{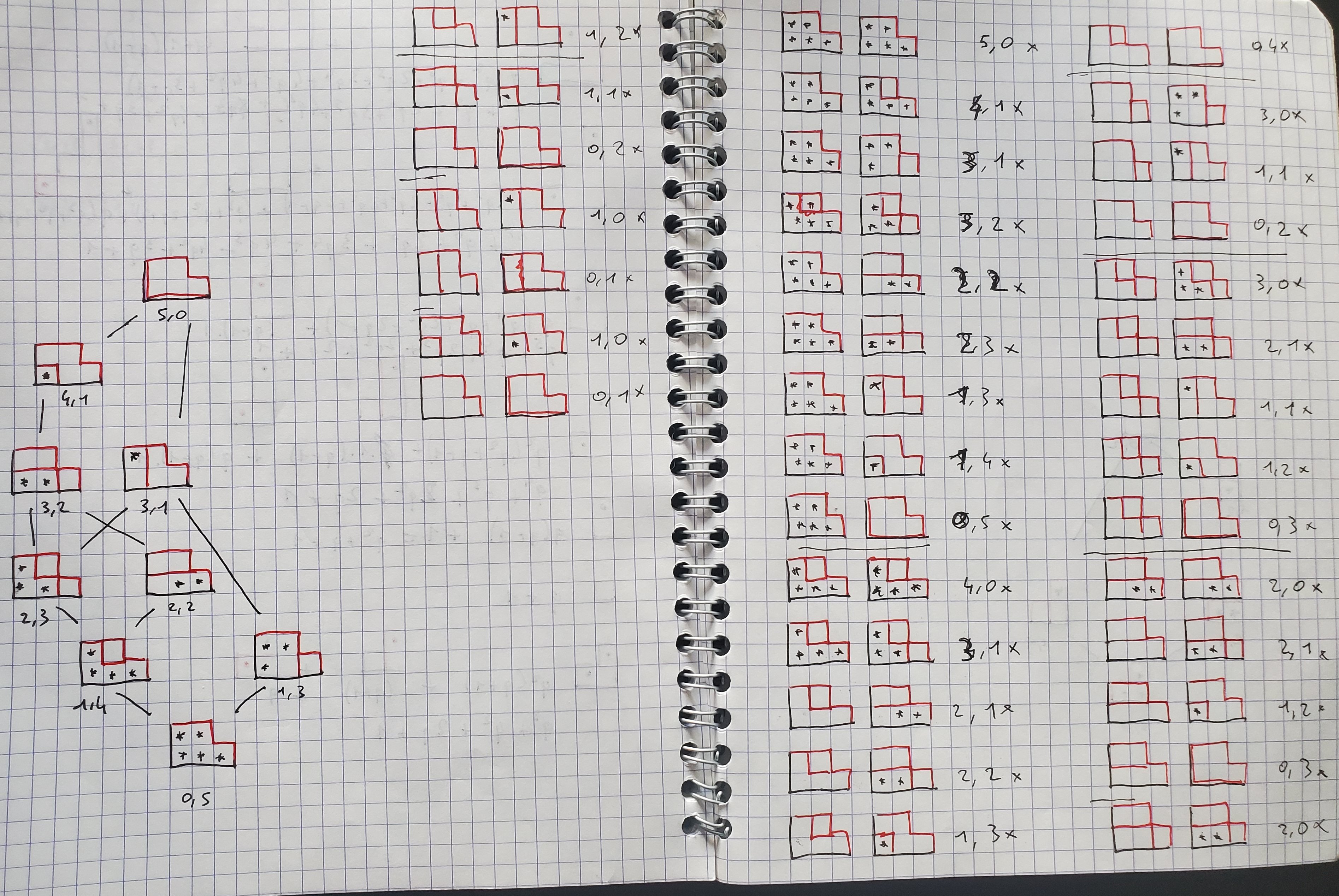}

\caption{A manual computation of $q,t$ symmetry on intervals of $\nu$-Tamari}
\label{fig:draft-nu-tam}
\end{figure}

Once I had read the paper, I wanted to study lattice structures on triangular Dyck paths to see if I could find a $q,t$ enumeration of intervals that would be symmetric and Schur positive (as in the classical and $m$-Tamari cases). As there was a general definition of Tamari in~\citeext{PRV17}, my first attempt was to look at it and see if it ``worked''. I manually computed the first non-trivial example (see Figure~\ref{fig:draft-nu-tam}) and noticed that, in this case, I had a $q,t$ symmetry.

Here, the manual computation helps me understand how the $\nu$-Tamari lattice translates into triangular partitions and gives me some initial insights on how the $q,t$ symmetry could work. Of course, the immediate next step is to implement this. Below are the methods I used to generate the $\nu$-Tamari lattice of triangular Dyck paths.

\begin{lstlisting}
    def path_tamari_rotate(self, line):
        L = self.skew_partition().row_lengths()
        p = list(self.path())
        v = L[line]
        p[line]-=1
        i = line-1
        while i >= 0 and L[i] > v:
            p[i]-=1
            i-=1
        return TriangularDyckPath(self.partition(), p)


    def path_tamari_up(self):
        for l,c in self.path().corners():
            yield self.path_tamari_rotate(l)
\end{lstlisting}

Once I had this, it quickly appeared that the $\nu$-Tamari lattice did \emph{not} give the desired $q,t$ interval enumeration. It is actually very common that the computer exploration is here to \emph{save us time} by rejecting false hypothesis early on. Not only does it tell us that some ideas are just wrong but it also exhibits interesting counter examples. In this case, the partition $(3,2,2,1)$ was the first truly interesting example where $\nu$-Tamari would fail to provide a symmetric $q,t$ enumeration of intervals. I could then concentrate my efforts on this particular partition. 

\begin{figure}[ht]

\includegraphics[height=10cm]{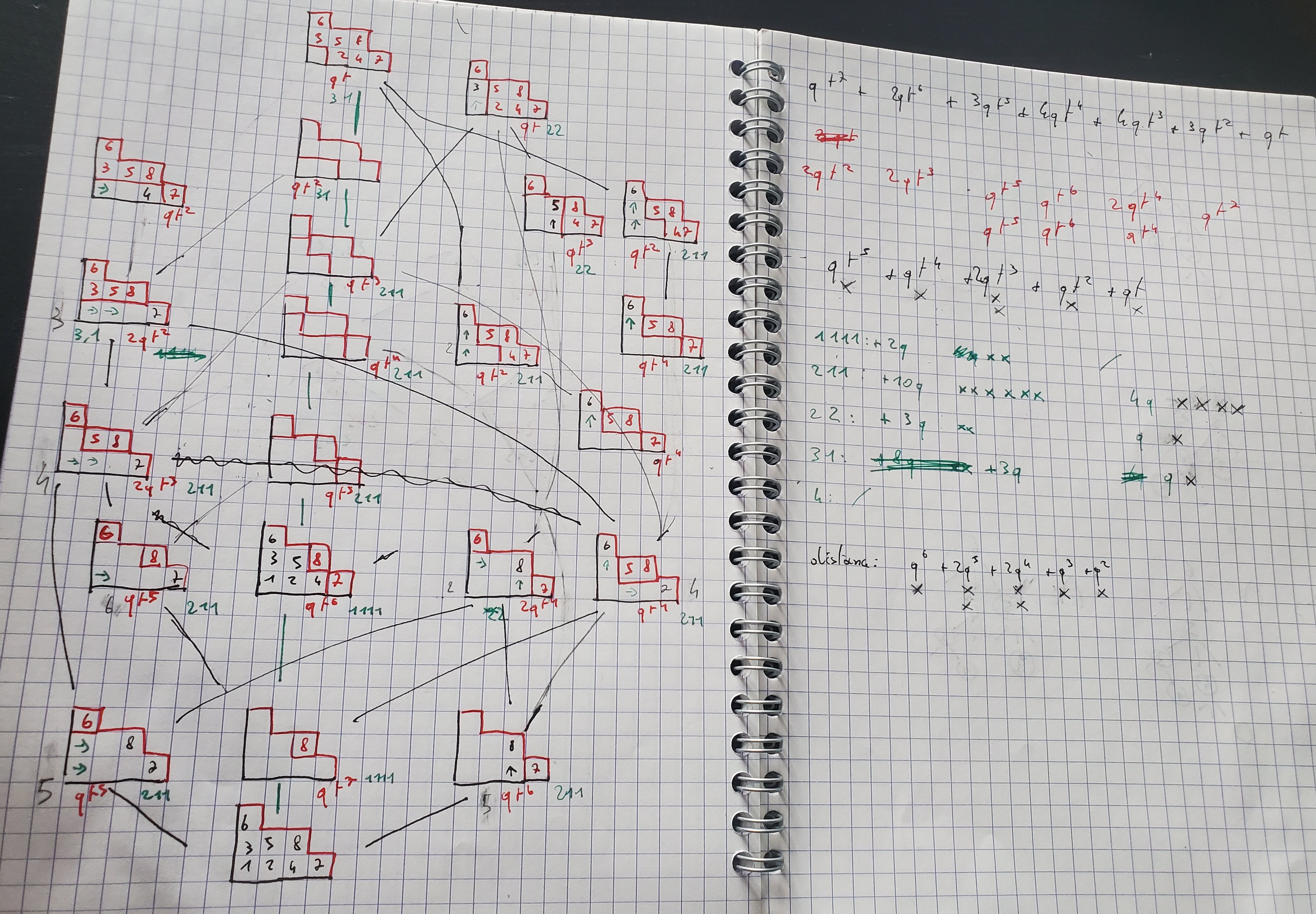}

\caption{Manual attempt to obtain the desired lattice structure on partition $(3,2,2,1)$.}
\label{fig:draft-lattice-3221}
\end{figure}

The partition $(3,2,2,1)$ has $23$ subpartitions and we expect to find $161$ intervals in the ``good'' lattice: $161$ is the result of a certain expansion of a Schur functions in $3$ variables $q,t,r$ and computing the polynomial for $q=t=r=1$. This is exactly a case which is ``in between'' manual and automatic exploration. Indeed, $23$ elements is rather small and manageable by hand while $161$ starts to be a bit too big. My research at this point illustrates that. I was literally looking for the right lattice by manually guessing what the correct cover relations might be using my intuition and some indication from the Schur computation (see Figure~\ref{fig:draft-lattice-3221}). I would then hard code the lattice inside my program and compute the $q,t$ enumeration of intervals.

Doing so, I was actually able to find a working structure for $(3,2,2,1)$ but could not deduce a general rule. But the multiple attempts has left me with a better knowledge of what works and what does not. 

\subsection{To find ideas and emit conjectures}

\begin{figure}[ht]
\center

\includegraphics[height=10cm]{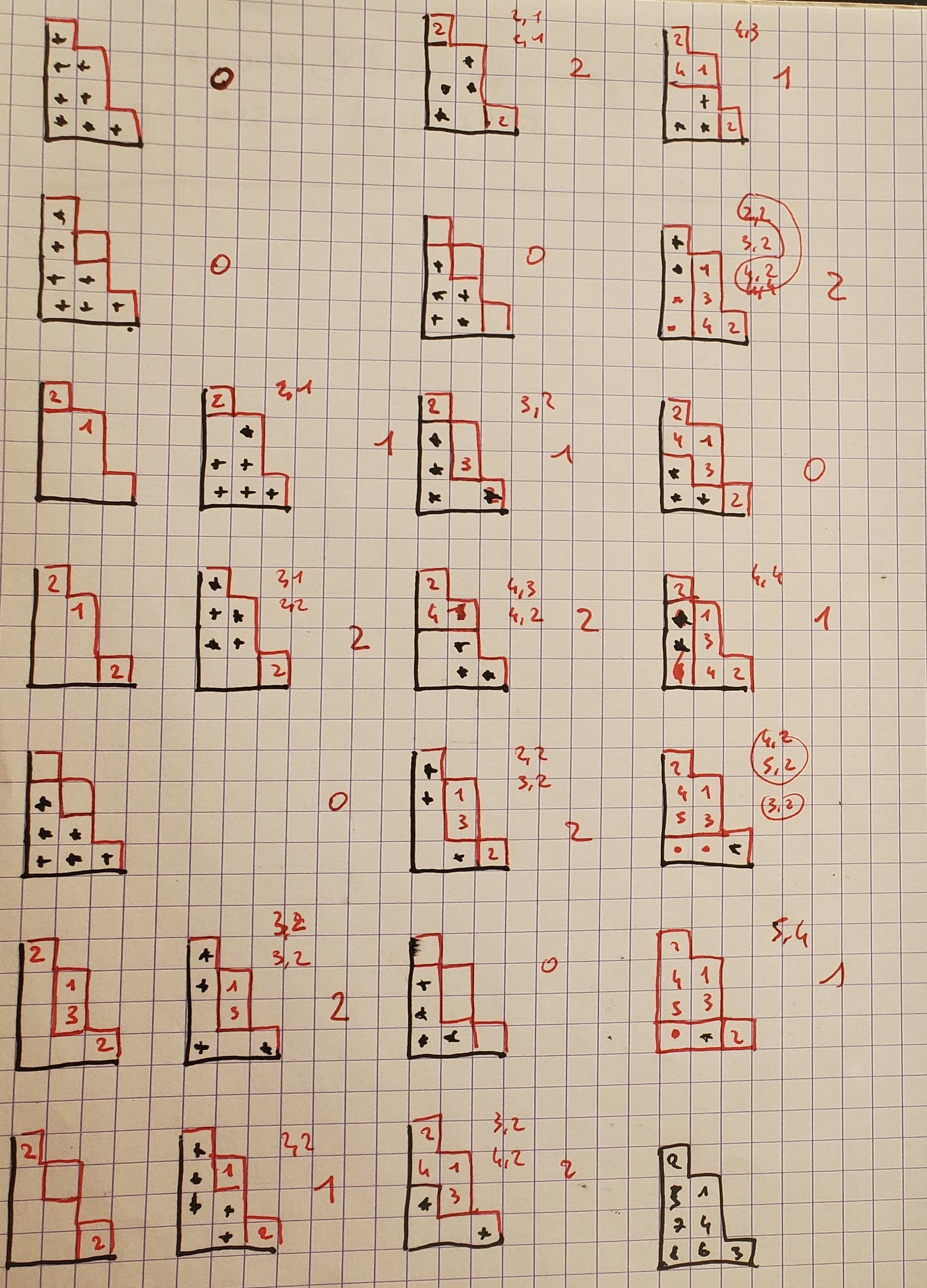}
\caption{Manual enumeration of triangular Dyck paths to define the tableau deficit.}
\label{fig:draf-deficit}
\end{figure}

Even though we are often wrong, we are thankfully sometimes right. There again, conjectures arise from a combination of manual and computer testing. By manually enumerating triangular Dyck paths along with their similar cells as in Figure~\ref{fig:draf-deficit}, I came up with the idea of the \defn{deficit} that I explain in Section~\ref{sec:deficit}. Only by computer exploration was I able to ``check'' that it was true with the following code.

\begin{lstlisting}
def test_deficit_cells(tdp):
    return all(tdp.is_sim_cell(i,j) != tdp.is_deficit_cell(i,j) \ 
    for i,j in tdp.path().cells())

def test_deficit_cells_partition(p):
    for tdp in p.triangular_dyck_paths():
        if not test_deficit_cells(tdp):
            return tdp
    return True

# tested 3..12
def test_all_partitions_deficit_cells(n):
    for p in TriangularPartitions(n):
        print(p)
        r = test_deficit_cells_partition(p)
        if r != True:
            return r
    return True
\end{lstlisting}

This was proved to be true by Lo\"ic Le Mogne in our paper~\citemeconf{LMP23} (this is an extended abstract, the complete proof will actually appear in our future complete paper).

Similarly, the final conjecture of~\citemeconf{LMP23} which I explain in Section~\ref{sec:triangular-intervals} has also been made possible by computer exploration. For example, the following function tests that the $\nu$-Tamari enumeration of intervals is Schur positive for all triangular partitions whose top-down tableau is \defn{sim-sym}.

\begin{lstlisting}
# tested 3..21
def test_schur_positive_PRVsimsym(n, q, t):
    for tr in TriangularPartitions(n):
        d = PRV_q_t_distribution(tr)
        if test_symmetry_distribution(d):
            print(tr)
            P = tr.path_tamari_lattice()
            d2 = PRV_poset_distribution(tr,P)
            pol = sum(d2[k]*q**k[0]*t**k[1] for k in d2)
            S = schur.from_polynomial(pol)
            if not all(c>0 for p,c in S if len(P) <= 2):
                return tr
    return True
\end{lstlisting}

The upper line is a comment which indicates which sizes have been tested. My research code is often full of testing functions with comments to keep track of positive / negative results. I also wrote the following function which tested that the Schur expansion I had corresponded to the actual Schur functions that had been sent to me by Bergeron.

\begin{lstlisting}
def test_PRV_polynomials(pols, q, t):
    wrong = []
    for key in pols:
        tr = TriangularPartition(key)
        d = PRV_q_t_distribution(tr)
        if test_symmetry_distribution(d):
            print(tr)
            P = tr.path_tamari_lattice()
            d2 = PRV_poset_distribution(tr,P)
            pol1 = sum(d2[k]*q**k[0]*t**k[1] for k in d2)
            pol2 = pols[key].expand(3)(q,t,1)
            if not pol1 == pol2:
                wrong.append(tr)
    return wrong
\end{lstlisting}

The interesting anecdote is that this function did not give a positive answer. Indeed, two of the Schur functions did not correspond. But it actually appeared that those Schur functions are very difficult to compute on the algebraic side. The functions I had received from Bergeron were actually missing some terms that he could not get at the time using his computational methods. I was able to \emph{guess} the missing terms using my conjecture. This not only strengthened my belief that the conjecture was true but also gave indications that the combinatorial computation (as inefficient as it is) would actually turn out to be more efficient than the algebraic computation which requires extensive computational power to reduce very big matrices.

\subsection{To solve and to prove}

As we saw, computer exploration helps us sort out positive and negative results before we even try to prove anything. This is of course useful at the investigation stage but not only. Indeed, once we have a result that \emph{seems} to be right and we are trying to prove it, we need to polish our statement as much as possible. Can we reduce the hypothesis? This other statement seems to \emph{just} be a rephrasing, is this the case? I believe my statement $A$ implies $B$, does the computer confirm that $B$ is true? Finding the route to the final proof, we might extract lemmas or other statements and the computer tests will be of great assistance. I consider that I have a certain pool of objects on which I can experiment to observe the results and guide me through my understanding\footnote{Note that I am not talking about formal proof assistants. This is another area altogether which I have not had time to explore.}.

\section{What to do with the code?}

My current file for research code on triangular partitions has more than 2500 lines, not including the code written by my student Lo\"ic Le Mogne. The code related to~\citemepre{CP22} and~\citemepre{CP23} is about 4500 lines. This contains implementations of the objects at stake and many test functions, failed attempts, leads, and so on. Very much like a draft paper, it is often badly written, with few comments and explanations. What is the value of this code? Does it have any and how can we use it?

\subsection{Contributions to open-source}

Computer exploration and experimental mathematics has developed with the general access to computer systems. As many mathematicians started to code, attempts to \emph{share} the produced code has arisen. I have been using the \Sage{} computer system since the beginning of my thesis. I believe in the open-source philosophy that the code we produce should be made accessible to all. Some objects are wildly used in combinatorics: permutations, Dyck paths, binary trees, and so on. They can be used in \Sage{} as the result of previous developments by other researchers. \Sage{} actually contains the basic grounds for most mathematical topics as can be read in~\citeext{SageBook}. 

I have myself contributed to \Sage{} whenever I felt that the code I had developed could be useful to others. My larger contribution is the implementation of Tamari interval-posets~\citemesoft{PonSage14} following my paper~\citeme{CP15}. Providing a direct implementation into \Sage{} was very useful in a research perspective: it helped Tamari interval-posets to become a major objects in the study of Tamari intervals as it could be used directly for computations. The initial class has been enriched by others, especially Frédéric Chapoton who works on this subject. Following my second paper~\citeme{Pon19}, we added together the \emph{rise-contact} involution which I had implemented during my investigations~\citemesoft{CPSage18}. Previously, I had also contributed to several methods on binary trees~\citemesoft{PonSage13}. These are just some examples, I have authored about~15 contributions and participated in about~50.

\subsection{Other good practices}

The road from writing research code for oneself and contributing to an open-source software (or develop a full independent software) is not easy. For one thing, much of the code we write for research does not belong to a general software like \Sage{}. A big portion of the code corresponds to failed attempts, or constructions that would eventually emerge in a different form in the final product. At this stage, it is not always easy to identify \emph{what} part of our code should belong the common shared knowledge of a software like \Sage{}. Besides, the code we write for research almost never follows the coding standard of a final \Sage{} contribution or of a software that could be distributed. There is often a long road of cleaning and documenting and testing. This takes a very long time that we, realistically, might not have. The result is that some useful, imperfect, code might never get shared.

Even if the quality is low, I believe that research code is worth sharing, especially when it leads to a published paper. Indeed, if I used this code for my own investigation, it could be useful for others. The objects used in the paper might turn out to be reused in another paper. Somebody else might decide to turn them into a proper software contribution. Or maybe, another researcher just needs to compute extra examples for their own purpose. This is why I have taken the habit to publish the research code along with the paper even if it is not a proper contribution to an existing software or a finished product. I first started a single \gls{github} repo where I would just post my research scripts organized by thematics along with some examples and presentations written with \gls{jupyter} notebooks~\citemesoft{PonSage18}. 

I provided an open license directly on the \gls{github} repo making it explicit that the code can be reused and, for example, integrated into another open-source software. To be fully efficient, one also needs to provide some \emph{environment} information about \emph{how} to run the code . I thus chose to include a \gls{docker} file which contains for example what \Sage{} version has been used to write the code. 

It became apparent that having only one repo for everything was not practical: I would often upgrade \Sage{} but not fix all my previous codes and the environment information would then be out of date. For my latest papers, I chose to create one repo per paper: this is what I did for~\citemesoft{PonSage22} as an annex to~\citemepre{Pon22}. I also have a specific code repo for this document~\citemesoft{PonSage23}. Finally, when I published the code related to~\citemepre{CP22}, I decided to create a specific archive of the code at the time of publication~\citemesoft{PonSage22b} using \gls{zenodo}. Indeed, as this code is still prone to be changed, it provides a link to a specific state of the software. 

As a summary, here is a list of codes related to my different papers:

\begin{itemize}
\item Code related to Tamari interval-posets (\citeme{CP15} and~\citeme{Pon19}) is available on my general repo~\citemesoft{PonSage18} with multiple demo notebooks. Some of it has been turned into \Sage{} contributions especially in~\citemesoft{PonSage14} and~\citemesoft{CPSage18}).
\item Code related to~\citeme{PP18, CPP19, PP20} with a demo notebook is available on my general repo~\citemesoft{CPSage18}.
\item Code related to~\citemepre{Pon22} with a demo notebook is available on a specific repo~\citemesoft{PonSage22}.
\item Code related to~\citemepre{CP22,CP23} with multiple demos is available on a specific repo and as a \gls{zenodo} archive~\citemesoft{PonSage22b}. 
\end{itemize}

Note that the code related to~\citemeconf{LMP23} that I mention in Section~\ref{sec:exp-process} is not yet published. Indeed, only an extended abstract has been published and the final paper is still in the making. Besides, Loïc Le Mogne also wrote some code and I would like to merge our contributions.

One reason for sharing this code is \emph{re-usability} which is why I make the effort of providing a running environment and often some demo notebooks that I edit as I write the corresponding papers. Speaking with historians and philosophers of mathematics, I also see another purpose. This ``unclean'' code is an equivalent of an early draft. It contains mistakes, failed attempts, ideas that sometimes were never fully developed. In other words, it contains some information about the process itself. Wether we share those codes or not, they hold an epistemic value and should at least be kept.

\chapter{Outreach}
\label{chap:outreach}

\chapcitation{Se vouloir libre, c'est aussi vouloir les autres libres.}{Simone de Beauvoir, Pour une morale de l'ambiguïté.}

In the past 10 years, I have been a promoter of the software \Sage{} and more generally of open-source development and free software. This has also led me to think more of my role as a scientist and the science community we build. In 2018, I gave the opening keynote at the programming event \gls{PyCon} Fr. I developed the idea that scientists and open source developers share a common ethic and motivation: we want to create knowledge and share it, we want contribute and to collaborate, we want our work to last to be re-used and improved.

Following this idea, it is important for me to take an active part in science communication and outreach at different levels. Within the math community, I work at promoting open-source software and programming knowledge in general. Within the scientific community at large, I have made presentations to explain the process of experimental mathematics. I also talk to general audiences and especially to middle and high schools to explain and share what research in mathematics and computer science mean. In this chapter, I give a general overview of this work.

\section{Advocate for open science}

Open science is the general idea that science productions should be made accessible to all: papers, data, software, etc. In mathematics and computer science, it concerns especially open-source software. Between 2015 and 2019, I have been part of the \gls{OpenDreamKit} project which was coordinated in my university Paris-Saclay and included 18 institutions and about 50 participants. We received significant funding that was used mostly for two things: hiring full paid developers to work on open-source software and organizing events within the different communities to foster collaboration and development.

I was the local coordinator for Paris-Saclay as well as the work package leader for ``Community Building, Training, Dissemination, Exploitation, and Outreach''. In the 4 years of the project, we were part of 110 events: either organized by members of the project or where they participated in a significant way. I was myself part of 23 events. For example, I co-organized the event ``Free Computational Mathematics at CIRM'' in 2019. The project also included 3 reviews in front of a European commission in which I took an active part presenting the activities of the project.  

My implication as a promoter of open science goes beyond my role in the project. I have organized several \gls{SageDays} workshop, even before the begining of \gls{OpenDreamKit}. In 2015, I organized a joint event at \gls{PyCon} US to help build relations between the \Sage{} community and the open-source community at large. At this occasion, I gave a presentation at the main conference to explain the concept of experimental mathematics~\citememisc{PonMisc15} to many non-mathematicians. 

I also participate in the effort to build a better scientific ecosystem, especially concerning scientific publications. I was an editor for the \emph{Journal of Open Source Software} between 2019 and 2022 to help better recognize the role of open-source development in scientific progress. I have also accepted a role as \emph{backend editor} for the new journal \emph{Combinatorial Theory} which is owned directly by its editorial board and charge no fees neither for readers nor writers.

\section{Talking about research to non researchers}

\begin{figure}[ht]
\center
\href{https://www.youtube.com/watch?v=RcXmhKF9ewo}{\includegraphics[height=5cm]{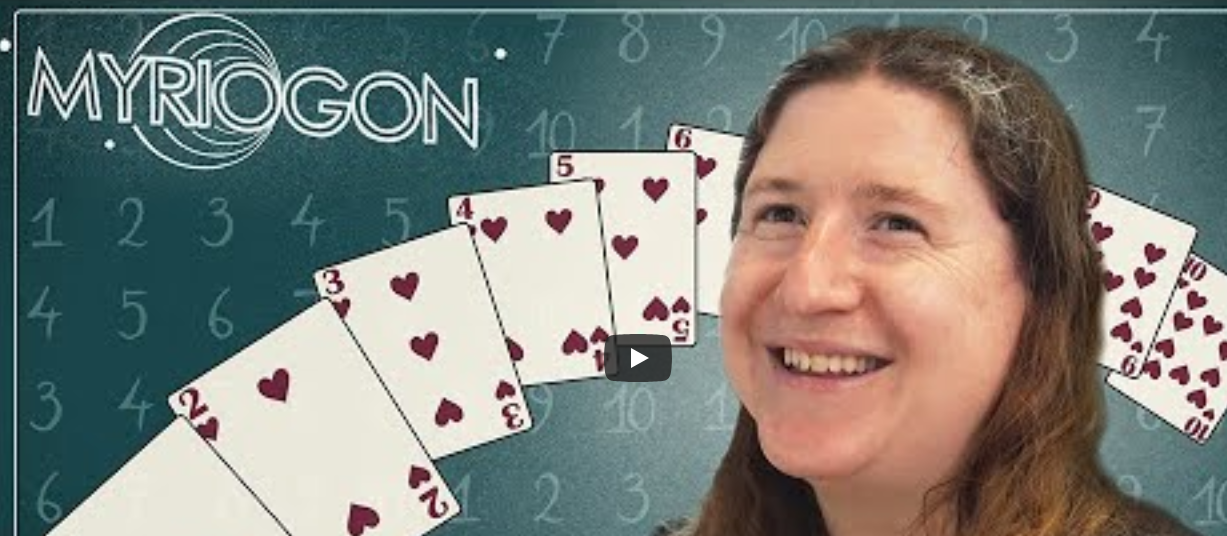}}

\caption{\emph{Des Cartes bien à leur place} for \emph{Le Myriogon} youtube channel}
\label{fig:myriogon}
\end{figure}

I also invest myself in promoting and explaining research to non researchers, especially in schools. I took part several times to the \gls{MathEnJean} project either by giving research talks to the schools, or by proposing research projects to pupils and meeting them throughout the year. Besides, I gave several general audience presentations and was invited to create an educational video with a professional team for the YouTube channel \emph{Le Myriogon}~\citememisc{Myr21}. I presented a simple bijection on permutations that I explained using playing cards as you can see on Figure~\ref{fig:myriogon}.

\section{Building a more inclusive community}

An open-source software is often free in the sense of ``free of charge'' (French word: \emph{``gratuit''}) and always free in the sense ``free to be reused and shared'' (French word: \emph{``libre''}). But this is not enough to make it equally accessible to all. This is also true for mathematics and computer science in general. Technical and societal barriers exist which affect certain demographics more than others. 

\begin{figure}[ht]
\center

\includegraphics[width=7cm]{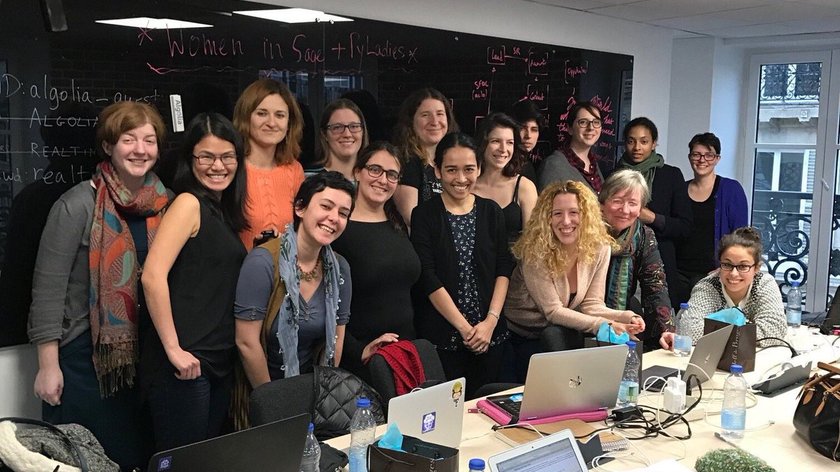}
\includegraphics[width=7cm]{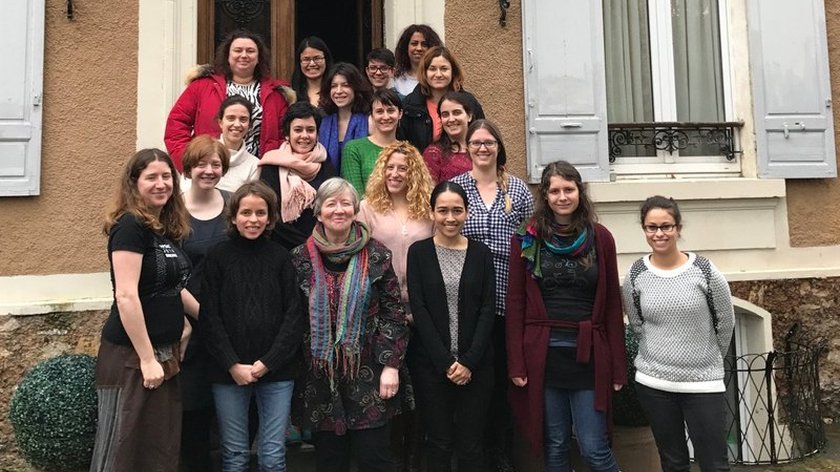}
\includegraphics[width=7cm]{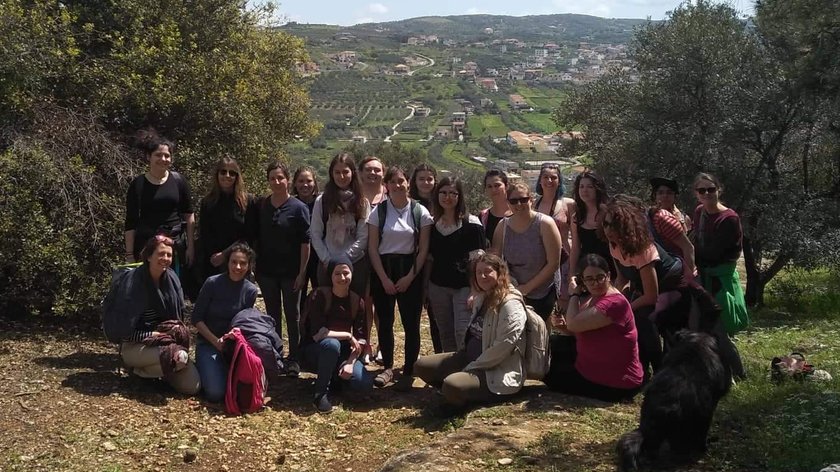}

\caption{Pictures taken at Women in Sage in 2017 and 2019}
\label{fig:womeninsage}
\end{figure}

Women especially suffer from negative stereotypes early on in relations to science and mathematics. This affects their own confidence and ability to solve problems~\citeext{DCL02}, the way they are perceived by others~\citeext{MDBGH12}, and their career choices~\citeext{NAV21}. Danica Savonick and Cathy N. Davidson have collected several studies on gender bias in academia and its combinations with other bias such as racism~\citeext{SD16}. As a result, women are largely underrepresented in mathematics and computer science~\citeext{BW18}. I have myself expressed the burden of \emph{being the only woman} in a blog post~\citememisc{PonMisc17b}. This has led me to take different actions. I have participated in many events such as the exhibition organized by \emph{Femmes et mathématiques}~\citeext{FetM22} to help make women in science more visible, especially to young girls. But I am also interested in building a better environment for current women scientists. For example, I run the \gls{PyLadies} parisian chapter between 2016 and 2020 where I organized weekly meetup with many women developers. I noticed that even at a very high level, women lack confidence in their ability to code and would be reluctant to join coding events such as \gls{SageDays} which often show a worst gender disparity than classical mathematics events. As part of the \gls{OpenDreamKit} project, I organized in 2017 and 2019 \emph{Women in Sage} events which were \gls{SageDays} workshops targeted at women. They were a great success. You can see some pictures on Figure~\ref{fig:womeninsage} and read the reports I published online~\citememisc{PonMisc17c, PonMisc19}.

\begin{figure}[h!]
\center

\includegraphics[height=5cm]{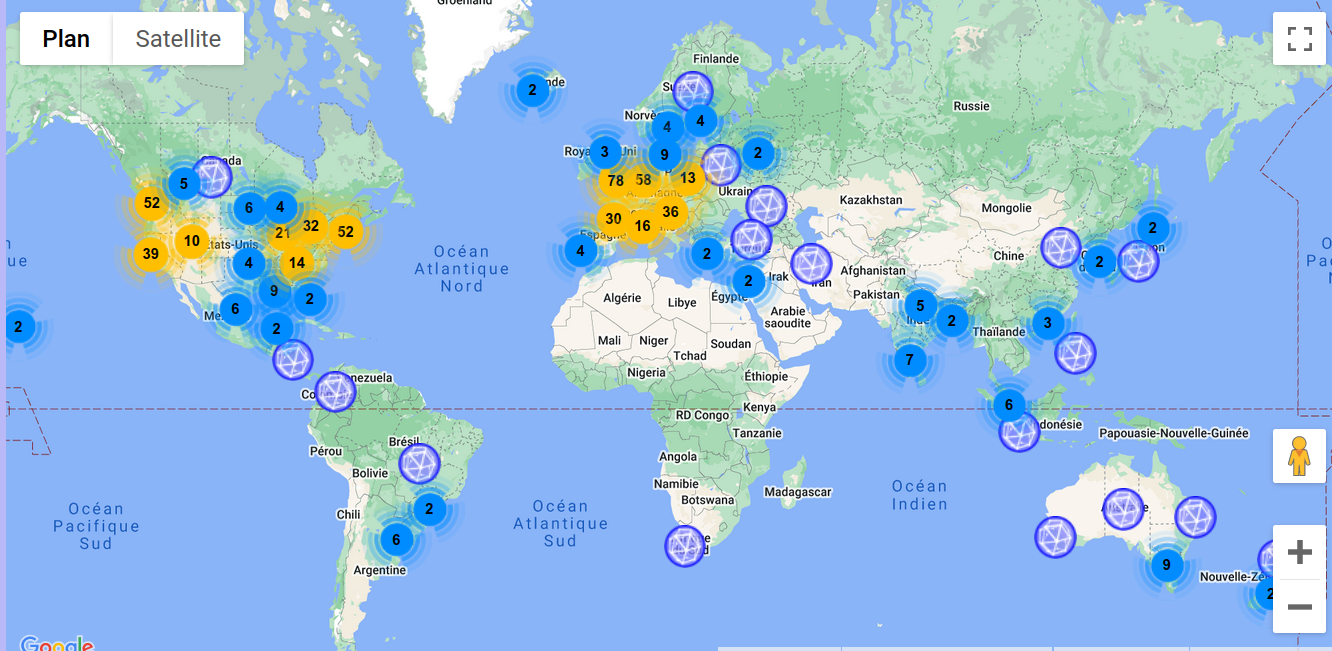}

\caption{The SageMath developer map taken from \url{https://www.sagemath.org/} in June 2023}
\label{fig:sagemap}
\end{figure}

Actual borders are physical barriers that stop people from doing research. Indeed, most resources are concentrated in ``western'' countries. In many parts of the world, students and researchers lack funding, basic infrastructure or even access to advanced education. This is visible in open-source software: see the big African gap on the \Sage{} developer map on Figure~\ref{fig:sagemap}. In 2019, two Nigerian women were supposed to join us at our \emph{Women in Sage} event organized in Archanes, Greece. Their trip was fully paid by the \gls{OpenDreamKit} project but they could not get their Visa in time. Nevertheless, the same year we had $3$ Nigerian researchers coming to our other event ``Free Computational Mathematics at CIRM''. They enjoyed the event very much and we decided to co-organize \gls{SageDays} in Nigeria the following summer. I did not attend myself as I was literally having a baby the same week. The project funded the trip of some \Sage{} instructors and relied on Ibadan university for local organization. It was a great success (see~\citeext{Bra19}) and also an indicator for the need of such events in West Africa as well as a lesson for us \Sage{} developers on how we can improve our software to allow an easier access.

In 2021, I co-organized and participated remotely to a \emph{Women in Sage} event in Senegal organized by Eliza Lorenzo. She was a participant to the \emph{Women in Sage} event I held in Paris in 2017. I had advertised the event to my contacts in Nigeria and some researchers from Ibadan were able to attend. This shows that these events create a good dynamic that can help West African researchers to reinforce their scientific network. Thanks to these events, I formed connections with two researchers: Ini Adinya from Nigeria and Olivia Nabawanda from Uganda with whom I try to build a scientific collaboration despite the systemic difficulties.

Thanks to \gls{OpenDreamKit}, I was also able to attend the event \gls{ECCO} (Colombia) in 2016 and 2018 as a \Sage{} instructor. These were extremely rewarding experiences (see~\citememisc{PonMisc16, PonMisc18}). The Colombian community in combinatorics has made great work in building an inclusive space and bringing Colombian and international students together. Following these events, I became a co-organiser for \gls{ECCO} 2022, especially in charge of the \gls{CIMPA} application and funding. I am now part of the \emph{ECCO Steering committee} whose role is to organize future events. Besides, \gls{ECCO} 2018 is where I met Daniel Tamayo Jiménez who is now finishing his PhD under my supervision and joined supervision of Vincent Pilaud.

\section{Political mathematics}

There is a general belief that mathematics is \emph{politically neutral}. Indeed, a theorem expresses a certain truth under some hypothesis, can it be racist or sexist? I do not believe this is the right question. Theorems are written by people in a certain context. Us, mathematicians, are part of a society which is, by nature, political. In~\citeext{God01, God03}, Roger Godement writes extensive preambles and notes to its books on Mathematical analysis. He questions the scientific responsibility of scientists in the creation of deadly weapons such as nuclear bombs. He also shows how political contexts and decisions impact what kind of science is considered \emph{good} or \emph{interesting} and thus developed.

In algebraic combinatorics, the long term applications of what we do, if they exist, may seem too far in the future to raise our concern. But this should not blind us in the illusion that what we do has no political implication. Defending fundamental science and the pursuit of knowledge for its own sake is already a political statement. Defending open-source software and open science is a political statement. Besides, there is the question of \emph{how} we do mathematics, \emph{who} do we teach it to, \emph{how} do we teach it, etc. Federico Ardila develops this idea in~\citeext{Ard20} and asks important questions about power.

Without any advanced statistical examination, I can say that I teach computer science mostly to young men and thus I am participating in a system which deprives young women of a useful knowledge in our society (I have chosen the gender disparity but, of course, many other inequalities are just as problematic). The academic system we are part of has some positive values and qualities and also some important setbacks. Of course, we cannot solve it all by ourselves. But we can choose to place our scientific activity in the general context of society and not see it as neutral. We can think about what we do and improve when we can. This is the choice I am always trying to make. This is also why I chose to include this last chapter. Not only does it relay some important activities I have run for the last 10 years but I consider these activities to be part of my scientific commitment just as much as my mathematical results. 

\bookmarksetup{startatroot}
\chapter*{Closing remarks}
\phantomsection
\addcontentsline{toc}{chapter}{Closing remarks}

\chapcitation{Truth is a matter of the imagination.}{Ursula K. Le Guin, The Left Hand of Darkness.}

I consider algebraic combinatorics to lie at the intersection of algorithmics and mathematics. My goal is to understand the inner structure of objects to exhibit the connection between topics that seem distant in nature. This is the case with the Tamari lattice and the Tamari intervals especially which I present in Chapters~\ref{chap:tamari-intervals},~\ref{chap:tam-stat}, and~\ref{chap:qt}. They appear in the context of planar maps as well as in representation theory. If I sometimes answer questions as in~\citeme{Pon19}, most of my papers actually \emph{create} new structures, like permutrees~\citeme{PP18}, or the $s$-weak order~\citemepre{CP22}. My hope is that those new objects deepen our understanding of the permutahedron and associahedron and create new points of views to explore.

My motivation comes from difficult open questions. I have mentioned the $q,t$-symmetry of Catalan numbers in Chapter~\ref{chap:qt}. To finish this document, I present two other open problems which have caught my attention recently. The Tamari lattice and its generalization the $\nu$-Tamari lattice can be understood in the more general framework of \defn{pipe dreams} and \defn{subword complexes}~\citeext{CLS14}. More precisely, pipe dreams are certain objects associated to permutations (each permutation gives a collection of pipe dreams). All pipe dreams of a given permutation are related through \defn{flips} and we can define a poset of increasing flips. For specific permutations, this gives the Tamari and $\nu$-Tamari lattices. In general, this is not a lattice but, in~\citeext{Rub12}, Rubey conjectures that it is indeed a lattice when flips are restricted to \defn{chute moves}. Beside, interesting geometrical structures appear which look like polyhedral complexes based on pure intervals like in our work on $s$-permutahedra~\citemepre{CP23}.

The second problem is related to the \defn{flip distance} of triangulations. In~\citeext{STT88}, Sleator, Tarjan, and Thurston prove a fundamental result on the associahedron. The diameter, \emph{i.e.}, the maximal rotation distance between two binary trees, is linear in the dimension $n$. The exact value was computed by Pournin in~\citeext{Pou14}. But it leaves open the following: what is the complexity of computing the minimal rotation distance between two given trees? This is the same as computing the minimal flip distance between two triangulations of a convex polygon. When the polygon is not convex, it is NP-complete~\citeext{OMP15} but the proof cannot be extended to the convex case. Besides, remember that the associahedron is \emph{in between} the permutahedron and the cube. In the permutahedron, the diameter is the maximal number of inversions, $\frac{n(n-1)}{2}$. In the cube, it is $n$. The minimal distance is also easy to compute in both those polytopes. Indeed, the distance between two permutations $\sigma$ and $\mu$ is the number of inversions of $\mu \sigma^{-1}$, and can be computed in $O(n^2)$. In the cube, the distance between two binary sequences is just the number of bits which are different. So it is computed in $O(n)$ where $n$ is the size of the sequence. It is interesting to notice these questions are easy in those two polytopes but difficult for the associahedron. What about permutreehedra and other generalized permutahedra? Can we use what we understand from the relations of these polytopes with the permutahedron and the cube to study those questions?

Solving such problems is a long term goal that fuels my current research. With each attempt, come new ideas that will flourish in unexpected manners, adding to the general understanding of the underlying structures. Besides, in addition to working on these questions myself, I believe that my role is to foster the scientific environment that makes progress possible. Encourage seminars, collaborations, workshops, open discussions, and freedom of investigation. This is the role I take especially as PhD advisor. With each of my students, I have tried to find the right balance between providing problems and questions and giving them the opportunity to express their own ideas and build their own research path. When I ask an open question, I do not always expect a final answer. I expect attempts and ideas. Then I am here to guide them through it, in a collaborative manner, to see where it could lead us. Sometimes, it eventually answers the original question, sometimes it opens new routes that we had not seen beforehand.

Someone told me an anecdote that they were once asked ``In your field, what do you think the next \emph{breakthrough} is going to be?''. I do not believe this is an interesting question. Breakthroughs do not happen because we chase after them (and their very nature makes them impossible to predict). They happen because we let our mind open to new ideas and keep inventing new truths.

\printglossaries 
\addcontentsline{toc}{chapter}{Glossary}

\backmatter
\lhead[\oldstylenums \thepage]{\rightmark}
\rhead[\leftmark]{\oldstylenums \thepage}

\lhead[\oldstylenums \thepage]{Bibliography}
\rhead[Biliography]{\oldstylenums \thepage}

\chapter*{Selected list of publications}
\phantomsection
\addcontentsline{toc}{chapter}{Selected list of publications}

\newcommand{\etalchar}[1]{$^{#1}$}

\chapter*{Other bibliographic references}
\phantomsection

\end{document}